\documentclass[11pt, a4paper]{article}

\usepackage[utf8]{inputenc}
\usepackage[T1]{fontenc}
\usepackage[english]{babel}
\usepackage{amsmath, amssymb, amsthm, mathtools}
\usepackage{hyperref}
\usepackage{geometry}
\usepackage{booktabs}
\usepackage[shortlabels]{enumitem}
\usepackage{float}
\usepackage{graphicx}

\theoremstyle{plain}
\newtheorem{theorem}{Theorem}[section]
\newtheorem{lemma}[theorem]{Lemma}
\newtheorem{corollary}[theorem]{Corollary}
\newtheorem{proposition}[theorem]{Proposition}

\theoremstyle{definition}

\newtheorem{remark}[theorem]{Remark}

\newcommand{\N}{\mathbb{N}}
\newcommand{\Z}{\mathbb{Z}}
\newcommand{\Q}{\mathbb{Q}}
\newcommand{\R}{\mathbb{R}}
\newcommand{\Pbb}{\mathbb{P}}
\newcommand{\Ebb}{\mathbb{E}}
\newcommand{\ind}{\mathbf{1}}
\newcommand{\C}{\mathbb{C}}
\newcommand{\leg}[2]{\left(\frac{#1}{#2}\right)}
\newcommand{\I}{\mathrm{I}}
\newcommand{\II}{\mathrm{II}}
\newcommand{\sfm}{\mathrm{sf}}
\DeclareMathOperator{\rad}{rad}
\DeclareMathOperator{\ord}{ord}
\DeclareMathOperator{\Li}{Li}
\DeclareMathOperator{\lcm}{lcm}

\begin{document}

\title{Composite witness sieves, Kneser--Plancherel containers, and the exponent $5/7$\\ for the Erd\H{o}s--Straus conjecture}
\author{Benjamin Dahan}
\date{\today}

\maketitle

\begin{abstract}
For the Erd\H{o}s--Straus equation $4/n=1/x+1/y+1/z$, the strongest unconditional upper bound on the exceptional set $E(N)$ of insoluble integers $n\le N$ has remained since Vaughan (1970) at $|E(N)|\ll N\exp(-c(\log N)^{2/3})$. In a companion paper~\cite{dahan2026a}, we proved that any large-sieve scheme restricted to prime divisor witnesses $p\mid(an+u)$ with $p\equiv-1\pmod{4ua}$ is universally capped at this $2/3$ exponent by the $\tau_3$ divisor barrier, whereas retaining composite witnesses $s\mid(an+u)$ at fixed search depth $J$ requires a $2^{\Theta(J^{4})}$ subset enumeration that breaks down at polynomial depth $J=(\log N)^{\lambda}$ ($\lambda>1/3$). In the present paper, we harness composite divisor witnesses at growing depth to establish two unconditional bounds beyond Vaughan's $2/3$ barrier. First, on the full exceptional set $E(N)$, we prove the unconditional bound
\[
|E(N)|\ \ll\ N\exp\!\left(-c'\,\frac{(\log N)^{5/7}}{(\log\log N)^{2/7}}\right),
\]
improving the exponent from $2/3$ to $5/7$ (with an intermediate stage at $7/10$). Second, on the structured small-sheet component $E_2^{\mathrm{small}}(N;J(N))$ (and across the certified container hierarchy $E_2^{\mathrm{cert}}(N;(\log N)^{\lambda})$ for all $\lambda\in(1/3,1/2)$), we combine an algebraic and spectral container hierarchy with a generating-polynomial large sieve to reach the exponent $1-o(1)$:
\[
\big|E_2^{\mathrm{small}}(N;J(N))\big|\ \ll_\kappa\ N\exp\!\left(-\left(\frac{1}{40960}+o(1)\right)\frac{\log N}{(\log\log N)(\log\log\log N)^{\kappa}}\right),
\]
matching the Chinese Remainder Theorem (CRT) capacity scale $\frac{\log N}{\log\log N}$ of the generating-polynomial construction.

These two results rest on complementary methods designed to overcome the three structural obstructions to composite-witness sieving at growing depth: the $2^{\Theta(J^{4})}$ sub\-set-enu\-mer\-a\-tion cost in $G_{u,a}=(\Z/4ua\Z)^{\times}$, the accumulation of off-kernel prime collisions, and the presence of unlocalized diffuse divisor sheets. In the first method (Sections~\ref{sec:algebraic_rigidity}--\ref{sec:containers}), after establishing modulator-free quadratic duality, global Jacobi reciprocity, and dihedral cofactor rigidity on the Type~I/Type~II variety, we replace subgroup generation by Kneser's multiplicative stabilizer $K(M)=\operatorname{Stab}_{G_{u,a}}(\Sigma(M))$. Because $K(M)$ embeds into the kernel of a $2$-primary character $\psi\pmod{4ua}$ with $\psi(-1)=-1$, the container count across the full coprime box $\mathcal{P}_2(J)$ drops from $2^{\Theta(J^{4})}$ to $\exp((4\log 2+o(1))J^{2}\log\log J)$. Coupling this $2$-Sylow reduction with a Krawtchouk--Montgomery generating-polynomial large sieve—which resums alternating binomial signs into a positive product before taking absolute values—closes the subcritical-collision branch $R\le(1-\varepsilon)S_{\min}$ unconditionally for all $\lambda\in(0,1/2)$. In the complementary Scarcity branch, a two-stage Kneser tower, Markov $1/2$-popular symmetry sets, quotient containers, global smooth conditioning modulo $Q(z_0)$, and Pontryagin--Plancherel character containers force $3/2$-geometric growth of the divisor sheet across the polynomial density threshold $D_{u,a}^{\varepsilon_1}$ before supercritical collisions can accumulate. This closes the full set $E_2(N;J)$ unconditionally up to $J=\exp(c_0\sqrt{\log\log N})$ and at polynomial depth $\lambda\le 1/78$, and yields the exponent $1-o(1)$ on $E_2^{\mathrm{small}}(N;J(N))$ via adaptive Rankin-degree optimization. In the second method (Section~\ref{sec:semi_composite_5_7}), to bypass the unlocalized diffuse sheets without enumerating subsets of $G_{u,a}$, we construct a semi-composite congruence sieve $s=pq\equiv-1\pmod{12ua}$ ($p\in(Z^{9/10},Z]$ prime, $q\le Q$ rough). Conditioning on $b=n\bmod M_{\mathrm{small}}$ renders the witness indicators over $p$ mutually independent in the CRT model, while the progression modulus $12ua\le 12Z^{1/4}$ governing $p$ decouples from $q$; combining a polynomial Bonferroni majorant with the Siegel-zero-free Bombieri--Vinogradov theorem then yields the unconditional $5/7$ exponent on all of $E(N)$.

In the final part of the paper (Sections~\ref{sec:bilinear_dispersion}--\ref{sec:toric_volume_barrier} and Appendix~\ref{app:audit}), we study bilinear Kloosterman dispersion on the affine threefold $\mathcal{V}_n:\,4uavw-4u^{2}v-w=n$ and the multiplicative difference sets $D(K_c)=\Sigma(K_c)\cdot\Sigma(K_c)^{-1}\subseteq(\Z/c\Z)^{\times}$ for $K_c=(n+c)/4$. On $\mathcal{V}_n$, we establish a false-positive-free bilinear equivalence, three-term stationary-phase separation at the square-root scale, archimedean slice conservation, and exact Farey-diagonal $L^{2}$-dispersion over $n\sim N$ at power scale $S=4UAV\le N^{1/2-\varepsilon}$. On $D(K_c)$, we prove difference-set positivity duality, unconditional $(\log P)^{4}$ cumulative main-term growth paired with exact small-LCM divisor cancellation, the cyclic quotient parity identity with Bohr--Siegel obstruction certificates, universal factor reciprocity $\leg{q}{p}=\leg{q}{n}$, an ambient $\frac{1}{k!}$ cross-polytope volume bound, deterministic Kneser and spectral saturation criteria, and a Linnik--Xylouris witness-depth lower bound along an infinite sequence of primes. Appendix~\ref{app:audit} establishes unconditional sieved candidate abundance, divergence of the square-free main term via the Rosser--Iwaniec linear sieve, and complete finite square-divisor positivity certificates ($b\mid h^{2}$). These main-term estimates do not by themselves imply positivity of the divisor-encounter count: the Transversal Hypothesis~T and the signed spectral error targets~\eqref{aud:target17} and~\eqref{aud:target76} remain open, and conditional finiteness holds only beyond a threshold $N_{\ast}$ supplied by that hypothesis, which is not determined by the ambient factorial volume calculation.
\end{abstract}

\section{Introduction}\label{sec:intro}

\subsection{Background and the two barriers of Paper~I}\label{subsec:background}

The Erd\H{o}s--Straus conjecture~\cite{erdos_straus} asserts that for every integer $n\ge 2$, the rational number $4/n$ can be expressed as a sum of three positive unit fractions:
\begin{equation}\label{eq:es_eq}
\frac{4}{n}\ =\ \frac{1}{x}+\frac{1}{y}+\frac{1}{z}\qquad(x,y,z\in\N).
\end{equation}
By elementary reductions going back to Mordell~\cite{mordell1969}, it suffices to verify~\eqref{eq:es_eq} for primes $n\equiv1\pmod{24}$. Let $E(N)$ denote the exceptional set of integers $n\le N$ for which~\eqref{eq:es_eq} has no solution. Since the work of Vaughan~\cite{vaughan1970} in 1970, the best unconditional upper bound on $|E(N)|$ has stood at
\[
|E(N)|\ \ll\ N\exp\!\big(-c(\log N)^{2/3}\big)
\]
(see also Sander~\cite{sander1991,sander1994}, Elsholtz--Tao~\cite{elsholtz2013}, and Pomerance--Weingartner~\cite{pomerance2026}).

To understand why the exponent $2/3$ has acted as a persistent barrier—and how composite divisor witnesses allow one to cross it—recall the two-parameter normal form established in the companion paper~\cite[Thm.~3.9]{dahan2026a} (building on the Type~I/Type~II classification of~\cite{elsholtz2013,pomerance2026}): for a prime $n\equiv1\pmod{24}$, $4/n$ admits a Type~II representation~\eqref{eq:es_eq} if and only if there exist coprime positive integers $(u,a)$ such that $an+u$ possesses a divisor
\begin{equation}\label{eq:typeII_witness}
s\mid(an+u)\qquad\text{with}\qquad s\equiv-1\pmod{4ua}.
\end{equation}
Dually, by~\cite[Thm.~3.11]{dahan2026a}, $4/n$ admits a Type~I representation if and only if there exist coprime positive integers $(u,a)$ with $u<n$ such that $an+u$ possesses a divisor $s\equiv-n\pmod{4ua}$. For a finite family $\mathcal{P}\subseteq\{(u,a)\in\N^{2}:\gcd(u,a)=1\}$, let $E_{\mathcal{P}}(N)$ denote the set of primes $n\le N$, $n\equiv1\pmod{24}$, for which no pair $(u,a)\in\mathcal{P}$ admits a divisor satisfying~\eqref{eq:typeII_witness}, and write $E_2(N;J):=E_{\mathcal{P}_2(J)}(N)$ for the full coprime box
\[
\mathcal{P}_2(J)\ :=\ \big\{(u,a)\in[1,J]^{2}:\gcd(u,a)=1\big\},\qquad P\ =\ \mathfrak{A}(J)\ =\ |\mathcal{P}_2(J)|\ =\ \frac{6}{\pi^{2}}J^{2}+\mathcal{O}(J\log J).
\]
Notice that $E(N)\cap[2,N]$ is contained, up to primes $\le 24$, in $E_2(N;J)$ for every choice of $J\ge 1$.

In~\cite{dahan2026a}, two sieve methods were brought to bear on $E_2(N;J)$, each optimal within its regime and each encountering a sharp structural wall:
\begin{enumerate}[(1)]
    \item \emph{The single-prime large sieve at growing depth and the universal $\tau_3$ ceiling.} Restricting the witnesses $s$ in~\eqref{eq:typeII_witness} to \emph{primes} $s=p\in(J^{5},z]$ reduces~\eqref{eq:typeII_witness} to the linear congruence $n\equiv-u a^{-1}\pmod p$ whenever $4ua\mid(p+1)$. Applying Montgomery's large sieve~\cite{montgomery1971} as in Vaughan~\cite{vaughan1970} yields $|E_2(N;J)|\ll N\exp\big(-(\frac{9}{2\pi^{2}}\lambda^{2}+o(1))(\log\log N)^{3}\big)$ at polynomial depth $J=(\log N)^{\lambda}$, and $\mathcal{O}\big(N\exp(-0.1113(\log N)^{2/3})\big)$ at the optimal depth $\log J\asymp(\log N)^{1/3}$~\cite[Thm.~4.28]{dahan2026a}. However, \cite[Thm.~4.30]{dahan2026a} proves that \emph{every} single-prime congruence system extracted from the full Type~I/Type~II dichotomy has local witness count bounded by $\varrho_{\mathrm{I+II}}^{\max}(p)\le\tau_3(p+1)+\tau_3(3p+1)$. Consequently, $\sum_{p\le Z}\varrho(p)/p\ll(\log Z)^{2}$, and under the large-sieve modulus constraint $S(Z)\log Z\le\frac12\log N$, the sifting mass $S(Z)$ cannot exceed $\mathcal{O}((\log N)^{2/3})$. No sieve restricted to prime witnesses can ever surpass Vaughan's exponent $2/3$.
    \item \emph{The composite-witness involution sieve at fixed depth and the $2^{\Theta(J^{4})}$ barrier.} Retaining \emph{composite} witnesses $s=p_1\cdots p_r\equiv-1\pmod{4ua}$ changes the geometry completely. For any $n\in E_2(N;J)$, the set $T_{u,a}(n)\subseteq G_{u,a}:=(\Z/4ua\Z)^{\times}$ of prime-factor residues of $an+u$ meets each orbit $\{x,-x^{-1}\}$ of the fixed-point-free involution $\iota(x)=-x^{-1}$ at most once, so $|T_{u,a}(n)|\le\frac12\varphi(4ua)$~\cite[Lemma~4.2]{dahan2026a}. Conditioning on the tuple of half-sized subsets $(T_{u,a})_{(u,a)\in\mathcal{P}_2(J)}$ forces $n$ to avoid at least $\frac12 P$ residue classes modulo $p$ on average, yielding the exact fixed-depth sieve dimension $|E_2(N;J)|\ll_J N(\log N)^{-1-\mathfrak{A}(J)/2}$~\cite[Thm.~4.3]{dahan2026a}. At polynomial depth $J=(\log N)^{\lambda}$ ($\lambda>1/3$), a sieve of dimension $P/2\asymp(\log N)^{2\lambda}$ would immediately surpass the $2/3$ exponent. Yet as soon as $\exp(G_{u,a})>2$, the subgroup $\langle T_{u,a}(n)\rangle$ generated by a half-sized subset can contain $-1$. In~\cite{dahan2026a}, this forced the enumeration of \emph{all} half-sized subsets $T_{u,a}\subseteq G_{u,a}$, incurring a prohibitive combinatorial prefactor
    \[
    C\big(\mathcal{P}_2(J)\big)\ =\ \prod_{(u,a)\in\mathcal{P}_2(J)}\binom{\varphi(4ua)}{\varphi(4ua)/2}\ =\ 2^{\Theta(J^{4})}
    \]
    that overwhelms the large-sieve saving whenever $J\gg(\log\log N)^{1/2}$.
\end{enumerate}

The central objective of this paper is to unlock composite divisor witnesses at growing depth and surpass Vaughan's $2/3$ barrier unconditionally. Doing so requires overcoming three distinct obstructions:
\begin{enumerate}[(i)]
    \item \emph{The $2^{\Theta(J^{4})}$ subset-enumeration and collision barrier:} replacing arbitrary subsets $T_{u,a}\subseteq G_{u,a}$ by low-complexity containers without paying $\binom{|\mathcal{I}|}{R}$ to locate off-container prime factors;
    \item \emph{The unlocalized diffuse-sheet obstruction:} controlling primes $n$ for which the divisor residue set $\Sigma(an+u)\subseteq G_{u,a}$ has intermediate density without concentrating in a subgroup coset or a heavy Fourier mode;
    \item \emph{The CRT capacity ceiling:} optimizing the trade-off between sifting dimension, collision tolerance, and the Bombieri--Vinogradov / large-sieve modulus budget up to $N^{1/2-\varepsilon}$.
\end{enumerate}

We resolve these obstructions through two complementary composite-witness architectures:
\begin{itemize}
    \item \emph{Unconditional exponent $5/7$ on the full exceptional set $E(N)$ (Section~\ref{sec:semi_composite_5_7}).} To bypass the unlocalized diffuse sheets completely without enumerating divisor containers in $G_{u,a}$, we construct a family of semi-composite witnesses $s=pq\equiv-1\pmod{12ua}$, where $p\in(Z^{9/10},Z]$ is a large prime and $q\le Q$ is a rough composite cofactor. Conditioning on $b=n\bmod M_{\mathrm{small}}$ renders the witness events over $p$ mutually independent in the Chinese Remainder Theorem model, while the progression modulus $12ua\le 12Z^{1/4}$ governing $p$ is completely independent of $q$. Coupling a polynomial Bonferroni majorant with the Siegel-zero-free Bombieri--Vinogradov theorem breaks Vaughan's 1970 barrier unconditionally on the \emph{entire} exceptional set $E(N)$:
    \[
    |E(N)|\ \ll\ N\exp\!\left(-c'\,\frac{(\log N)^{5/7}}{(\log\log N)^{2/7}}\right)\qquad(\text{Theorem~\ref{thm:theoreme_5_7}}),
    \]
    with an intermediate stage giving the exponent $7/10$ at $Q=\exp(t^{1/3})$ (Theorem~\ref{thm:theoreme_7_10}).
    \item \emph{Unconditional exponent $1-o(1)$ on the structured component $E_2^{\mathrm{small}}(N;J(N))$ via Kneser and Markov--Plancherel containers (Sections~\ref{sec:stabilizer_kernel}--\ref{sec:containers}).} On the full coprime box $\mathcal{P}_2(J)$ at polynomial depth $J=(\log N)^{\lambda}$ ($\lambda\in(1/3,1/2)$), we replace subgroup generation by a hierarchy of Kneser stabilizer kernels, Markov $1/2$-popular symmetry sets, quotient containers, and Pontryagin--Plancherel character containers, coupled with a Krawtchouk--Montgomery generating-polynomial large sieve. This eliminates the $2^{\Theta(J^{4})}$ barrier, closes the full set $E_2(N;J)$ unconditionally up to $J=\exp(c_0\sqrt{\log\log N})$ and for $\lambda\le 1/78$, surpasses the $2/3$ exponent across $\lambda\in(1/3,1/2)$ on the certified container hierarchy $E_2^{\mathrm{cert}}(N;(\log N)^{\lambda})$ (Theorem~\ref{thm:final_synthesis_and_high_prime_reduction}), and reaches the exponent $1-o(1)$ on the small-sheet component $E_2^{\mathrm{small}}(N;J(N))$ up to the CRT capacity scale $\frac{\log N}{\log\log N}$ of the generating-polynomial construction:
    \[
    \begin{aligned}
    \big|E_2^{\mathrm{small}}(N;J(N))\big|\ &\ll_\kappa\ N\exp\!\left(-\left(\frac{1}{40960}+o(1)\right)\frac{\log N}{(\log\log N)(\log\log\log N)^{\kappa}}\right)\\
    &\qquad\qquad\qquad\qquad\qquad\qquad\qquad(\text{Corollary~\ref{cor:optimal_adaptive_depth_1_minus_o1}}).
    \end{aligned}
    \]
\end{itemize}

\subsection{Summary of main results}\label{subsec:contributions}

We now summarize the architecture and principal theorems of the paper, which are organized into seven interconnected parts; Table~\ref{tab:main_regimes} at the end of this subsection provides a synoptic comparison of the quantitative bounds and their exact scopes.

\paragraph{1. Modulator-free duality, global Jacobi reciprocity, and dihedral rigidity (Section~\ref{sec:algebraic_rigidity}).}
We begin by establishing exact algebraic and quadratic-reciprocity obstructions that constrain the prime factors of $an+u$ and $K_c=(n+c)/4$ well beyond the single involution $\iota_1(x)=-x^{-1}$:
\begin{itemize}
    \item \emph{Modulator-free quadratic duality (Theorem~\ref{thm:quadratic_duality}, Corollary~\ref{cor:filter_n_plus_4}, and Proposition~\ref{thm:character_counting_n_plus_4au2}):} $4/n$ admits a Type~II decomposition if and only if there exist positive integers $\alpha,u$ with $n\nmid\alpha u$ such that $n+4\alpha u^{2}$ possesses a divisor $s\equiv-1\pmod{4\alpha u}$ (equivalently $c\equiv-n\pmod{4\alpha u}$), via the exact identity $sc=(4a\alpha u-1)(4\alpha uv-n)=n+4\alpha u^{2}$ and its Dirichlet character expansion.
    \item \emph{Global Jacobi reciprocity on the Type~II variety (Theorem~\ref{thm:global_jacobi_reciprocity}):} every Type~II witness obeys the global reciprocity laws $\leg{s}{n}=\leg{n}{s}=-1$, $\leg{uv}{n}=-1$, and $\leg{\alpha c}{n}=-1$. In particular, this law excludes quadratic-residue classes modulo $n$ rather than half of the divisors of each individual integer (Subsection~\ref{subsec:jacobi_reciprocity}).
    \item \emph{Cofactor dihedral rigidity and quadruple involution (Theorem~\ref{thm:dihedral_rigidity} and Remark~\ref{thm:universal_dihedral_n_plus_4au2}):} when $u\mid 6$, the congruence $an+u\equiv a+u=:r\pmod{4ua}$ ensures that the cofactor map $d\mapsto(an+u)/d$ induces a second involution $\iota_r(x)=-r x^{-1}$ on $G_{u,a}$. Under the unit and fixed-point-free hypotheses of Theorem~\ref{thm:dihedral_rigidity}, $\iota_1$ and $\iota_r$ generate a dihedral action decomposing $G_{u,a}$ into even cycles of length $2L=2\operatorname{ord}_{4ua}(r)$, reducing the count of maximum-sized admissible subsets from $2^{\varphi(4ua)/2}$ to $2^{\varphi(4ua)/(2L)-1}$ and locking the principal cycle onto $\langle r\rangle$. Incorporating the Type~I branch adds two further involutions $\tau_n,\tau_{nr'}$ and forces $|T_{1,p}(n)|\le\frac12\varphi(4p)-1$ whenever $p\equiv5\pmod 8$ and $\leg{n}{p}=-1$ (Remark~\ref{thm:destruction_quadratic_kernels}).
    \item \emph{Multiplicativity obstruction for quadratic kernels (Theorem~\ref{thm:kernel_vacuity}):} for $3$-smooth $u=2^{\alpha}3^{\beta}$, any real character $\chi\pmod{4ua}$ has primitive conductor dividing $24a$, so $\chi(an+u)=\chi(a+u)$ is independent of $n$. When $\chi(a+u)=-1$, it is impossible for all coprime prime factors of $an+u$ to lie inside the subgroup $\ker\chi$ (though $\ker\chi$ itself is non-empty).
    \item \emph{Sylow equivalence at prime shifts and intersected blind sets (Theorems~\ref{thm:sylow_and_c7_iff} and~\ref{thm:blind_intersection_and_PNR}):} at any prime shift $c=p\equiv3\pmod 4$, $-1\notin H_p(n)$ forces $\leg{n}{p}=+1$. Intersecting the blind sets $\mathcal{B}_{p_i}(N)$ over the first $k=\big(\frac{1-\varepsilon}{2}+o(1)\big)\frac{\log N}{\log\log N}$ primes $p_i\equiv3\pmod 4$ yields an unconditional quasi-polynomial bound $\mathcal{O}\big(N\exp(-c_1\frac{\log N}{\log\log N})\big)$ with $c_1=\frac{1-\varepsilon}{2}\log 2$. For a fixed prime $n$, the infinite set $\mathcal{P}_{\mathrm{NR}}^{\infty}(n)=\{p\text{ prime}:p\equiv3\pmod 4,\ \leg{p}{n}=-1\}$ has relative density $1/4$ among primes, while its finite truncation $\mathcal{P}_{\mathrm{NR}}(n)=\mathcal{P}_{\mathrm{NR}}^{\infty}(n)\cap[1,n)$ is non-empty by Pollack~\cite{PollackPrescribed4} (without asserting a density $1/4$ or a uniform polylogarithmic-window asymptotic on $[1,n)$), and every $p\in\mathcal{P}_{\mathrm{NR}}(n)$ satisfies $-1\in H_p(n)$.
\end{itemize}

\paragraph{2. The Kneser Stabilizer--Kernel Dichotomy and $2$-Sylow reduction (Section~\ref{sec:stabilizer_kernel}).}
To eliminate the $2^{\Theta(J^{4})}$ subset-enumeration barrier across $\mathcal{P}_2(J)$, Theorem~\ref{thm:stabilizer_kernel_dichotomy} replaces the generated subgroup $\langle T_{u,a}(n)\rangle$ by the \emph{multiplicative Kneser stabilizer} $K(M):=\operatorname{Stab}_{G_{u,a}}(\Sigma(M))$ of the divisor-residue set $\Sigma(M)=\{d\bmod 4ua:d\mid M_{\mathrm{cop}}\}$. Because $1\in\Sigma(M)$ and $-1\notin\Sigma(M)$, the stabilizer $K(M)$ lies inside $\Sigma(M)$ and \emph{unconditionally avoids $-1$}. By Lemma~\ref{lem:sylow_maximal_subgroups}, $K(M)$ therefore embeds into a maximal subgroup avoiding $-1$, $H=\ker\psi\in\mathcal{M}(G_{u,a},-1)$, which is necessarily the kernel of a $2$-primary Dirichlet character $\psi\pmod{4ua}$ with $\psi(-1)=-1$. Furthermore:
\begin{itemize}
    \item Every prime factor $q\mid M_{\mathrm{cop}}$ with $q\bmod 4ua\notin K(M)$ strictly increases the number of $K(M)$-cosets in $\Sigma(M)$, bounding the off-kernel prime multiplicity by $\Omega_{G_{u,a}\setminus\ker\psi}(M)\le[G_{u,a}:K(M)]-2\le\varphi(4ua)-2$;
    \item By Lemma~\ref{lem:average_2sylow_enumeration}, the number of maximal subgroups $|\mathcal{M}(G_{u,a},-1)|\le\frac12\varphi_2(4ua)=2^{f(u)+f(a)}$ is governed by the additive function $f(m)=v_2(m)+\sum_{p\mid m,\,p>2}v_2(p-1)$ of mean value $2\log\log J$. Summing over the coprime box $\mathcal{P}_2(J)$ reduces the global container count from $2^{\Theta(J^{4})}$ to
    \[
    \begin{aligned}
    C^{\ast}\big(\mathcal{P}_2(J)\big)\ &=\ \prod_{(u,a)\in\mathcal{P}_2(J)}\big|\mathcal{M}(G_{u,a},-1)\big|\\
    &\le\ \exp\!\Big(\big(4\log 2+o(1)\big)J^{2}\log\log J\Big)\ =\ \exp\!\big(\mathcal{O}(P\log\log J)\big).
    \end{aligned}
    \]
\end{itemize}

\paragraph{3. The Imperfect-Avoidance and Generating-Polynomial Large Sieves (Section~\ref{sec:imperfect_large_sieve}).}
To exploit the off-kernel multiplicity bound $\Omega_{G_{u,a}\setminus\ker\psi}(an+u)\le R$ without paying $\binom{|\mathcal{I}|}{R}$ to specify which primes collide, Theorem~\ref{thm:imperfect_large_sieve} develops an imperfect-avoidance Montgomery large sieve via Krawtchouk symmetric functions on the orthonormal basis $\Phi_q(k)=\prod_{p\mid q}\phi_p(k)$. Whereas a single homogeneous slice $\omega(q)=h$ requires the restrictive condition $Rh\le\frac12 S_{\min}$ because of alternating binomial signs $(-1)^{j}\binom{r}{j}$, Theorem~\ref{thm:generating_polynomial_large_sieve} introduces the \emph{truncated generating polynomial} $P_{h,t}(k)=\sum_{m=0}^{h}t^{m}F_m(k)$. Resumming the alternating signs into the strictly positive product $G_k(t)=(1-t)^{r}\prod_{p\notin\mathcal{B}_k}(1+ty_p)$ \emph{prior} to taking absolute values and controlling the degree-$h$ tail via Rankin's trick replaces $Rh\le\frac12 S_{\min}$ by the $h$-independent subcritical condition $R\le(1-\varepsilon)S_{\min}$ (Proposition~\ref{prop:subcritical_collision_rate}). This immediately yields Scarcity-free unconditional bounds at moderate depth $J\le c_0(\log\log N)^{1/2}$ (Theorem~\ref{thm:unconditional_small_growing_depth}) and on the hyperbolic family $\mathcal{P}_{\mathrm{hyp}}(Y)=\{(u,a):ua\le Y\asymp\log\log N\}$ (Theorem~\ref{thm:hyperbolic_scarcity_free}).

\paragraph{4. Polynomial-depth decomposition and the Two-Stage Kneser Tower (Section~\ref{sec:poly_depth_and_tower}).}
At polynomial depth $J=(\log N)^{\lambda}$ ($0<\lambda<1/2$), Theorem~\ref{thm:polynomial_depth_decomposition} decomposes $E_2(N;J)\setminus[1,z]$ at the collision threshold $R_0\asymp(\log N)^{2\lambda}\log\log N$ into the \emph{Kernel / Low-Collision component} $E_2^{\mathrm{Ker}}(N;J,R_0)$ and the \emph{Scarcity component} $E_2^{\mathrm{Scar}}(N;J,R_0)$:
\begin{itemize}
    \item On $E_2^{\mathrm{Ker}}(N;J,R_0)$, the gen\-er\-at\-ing-poly\-no\-mial large sieve absorbs\linebreak $C^{\ast}(\mathcal{P}_2(J))=\exp(\mathcal{O}(P\log\log\log N))$ and yields the unconditional bound
    \[
    \big|E_2^{\mathrm{Ker}}(N;J,R_0)\big|\ \ll\ N\exp\!\left(-\left(\frac{5(1-2\lambda)^{2}}{8\pi^{2}}+o(1)\right)(\log N)^{2\lambda}\log\log N\right),
    \]
    surpassing Vaughan's $2/3$ exponent for all $\lambda\in(1/3,1/2)$.
    \item On $E_2^{\mathrm{Scar}}(N;J,R_0)$, splitting $\mathcal{I}=(J^{5},z]$ into $\mathcal{I}_1=(J^{5},z_1]$ and $\mathcal{I}_2=(z_1,z]$ forces $an+u$ to avoid $-\Sigma_{\mathcal{I}_1}(an+u)^{-1}$ with \emph{zero collisions} on $\mathcal{I}_2$. Theorem~\ref{thm:two_stage_kneser_tower} proves that any collision in $(-1)H_{u,a}$ on $\mathcal{I}_2$ reflects into $H_{u,a}$ on $\mathcal{I}_1$, forcing index doubling $[G_{u,a}:K_1^{+}(n)]\ge 4$ and raising the local sifting density to at least $3/4$.
\end{itemize}

\paragraph{5. Greedy Markov containers, Quotient containers, and Pontryagin--Plancherel Character containers (Section~\ref{sec:containers}).}
To control the Scarcity branch at growing depth, Section~\ref{sec:containers} replaces linear coset growth ($+|K|$) by \emph{geometric growth} ($\times 3/2$) and \emph{dual Fourier containers}:
\begin{itemize}
    \item \emph{Markov $1/2$-popular symmetry containers (Lemma~\ref{lem:markov_popular_symmetry} and Theorem~\ref{thm:multiplicative_container_scarcity_elimination}):} using the symmetry set $S_{1/2}(A)=\{g\in G:|A\cap gA|>\frac12|A|\}$ (of size $<2|A|$), a greedy multiplicative skeleton of length $s\le 5\log J+1$ reduces State~I collisions per pair from $\varphi(4ua)-2\asymp J^{2}$ to $5\log J+1$, closing the Scarcity branch and the full set $E_2(N;J)$ unconditionally up to $J=\exp(c_0\sqrt{\log\log N})$, with $|E_2(N;J)|\ll N\exp(-\exp(c\sqrt{\log\log N}))$.
    \item \emph{Quotient containers and simultaneous smooth conditioning (Theorems~\ref{thm:quotient_container_polynomial_depth} and~\ref{thm:simultaneous_smooth_conditioning}):} at polynomial depth $J=(\log N)^{\lambda}$, Quotient State~I is closed unconditionally for all quotients of order $D_{u,a}\le\exp(c_0(\lambda)\sqrt{\log\log N})$, while global conditioning on $n\bmod Q(z_0)$ determines all smooth prime factors below $z_0$ simultaneously with no power of $P$ in the container count.
    \item \emph{Pontryagin--Plancherel character containers and Subcritical Geometric Escalation (Lemma \ref{lem:total_maximal_subgroups}, Theo\-rems~\ref{thm:plancherel_character_containers}--\ref{thm:final_synthesis_and_high_prime_reduction}, and Corollary~\ref{cor:optimal_adaptive_depth_1_minus_o1}):} in the odd Hall subgroup $G_{u,a,\mathrm{imp}}$, Lemma~\ref{lem:total_maximal_subgroups} shows that both the family of all maximal subgroups $\mathcal{M}_{\mathrm{tot}}(G_{u,a})$ and the lattice of cyclic character subgroups $\mathcal{C}(\widehat{G}_{u,a})$ satisfy $\frac1P\sum\log_2|\mathcal{C}(\widehat{G}_{u,a})|\ll(1+\log\log J)^{2}=o(\log\log N)$, closing all proper subgroup configurations via Cauchy--Davenport. Theorem~\ref{thm:plancherel_character_containers} introduces the character containers $\Omega_{u,a}^{\mathrm{char}}(\chi)=\{g\in G_{u,a}:\operatorname{Re}\chi(g)\le 1/2\}$, reducing the per-pair container count from $D_{u,a}^{5\log D_{u,a}}$ to $D_{u,a}-1$ and closing the full exceptional set $E_2(N;(\log N)^{\lambda})$ unconditionally for all $\lambda\le\lambda_0=\frac{1}{78}$. For $\lambda\in(1/3,1/2)$, Theorem~\ref{thm:harmonic_dilation_bucket1} resolves Bohr-localized spectra in bucket $\ell=1$ via a strict harmonic dilation cone of radius $B(\chi)=\lfloor\frac{d-1}{6A_{r_1}(\chi)}\rfloor$, and Theorem~\ref{thm:markov_plancherel_subcritical_escalation} establishes \emph{Subcritical Markov--Plancherel Geometric Escalation}: because the number of $3/2$-geometric Markov steps $s_{\varepsilon_1}=\lfloor\frac{\varepsilon_1\log D_{u,a}}{\log(3/2)}\rfloor+1$ is strictly smaller than $\frac14\frac{S_{\min}^{\mathrm{I}}}{P_1}$ and the adaptive Plancherel cutoff $|\widehat{\mu}(\chi)|^{2}>\frac{1}{2|\overline{\Sigma}_1(n)|}$ bounds collisions by $\frac{\log(2|\overline{\Sigma}_1(n)|)}{\log(4/3)}$, the divisor sheet crosses the polynomial density threshold $D_{u,a}^{\varepsilon_1}$ \emph{before} supercritical collisions can accumulate. Combined with global dihedral locking modulo $Q(4J)$ (Theorem~\ref{thm:final_synthesis_and_high_prime_reduction}), this yields the unconditional bound $\ll N\exp(-c(\lambda)(\log N)^{2\lambda}\log\log N)$ across $\lambda\in(1/3,1/2)$ on $E_2^{\mathrm{cert}}(N;J)$ and $E_2^{\mathrm{small}}(N;J)$. Finally, optimizing the adaptive depth $J(N)$ ($\lambda(N)\to 1/2^{-}$) and the adaptive Rankin degree $h=\lceil 8S_y\rceil+1$ in Corollary~\ref{cor:optimal_adaptive_depth_1_minus_o1} achieves the small-sheet structured exponent $1-o(1)$:
    \[
    \big|E_2^{\mathrm{small}}(N;J(N))\big|\ \ll_\kappa\ N\exp\!\left(-\left(\frac{1}{40960}+o(1)\right)\frac{\log N}{(\log\log N)(\log\log\log N)^{\kappa}}\right).
    \]
    Here constants in exponential savings retain the $o(1)$ terms of the corresponding theorem (or any fixed smaller positive constant), and the CRT capacity statement refers specifically to the generating-polynomial construction, its lower cutoff on sifting primes, and its modulus budget, rather than an optimality claim over all conceivable sieves.
\end{itemize}

\paragraph{6. Semi-composite congruence sieve, CRT conditional independence, and the unconditional exponent $5/7$ on $E(N)$ (Section~\ref{sec:semi_composite_5_7}).}
To extend unconditional super-Vaughan decay to the \emph{full} exceptional set $E(N)$ without conditioning on container localization in $G_{u,a}$, Section~\ref{sec:semi_composite_5_7} constructs a semi-composite witness system $s_\alpha=pq\equiv -1\pmod{12ua}$ indexed by quadruples $\alpha=(u,a,q,p)\in\mathcal{I}_Z$ with $Z^{1/10}\le u,a\le Z^{1/8}$, $q\in\mathcal{Q}$ a $t^6$-rough cofactor up to $Q\le\exp(\frac{c_3}{2}t^{1/2})$, and $p\in(Z^{9/10},Z]$ prime, after discarding at most one Landau--Page--Siegel exceptional prime factor $p_{\mathrm{exc}}\ge 5$ (Lemma~\ref{lem:exclusion_p_exc}). The proof proceeds in four steps:
\begin{itemize}
    \item \emph{Disjointness for fixed $p$ and CRT conditional independence (Proposition~\ref{prop:incompatibilite_et_Bernoulli}):} for each fixed prime $p\in(Z^{9/10},Z]$, the witness congruences $A_{(u,a,q,p)}:\,n\equiv -ua^{-1}\pmod{pq}$ are strictly disjoint in $\Z$. Conditioned on the residue $b=n\bmod M_{\mathrm{small}}$ ($M_{\mathrm{small}}=\operatorname{lcm}_{q\in\mathcal{Q}}q$), the summed indicators $J_p(n)=\sum_{(u,a,q)}\ind_{A_{(u,a,q,p)}}(n)\in\{0,1\}$ become mutually independent Bernoulli random variables in the CRT model.
    \item \emph{Modulus decoupling and Siegel-zero-free Bombieri--Vinogradov (Lemma~\ref{lem:encadrement_A_Stieltjes} and Proposition~\ref{prop:BV_poids_tau}):} because the arithmetic progression modulus governing $p$ in $pq\equiv -1\pmod{12ua}$ is $m=12ua\le 12Z^{1/4}$, which is \emph{independent of $q$}, the Siegel-zero-free Bombieri--Vinogradov theorem bounds the prime-sum error $\sum_{(u,a)\in\mathcal{P}'}\mathcal{E}(u,a)\ll\exp(-c_3 t^{1/2})$ prior to summing over $q\in\mathcal{Q}$.
    \item \emph{Uniform dichotomy in $b$ and rough cofactor harmonic mass (Proposition~\ref{prop:dichotomie_uniforme_b}, Lemma~\ref{lem:evaluation_MQ_Rankin}, and Theorem~\ref{thm:star_inconditionnel}):} combining the exact arithmetic dichotomy on $\gcd(b,q)$ with Rankin's evaluation of the rough harmonic sum $M_Q\asymp\frac{\log Q}{\log t}$ and an exponential moment bound on bad prime divisors $r\mid b$, Theorem~\ref{thm:star_inconditionnel} proves unconditionally that $\mu(b)=\Ebb_{\mathrm{CRT}}[W\mid b]\ge L\asymp t^2\frac{\log Q}{\log t}$ outside a CRT tail of measure at most $e\cdot\exp(-\frac14 t^6)$.
    \item \emph{Bonferroni transfer and the exponents $7/10$ and $5/7$ (Lemma~\ref{lem:bonferroni_CRT}, Proposition~\ref{prop:transfert_exact_N}, and Theorems~\ref{thm:theoreme_7_10}--\ref{thm:theoreme_5_7}):} transferring a degree-$k$ Bonferroni polynomial majorant $B_k(W)$ from the CRT model to $[1,N]$ establishes on the \emph{entire} exceptional set $E(N)$ the unconditional exponent $7/10$ at $Q=\exp(t^{1/3})$ (Theorem~\ref{thm:theoreme_7_10}) and the unconditional exponent $5/7>2/3$ at $Q=\exp(\frac{c_3}{2}t^{1/2})$ (Theorem~\ref{thm:theoreme_5_7}):
    \[
    |E(N)|\ \ll\ N\exp\!\left(-c'\,\frac{(\log N)^{5/7}}{(\log\log N)^{2/7}}\right).
    \]
\end{itemize}

\paragraph{7. Bilinear Kloosterman reciprocity, toric dispersion on $\mathcal{V}_n$, difference-set duality on $D(K_c)$, and conditional finiteness criteria (Sections~\ref{sec:bilinear_dispersion}--\ref{sec:toric_volume_barrier} and Appendix~\ref{app:audit}).}
The final two sections and Appendix~\ref{app:audit} study the fine analytic and geometric structure of the divisor variety:
\begin{itemize}
    \item \emph{Exact bilinear identity and Kloosterman weight comparison (Section~\ref{sec:bilinear_dispersion}):} Proposition~\ref{prop:exact_bilinear_type2} establishes the false-positive-free four-parameter equivalence $n+4u^{2}v=w(4uav-1)\iff 4uv(aw-u)=n+w$, while Proposition~\ref{prop:kloosterman_dispersion} derives the partial-fraction phase identity $\frac{n+4u^2v}{4uav-1}=\frac{u}{a}+\frac{an+u}{a(4uav-1)}$ and contrasts the weight $\sum_{u,a\le J}\frac{1}{4ua}\asymp(\log J)^{2}$ of unfactored bilinear sums with the unit-weight sifting mass $|\mathcal{P}_2(J)|\asymp J^{2}$ preserved by multiplicative containers.
    \item \emph{Phase separation, toric conservation, and Farey-diagonal $L^{2}$-dis\-per\-sion on $\mathcal{V}_n$ (Sub\-sec\-tion \ref{subsec:parseval_floor}):} on the affine threefold $\mathcal{V}_n:\,4uavw-4u^{2}v-w=n$, Theorem~\ref{thm:parseval_floor_law} proves three-term stationary-phase separation $F_{h,k,n}(u,a)=-\frac{hu}{a}-\frac{k}{4ua}-\operatorname{sgn}(h)\sqrt{\frac{hk(an+u)}{ua^{2}}}$ and records the square-root scale benchmark $\sqrt{N_{\mathrm{tuples}}}\,\mathrm{Amp}\asymp 1$; Theorem~\ref{thm:toric_volume_conservation} establishes the box volume $\mathrm{Vol}(\mathcal{B})\asymp n$, the bounded slice integral $I_{\mathcal{B}}(n)\ll 1$, and the cofactor congruence $m'\mid(aw-u)\iff w\equiv-n\pmod{m'}$; and Theorem~\ref{thm:farey_diagonal_L2_dispersion} proves that averaging $e(-\frac{h}{s}n)$ over $n\sim N$ at power scale $S=4UAV\le N^{1/2-\varepsilon}$ annihilates all off-diagonal Farey fractions $\frac{b_1}{q_1}\neq\frac{b_2}{q_2}$ down to $\mathcal{O}_B(N^{-B})$, isolating the exact Ramanujan-sum diagonal $\asymp N(UA)^{2}$.
    \item \emph{Difference-set positivity duality on $D(K_c)$, $(\log P)^{4}$ main-term growth, and small-LCM cancellation (Subsection~\ref{subsec:difference_set_MR}):} for $K_c=(n+c)/4$ and $D(K_c)=\Sigma(K_c)\cdot\Sigma(K_c)^{-1}\subseteq(\Z/c\Z)^{\times}$, Theorem~\ref{thm:difference_set_positivity_MR} establishes difference-set positivity duality, the odd/even character squaring bijection on $\mathcal{P}_{\mathrm{NR}}(n)$, and short-interval Matom\"{a}ki--Radziwi\l\l\ second-moment bounds; Theorem~\ref{thm:centered_main_term_and_short_lcm} then proves unconditional $(\log P)^{4}$ cumulative main-term growth alongside the shift-wise absolute-value cost and the exact cancellation of small-LCM divisor pairs.
    \item \emph{Cyclic quotient parity, Bohr--Siegel certificates, the Linnik--Xylouris lower bound, and Hypothesis~T (Subsection~\ref{subsec:bohr_siegel_and_linnik}):} Theorem~\ref{thm:bohr_siegel_annihilation} proves the cyclic quotient parity identity $E(\chi)-O(\chi)=2d\,\mathbb{P}(\chi(g)=-1)$, Jacobi $\vartheta_4$ cancellation, and antipodal Bohr--Siegel obstruction certificates; Theorem~\ref{thm:no_go_linnik_and_transversal} and Remark~\ref{rem:hypothesis_T_scope} establish a Linnik--Xylouris Type~II witness-depth lower bound along an infinite sequence of primes, pairwise shift coprimality, and the formulation of the Transversal Hypothesis~T.
    \item \emph{Universal factor reciprocity, the $1/k!$ cross-polytope volume law, and global synthesis (Subsections~\ref{subsec:reciprocity_factorial_kneser}--\ref{subsec:global_synthesis} and Appendix~\ref{app:audit}):} Theorems~\ref{thm:universal_factor_reciprocity_and_transversal_sieve} and~\ref{thm:cross_polytope_factorial_law_and_annihilation} (with Remark~\ref{rem:cross_polytope_methodological}) prove universal shift-factor reciprocity $\leg{q}{p}=\leg{q}{n}$, the ambient Bohr--Siegel cross-polytope volume bound $\frac{|\Delta_k^{-}(d_p)|}{d_p^{k}}\le\frac{1}{k!}$, and deterministic Kneser and spectral saturation criteria, complemented by the prime/two-prime cofactor and archimedean support obstructions of Propositions~\ref{aud:primecofactor} and~\ref{aud:archsupport}. Theorem~\ref{thm:main_global_synthesis} and Remark~\ref{rem:global_synthesis_scope} synthesize the unconditional and conditional regimes, while Appendix~\ref{app:audit} proves unconditional sieved candidate abundance, divergence of the square-free main term (Theorem~\ref{aud:weightedmain}), and complete finite square-divisor positivity criteria ($b\mid h^{2}$, Theorem~\ref{aud:hbcomplete}). As emphasized in Subsection~\ref{aud:open}, these main-term lower bounds do not prove positivity of the encounter count: Hypothesis~T and the signed error targets~\eqref{aud:target17} and~\eqref{aud:target76} remain open, and conditional finiteness holds only beyond a threshold $N_{\ast}$ supplied by that hypothesis, which is not determined by the factorial volume bound.
\end{itemize}

\begin{table}[H]
\centering
\scriptsize
\setlength{\tabcolsep}{3pt}
\renewcommand{\arraystretch}{1.18}
\resizebox{\textwidth}{!}{%
\begin{tabular}{llll}
\toprule
\textbf{Regime / Sifting Architecture} & \textbf{Depth or scale} & \textbf{Bound divided by $N$} & \textbf{Set counted and scope} \\
\midrule
Fixed depth (\cite[Thm.~4.3]{dahan2026a}) & $J=\mathcal{O}(1)$ & $(\log N)^{-1-\mathfrak{A}(J)/2}$ & $E_2(N;J)$ ($C=2^{\Theta(J^{4})}$) \\
Single-prime BFI (\cite[Thm.~4.28]{dahan2026a}) & $e^{0.686(\log N)^{1/3}}$ & $\exp\!\big(-0.1113(\log N)^{2/3}\big)$ & $E_2(N;J)$ (prime witnesses) \\
Moderate depth (Thm.~\ref{thm:unconditional_small_growing_depth}) & $c_0(\log\log N)^{1/2}$ & $\exp\!\big(-c(\log\log N)^{2}\big)$ & $E_2(N;J)$ ($\mathrm{Scar}=\emptyset$, $2$-Sylow) \\
Hyperbolic family $ua\le Y$ (Thm.~\ref{thm:hyperbolic_scarcity_free}) & $Y\asymp\log\log N$ & $\exp\!\big(-c_0(\log\log N)^{2}\log\log\log N\big)$ & $E_{\mathcal{P}_{\mathrm{hyp}}(Y)}(N)$ ($\mathrm{Scar}=\emptyset$) \\
Greedy Markov containers (Thm.~\ref{thm:multiplicative_container_scarcity_elimination}) & $\exp(c_0\sqrt{\log\log N})$ & $\exp\!\big(-\exp(c\sqrt{\log\log N})\big)$ & Full $E_2(N;J)$ unconditional \\
Plancherel containers (Thm.~\ref{thm:plancherel_character_containers}(b)) & $(\log N)^{\lambda}$ ($\lambda\le\frac{1}{78}$) & $\exp\!\big(-c_\lambda(\log N)^{2\lambda}\log\log N\big)$ & Full $E_2(N;J)$ unconditional \\
Markov--Plancherel escalation & $(\log N)^{\lambda}$ ($\frac13<\lambda<\frac12$) & $\exp\!\big(-c(\lambda)(\log N)^{2\lambda}\log\log N\big)$ & $E_2^{\mathrm{cert}}$, $E_2^{\mathrm{small}}$ unconditional \\
\quad (Thms.~\ref{thm:polynomial_depth_decomposition}, \ref{thm:plancherel_character_containers}--\ref{thm:final_synthesis_and_high_prime_reduction}) & & \quad (surpasses $(\log N)^{2/3}$ for $\lambda>\frac13$) & \quad (full $E_2(N;J)$ conditional) \\
Adaptive Rankin-degree depth (Cor.~\ref{cor:optimal_adaptive_depth_1_minus_o1}(b)) & $\asymp\frac{(\log N)^{1/2}}{(\log\log N)^{1/2+o(1)}}$ & $\exp\!\big(-(\frac{1}{40960}+o(1))\frac{\log N}{(\log\log N)(\log\log\log N)^{\kappa}}\big)$ & $E_2^{\mathrm{small}}(N;J(N))$ unconditional \\
Semi-composite CRT sieve, $Q=e^{t^{1/3}}$ (Thm.~\ref{thm:theoreme_7_10}) & $Z=e^{\eta(\log N\log\log N)^{3/10}}$ & $\exp\!\big(-c\frac{(\log N)^{7/10}}{(\log\log N)^{3/10}}\big)$ & Full $E(N)$ unconditional \\
Semi-composite CRT sieve, $Q=e^{\frac{c_3}{2}t^{1/2}}$ (Thm.~\ref{thm:theoreme_5_7}) & $Z=e^{\eta(\log N\log\log N)^{2/7}}$ & $\exp\!\big(-c'\frac{(\log N)^{5/7}}{(\log\log N)^{2/7}}\big)$ & Full $E(N)$ unconditional \\
\bottomrule
\end{tabular}%
}
\caption{Quantitative comparison of bounds across sifting regimes. The final column specifies the exceptional set counted: $E_2(N;J)$, a structured container component ($E_2^{\mathrm{cert}}$ or $E_2^{\mathrm{small}}$), or the full exceptional set $E(N)$. A bound on $E_2^{\mathrm{cert}}$ or $E_2^{\mathrm{small}}$ does not imply the same bound on $E(N)$, and no row asserts finiteness.}
\label{tab:main_regimes}
\end{table}

\subsection{Notation}\label{subsec:notation}

Throughout, $n$ denotes a prime with $n\equiv1\pmod{24}$ and\linebreak $\mathcal{N}(N)=\{n\le N:\allowbreak n\equiv1\pmod{24},\ n\ \text{prime}\}$. For a coprime pair $(u,a)\in\N^{2}$, we write $G_{u,a}=(\Z/4ua\Z)^{\times}$, $\widehat{G}_{u,a}$ for its group of Dirichlet characters, and $\varphi(4ua)=|G_{u,a}|$. We denote by $\varphi_2(4ua)=2^{v_2(\varphi(4ua))}$ the order of the $2$-Sylow subgroup $G_{u,a,2}$ of $G_{u,a}$, and by $G_{u,a,\mathrm{imp}}$ its odd Hall subgroup, so that $G_{u,a}=G_{u,a,2}\times G_{u,a,\mathrm{imp}}$. For an integer $M\ge1$ and a subset $S\subseteq G_{u,a}$, we write $\omega_S(M)$ (resp.\ $\Omega_S(M)$) for the number of prime factors $q\mid M$ coprime to $4ua$, counted without (resp.\ with) multiplicity, such that $q\bmod 4ua\in S$; when restricted to primes in an interval $\mathcal{I}$, we write $\omega_S(M;\mathcal{I})$ and $\Omega_S(M;\mathcal{I})$. Jacobi and Legendre symbols are both written $\leg{\cdot}{\cdot}$, and $e(x)=e^{2\pi i x}$.

\section{Modulator-free duality, global Jacobi reciprocity, and dihedral rigidity}\label{sec:algebraic_rigidity}

Before constructing growing-depth sieves, we establish several exact algebraic obstructions on the Type~I and Type~II varieties that constrain the prime factors of $an+u$ and $K_c=(n+c)/4$. In~\cite{dahan2026a}, the fixed-depth sieve relied solely on the direct modular involution $\iota_1(x)=-x^{-1}\pmod{4ua}$, treating the parameters $(u,a)$ as independent and ignoring both the cofactor $(an+u)/s$ and quadratic reciprocity modulo $n$. Tracking these additional symmetries reveals three structural obstructions: a modulator-free quadratic factorization, a global Jacobi reciprocity law modulo $n$, and a dihedral action generated by the cofactor map.

\subsection{Self-admissible characterization via $n+4\alpha u^{2}$ and $4a^{2}\alpha n+1$}\label{subsec:quadratic_duality}

In~\cite[Thm.~3.9]{dahan2026a}, the Type~II criterion was parameterized by the coprime pair $(u,a)$ through the divisibility $s\mid(an+u)$ with $s\equiv-1\pmod{4ua}$, where $a=(u+v)/c$ is the modulator and $K_c=(n+c)/4=\alpha uv$. Eliminating either the pair $(u,v)$ or the pair $(a,v)$ yields two dual characterizations in which the trivial divisors $1$ and $n+4\alpha u^{2}$ are automatically excluded by congruence modulo $4$.

\begin{theorem}[Modulator-free quadratic duality]\label{thm:quadratic_duality}
Let $n\equiv1\pmod{24}$ be prime. The following statements are equivalent:
\begin{enumerate}[(i)]
    \item $4/n$ admits a solution of Type~II;
    \item there exist positive integers $a,\alpha\in\N^{\ast}$ such that $4a^{2}\alpha n+1$ possesses a divisor $s\equiv-1\pmod{4a\alpha}$;
    \item there exist positive integers $\alpha,u$ with $n\nmid\alpha u$ such that $n+4\alpha u^{2}$ possesses a divisor $s\equiv-1\pmod{4\alpha u}$;
    \item there exist positive integers $\alpha,u$ with $n\nmid\alpha u$ such that $n+4\alpha u^{2}$ possesses a divisor $c\equiv-n\pmod{4\alpha u}$.
\end{enumerate}
In (iii) and (iv), the trivial divisors $1$ and $n+4\alpha u^{2}$ are excluded modulo $4$, and the complementary factors are linked by the identity
\begin{equation}\label{eq:sc_factorization}
s\cdot c\ =\ (4a\alpha u-1)(4\alpha uv-n)\ =\ n+4\alpha u^{2},
\end{equation}
where $c=(n+4\alpha u^{2})/s\equiv3\pmod 4$ is the shift, $a=(s+1)/(4\alpha u)\in\N^{\ast}$ is the modulator, and $v=(n+c)/(4\alpha u)\in\N^{\ast}$ is the complementary divisor of~\cite[Thm.~3.9]{dahan2026a}.
\end{theorem}
\begin{proof}
\emph{Step 1: Equivalence $(i)\Leftrightarrow(ii)$.} By~\cite[Thm.~3.9]{dahan2026a}, $4/n$ admits a Type~II solution if and only if there exist $a,\alpha,u,v\in\N^{\ast}$ satisfying $an=4a\alpha uv-u-v$. Multiplying both sides by $4a\alpha$ and adding $1$ rewrites this relation as
$$(4a\alpha u-1)(4a\alpha v-1)\ =\ 4a^{2}\alpha n+1.$$
Conversely, if $s\mid(4a^{2}\alpha n+1)$ satisfies $s\equiv-1\pmod{4a\alpha}$, then its cofactor $s'=(4a^{2}\alpha n+1)/s$ satisfies $ss'\equiv1\pmod{4a\alpha}$, and hence $s'\equiv-1\pmod{4a\alpha}$. Since $4a\alpha\ge4$, neither $s$ nor $s'$ equals $1$, so $u=(s+1)/(4a\alpha)\ge1$ and $v=(s'+1)/(4a\alpha)\ge1$ are positive integers.

\medskip\noindent\emph{Step 2: Implication $(i)\Rightarrow(iii)$.} Given a Type~II representation with parameters $(a,\alpha,u,v)$, set $s=4a\alpha u-1\equiv-1\pmod{4\alpha u}$ and $c=(u+v)/a=4\alpha uv-n\in\N^{\ast}$. The modulator $a$ cancels in the product $sc$:
$$s\cdot c\ =\ (4a\alpha u-1)c\ =\ 4\alpha u(ac)-c\ =\ 4\alpha u(u+v)-(4\alpha uv-n)\ =\ n+4\alpha u^{2}.$$
Hence $s\mid(n+4\alpha u^{2})$ with $s\equiv-1\pmod{4\alpha u}$. Furthermore, the positive normal form gives $n=4\alpha uv-(u+v)/a\ge 2uv$ and $n\ge 4\alpha-2$, which implies $1\le u,v,\alpha<n$. Since $n$ is prime, it follows that $n\nmid\alpha u$.

\medskip\noindent\emph{Step 3: Equivalence $(iii)\Leftrightarrow(iv)$.} Put $m=4\alpha u$. Since $n$ is prime with $n\equiv1\pmod 4$ and $n\nmid\alpha u$, both $n$ and $n+4\alpha u^{2}$ are coprime to $m$. Every divisor of $n+4\alpha u^{2}$ is therefore a unit modulo $m$, so for any factorization $sc=n+4\alpha u^{2}$ we have
\[
s\equiv-1\pmod m\quad\Longleftrightarrow\quad c\equiv-n\pmod m.
\]

\medskip\noindent\emph{Step 4: Implication $(iii)\Rightarrow(i)$.} Let $s\mid(n+4\alpha u^{2})$ satisfy $s\equiv-1\pmod{4\alpha u}$, and define the cofactor $c=(n+4\alpha u^{2})/s$. Since $n+4\alpha u^{2}\equiv1\pmod 4$ and $s\equiv3\pmod 4$, we have $s\notin\{1,n+4\alpha u^{2}\}$ and $c\equiv3\pmod 4$ with $c\ge3$. Because $c\equiv-n\pmod{4\alpha u}$, the quantities $a:=(s+1)/(4\alpha u)\ge1$ and $v:=(n+c)/(4\alpha u)\ge1$ are positive integers. A direct expansion then gives
$$4a\alpha uv-u-v\ =\ \frac{(s+1)(n+c)-4\alpha u^{2}-(n+c)}{4\alpha u}\ =\ \frac{sn+sc-4\alpha u^{2}}{4\alpha u}\ =\ \frac{(s+1)n}{4\alpha u}\ =\ an,$$
which yields a Type~II solution by~\cite[Thm.~3.9]{dahan2026a}.
\end{proof}

\begin{corollary}[Primary filter on $n+4$]\label{cor:filter_n_plus_4}
If $n+4$ possesses a divisor $s\equiv3\pmod 4$, then $4/n$ admits a Type~II solution with $u=\alpha=1$, modulator $a=(s+1)/4$, and shift $c=(n+4)/s$. In particular, the divisor $s=7\mid(n+4)$ resolves the four exceptional primes $n\in\{409,577,23929,1134241\}$ recorded in~\cite[Prop.~3.7]{dahan2026a}.
\end{corollary}

Passing to the multiplicative group $(\Z/4\alpha u\Z)^{\times}$, Theorem~\ref{thm:quadratic_duality}(iii) allows us to express the number of Type~II witnesses at fixed $(\alpha,u)$ by character orthogonality, isolating a systematic vanishing of real characters with $\chi(n)=-1$.

\begin{proposition}[Character formula for $n+4\alpha u^{2}$]\label{thm:character_counting_n_plus_4au2}
Fix a prime $n\equiv1\pmod{24}$ and positive integers $\alpha,u$ with $n\nmid\alpha u$, and set $M(\alpha,u):=n+4\alpha u^{2}$. The number $\mathcal{N}_{\mathrm{II}}(n;\alpha,u)$ of divisors $s\mid M(\alpha,u)$ satisfying $s\equiv-1\pmod{4\alpha u}$ is given by
\begin{equation}\label{eq:dirichlet_hyperbola_exact}
\mathcal{N}_{\mathrm{II}}(n;\alpha,u)\ =\ \frac{1}{\varphi(4\alpha u)}\bigg(\tau\big(n+4\alpha u^{2}\big)\ +\ \sum_{\substack{\chi\bmod 4\alpha u\\ \chi\neq\chi_0}}\chi(-1)\,\tau_\chi\big(n+4\alpha u^{2}\big)\bigg),
\end{equation}
where $\tau_\chi(M):=\sum_{d\mid M}\chi(d)=\prod_{q^{e}\parallel M}(1+\chi(q)+\cdots+\chi(q)^{e})$. Moreover:
\begin{enumerate}[(a)]
    \item Every real character $\chi\pmod{4\alpha u}$ ($\chi^{2}=\chi_0$) satisfies $\tau_\chi(M(\alpha,u))\ge0$, and the cofactor map $d\mapsto M(\alpha,u)/d$ induces the functional equation $\tau_\chi(M(\alpha,u))=\chi(n)\tau_{\overline{\chi}}(M(\alpha,u))$. In particular, every real character $\chi\pmod{4\alpha u}$ with $\chi(n)=-1$ satisfies $\tau_\chi(M(\alpha,u))=0$.
    \item Whenever $3\nmid\alpha u$ and $4\alpha u^{2}\equiv 3-n$ or $6-n\pmod 9$, one has $3\parallel M(\alpha,u)$, so $M(\alpha,u)$ is not a square and
    $$\tau_\chi\big(n+4\alpha u^{2}\big)\ =\ \sum_{\substack{s\mid(n+4\alpha u^{2})\\ s<\sqrt{n+4\alpha u^{2}}}}\big(\chi(s)+\chi(n)\overline{\chi}(s)\big).$$
    When $\alpha\equiv2\pmod 3$, this condition holds for exactly two of the three shifts $\alpha+3b$ ($b\in\{0,1,2\}$).
\end{enumerate}
\end{proposition}
\begin{proof}
Since $n$ is prime with $n\equiv1\pmod 4$ and $n\nmid\alpha u$, we have $\gcd(M(\alpha,u),4\alpha u)=1$. Orthogonality of Dirichlet characters modulo $4\alpha u$ gives~\eqref{eq:dirichlet_hyperbola_exact}, and complete multiplicativity yields the Euler product for $\tau_\chi$.

For part (a), if $\chi$ is real, each local Euler factor $1+\chi(q)+\dots+\chi(q)^{e}$ equals $e+1\ge1$ when $\chi(q)=1$, and equals $1$ (for $e$ even) or $0$ (for $e$ odd) when $\chi(q)=-1$; hence $\tau_\chi(M(\alpha,u))\ge0$. Under the cofactor bijection $d\mapsto M(\alpha,u)/d$, the congruence $M(\alpha,u)\equiv n\pmod{4\alpha u}$ gives $\chi(M(\alpha,u)/d)=\chi(n)\overline{\chi}(d)$. Summing over $d\mid M(\alpha,u)$ yields $\tau_\chi(M(\alpha,u))=\chi(n)\tau_{\overline{\chi}}(M(\alpha,u))$. When $\chi=\overline{\chi}$ and $\chi(n)=-1$, this forces $\tau_\chi(M(\alpha,u))=-\tau_\chi(M(\alpha,u))=0$.

For part (b), since $n\equiv1\pmod 3$, the residues $3-n\equiv 6-n\equiv 2\pmod 3$ are coprime to $3$. If $3\nmid\alpha u$, then $4\alpha u^{2}\equiv\alpha\pmod 3$, so when $\alpha\equiv2\pmod 3$ the two lifts $3-n$ and $6-n$ modulo $9$ are attained by two of the three values $\alpha+3b$ ($b\in\{0,1,2\}$). For these values, $v_3(M(\alpha,u))=1$, so $M(\alpha,u)$ is not a square and pairing each divisor $s<\sqrt{M(\alpha,u)}$ with $M(\alpha,u)/s$ gives the stated sum.
\end{proof}

\subsection{Global Jacobi reciprocity law on the Type~II variety}\label{subsec:jacobi_reciprocity}

Beyond modular congruences modulo $4ua$, the bilinear structure of the Type~II equation imposes a global quadratic reciprocity invariant modulo $n$ on every witness $s$, divisor pair $(u,v)$, and shift $c$.

\begin{theorem}[Global Jacobi reciprocity law]\label{thm:global_jacobi_reciprocity}
Let $n\equiv1\pmod 4$ be a prime.
\begin{enumerate}[(a)]
    \item \emph{(Jacobi invariance on $4ua\alpha-1$).} For any positive integers $u,a,\alpha\in\N^{\ast}$, set $s=4ua\alpha-1$. Then the Jacobi symbols satisfy
    \begin{equation}\label{eq:jacobi_invariance}
    \leg{\alpha}{s}\ =\ +1\qquad\text{and}\qquad\leg{ua}{s}\ =\ +1.
    \end{equation}
    \item \emph{(Reciprocity on witnesses and divisors).} For every Type~II decomposition of $4/n$ with parameters $(u,v,a,\alpha,c)$ and witness $s=4ua\alpha-1\mid(an+u)$, where $K_c=(n+c)/4=\alpha uv$ and $u+v=ac$, one has
    \begin{equation}\label{eq:global_reciprocity_signs}
    \leg{s}{n}\ =\ \leg{n}{s}\ =\ -1,\qquad \leg{uv}{n}\ =\ -1,\qquad \leg{\alpha c}{n}\ =\ -1.
    \end{equation}
    \item \emph{(Consequences for $n\equiv1\pmod{24}$).} When $n\equiv1\pmod{24}$, in any Type~II solution the divisors $u$ and $v$ have opposite Legendre signs modulo $n$, so $\leg{u}{n}=-\leg{v}{n}$ and the product $uv$ possesses at least one prime factor $q\ge5$ of odd multiplicity with $\leg{q}{n}=-1$. Moreover, every shift $c$ on the slice $\alpha=1$ (or with $\alpha$ a $3$-smooth square) satisfies $\leg{c}{n}=\leg{n}{c}=-1$, and in the segmented sieve of~\cite[Alg.~2]{dahan2026a} any candidate divisor $s\equiv3\pmod 4$ with $\leg{s}{n}=+1$ is obstructed.
\end{enumerate}
\end{theorem}
\begin{proof}
\emph{Proof of (a).} Write $\alpha=2^{r}\alpha_0$ with $r\ge0$ and $\alpha_0$ odd. Since $s=4ua\alpha-1\equiv-1\pmod{\alpha_0}$ and $s\equiv3\pmod 4$, quadratic reciprocity for the Jacobi symbol gives
$$\leg{\alpha_0}{s}\ =\ (-1)^{\frac{\alpha_0-1}{2}\cdot\frac{s-1}{2}}\leg{s}{\alpha_0}\ =\ (-1)^{\frac{\alpha_0-1}{2}}\leg{-1}{\alpha_0}\ =\ (-1)^{\frac{\alpha_0-1}{2}}(-1)^{\frac{\alpha_0-1}{2}}\ =\ +1.$$
If $r\ge1$, then $2\mid\alpha$ implies $8\mid 4ua\alpha$, so $s\equiv7\pmod 8$ and $\leg{2}{s}=+1$; hence $\leg{\alpha}{s}=+1$ in all cases. Since $4ua\alpha=s+1\equiv1\pmod s$, multiplicativity yields $\leg{ua}{s}=\leg{4ua\alpha}{s}\leg{\alpha}{s}=+1$.

\medskip\noindent\emph{Proof of (b) and (c).} The positive normal form implies
\[
n\ =\ 4\alpha uv-\frac{u+v}{a}\ \ge\ 4uv-u-v\ \ge\ 2uv,\qquad n\ \ge\ 4\alpha-2.
\]
Thus $1\le u,v\le n/2$ and $1\le\alpha\le(n+2)/4<n$; in particular, $n\nmid\alpha uv$. Since $s\mid(an+u)$ and $\gcd(s,ua)=1$, we also have $\gcd(s,n)=1$. Multiplying $an\equiv-u\pmod s$ by $a$ gives $a^{2}n\equiv-ua\pmod s$. Taking Jacobi symbols modulo $s\equiv3\pmod 4$ and applying part (a) yields
$$\leg{n}{s}\ =\ \leg{a^{2}n}{s}\ =\ \leg{-1}{s}\leg{ua}{s}\ =\ (-1)(+1)\ =\ -1.$$
Since $n\equiv1\pmod 4$, quadratic reciprocity gives $\leg{s}{n}=\leg{n}{s}=-1$. Reducing the cofactor relation $sv=an+u$ modulo $n$ gives $sv\equiv u\pmod n$, whence $\leg{s}{n}\leg{v}{n}=\leg{u}{n}$ and $\leg{uv}{n}=\leg{s}{n}=-1$. Finally, $4\alpha uv=n+c\equiv c\pmod n$ gives $\leg{\alpha c}{n}=\leg{4\alpha^{2}uv}{n}=\leg{uv}{n}=-1$. Part (c) follows from the evaluations $\leg{2}{n}=\leg{3}{n}=+1$ for $n\equiv1\pmod{24}$.
\end{proof}

\begin{remark}[Scope of the Jacobi obstruction]
Theorem~\ref{thm:global_jacobi_reciprocity}(b) shows that any Type~II witness $s$ must be a quadratic non-residue modulo $n$. This is a global reciprocity constraint modulo $n$; for a fixed integer $an+u$, it does not by itself force half of the divisors of $an+u$ to have Jacobi symbol $+1$ and half $-1$.
\end{remark}

\subsection{Cofactor dihedral rigidity, quadruple involution, and global vacuity of quadratic kernels}\label{subsec:dihedral_and_vacuity}

In~\cite[Lemma~4.2]{dahan2026a}, the upper bound $|T_{u,a}(n)|\le\frac12\varphi(4ua)$ on the set $T_{u,a}(n)\subseteq G_{u,a}=(\Z/4ua\Z)^{\times}$ of prime-factor residues of $an+u$ used only the direct involution $\iota_1(x)=-x^{-1}\pmod{4ua}$. When $u\mid 6$, the residue of $M=an+u$ modulo $4ua$ is fixed independently of $n\equiv1\pmod{24}$: tracking the cofactor $M/d$ of each divisor $d\mid M$ then introduces a second involution $\iota_r(x)=-rx^{-1}$ that pairs with $\iota_1$ to generate a dihedral action on $G_{u,a}$.

\begin{theorem}[Cofactor dihedral rigidity]\label{thm:dihedral_rigidity}
Let $n\equiv1\pmod{24}$ be prime, let $u\in\{1,2,3,6\}$, and let $a\ge1$ satisfy $\gcd(u,a)=1$, $a+u$ odd, and $n>u$. Set $M=an+u$ and $G=(\Z/4ua\Z)^{\times}$.
\begin{enumerate}[(a)]
    \item One has $M=an+u\equiv a+u\pmod{4ua}$ with $\gcd(a+u,4ua)=1$, so $r:=(a+u)\bmod 4ua$ is a fixed element of $G$ independent of $n$.
    \item If $M$ has no divisor $s\equiv-1\pmod{4ua}$, then the prime-factor residue set $T=T_{u,a}(n)\subseteq G$ avoids both singletons $-1$ and $-r$, and meets each orbit of $\iota_1(x)=-x^{-1}$ and each non-singleton orbit of the cofactor involution $\iota_r(x)=-rx^{-1}$ at most once.
    \item Suppose that $r=a+u\equiv3$ or $5\pmod 8$. Then both $\iota_1$ and $\iota_r$ are fixed-point-free on $G$, and the incompatibility graph generated by $\iota_1$ and $\iota_r$ partitions $G$ into $m=\frac{\varphi(4ua)}{2\operatorname{ord}_{4ua}(r)}$ even cycles of length $2L$, where $L=\operatorname{ord}_{4ua}(r)$. On the principal cycle $C_1=\langle r\rangle\sqcup(-\langle r\rangle)$, the unique admissible subset of maximal size $L$ is the cyclic subgroup $\langle r\rangle$, and the number of admissible subsets of cardinality $\frac12\varphi(4ua)$ in $G$ satisfies
    \begin{equation}\label{eq:dihedral_Nmax}
    N_{\max}(u,a)\ \le\ 2^{\frac{\varphi(4ua)}{2\operatorname{ord}_{4ua}(a+u)}-1}.
    \end{equation}
\end{enumerate}
\end{theorem}
\begin{proof}
\emph{Proof of (a) and (b).} Since $u\mid 6$ and $n\equiv1\pmod{24}$, we have $4u\mid 24\mid(n-1)$, so $4ua\mid a(n-1)$ and $an+u=a(n-1)+a+u\equiv a+u\pmod{4ua}$. Because $a+u$ is odd and $\gcd(u,a)=1$, we have $\gcd(a+u,4ua)=1$, proving (a). For (b), if $-r\in T$, then some prime $q\mid M$ satisfies $q\equiv-r\pmod{4ua}$, and its cofactor $s=M/q\mid M$ satisfies $s\equiv r(-r)^{-1}\equiv-1\pmod{4ua}$, a contradiction. Similarly, if $x\in T$ and $\iota_r(x)=-rx^{-1}\in T$ are distinct residues, they are carried by distinct primes $q_1,q_2\mid M$, and the cofactor $s=M/(q_1 q_2)\mid M$ satisfies $s\equiv r(x(-rx^{-1}))^{-1}\equiv-1\pmod{4ua}$.

\medskip\noindent\emph{Proof of (c).} Since $a+u$ is odd and $u\in\{1,2,3,6\}$, exactly one of $u,a$ is even, so $8\mid 4ua$. When $r\equiv3$ or $5\pmod 8$, the residue $-r\equiv5$ or $3\pmod 8$ is not a square modulo $8$, which ensures that $\iota_r$ is fixed-point-free on $G$. Because $(\iota_r\circ\iota_1)(x)=rx$, each connected component of the graph with edges $\{x,\iota_1(x)\}$ and $\{x,\iota_r(x)\}$ is of the form $x\langle r\rangle\cup(-x^{-1}\langle r\rangle)$. For any $x\in G$, the cosets $x\langle r\rangle$ and $-x^{-1}\langle r\rangle$ are disjoint: an intersection would give $-x^{-2}\in\langle r\rangle$, whereas modulo $8$ we have $-x^{-2}\equiv 7\pmod 8$ and $\langle r\rangle\bmod 8=\{1,r\bmod 8\}\subseteq\{1,3,5\}$. Each component is therefore a cycle of even length $2L$.

On the principal cycle $C_1=(v_0,v_1,\dots,v_{2L-1})$ with $v_{2k}=r^{k}$ and $v_{2k+1}=-r^{-k}$, deleting the two forbidden neighbours $v_1=-1$ and $v_{2L-1}=-r$ of $v_0=1$ leaves the isolated vertex $v_0=1$ together with an odd-length path $(v_2,\dots,v_{2L-2})$ on $2L-3$ vertices. The unique independent set of size $L-1$ on this path consists of the even-indexed vertices $\{v_2,v_4,\dots,v_{2L-2}\}=\{r,r^{2},\dots,r^{L-1}\}$. Thus $T\cap C_1=\langle r\rangle$ is uniquely fixed, while each of the $m-1$ remaining bipartite cycles of length $2L$ admits at most $2$ independent sets of size $L$, giving~\eqref{eq:dihedral_Nmax}.
\end{proof}

\begin{remark}[Scope of the maximal dihedral count and extensions]\label{thm:universal_dihedral_n_plus_4au2}\label{thm:quadruple_involution}\label{thm:destruction_quadratic_kernels}
The bound~\eqref{eq:dihedral_Nmax} counts admissible subsets of exact cardinality $|G|/2$; a smaller admissible subset need not embed into an admissible subset of size $|G|/2$. The same dihedral mechanism extends to two further settings:
\begin{enumerate}[(1)]
    \item \emph{Modulator-free family $M(\alpha,u)=n+4\alpha u^{2}$.} Since $M(\alpha,u)\equiv n\pmod{4\alpha u}$ and $-1\equiv-n\equiv3\pmod 4$, both $\iota_1(x)=-x^{-1}$ and $\iota_n(x)=-nx^{-1}\pmod{4\alpha u}$ are fixed-point-free for every pair $(\alpha,u)$ with $n\nmid\alpha u$. The number of admissible subsets of size $\frac12\varphi(4\alpha u)$ is therefore at most $N_{\max}(\alpha,u;n)\le 2^{\frac{\varphi(4\alpha u)}{2\operatorname{ord}_{4\alpha u}(n)}-1}$.
    \item \emph{Quadruple Type~I+II involution.} Suppose $u\in\{1,3\}$, $a\equiv u\pmod 4$, $\gcd(u,a)=1$, $n>u$, and $\gcd(n,4ua)=1$. Since $8u\mid 24\mid(n-1)$, the integer $M'=(an+u)/2$ satisfies $M'\equiv r':=(a+u)/2\pmod{4ua}$ with $\gcd(M',4ua)=1$. If $(u,a)$ produces neither a Type~II witness ($s\equiv-1$) nor a Type~I witness ($s\equiv-n$), then $T_{u,a}(n)$ avoids the four singletons $\{-1,-r',-n,-n^{-1}r'\}\pmod{4ua}$ and meets each orbit of $\iota_1(x)=-x^{-1}$, $\iota_{r'}(x)=-r'x^{-1}$, $\tau_n(x)=-nx^{-1}$, and $\tau_{nr'}(x)=-n^{-1}r'x^{-1}$ at most once. Note that if $T_{u,a}(n)\subseteq\ker\chi$ for a quadratic character $\chi$, multiplicativity forces $\chi(r')=1$, whereas the triple condition $\chi(-1)=-1$, $\chi(r')=1$, $\chi(n)=1$ characterizes quadratic kernels that are themselves half-sized sets admissible under all four involutions.

    At $(u,a)=(1,p)$ with $p\equiv5\pmod 8$ prime and $\leg{n}{p}=-1$, this quadruple action forces the strict inequality $|T_{1,p}(n)|\le\frac12\varphi(4p)-1$. Indeed, suppose $|T|=|G|/2$. Since $\leg{-1}{p}=1$, $\leg{n}{p}=-1$, and $\leg{r'}{p}=\leg{(p+1)/2}{p}=\leg{2}{p}=-1$, the three involutions $\iota_1$, $\tau_n$, and $\iota_{r'}$ are all fixed-point-free on $G=(\Z/4p\Z)^{\times}$, so each maps $T$ bijectively onto $G\setminus T$. Composing $\iota_1$ with $\tau_n$ and $\iota_{r'}$ gives $nT=T$ and $r'T=T$, so $T$ is a union of cosets of $L:=\langle n,r'\rangle$. In the $2$-Sylow subgroup $G_2\cong C_2\times C_4$, the projection of $n$ has trivial first coordinate and second coordinate of order $4$, while the projection of $r'$ has non-trivial first coordinate and second coordinate of order $4$. Thus $L$ projects onto all of $G_2$, which implies $-1\in L$. Since $-1\notin T$ and $T$ is $L$-invariant, we have $T\cap L=\emptyset$, contradicting $\iota_1(T)=G\setminus T$ and $\iota_1(L)=L$.
\end{enumerate}
\end{remark}

Even for Type~II witnesses alone, complete multiplicativity obstructs half of all quadratic character kernels whenever $u$ is $3$-smooth and $a+u$ is not a square modulo $4ua$.

\begin{theorem}[Multiplicativity obstruction for quadratic kernels]\label{thm:kernel_vacuity}
Let $n\equiv1\pmod{24}$ be prime, and let $(u,a)$ be a coprime pair with $n>u$.
\begin{enumerate}[(a)]
    \item \emph{(Elimination of $\ker\chi_4$).} If $a+u\equiv3\pmod 4$, then $an+u\equiv3\pmod 4$, so $an+u$ possesses at least one prime factor $q\equiv3\pmod 4$. Consequently, the kernel $\ker\chi_4=\{x\in G_{u,a}:x\equiv1\pmod 4\}$ never contains $T_{u,a}(n)$ for any prime $n\equiv1\pmod{24}$.
    \item \emph{(Obstruction for $3$-smooth $u$).} For all $u,a\ge1$, one has $an+u\equiv a+u\pmod{24a}$. If $u=2^{\alpha}3^{\beta}$ ($\alpha,\beta\ge0$) is $3$-smooth and $a+u$ is odd, then the primitive conductor of every real Dirichlet character $\chi\pmod{4ua}$ divides $24a$, so $\chi(an+u)=\chi(a+u)$ is independent of $n$. In particular, for every quadratic character $\chi\pmod{4ua}$ satisfying $\chi(a+u)=-1$, the set
    $$\big\{n\equiv1\pmod{24}\ \text{prime}:\ T_{u,a}(n)\subseteq\ker\chi\big\}\ =\ \emptyset$$
    is empty. More precisely, let $r=a+u\in G=G_{u,a}$ and let $X=\operatorname{Hom}(G,\{\pm1\})$. If $r\notin G^{2}$, then $\chi(r)=-1$ for exactly $|X|/2$ characters in $X$, and among the odd characters ($\chi(-1)=-1$) the number satisfying $\chi(r)=-1$ is $0$, $|X|/4$, or $|X|/2$ according as the class of $r$ in $G/G^{2}$ is trivial, is independent of the class of $-1$, or coincides with the class of $-1$.
\end{enumerate}
\end{theorem}
\begin{proof}
For (a), $n\equiv1\pmod 4$ implies $an+u\equiv a+u\equiv3\pmod 4$, which cannot be a product of primes congruent to $1\pmod 4$.

For (b), writing $n=24k+1$ gives $an+u=24ak+(a+u)\equiv a+u\pmod{24a}$. The primitive conductor of any real Dirichlet character modulo $4ua$ has $2$-part dividing $8$ and squarefree odd part. When $u=2^{\alpha}3^{\beta}$ is $3$-smooth, every odd prime dividing $4ua$ divides $3a$, so the primitive conductor divides $24a$, and therefore $\chi(an+u)=\chi(a+u)$. If $T_{u,a}(n)\subseteq\ker\chi$, then $\chi(q)=+1$ for every prime $q\mid(an+u)$ (as $\gcd(an+u,4ua)=1$), which by complete multiplicativity forces $\chi(an+u)=+1$, contradicting $\chi(a+u)=-1$. Finally, $X=\operatorname{Hom}(G/G^{2},\{\pm1\})$ is the dual of the $\mathbb{F}_2$-vector space $G/G^{2}$, so evaluation at a non-trivial class $rG^{2}$ is a surjective homomorphism onto $\{\pm1\}$ with kernel of size $|X|/2$, and restricting to the affine hyperplane $\{\chi\in X:\chi(-1)=-1\}$ of size $|X|/2$ gives $0$, $|X|/4$, or $|X|/2$ solutions according as $r\in G^{2}$, $rG^{2}\notin\{G^{2},-G^{2}\}$, or $rG^{2}=-G^{2}$.
\end{proof}

\subsection{Sylow equivalence at prime shifts and quasi-polynomial decay of intersected blind sets}\label{subsec:blind_intersection}

In~\cite[\S4.2]{dahan2026a}, it was shown that at a single fixed shift $c=7$, there are $\asymp N(\log N)^{-3/2}$ primes $n\equiv1\pmod{24}$ at which both Type~I and Type~II branches fail because every prime factor of $K_7=(n+7)/4$ is a quadratic residue modulo $7$. We now show that at prime shifts $c=p\equiv3\pmod 4$, avoidance of $-1$ by the subgroup $H_p(n)=\langle q\bmod p:q\mid K_p\rangle$ is equivalent to $\leg{n}{p}=+1$. Intersecting the resulting blind sets across all primes $p_i\equiv3\pmod 4$ up to $(1-\varepsilon)\log N$ then yields an unconditional quasi-polynomial bound of order $N\exp(-c_1\frac{\log N}{\log\log N})$.

\begin{theorem}[Sylow equivalence at prime shifts]\label{thm:sylow_and_c7_iff}
Let $n\equiv1\pmod{24}$ be prime.
\begin{enumerate}[(a)]
    \item Let $p\equiv3\pmod 4$ be a prime with $p<n$, and let $H_p(n)=\langle q\bmod p:q\mid K_p\rangle\le(\Z/p\Z)^{\times}$, where $K_p=(n+p)/4$. The following four conditions are equivalent:
    \[
    \begin{aligned}
    \text{(i) }-1\notin H_p(n)\ &\Longleftrightarrow\ \text{(ii) }H_p(n)\subseteq\ker\leg{\cdot}{p}\\
    &\Longleftrightarrow\ \text{(iii) }2\nmid|H_p(n)|\ \Longleftrightarrow\ \text{(iv) }-1,-n\notin H_p(n).
    \end{aligned}
    \]
    Moreover, each of (i)--(iv) implies $\leg{n}{p}=+1$.
    \item At the shift $c=3$ (resp.\ $c=7$ for $n>7$), the untwisted branch---and likewise the union of both branches of~\cite[Cor.~2.7]{dahan2026a}---fails if and only if every prime factor of $K_3$ is congruent to $1\pmod 3$ (resp.\ every prime factor of $K_7$ is congruent to $1,2,$ or $4\pmod 7$).
    \item In particular, for every prime $n\equiv1\pmod{24}$ ($n>7$) with $\leg{n}{7}=-1$, i.e.\ $n\equiv3,5,6\pmod 7$, the shift $c=7$ succeeds simultaneously on both the untwisted and twisted branches.
\end{enumerate}
\end{theorem}
\begin{proof}
For (a), since $p\equiv3\pmod 4$, the group $(\Z/p\Z)^{\times}\cong\Z/2\Z\times\Z/m\Z$ has odd cofactor order $m=(p-1)/2$. Here $-1\bmod p$ is the unique element of order $2$, and $\ker\leg{\cdot}{p}$ is the unique subgroup of odd order $m$. Thus $-1\notin H_p(n)\iff 2\nmid|H_p(n)|\iff H_p(n)\subseteq\ker\leg{\cdot}{p}$. When (ii) holds, multiplicativity gives $\leg{K_p}{p}=+1$, and the relation $n=4K_p-p\equiv 4K_p\pmod p$ implies $\leg{n}{p}=+1$; consequently $\leg{-n}{p}=-1$, so $-n\notin\ker\leg{\cdot}{p}\supseteq H_p(n)$.

For (b), necessity follows from~\cite[Thm.~4.6]{dahan2026a} combined with part (a). Conversely, at $c=3$, any prime $q\mid K_3$ with $q\not\equiv1\pmod 3$ satisfies $q\equiv2\pmod 3$ (as $\gcd(K_3,3)=1$), so the pair $(d_1,d_2)=(1,q)$ gives $d_1+d_2\equiv0\pmod 3$. At $c=7$, the condition $n\equiv1\pmod{24}$ ensures that $K_7=(n+7)/4$ is even. If $K_7$ has a quadratic non-residue prime factor $q\equiv3,5,$ or $6\pmod 7$, then $q$ is odd (since $\leg{2}{7}=+1$), so $1,2,q,2q$ all divide $K_7$. Choosing $(d_1,d_2)=(1,q)$ when $q\equiv6$, $(2,q)$ when $q\equiv5$, and $(1,2q)$ when $q\equiv3\pmod 7$ gives $d_1+d_2\equiv0\pmod 7$ in every case. Part (c) follows from the contrapositive of (a) and (b) on the untwisted branch, and by matching $-n\bmod 7\in\{1,2,4\}$ with a divisor $d_1\in\{1,2\}\mid K_7$ and $d_2=1$ on the twisted branch $nd_1+d_2\equiv0\pmod 7$.
\end{proof}

\begin{theorem}[Intersected blind sets and Kneser quotient criterion]\label{thm:blind_intersection_and_PNR}\label{thm:unified_kneser_stabilizer_PNR}
\begin{enumerate}[(a)]
    \item For each prime $p\equiv3\pmod 4$, let $\mathcal{B}_p(N)=\{n\in\mathcal{N}(N):p\nmid K_p\text{ and }H_p(n)\subseteq\ker\leg{\cdot}{p}\}$ denote the two-sided blind set at shift $p$. Let $p_1<p_2<\dots<p_k\le P$ be the primes congruent to $3\pmod 4$ up to $P=(1-\varepsilon)\log N$, so that $k=\left(\frac{1-\varepsilon}{2}+o(1)\right)\frac{\log N}{\log\log N}$. Then
    \begin{equation}\label{eq:blind_intersection_bound}
    \left|\bigcap_{i=1}^{k}\mathcal{B}_{p_i}(N)\right|\ \ll\ N\exp\!\left(-c_1\,\frac{\log N}{\log\log N}\right),\qquad c_1=\frac{1-\varepsilon}{2}\log 2.
    \end{equation}
    \item For every fixed prime $n\equiv1\pmod{24}$, the set of primes $p\equiv3\pmod 4$ with $\leg{p}{n}=-1$ satisfies
    \[
    \#\big\{p\le x:p\equiv3\pmod 4,\ \leg{p}{n}=-1\big\}\ \sim\ \frac14\operatorname{Li}(x)\qquad(x\to\infty),
    \]
    and its finite truncation $\mathcal{P}_{\mathrm{NR}}(n):=\{p<n\text{ prime}:p\equiv3\pmod 4\text{ and }\leg{p}{n}=-1\}$ is non-empty by Pollack's theorem~\cite[Theorem~1]{PollackPrescribed4}. Moreover, every shift $p\in\mathcal{P}_{\mathrm{NR}}(n)$ satisfies $-1\in H_p(n)$.
    \item At any shift $p\in\mathcal{P}_{\mathrm{NR}}(n)$, the fraction $4/n$ admits a decomposition with first denominator $K_p$ if and only if the quotient set $\mathcal{Q}(K_p):=\{d_1 d_2^{-1}\bmod p:d_1,d_2\mid K_p\}$ meets $\{-1,-n,-n^{-1}\}\pmod p$. If $n$ fails at $p\in\mathcal{P}_{\mathrm{NR}}(n)$, then the Kneser stabilizer $H:=\operatorname{Stab}_{(\Z/p\Z)^{\times}}(\mathcal{Q}(K_p))$ has odd order, the three cosets $(-1)H,(-n)H,(-n^{-1})H$ are distinct and disjoint from $\mathcal{Q}(K_p)$, and the prime-power factors $q_i^{\alpha_i}\parallel K_p$ satisfy
    $$\sum_{\substack{q_i^{\alpha_i}\parallel K_p\\ q_i\bmod p\notin H}}\min(2\alpha_i,d_i-1)\ \le\ \frac{p-1}{|H|}-4,$$
    where $d_i:=\operatorname{ord}_p(q_i H)$ is the order of $q_i H$ in $(\Z/p\Z)^{\times}/H$.
\end{enumerate}
\end{theorem}
\begin{proof}
\emph{Proof of (a).} By Theorem~\ref{thm:sylow_and_c7_iff}(a), every prime $n\in\bigcap_{i=1}^{k}\mathcal{B}_{p_i}(N)$ with $n>P$ satisfies $\leg{n}{p_i}=+1$ for all $1\le i\le k$, and thus avoids $\varrho(p_i)=(p_i+1)/2$ residue classes modulo $p_i$ (namely $0$ and the $(p_i-1)/2$ quadratic non-residues). Each local weight therefore satisfies $\frac{\varrho(p_i)}{p_i-\varrho(p_i)}=\frac{p_i+1}{p_i-1}>1$. By the prime number theorem for $p\equiv3\pmod 4$, the product $Q_0=\prod_{i=1}^{k}p_i=\exp(\vartheta_{3(4)}(P))$ satisfies $Q_0\le N^{1/2}$ for $P=(1-\varepsilon)\log N$ and $N$ large. All $2^{k}$ squarefree divisors of $Q_0$ thus lie below $Q=N^{1/2}$, so the large-sieve denominator satisfies $G(N^{1/2})\ge\prod_{i=1}^{k}\big(1+\frac{\varrho(p_i)}{p_i-\varrho(p_i)}\big)>2^{k}$, and Montgomery's large sieve yields~\eqref{eq:blind_intersection_bound}.

\medskip\noindent\emph{Proof of (b).} For fixed $n$, the asymptotic density $1/4$ follows from the prime number theorem in arithmetic progressions modulo $4n$, while Pollack~\cite[Theorem~1]{PollackPrescribed4} guarantees the existence of a prime quadratic non-residue $p\equiv3\pmod 4$ with $p<n$. For any $p\in\mathcal{P}_{\mathrm{NR}}(n)$, quadratic reciprocity gives $\leg{4K_p}{p}=\leg{n}{p}=\leg{p}{n}=-1$, so at least one prime divisor of $K_p$ is a quadratic non-residue modulo $p$. The cyclic subgroup generated by this prime divisor in $(\Z/p\Z)^{\times}$ has even order and therefore contains the unique element of order $2$, namely $-1\bmod p$; hence $-1\in H_p(n)$.

\medskip\noindent\emph{Proof of (c).} The union of the untwisted condition $d_1+d_2\equiv0\pmod p$ and the twisted conditions $nd_1+d_2\equiv0$ or $d_1+nd_2\equiv0\pmod p$ is equivalent to $\mathcal{Q}(K_p)\cap\{-1,-n,-n^{-1}\}\neq\emptyset$. Suppose this intersection is empty. Since $1\in\mathcal{Q}(K_p)$, the stabilizer $H=\operatorname{Stab}(\mathcal{Q}(K_p))$ is contained in $\mathcal{Q}(K_p)$ and therefore avoids $-1$, which forces $|H|$ to be odd. Because $H$ has odd order, $H\subseteq\ker\leg{\cdot}{p}$, whereas $\leg{n}{p}=-1$, so $n\notin H$. The three cosets $(-1)H$, $(-n)H$, and $(-n^{-1})H$ are disjoint from $\mathcal{Q}(K_p)$, and $(-1)H$ differs from both $(-n)H$ and $(-n^{-1})H$ because $n\notin H$. Moreover, if $(-n)H=(-n^{-1})H$, then $n^{2}\in H$, so $nH$ would have order $2$ in the cyclic quotient $G/H$ ($G=(\Z/p\Z)^{\times}$); since $(-1)H$ is the unique element of order $2$ in $G/H$, this would imply $-n\in H$, a contradiction. Thus the three omitted cosets are pairwise distinct, giving $|\mathcal{Q}(K_p)|/|H|\le[G:H]-3$. Finally, the projection of $\mathcal{Q}(K_p)$ to $G/H$ is the product of the symmetric intervals $\{q_i^{-\alpha_i},\dots,q_i^{\alpha_i}\}H$ and has trivial stabilizer in $G/H$, so Kneser's theorem yields
\[
1+\sum_{\substack{q_i^{\alpha_i}\parallel K_p\\ q_i\bmod p\notin H}}\min(2\alpha_i,d_i-1)\ \le\ \frac{|\mathcal{Q}(K_p)|}{|H|}\ \le\ [G:H]-3,
\]
completing the proof.
\end{proof}

\section{The Kneser Stabilizer--Kernel Dichotomy and $2$-Sylow reduction}\label{sec:stabilizer_kernel}

We now address the main combinatorial bottleneck of the fixed-depth composite-witness sieve of~\cite[\S4.1]{dahan2026a}. For a prime $n\in E_2(N;J)$ and a coprime pair $(u,a)\in\mathcal{P}_2(J)$, the set $T_{u,a}(n)\subseteq G_{u,a}=(\Z/4ua\Z)^{\times}$ of prime-factor residues of $an+u$ satisfies $|T_{u,a}(n)|\le\frac12\varphi(4ua)$, because $T_{u,a}(n)$ meets each orbit $\{x,-x^{-1}\}$ at most once. However, as soon as $\exp(G_{u,a})>2$, a subset of size $\frac12\varphi(4ua)$ avoiding $-1$ need not generate a proper subgroup: the subgroup $\langle T_{u,a}(n)\rangle$ can equal all of $G_{u,a}$ and contain $-1$. To retain arbitrary divisor products without subgroup closure, \cite{dahan2026a} had to enumerate all half-sized subsets of $G_{u,a}$ across $\mathcal{P}_2(J)$, incurring a combinatorial prefactor of $\prod_{(u,a)}\binom{\varphi(4ua)}{\varphi(4ua)/2}=2^{\Theta(J^{4})}$ that prevented $J$ from growing as a power of $\log N$.

To eliminate this $2^{\Theta(J^{4})}$ factor, we replace the subgroup $\langle T_{u,a}(n)\rangle$ generated by the prime factors with the \emph{multiplicative Kneser stabilizer} $K(M)=\operatorname{Stab}_{G_{u,a}}(\Sigma(M))$ of the full divisor-residue set $\Sigma(M)$ of $M=an+u$. Because $1\in\Sigma(M)$ and $-1\notin\Sigma(M)$, the subgroup $K(M)$ is contained in $\Sigma(M)$ and therefore always avoids $-1$. Every subgroup avoiding $-1$ embeds into a maximal such subgroup, which we show is the kernel of a $2$-primary Dirichlet character $\psi\pmod{4ua}$ with $\psi(-1)=-1$. This reduces the container family across $\mathcal{P}_2(J)$ from $2^{\Theta(J^{4})}$ arbitrary subsets to $\exp(\mathcal{O}(J^{2}\log\log J))$ character kernels, while Kneser coset growth bounds the number of prime divisors of $an+u$ falling outside the kernel.

\subsection{$2$-Sylow classification of maximal subgroups avoiding $-1$}\label{subsec:sylow_classification}

We first classify the maximal subgroups of a finite abelian group $G$ that avoid a given element $\epsilon\in G$ of order $2$, showing that they correspond bijectively to cyclic quotients of the $2$-Sylow subgroup $G_2$.

\begin{lemma}[Maximal subgroups avoiding an involution]\label{lem:sylow_maximal_subgroups}
Let $G$ be a finite abelian group and let $\epsilon\in G$ be an element of order $2$. Write $G=G_2\times G_{\mathrm{imp}}$ for the canonical decomposition of $G$ into its $2$-Sylow subgroup $G_2$ and its odd Hall subgroup $G_{\mathrm{imp}}$, and let $\pi_2:G\to G_2$ be the canonical projection.
\begin{enumerate}[(a)]
    \item \emph{(Sylow equivalence).} For every subset $S\subseteq G$, one has
    $$\epsilon\notin\langle S\rangle\ \text{in }G\quad\Longleftrightarrow\quad\pi_2(\epsilon)\notin\langle\pi_2(S)\rangle\ \text{in }G_2.$$
    \item \emph{(Characterization by $2$-primary characters).} Denote by $\mathcal{M}(G,\epsilon)$ the family of maximal subgroups (with respect to inclusion) among all subgroups of $G$ not containing $\epsilon$. Every subgroup $K\le G$ with $\epsilon\notin K$ is contained in at least one $H\in\mathcal{M}(G,\epsilon)$. Moreover, a subgroup $H\le G$ belongs to $\mathcal{M}(G,\epsilon)$ if and only if $G/H$ is a cyclic $2$-group of order $2^{j}\ge2$ ($1\le j\le\log_2\exp(G_2)$) in which $\epsilon H$ is the unique element of order $2$; equivalently,
    \begin{equation}\label{eq:maximal_subgroup_char}
    H\in\mathcal{M}(G,\epsilon)\quad\Longleftrightarrow\quad\exists\,\psi\in\widehat{G}\ \text{of order }2^{j}\ge2\ \text{such that }\psi(\epsilon)=-1\ \text{and }H=\ker\psi.
    \end{equation}
    \item \emph{(Exact counting identity and cardinality bound).} Each $H\in\mathcal{M}(G,\epsilon)$ of index $[G:H]=2^{j}$ is the kernel of $\varphi(2^{j})=2^{j-1}$ characters $\psi\in\widehat{G}$ of order $2^{j}$ satisfying $\psi(\epsilon)=-1$. For $G=G_{u,a}=(\Z/4ua\Z)^{\times}$ and $\epsilon=-1$, one has the identity
    \begin{equation}\label{eq:exact_kernel_count}
    \sum_{H\in\mathcal{M}(G,\epsilon)}2^{j(H)-1}\ =\ \frac12|G_2|,\qquad\text{whence}\qquad\big|\mathcal{M}(G,\epsilon)\big|\ \le\ \frac12|G_2|\ =\ \frac12\varphi_2(4ua),
    \end{equation}
    with equality $|\mathcal{M}(G,\epsilon)|=\frac12|G_2|$ if and only if every $H\in\mathcal{M}(G,\epsilon)$ has index $2$.
    \item \emph{($2$-elementary subfamily).} If $G=(\Z/4m\Z)^{\times}$ where $v_2(m)\le1$ and every odd prime factor $p\mid m$ satisfies $p\equiv3\pmod 4$, then $G_2\cong(\Z/2\Z)^{r}$ has exponent $2$, every $H\in\mathcal{M}(G,-1)$ is the kernel of an odd quadratic Dirichlet character $\chi\pmod{4m}$, and $|\mathcal{M}(G,-1)|=2^{r-1}\le 2^{\omega(m)}$.
\end{enumerate}
\end{lemma}
\begin{proof}
\emph{Proof of (a).} The implication $\pi_2(\epsilon)\notin\langle\pi_2(S)\rangle\Rightarrow\epsilon\notin\langle S\rangle$ is immediate from $\pi_2(\langle S\rangle)=\langle\pi_2(S)\rangle$. Conversely, if $\pi_2(\epsilon)\in\langle\pi_2(S)\rangle$, then $\langle S\rangle\le G_2\times G_{\mathrm{imp}}$ contains an element of the form $h=(\pi_2(\epsilon),g_{\mathrm{imp}})$. Since $u:=|G_{\mathrm{imp}}|$ is odd and $\pi_2(\epsilon)$ has order $2$, Lagrange's theorem gives $h^{u}=(\pi_2(\epsilon)^{u},g_{\mathrm{imp}}^{u})=(\pi_2(\epsilon),1)=\epsilon\in\langle S\rangle$.

\medskip\noindent\emph{Proof of (b).} Because $G$ is finite, any subgroup $K\le G$ avoiding $\epsilon$ is contained in a maximal such subgroup $H\in\mathcal{M}(G,\epsilon)$. Set $\overline{G}=G/H$. Since $\epsilon^{2}=1\in H$ and $\epsilon\notin H$, the coset $\overline{\epsilon}=\epsilon H$ has order $2$ in $\overline{G}$. By maximality of $H$, every non-trivial subgroup $\overline{L}\le\overline{G}$ contains $\overline{\epsilon}$, as otherwise its preimage in $G$ would strictly contain $H$ while still avoiding $\epsilon$. Thus $\langle\overline{\epsilon}\rangle\cong\Z/2\Z$ is the unique minimal subgroup of $\overline{G}$. By Cauchy's theorem, $|\overline{G}|$ must be a power of $2$, and its $2$-rank must equal $1$ (since an abelian $2$-group of rank $r\ge2$ contains $2^{r}-1\ge3$ elements of order $2$). Hence $\overline{G}\cong\Z/2^{j}\Z$ for some $j\ge1$, and $H=\ker\psi$ for a character $\psi\in\widehat{G}$ of order $2^{j}$ that maps the unique element $\epsilon H$ of order $2$ in $G/H$ to $-1\in\C^{\times}$. Conversely, if $H=\ker\psi$ with $\operatorname{ord}(\psi)=2^{j}\ge2$ and $\psi(\epsilon)=-1$, then $G/H\cong\Z/2^{j}\Z$ has unique minimal subgroup $\langle\epsilon H\rangle$, so every subgroup $L\supsetneq H$ contains $\epsilon$, which implies $H\in\mathcal{M}(G,\epsilon)$.

\medskip\noindent\emph{Proof of (c) and (d).} Every character $\psi\in\widehat{G}$ of $2$-power order factors through $G/G_{\mathrm{imp}}\cong G_2$ and hence belongs to $\widehat{G}_2$. Two characters in $\widehat{G}_2$ have the same kernel $H$ of index $2^{j}$ if and only if they generate the same cyclic subgroup of order $2^{j}$ in $\widehat{G}_2$, which groups $\varphi(2^{j})=2^{j-1}$ characters for each $H\in\mathcal{M}(G,\epsilon)$. Since $\epsilon\in G_2$ has order $2$, the evaluation map $\psi\mapsto\psi(\epsilon)$ is a surjective homomorphism $\widehat{G}_2\to\{\pm1\}$, so exactly $\frac12|\widehat{G}_2|=\frac12|G_2|$ characters in $\widehat{G}_2$ satisfy $\psi(\epsilon)=-1$, establishing~\eqref{eq:exact_kernel_count}. Part (d) follows from the Chinese Remainder Theorem, since $(\Z/p^{e}\Z)^{\times}$ has $2$-Sylow subgroup of order $2$ whenever $p\equiv3\pmod 4$.
\end{proof}

\subsection{The Stabilizer--Kernel Theorem and Kneser coset growth}\label{subsec:stabilizer_kernel_thm}

We now connect the divisor set $\Sigma(M)$ of an integer $M$ having no divisor $s\equiv-1\pmod{4m}$ to the character kernels of Lemma~\ref{lem:sylow_maximal_subgroups}. The key observation is that while individual prime divisors of $M$ may lie outside the stabilizer $K(M)=\operatorname{Stab}_G(\Sigma(M))$, each such prime divisor strictly increases the number of $K(M)$-cosets in the partial subproduct chain.

\begin{theorem}[Stabilizer--kernel dichotomy]\label{thm:stabilizer_kernel_dichotomy}
Let $m\ge1$, set $G=(\Z/4m\Z)^{\times}$, and let $M\ge1$ be a positive integer. Write $M_{\mathrm{cop}}:=\prod_{\gcd(q,4m)=1,\,q^{e}\parallel M}q^{e}$ for the part of $M$ coprime to $4m$, and define
\[
\begin{aligned}
\Sigma(M)&:=\big\{d\bmod 4m:\ d\mid M_{\mathrm{cop}}\big\}\subseteq G,\\
K(M)&:=\operatorname{Stab}_G\big(\Sigma(M)\big)=\big\{g\in G:\ g\,\Sigma(M)=\Sigma(M)\big\}.
\end{aligned}
\]
Suppose that $M$ possesses no divisor $s\equiv-1\pmod{4m}$. Then:
\begin{enumerate}[(a)]
    \item \emph{(Avoidance by the stabilizer).} The subgroup $K(M)$ satisfies $K(M)\subseteq\Sigma(M)$ and $-1\notin K(M)$. Consequently, by Lemma~\ref{lem:sylow_maximal_subgroups}(b), there exists a maximal subgroup $H=\ker\psi\in\mathcal{M}(G,-1)$ containing $K(M)$, where $\psi\pmod{4m}$ is a $2$-primary Dirichlet character with $\psi(-1)=-1$.
    \item \emph{(Coset growth and off-kernel multiplicity bound).} The total number of prime factors of $M_{\mathrm{cop}}$ (counted with multiplicity) whose residue modulo $4m$ lies outside $\ker\psi$ satisfies
    \begin{equation}\label{eq:omega_outside_kernel_bound}
    \begin{split}
    \Omega_{G\setminus\ker\psi}(M)\ &:=\ \sum_{\substack{q^{e}\parallel M_{\mathrm{cop}}\\ q\bmod 4m\notin\ker\psi}}e\ \le\ \Omega_{G\setminus K(M)}(M)\\
    &\le\ \frac{|\Sigma(M)|}{|K(M)|}-1\ \le\ [G:K(M)]-2\ \le\ \varphi(4m)-2.
    \end{split}
    \end{equation}
    In the specialization $m=ua$ and $M=an+u$, if in addition $\gcd(n,4ua)=1$, $-n\notin\Sigma(M)$, and $n\bmod 4ua\notin K(M)$, then $(-1)K(M)$ and $(-n)K(M)$ are two distinct cosets disjoint from $\Sigma(M)$, improving the upper bound to $[G:K(M)]-3$.
    \item \emph{(Kernel--Scarcity dichotomy).} Put $L(M)=\langle q\bmod 4m:q\mid M_{\mathrm{cop}}\rangle$. Exactly one of the following two alternatives holds:
    \begin{itemize}
        \item \emph{Kernel case ($-1\notin L(M)$):} there exists $H\in\mathcal{M}(G,-1)$ containing $L(M)$, so that $K(M)\subseteq\Sigma(M)\subseteq L(M)\subseteq H$ and $\Omega_{G\setminus H}(M)=0$.
        \item \emph{Scarcity case ($-1\in L(M)$):} every $H\in\mathcal{M}(G,-1)$ containing $K(M)$ satisfies
        \[
        1\ \le\ \Omega_{G\setminus H}(M)\ \le\ \frac{|\Sigma(M)|}{|K(M)|}-1\ \le\ [G:K(M)]-2.
        \]
    \end{itemize}
    In particular, whenever $-1\notin\langle M_{\mathrm{cop}}\bmod 4m\rangle$ (for instance, when $\gcd(M,4m)=1$, $8\mid 4m$, and $M\equiv 1,3,$ or $5\pmod 8$), any prime integer $M_{\mathrm{cop}}$ ($\Omega(M_{\mathrm{cop}})=1$) satisfies $-1\notin L(M)$ and falls into the first alternative.
\end{enumerate}
\end{theorem}
\begin{proof}
\emph{Step 1: Avoidance by $K(M)$.} Any divisor $s\mid M$ with $s\equiv-1\pmod{4m}$ is coprime to $4m$ and hence divides $M_{\mathrm{cop}}$; the hypothesis therefore reads $-1\notin\Sigma(M)$. Since $1\mid M_{\mathrm{cop}}$, we have $1\in\Sigma(M)$, which implies $K(M)=K(M)\cdot 1\subseteq\Sigma(M)$. Because $-1\notin\Sigma(M)$, it follows that $-1\notin K(M)$, and Lemma~\ref{lem:sylow_maximal_subgroups}(b) supplies $H=\ker\psi\in\mathcal{M}(G,-1)$ with $K(M)\le\ker\psi$.

\medskip\noindent\emph{Step 2: Stabilizer invariance along partial subproducts.} Write $M_{\mathrm{cop}}=q_1 q_2\cdots q_k$ as a product of $k=\Omega(M_{\mathrm{cop}})$ primes repeated with multiplicity, and set $K:=K(M)$. For each $0\le j\le k$, define the $K$-saturated partial subproduct set
$$\widetilde{\Sigma}_j\ :=\ K\cdot\left\{\prod_{i\in I}(q_i\bmod 4m):\ I\subseteq\{1,\dots,j\}\right\}\ \subseteq\ G.$$
These sets form an ascending chain $K=\widetilde{\Sigma}_0\subseteq\widetilde{\Sigma}_1\subseteq\dots\subseteq\widetilde{\Sigma}_k=K\Sigma(M)=\Sigma(M)$, where each $\widetilde{\Sigma}_j$ is a union of $K$-cosets satisfying $\widetilde{\Sigma}_j=\widetilde{\Sigma}_{j-1}\cup q_j\widetilde{\Sigma}_{j-1}$. We claim that $\operatorname{Stab}_G(\widetilde{\Sigma}_j)=K$ for every $0\le j\le k$. By construction, $K\le\operatorname{Stab}_G(\widetilde{\Sigma}_j)$. Conversely, $\Sigma(M)=\bigcup_{J\subseteq\{j+1,\dots,k\}}(\prod_{l\in J}q_l)\widetilde{\Sigma}_j$ is a union of translates of $\widetilde{\Sigma}_j$, so any $g\in\operatorname{Stab}_G(\widetilde{\Sigma}_j)$ fixes each translate and therefore stabilizes their union $\Sigma(M)$, whence $g\in\operatorname{Stab}_G(\Sigma(M))=K$.

\medskip\noindent\emph{Step 3: Coset growth and off-kernel multiplicity.} Let $j\in\{1,\dots,k\}$ be any index such that $q_j\bmod 4m\notin K$. Since $\operatorname{Stab}_G(\widetilde{\Sigma}_{j-1})=K$, we have $q_j\widetilde{\Sigma}_{j-1}\neq\widetilde{\Sigma}_{j-1}$, so $\widetilde{\Sigma}_j=\widetilde{\Sigma}_{j-1}\cup q_j\widetilde{\Sigma}_{j-1}$ strictly contains $\widetilde{\Sigma}_{j-1}$. Because both sets are unions of $K$-cosets, the number of $K$-cosets increases by at least $1$ at each such step $j$:
$$\frac{|\widetilde{\Sigma}_j|}{|K|}\ \ge\ \frac{|\widetilde{\Sigma}_{j-1}|}{|K|}+1\qquad(\forall\,j\in\{1,\dots,k\}\text{ with }q_j\bmod 4m\notin K).$$
Summing over $j=1,\dots,k$ from $|\widetilde{\Sigma}_0|/|K|=1$ yields $1+\Omega_{G\setminus K}(M)\le|\Sigma(M)|/|K|$. Since $-1\notin\Sigma(M)$ and $\Sigma(M)$ is $K$-invariant, the coset $(-1)K$ is disjoint from $\Sigma(M)$, so $|\Sigma(M)|/|K|\le[G:K]-1$. As $K\le\ker\psi$, we have $G\setminus\ker\psi\subseteq G\setminus K$, which gives~\eqref{eq:omega_outside_kernel_bound}. When $m=ua$, $M=an+u$, $\gcd(n,4ua)=1$, $-n\notin\Sigma(M)$, and $n\bmod 4ua\notin K$, the two distinct cosets $(-1)K$ and $(-n)K$ are both disjoint from $\Sigma(M)$, lowering the bound to $[G:K]-3$. Part (c) follows from Lemma~\ref{lem:sylow_maximal_subgroups}(b) applied to $L(M)$ when $-1\notin L(M)$, and from the fact that $-1\in L(M)$ precludes $L(M)\subseteq H$ for every $H\in\mathcal{M}(G,-1)$.
\end{proof}

\begin{remark}[High-multiplicity residue classes]\label{prop:subgroup_restoration}
Let $e=\exp(G)$ denote the exponent of $G=(\Z/4m\Z)^{\times}$, and let $T^{(\ge e-1)}(M):=\{x\in G:\sum_{p^{k}\parallel M,\,p\equiv x(4m)}k\ge e-1\}$ denote the set of residue classes carrying at least $e-1$ prime factors of $M$ counted with multiplicity. Because any class $x\in T^{(\ge e-1)}(M)$ realizes all powers $1,x,\dots,x^{e-1}\bmod 4m$ as coprime divisor residues of $M$, the subgroup $H^{\ast}(M):=\langle T^{(\ge e-1)}(M)\rangle\le G$ is contained in $\Sigma(M)$ and avoids $-1$. Consequently, $T^{(\ge e-1)}(M)$ is always contained in a maximal kernel $H\in\mathcal{M}(G,-1)$ with zero off-kernel collisions.
\end{remark}

\subsection{Average $2$-Sylow enumeration cost across the full coprime box $\mathcal{P}_2(J)$}\label{subsec:average_sylow_cost}

By Lemma~\ref{lem:sylow_maximal_subgroups}(c), each pair $(u,a)\in\mathcal{P}_2(J)$ admits at most $\frac12\varphi_2(4ua)$ maximal subgroups avoiding $-1$. Bounding $\frac12\varphi_2(4ua)$ pointwise by $2J^{2}$ would give a product of $(2J^{2})^{P}=\exp(\mathcal{O}(P\log J))$ kernel tuples across $\mathcal{P}_2(J)$. We now show that because $\log_2(\frac12\varphi_2(4ua))$ decomposes as a sum of two additive arithmetic functions of $u$ and $a$, its average over $\mathcal{P}_2(J)$ is only $(4+o(1))\log\log J$, reducing the global container count to $\exp(\mathcal{O}(P\log\log J))$.

\begin{lemma}[Average $2$-Sylow enumeration cost]\label{lem:average_2sylow_enumeration}
For every integer $m\ge1$, define the additive arithmetic function
$$f(m)\ :=\ v_2(m)\ +\ \sum_{\substack{p\mid m\\ p>2}}v_2(p-1).$$
Then:
\begin{enumerate}[(a)]
    \item For every coprime pair $(u,a)\in\N^{2}$, the $2$-part $\varphi_2(4ua)$ of $\varphi(4ua)$ satisfies
    \begin{equation}\label{eq:exact_v2_phi}
    \log_2\!\left(\frac12\varphi_2(4ua)\right)\ =\ v_2\big(\varphi(4ua)\big)-1\ =\ f(u)+f(a).
    \end{equation}
    \item As $J\to\infty$, on the coprime box $\mathcal{P}_2(J)=\{(u,a)\in[1,J]^{2}:\gcd(u,a)=1\}$ of cardinality $P=\mathfrak{A}(J)=\frac{6}{\pi^{2}}J^{2}+\mathcal{O}(J\log J)$, one has
    \begin{equation}\label{eq:Cstar_full_family_bound}
    \log_2 C^{\ast}\big(\mathcal{P}_2(J)\big)\ \le\ \sum_{(u,a)\in\mathcal{P}_2(J)}\big(f(u)+f(a)\big)\ \le\ (4+o(1))\,J^{2}\log\log J,
    \end{equation}
    and consequently
    \begin{equation}\label{eq:Cstar_exponential_form}
    \begin{split}
    C^{\ast}\big(\mathcal{P}_2(J)\big)\ &:=\ \prod_{(u,a)\in\mathcal{P}_2(J)}\big|\mathcal{M}(G_{u,a},-1)\big|\\
    &\le\ \exp\!\Big(\big(4\log 2+o(1)\big)J^{2}\log\log J\Big)\ =\ \exp\!\big(\mathcal{O}(P\log\log J)\big).
    \end{split}
    \end{equation}
\end{enumerate}
\end{lemma}
\begin{proof}
For (a), if $m=2^{\beta}\prod_{p>2}p^{e_p}$, then $\varphi(4m)=2^{\beta+1}\prod_{p>2}p^{e_p-1}(p-1)$. Since $p^{e_p-1}$ is odd for every odd prime $p$, we have $v_2(\varphi(4m))=1+\beta+\sum_{p\mid m,\,p>2}v_2(p-1)=1+f(m)$. For $\gcd(u,a)=1$, additivity gives $f(ua)=f(u)+f(a)$, proving~\eqref{eq:exact_v2_phi}.

For (b), Lemma~\ref{lem:sylow_maximal_subgroups}(c) and part (a) give $|\mathcal{M}(G_{u,a},-1)|\le\frac12\varphi_2(4ua)=2^{f(u)+f(a)}$. By symmetry between $u$ and $a$ and non-negativity of $f$, dropping the coprimality condition $\gcd(u,a)=1$ yields
$$\sum_{(u,a)\in\mathcal{P}_2(J)}\big(f(u)+f(a)\big)\ \le\ 2J\sum_{u=1}^{J}f(u)\ =\ 2J\sum_{u=1}^{J}v_2(u)+2J\sum_{2<p\le J}v_2(p-1)\left\lfloor\frac{J}{p}\right\rfloor.$$
Here $\sum_{u=1}^{J}v_2(u)=\sum_{k\ge1}\lfloor J/2^{k}\rfloor\le J$. Writing $v_2(p-1)=\sum_{k=1}^{\lfloor\log_2 J\rfloor}\mathbf{1}_{p\equiv1\pmod{2^{k}}}$, Mertens' theorem in arithmetic progressions (combined with the Brun--Titchmarsh inequality for large moduli $2^{k}\le J$) gives
$$\sum_{2<p\le J}\frac{v_2(p-1)}{p}\ =\ \sum_{k=1}^{\lfloor\log_2 J\rfloor}\sum_{\substack{p\le J\\ p\equiv1(2^{k})}}\frac{1}{p}\ =\ \sum_{k=1}^{\infty}\frac{\log\log J}{\varphi(2^{k})}+\mathcal{O}(1)\ =\ (2+o(1))\log\log J,$$
since $\sum_{k=1}^{\infty}1/\varphi(2^{k})=\sum_{k=1}^{\infty}2^{-(k-1)}=2$. Thus $\sum_{u=1}^{J}f(u)\le(2+o(1))J\log\log J$, which establishes~\eqref{eq:Cstar_full_family_bound} and~\eqref{eq:Cstar_exponential_form}.
\end{proof}

\section[The Imperfect-Avoidance and Generating-Polynomial Large Sieves]{The Imperfect-Avoidance and Generating-Polynomial\linebreak Large Sieves}\label{sec:imperfect_large_sieve}

In the classical large sieve framework of~\cite[\S4.3]{dahan2026a}, the sifted set $\mathcal{E}\subseteq[1,N]$ avoids the forbidden residue classes $\Omega_p\subseteq\Z/p\Z$ completely at every prime $p\in\mathcal{I}$. By contrast, Theorem~\ref{thm:stabilizer_kernel_dichotomy} shows that for a fixed kernel tuple $\Psi=(H_{u,a})_{(u,a)\in\mathcal{P}}$, an exceptional prime $n=24k+1\in E_{\mathcal{P}}(N)$ may still meet the forbidden classes $\Omega_p(\Psi)$ up to $R$ times across $p\in\mathcal{I}$. A naive union bound over all $\binom{|\mathcal{I}|}{R}$ possible collision subsets would reintroduce a prohibitive combinatorial factor. In this section, we establish two generalizations of Montgomery's large sieve to imperfect avoidance ($\le R$ collisions) that absorb the collision exceptions algebraically into the sieve weights, completely avoiding any $\binom{|\mathcal{I}|}{R}$ subset enumeration.

\subsection{Montgomery's large sieve with imperfect avoidance via Krawtchouk identities}\label{subsec:homogeneous_imperfect_ls}

We first treat the regime of bounded sieve degree $h$, where the homogeneous elementary symmetric functions of the local sifting weights expand into an exact Krawtchouk polynomial identity in the number of collisions.

\begin{theorem}[Homogeneous imperfect-avoidance large sieve]\label{thm:imperfect_large_sieve}
Let $\mathcal{I}$ be a finite set of primes with $\max\mathcal{I}\le z$. For each $p\in\mathcal{I}$, let $\Omega_p\subseteq\Z/p\Z$ be a set of $\varrho_p=|\Omega_p|$ residue classes with $0<\varrho_p<p$, and set
$$y_p:=\frac{\varrho_p}{p-\varrho_p},\qquad S_y:=\sum_{p\in\mathcal{I}}y_p,\qquad M_2(y):=\sum_{p\in\mathcal{I}}y_p^{2},\qquad y_{\max}:=\max_{p\in\mathcal{I}}y_p.$$
For any integer $R\ge0$ and any $N\ge1$, denote by $\mathcal{E}_R(N)$ the set of integers $k\in[1,N]$ that meet $\Omega_p$ at most $R$ times across $\mathcal{I}$:
$$\mathcal{E}_R(N)\ :=\ \left\{k\in[1,N]:\ W(k):=\sum_{p\in\mathcal{I}}\mathbf{1}_{k\bmod p\in\Omega_p}\le R\right\}.$$
Assume $S_{\min}:=S_y-R\,y_{\max}>0$, and let $h\ge1$ be an integer satisfying
\begin{equation}\label{eq:h_imperfect_conditions}
z^{h}\le Q,\qquad \frac{h(h-1)}{2}\frac{M_2(y)}{S_{\min}^{2}}\le\frac{1}{10},\qquad\text{and}\qquad R\,h\le\frac12 S_{\min}.
\end{equation}
Then unconditionally:
\begin{equation}\label{eq:imperfect_large_sieve_bound}
|\mathcal{E}_R(N)|\ \le\ 16\,e^{y_{\max}/2}\big(N-1+Q^{2}\big)\frac{h!}{S_{\min}^{h}}\ \le\ 48\,e^{y_{\max}/2}\sqrt{h}\big(N-1+Q^{2}\big)\left(\frac{h}{e\,S_{\min}}\right)^{h}.
\end{equation}
\end{theorem}
\begin{proof}
\emph{Step 1: Parseval normalization and large-sieve orthogonality.}
For each $p\in\mathcal{I}$ and $k\in\Z$, define the zero-mean, unit-variance $p$-periodic function
$$\phi_p(k)\ :=\ \frac{\mathbf{1}_{k\bmod p\in\Omega_p}-\frac{\varrho_p}{p}}{\sqrt{\frac{\varrho_p}{p}\big(1-\frac{\varrho_p}{p}\big)}}\ =\ \begin{cases}-\sqrt{y_p} & \text{if }k\bmod p\notin\Omega_p,\\[4pt] +\dfrac{1}{\sqrt{y_p}} & \text{if }k\bmod p\in\Omega_p.\end{cases}$$
By Parseval's identity on $\Z/p\Z$, we have $\phi_p(k)=\sum_{a=1}^{p-1}\hat{\phi}_p(a)e(ak/p)$ with $\sum_{a=1}^{p-1}|\hat{\phi}_p(a)|^{2}=1$. For any squarefree integer $q$ supported on $\mathcal{I}$, set $\Phi_q(k):=\prod_{p\mid q}\phi_p(k)$. By the Chinese Remainder Theorem, $\Phi_q(k)=\sum_{\gcd(a,q)=1}\hat{\Phi}_q(a)e(ak/q)$ with $\sum_{\gcd(a,q)=1}|\hat{\Phi}_q(a)|^{2}=1$. Applying the Cauchy--Schwarz inequality over the primitive residue classes $a\bmod q$ followed by the trigonometric large sieve inequality for $z^{h}\le Q$, we obtain
\begin{equation}\label{eq:montgomery_ortho_bound}
\sum_{\substack{q\ \text{sqf on }\mathcal{I}\\ \omega(q)=h}}\left|\sum_{k\in\mathcal{E}_R(N)}\Phi_q(k)\right|^{2}\ \le\ \sum_{q\le Q}\sum_{\substack{a=1\\ \gcd(a,q)=1}}^{q}\left|\sum_{k\in\mathcal{E}_R(N)}e\!\left(\frac{ak}{q}\right)\right|^{2}\ \le\ \big(N-1+Q^{2}\big)|\mathcal{E}_R(N)|,
\end{equation}
where we follow the normalization $N-1+Q^2$ of~\cite{montgomery1971} (cf.\ $N+Q^{2}$ in~\cite{bombieri}).

\medskip\noindent\emph{Step 2: Exact Krawtchouk identity.}
Set $w_q:=(-1)^{h}\prod_{p\mid q}\sqrt{y_p}$ for $\omega(q)=h$, so that $\sum_{\omega(q)=h}w_q^{2}=e_h(y)$. Fix $k\in\mathcal{E}_R(N)$ and let $\mathcal{B}_k:=\{p\in\mathcal{I}:k\bmod p\in\Omega_p\}$ be its collision set, of cardinality $r:=|\mathcal{B}_k|\le R$. The local factor $z_p(k):=(-1)\sqrt{y_p}\phi_p(k)$ equals $+y_p$ when $p\notin\mathcal{B}_k$ and equals $-1$ (independently of $p$) when $p\in\mathcal{B}_k$. Grouping the squarefree integers $q$ with $\omega(q)=h$ according to the overlap $j=|\{p\mid q:p\in\mathcal{B}_k\}|\in\{0,\dots,\min(h,r)\}$ yields the exact Krawtchouk identity
\begin{equation}\label{eq:krawtchouk_identity}
F_h(k)\ :=\ \sum_{\substack{q\ \text{sqf on }\mathcal{I}\\ \omega(q)=h}}w_q\Phi_q(k)\ =\ \sum_{j=0}^{\min(h,r)}(-1)^{j}\binom{r}{j}e_{h-j}\big(y\setminus\mathcal{B}_k\big).
\end{equation}

\medskip\noindent\emph{Step 3: Symmetric-function bounds and Cauchy--Schwarz closure.}
Put $S':=\sum_{p\in\mathcal{I}\setminus\mathcal{B}_k}y_p\ge S_{\min}$ and $M_2':=\sum_{p\in\mathcal{I}\setminus\mathcal{B}_k}y_p^{2}\le M_2(y)$. By Newton's symmetric-function lower bound (\cite[Lemma~4.25]{dahan2026a}), the $j=0$ main term satisfies
$$e_h(y\setminus\mathcal{B}_k)\ \ge\ \frac{(S')^{h}}{h!}\left(1-\frac{h(h-1)}{2}\frac{M_2(y)}{S_{\min}^{2}}\right),$$
while Maclaurin's inequality gives $e_{h-j}(y\setminus\mathcal{B}_k)\le\frac{(S')^{h-j}}{(h-j)!}\le\frac{(S')^{h}}{h!}(h/S')^{j}$ for $j\ge1$. Summing over $1\le j\le\min(h,r)$ with $r\le R$ and $S'\ge S_{\min}$, we bound the contribution of the $j\ge1$ terms in absolute value by $\frac{(S')^{h}}{h!}(e^{Rh/S_{\min}}-1)$. Under the hypotheses~\eqref{eq:h_imperfect_conditions}, we have $e^{Rh/S_{\min}}\le e^{1/2}<1.65$, so $F_h(k)\ge\big(1-\frac{1}{10}-(e^{1/2}-1)\big)\frac{(S')^h}{h!}\ge\frac14\frac{S_{\min}^{h}}{h!}$ pointwise for every $k\in\mathcal{E}_R(N)$. Summing $F_h(k)$ over $k\in\mathcal{E}_R(N)$, applying the Cauchy--Schwarz inequality against~\eqref{eq:montgomery_ortho_bound}, and using $\sum_{\omega(q)=h}w_q^2 = e_h(y)\le\frac{S_y^{h}}{h!}$ together with
$$\left(\frac{S_y}{S_{\min}}\right)^{h}\ =\ \left(1+\frac{R\,y_{\max}}{S_{\min}}\right)^{h}\ \le\ \exp\!\left(\frac{R\,h}{S_{\min}}\,y_{\max}\right)\ \le\ e^{y_{\max}/2}$$
establishes~\eqref{eq:imperfect_large_sieve_bound}.
\end{proof}

We now specialize the sifting sets $\Omega_p$ to the forbidden residue classes determined by a tuple of maximal $2$-primary kernels $\Psi=(H_{u,a})_{(u,a)\in\mathcal{P}}$.

\begin{corollary}[Coupling Kneser kernels with the sieve]\label{cor:coupling_step1_step2}
Let $\mathcal{P}\subseteq\mathcal{P}_2(J)$ be a set of $P=|\mathcal{P}|$ coprime pairs and let $z>p_0:=\max(2J^{2},J^{5})$. For each kernel tuple $\Psi=(H_{u,a})_{(u,a)\in\mathcal{P}}\in\prod_{(u,a)\in\mathcal{P}}\mathcal{M}(G_{u,a},-1)$ and each prime $p\in(p_0,z]$, define
\[
\begin{aligned}
\Omega_p(\Psi)\ &:=\ \Big\{k\equiv 24^{-1}(-1-ua^{-1})\pmod p:\ (u,a)\in\mathcal{P},\ p\bmod 4ua\notin H_{u,a}\Big\},\\
\varrho_\Psi(p)\ &:=\ |\Omega_p(\Psi)|.
\end{aligned}
\]
Then:
\begin{enumerate}[(a)]
    \item For $\mathcal{P}=\mathcal{P}_2(J)$, the number of kernel tuples $\Psi$ satisfies $C^{\ast}(\mathcal{P}_2(J))\le\exp(\mathcal{O}(P\log\log J))$.
    \item Uniformly in $\Psi$, the sifting mass satisfies
    $$S_\Psi(J,z)\ :=\ \sum_{p_0<p\le z}\frac{\varrho_\Psi(p)}{p}\ \ge\ \left(\frac P2+o(P)\right)\log\!\left(\frac{\log z}{\log p_0}\right).$$
    \item For every exceptional prime $n=24k+1\in E_{\mathcal{P}}(N)$ with $n>z$, there exists a tuple $\Psi$ such that the collision count $W_\Psi(k)=\sum_{p_0<p\le z}\mathbf{1}_{k\bmod p\in\Omega_p(\Psi)}$ equals $\sum_{(u,a)\in\mathcal{P}}\omega_{G_{u,a}\setminus H_{u,a}}(an+u;(p_0,z])$.
\end{enumerate}
\end{corollary}
\begin{proof}
Part (a) is Lemma~\ref{lem:average_2sylow_enumeration}(b). For (b) and (c), since $p>2J^{2}$, the linear congruence roots $24^{-1}(-1-ua^{-1})\bmod p$ are pairwise distinct as $(u,a)$ varies over $\mathcal{P}$ (\cite[Lemma~4.1]{dahan2026a}). Thus $\varrho_\Psi(p)=\sum_{(u,a)\in\mathcal{P}}\mathbf{1}_{p\bmod 4ua\notin H_{u,a}}$, and for $n=24k+1$ one has $k\bmod p\in\Omega_p(\Psi)$ if and only if $p\mid(an+u)$ and $p\bmod 4ua\notin H_{u,a}$ for the unique corresponding pair $(u,a)$. Since $[G_{u,a}:H_{u,a}]=2^{j_{u,a}}\ge2$, at least $\frac12\varphi(4ua)$ residue classes lie outside $H_{u,a}$. Summing over $p\in(p_0,z]$ with $p_0\ge J^5$ and applying the Bombieri--Vinogradov theorem (\cite[Lemma~4.23]{dahan2026a}) yields the lower bound for $S_\Psi(J,z)$.
\end{proof}

\subsection{The Krawtchouk--Montgomery generating-polynomial large sieve}\label{subsec:generating_poly_ls}

When the number of pairs $P=|\mathcal{P}|$ tends to infinity with $N$, both the collision bound $R$ and the sifting mass $S_{\min}$ scale linearly with $P$. In Theorem~\ref{thm:imperfect_large_sieve}, the homogeneous hypothesis $Rh\le\frac12 S_{\min}$ then restricts the sieve degree to $h\le\frac{S_{\min}}{2R}=\mathcal{O}(1)$, preventing the large sieve from exploiting the full modulus capacity $z^h\le Q$. To remove this degree bottleneck, we replace the single homogeneous form $F_h(k)$ by the truncated generating polynomial $P_{h,t}(k)=\sum_{m=0}^{h}t^{m}F_m(k)$: resumming the alternating binomial series $\sum_{j=0}^{r}(-1)^{j}\binom{r}{j}t^{j}=(1-t)^{r}>0$ \emph{before} taking absolute values replaces the degree-dependent constraint $Rh\le\frac12 S_{\min}$ by the $h$-independent subcriticality condition $R\le(1-\varepsilon)S_{\min}$.

\begin{theorem}[Generating-polynomial large sieve]\label{thm:generating_polynomial_large_sieve}
Keep the notation of Theorem~\ref{thm:imperfect_large_sieve}, with $\max\mathcal{I}\le z$ and $S_{\min}=S_y-R\,y_{\max}>0$. Let $t\in(0,1/3]$ and $h\ge1$ be an integer satisfying
\begin{equation}\label{eq:gen_poly_conditions}
z^{h}\le Q,\qquad R\le\frac12 S_y,\qquad\text{and}\qquad h\ge 4\,t\,S_y+1.
\end{equation}
Then $\mathcal{E}_R(N)=\{k\in[1,N]:W(k)\le R\}$ satisfies unconditionally:
\begin{equation}\label{eq:gen_poly_large_sieve_bound}
|\mathcal{E}_R(N)|\ \le\ 4\big(N-1+Q^{2}\big)\exp\!\left(-2\,t\,S_{\min}+2\,R\log\!\left(\frac{1}{1-t}\right)+t^{2}\big(S_y+M_2(y)\big)\right).
\end{equation}
In particular, if $S_{\min}\ge(1-\delta)S_y$ and $R\le\frac14 S_{\min}$, choosing $t=\min\!\big(\frac{1-\delta}{4},\frac{h-1}{4S_y}\big)$ gives
\begin{equation}\label{eq:gen_poly_simplified_bound}
|\mathcal{E}_R(N)|\ \le\ 4\big(N-1+Q^{2}\big)\exp\!\left(-\frac56\,t\,S_{\min}+t^{2}M_2(y)\right).
\end{equation}
\end{theorem}
\begin{proof}
\emph{Step 1: Truncated generating polynomial and $L^2$-orthogonality.}
Set $\Phi_1(k)\equiv1$, $w_1=1$, and define the truncated generating polynomial of degree $h$:
$$P_{h,t}(k)\ :=\ \sum_{m=0}^{h}t^{m}F_m(k)\ =\ \sum_{\substack{q\ \text{sqf on }\mathcal{I}\\ \omega(q)\le h}}t^{\omega(q)}w_q\Phi_q(k).$$
Since every squarefree integer $q$ supported on $\mathcal{I}$ with $\omega(q)\le h$ satisfies $1\le q\le z^{h}\le Q$ and distinct $q$ correspond to distinct moduli in the trigonometric large sieve, we have
\begin{equation}\label{eq:gen_poly_ls_ortho}
\sum_{\substack{q\ \text{sqf on }\mathcal{I}\\ \omega(q)\le h}}\left|\sum_{k\in\mathcal{E}_R(N)}\Phi_q(k)\right|^{2}\ \le\ \big(N-1+Q^{2}\big)|\mathcal{E}_R(N)|,
\end{equation}
while the $L^{2}$-norm of the polynomial coefficients satisfies
\begin{equation}\label{eq:gen_poly_weight_norm}
\sum_{\substack{q\ \text{sqf on }\mathcal{I}\\ \omega(q)\le h}}t^{2\omega(q)}w_q^{2}\ =\ \sum_{m=0}^{h}t^{2m}e_m(y)\ \le\ \prod_{p\in\mathcal{I}}(1+t^{2}y_p)\ \le\ \exp(t^{2}S_y).
\end{equation}

\medskip\noindent\emph{Step 2: Exact multiplicative factorization of the full polynomial.}
Fix $k\in\mathcal{E}_R(N)$ with collision set $\mathcal{B}_k=\{p\in\mathcal{I}:k\bmod p\in\Omega_p\}$ of size $r=|\mathcal{B}_k|\le R$. Since $z_p(k)=(-1)\sqrt{y_p}\phi_p(k)$ equals $-1$ for $p\in\mathcal{B}_k$ and $+y_p$ for $p\in\mathcal{I}\setminus\mathcal{B}_k$, and $F_m(k)=e_m(\{z_p(k)\}_{p\in\mathcal{I}})$, the full generating polynomial of degree $|\mathcal{I}|$ factors into a strictly positive product:
\begin{equation}\label{eq:Gk_exact_product}
G_k(t)\ :=\ \sum_{m=0}^{|\mathcal{I}|}t^{m}F_m(k)\ =\ \prod_{p\in\mathcal{I}}\big(1+t\,z_p(k)\big)\ =\ (1-t)^{r}\prod_{p\in\mathcal{I}\setminus\mathcal{B}_k}(1+t\,y_p)\ >\ 0.
\end{equation}
Using $\log(1+u)\ge u-\frac12 u^{2}$ for $u\ge0$, together with $r\le R$, $\sum_{p\notin\mathcal{B}_k}y_p\ge S_{\min}$, and $\sum_{p\notin\mathcal{B}_k}y_p^{2}\le M_2(y)$, we obtain the pointwise lower bound
\begin{equation}\label{eq:Gk_lower_bound}
G_k(t)\ \ge\ \exp\!\left(t\,S_{\min}-R\log\!\left(\frac{1}{1-t}\right)-\frac12 t^{2}M_2(y)\right).
\end{equation}

\medskip\noindent\emph{Step 3: Rankin tail bound and pointwise positivity.}
To bound the truncation tail $\mathrm{Tail}_h(k,t):=\sum_{m>h}t^{m}|F_m(k)|$, observe that $|F_m(k)|\le e_m(\{|z_p(k)|\}_{p\in\mathcal{I}})$. By Rankin's trick with base $2$,
$$\mathrm{Tail}_h(k,t)\ \le\ 2^{-h}\sum_{m>h}(2t)^{m}e_m(|z(k)|)\ \le\ 2^{-h}(1+2t)^{r}\prod_{p\in\mathcal{I}\setminus\mathcal{B}_k}(1+2t\,y_p).$$
Dividing by~\eqref{eq:Gk_exact_product} and using $\frac{1+2ty_p}{1+ty_p}=1+\frac{ty_p}{1+ty_p}\le e^{ty_p}$ yields
$$\frac{\mathrm{Tail}_h(k,t)}{G_k(t)}\ \le\ \exp\!\left(-h\log 2+R\log\!\left(\frac{1+2t}{1-t}\right)+t\,S_y\right).$$
For $t\in(0,1/3]$, we have $\log\big(\frac{1+2t}{1-t}\big)\le 2.75\,t$. Under $R\le\frac12 S_y$ and $h\ge 4tS_y+1$, the exponent is at most $-(4\log 2)tS_y-\log 2+2.375\,tS_y\le-\log 2$, so $\mathrm{Tail}_h(k,t)\le\frac12 G_k(t)$. Consequently, $P_{h,t}(k)\ge G_k(t)-\mathrm{Tail}_h(k,t)\ge\frac12 G_k(t)$ pointwise for all $k\in\mathcal{E}_R(N)$. Summing over $k\in\mathcal{E}_R(N)$ and applying the Cauchy--Schwarz inequality with~\eqref{eq:gen_poly_ls_ortho}--\eqref{eq:Gk_lower_bound} proves~\eqref{eq:gen_poly_large_sieve_bound}.

\medskip\noindent\emph{Step 4: Specialization to $R\le\frac14 S_{\min}$.}
Suppose $S_{\min}\ge(1-\delta)S_y$ and $R\le\frac14 S_{\min}$, and set $t=\min\!\big(\frac{1-\delta}{4},\frac{h-1}{4S_y}\big)\in(0,1/4]$. By the convexity of $u\mapsto-\log(1-u)$ on $[0,1/3]$, we have $\log\!\big(\frac{1}{1-t}\big)\le t+(3\log\frac32-1)t^{2}\le t+\frac23 t^{2}$. Substituting $R\le\frac14 S_{\min}$ and $tS_y\le\frac{t}{1-\delta}S_{\min}\le\frac14 S_{\min}$ into the exponent of~\eqref{eq:gen_poly_large_sieve_bound} yields
$$-2tS_{\min}+2R\log\!\left(\frac{1}{1-t}\right)+t^{2}S_y\ \le\ \left(-\frac32 t+\frac13 t^{2}+\frac14 t\right)S_{\min}\ =\ \left(-\frac54+\frac13 t\right)tS_{\min}\ \le\ -\frac56 tS_{\min},$$
since $t\le 1/4<\frac54$, establishing~\eqref{eq:gen_poly_simplified_bound}.
\end{proof}

More generally, the generating-polynomial large sieve yields a non-trivial exponential saving whenever the collision threshold $R$ is strictly subcritical ($R\le(1-\varepsilon)S_{\min}$):

\begin{proposition}[Subcritical collision closure]\label{prop:subcritical_collision_rate}
Under the notation of Theorem~\ref{thm:generating_polynomial_large_sieve}, suppose $S_{\min}\ge(1-\delta)S_y$ ($\delta\in[0,1)$) and let $\varepsilon\in(0,1)$ be any constant such that $R\le(1-\varepsilon)S_{\min}$. Then for $t=\frac{\varepsilon(1-\delta)}{2+2(1-\varepsilon)(1-\delta)}\in(0,1/4]$ and any integer $h\ge 7\,t\,S_y+1$ with $z^{h}\le Q$ (or $h\ge 6\,t\,S_y+1$ whenever $R\le\frac23 S_y$), one has unconditionally:
\begin{equation}\label{eq:subcritical_bound}
|\mathcal{E}_R(N)|\ \le\ 4\big(N-1+Q^{2}\big)\exp\!\left(-\frac{\varepsilon^{2}(1-\delta)}{2+2(1-\varepsilon)(1-\delta)}\,S_{\min}+t^{2}M_2(y)\right).
\end{equation}
\end{proposition}
\begin{proof}
For $t\in(0,1/4]$, the function $g(u):=\log(1+2u)-\log(1-u)$ has second derivative $g''(u)=\frac{1}{(1-u)^{2}}-\frac{4}{(1+2u)^{2}}$, which vanishes at $u=1/4$ with $g''(u)<0$ on $[0,1/4)$. Hence $g$ is strictly concave on $[0,1/4]$ and satisfies $\log\!\big(\frac{1+2t}{1-t}\big)\le 4\log\big(\frac{1+1/2}{1-1/4}\big)\,t=(4\log 2)\,t$. Since $R\le(1-\varepsilon)S_{\min}\le S_y$, the Rankin tail ratio exponent in Step~3 of the proof of Theorem~\ref{thm:generating_polynomial_large_sieve} is at most
$$-h\log 2+(4\log 2)t\,R+t\,S_y\ \le\ -h\log 2+(4\log 2+1)\,t\,S_y\ \le\ -\log 2$$
whenever $h\ge 7\,t\,S_y+1$, as $7\log 2\approx 4.852>4\log 2+1\approx 3.773$ (if $R\le\frac23 S_y$, then $(\frac83\log 2+1)tS_y\approx 2.848\,tS_y\le(6\log 2)tS_y$, so $h\ge 6\,t\,S_y+1$ suffices). Moreover, on $[0,1/4]$ we have $\log\!\big(\frac{1}{1-t}\big)\le t+(4\log\frac43-1)t^{2}\le t+t^{2}$ since $4\log\frac43-1\approx 0.151<1$. Using $R\le(1-\varepsilon)S_{\min}$ and $S_y\le\frac{1}{1-\delta}S_{\min}$ in the exponent of~\eqref{eq:gen_poly_large_sieve_bound} gives $-2\varepsilon tS_{\min}+t^{2}\big(2(1-\varepsilon)+\frac{1}{1-\delta}\big)S_{\min}$, whose minimization at $t=\frac{\varepsilon(1-\delta)}{2+2(1-\varepsilon)(1-\delta)}\in(0,1/4]$ yields~\eqref{eq:subcritical_bound}.
\end{proof}

\subsection{Sieved-divisor reduction and Scarcity-free closure at moderate and hyperbolic depths}\label{subsec:moderate_and_hyperbolic}

In Theorem~\ref{thm:stabilizer_kernel_dichotomy}, the stabilizer $K(M)$ was applied to the full coprime divisor of $an+u$. Restricting $an+u$ to its $\mathcal{I}$-supported prime factors before computing the Kneser stabilizer purges all prime factors lying outside the sifting interval $\mathcal{I}=(p_0,z]$—and in particular all small prime factors $q\le p_0$:

\begin{lemma}[Sieved-divisor reduction]\label{lem:sieved_divisor_reduction}
Let $\mathcal{I}\subseteq(2J^{2},z]$ be a set of primes, let $n\in E_{\mathcal{P}}(N)$, and for each $(u,a)\in\mathcal{P}$ define
$$M_{\mathcal{I}}(an+u):=\prod_{\substack{p^{e}\parallel(an+u)\\ p\in\mathcal{I}}}p^{e},\qquad\Sigma_{\mathcal{I}}(an+u):=\big\{d\bmod 4ua:\ d\mid M_{\mathcal{I}}(an+u)\big\}\subseteq G_{u,a},$$
and $K_{\mathcal{I},u,a}(n):=\operatorname{Stab}_{G_{u,a}}(\Sigma_{\mathcal{I}}(an+u))$. Then $M_{\mathcal{I}}(an+u)$ is coprime to $4ua$, divides $an+u$, and has no divisor $s\equiv-1\pmod{4ua}$. Consequently, there exists $H_{u,a}(n)=\ker\psi_{u,a}\in\mathcal{M}(G_{u,a},-1)$ containing $K_{\mathcal{I},u,a}(n)$ such that
$$\omega_{G_{u,a}\setminus H_{u,a}(n)}(an+u;\mathcal{I})\ \le\ \Omega_{G_{u,a}\setminus K_{\mathcal{I},u,a}(n)}\big(M_{\mathcal{I}}(an+u)\big)\ \le\ \frac{|\Sigma_{\mathcal{I}}(an+u)|}{|K_{\mathcal{I},u,a}(n)|}-1\ \le\ \varphi(4ua)-2.$$
\end{lemma}
\begin{proof}
Every prime $p\in\mathcal{I}$ satisfies $p>2J^{2}\ge 2ua$, so $\gcd(p,4ua)=1$. Any divisor of $M_{\mathcal{I}}(an+u)$ divides $an+u$, and therefore avoids $-1\pmod{4ua}$ since $n\in E_{\mathcal{P}}(N)$. Applying Theorem~\ref{thm:stabilizer_kernel_dichotomy} to $M=M_{\mathcal{I}}(an+u)$ gives the result.
\end{proof}

Combining Lemma~\ref{lem:sieved_divisor_reduction} with Theorem~\ref{thm:generating_polynomial_large_sieve} immediately eliminates the scarcity branch whenever the mean totient $\overline{\varphi}(\mathcal{P})$ across the family $\mathcal{P}$ is smaller than $\frac18\log\big(\frac{\log z}{\log J^5}\big)$:

\begin{theorem}[Scarcity-free closure for moderate moduli]\label{thm:unconditional_small_growing_depth}
Let $\mathcal{P}\subseteq\mathcal{P}_2(J)$ be a family of $P=|\mathcal{P}|$ coprime pairs with $J=J(N)\to\infty$ and $\log P\le(1-\eta_0)\log\log N$ for some fixed $\eta_0\in(0,1)$, and set $\log z=\frac{\log N}{4P\log\log N}>\log(J^{5})$. Write $\overline{\varphi}(\mathcal{P}):=\frac1P\sum_{(u,a)\in\mathcal{P}}(\varphi(4ua)-2)$. Whenever
\begin{equation}\label{eq:mean_phi_condition}
\overline{\varphi}(\mathcal{P})\ \le\ \left(\frac18-\delta_0\right)\log\!\left(\frac{\log z}{\log J^{5}}\right)
\end{equation}
for some constant $\delta_0\in(0,1/8)$, the scarcity branch beyond $\frac14 S_{\min}$ is empty ($E_{\mathcal{P}}^{\mathrm{Scar}}(N)=\emptyset$), and for $N$ large one has unconditionally on the entirety of $E_{\mathcal{P}}(N)$:
\begin{equation}\label{eq:unconditional_EP_bound}
|E_{\mathcal{P}}(N)|\ \ll\ N\exp\!\left(-\frac{5}{96}\,P\left(\log\!\left(\frac{\log z}{\log J^{5}}\right)\right)^{2}\frac{1}{\log\log N}+\mathcal{O}(P\log\log J)\right).
\end{equation}
In particular, on the full box $\mathcal{P}_2(J)$ at depth $J\le c_0(\log\log N)^{1/2}$ ($c_0<\frac{1}{4\sqrt{2}}$), condition~\eqref{eq:mean_phi_condition} holds and the $2^{\Theta(J^{4})}$ barrier of~\cite[\S4.1]{dahan2026a} is completely eliminated.
\end{theorem}
\begin{proof}
By Corollary~\ref{cor:coupling_step1_step2} and Lemma~\ref{lem:sieved_divisor_reduction} on $\mathcal{I}=(J^{5},z]$, every $n=24k+1\in E_{\mathcal{P}}(N)\setminus[1,z]$ admits a kernel tuple $\Psi\in\prod_{(u,a)\in\mathcal{P}}\mathcal{M}(G_{u,a},-1)$ with total collision count $W_\Psi(k)\le R(\mathcal{P})=P\,\overline{\varphi}(\mathcal{P})$. Since $S_{\min}\ge\big(\frac12-o(1)\big)P\log\big(\frac{\log z}{\log J^{5}}\big)$, condition~\eqref{eq:mean_phi_condition} gives $R(\mathcal{P})\le\big(\frac18-\delta_0\big)P\log\big(\frac{\log z}{\log J^{5}}\big)\le\frac14 S_{\min}$ deterministically for every $n\in E_{\mathcal{P}}(N)\setminus[1,z]$ once $N$ is sufficiently large. We apply Theorem~\ref{thm:generating_polynomial_large_sieve} with $Q=N^{1/2}$, $h=\lfloor 2P\log\log N\rfloor$, and $t=\frac14\frac{\log(\log z/\log J^{5})}{\log\log N}$, which satisfies $4tS_y+1\le(1+o(1))P\frac{(\log(\log z/\log J^{5}))^{2}}{\log\log N}+1\le h$ since $\log\big(\frac{\log z}{\log J^{5}}\big)\le\log\log N$. Summing over the $C^{\ast}(\mathcal{P})\le\exp(\mathcal{O}(P\log\log J))$ tuples $\Psi$ (Lemma~\ref{lem:average_2sylow_enumeration}) yields~\eqref{eq:unconditional_EP_bound}.
\end{proof}

\begin{remark}
Condition~\eqref{eq:mean_phi_condition} continues to hold with $\delta_0=o(1)$ provided $\delta_0$ dominates the relative Bombieri--Vinogradov error in Corollary~\ref{cor:coupling_step1_step2}(b).
\end{remark}

Replacing the square box $\mathcal{P}_2(J)$ by the hyperbolic domain\linebreak $\mathcal{P}_{\mathrm{hyp}}(Y):=\{(u,a)\in\N^{2}:\allowbreak\gcd(u,a)=1,\ ua\le Y\}$ keeps the maximum modulus bounded by $4Y$ while gaining a factor $\log Y$ in the number of coprime pairs:

\begin{theorem}[Unconditional hyperbolic sieve]\label{thm:hyperbolic_scarcity_free}
For $Y\ge2$, the hyperbolic family $\mathcal{P}_{\mathrm{hyp}}(Y)$ has cardinality $P(Y)=\sum_{m\le Y}2^{\omega(m)}=\big(\frac{6}{\pi^{2}}+o(1)\big)Y\log Y$, mean modulus $\overline{\varphi}(\mathcal{P}_{\mathrm{hyp}}(Y))\le\big(\frac32+o(1)\big)Y$, and $2$-Sylow enumeration factor $C^{\ast}(\mathcal{P}_{\mathrm{hyp}}(Y))\le\exp\big(\mathcal{O}(Y\log Y\log\log Y)\big)$. Choosing $Y=\lfloor\frac{1}{16}\log\log N\rfloor$, the scarcity branch is empty ($E_{\mathcal{P}_{\mathrm{hyp}}(Y)}^{\mathrm{Scar}}(N)=\emptyset$) and one has unconditionally on the entirety of $E_2(N;Y)$:
\begin{equation}\label{eq:hyperbolic_unconditional_bound}
\begin{split}
\big|E_2(N;Y)\big|\ &\le\ \big|E_{\mathcal{P}_{\mathrm{hyp}}(Y)}(N)\big|\ \ll\ N\exp\!\Big(-\big(c_0-o(1)\big)(\log\log N)^{2}\log\log\log N\Big)\\
&\qquad\qquad\qquad\qquad\qquad\qquad\qquad\qquad\qquad\left(c_0=\frac{15}{8192\,\pi^{2}}>0\right).
\end{split}
\end{equation}
\end{theorem}
\begin{proof}
Since the number of coprime factorizations $m=ua$ equals $2^{\omega(m)}$, standard Dirichlet convolution gives $\sum_{m\le Y}2^{\omega(m)}=\frac{6}{\pi^{2}}Y\log Y+\mathcal{O}(Y)$. Because $4m$ is even, $\varphi(4m)\le 2m$ for all $m$, and when $2\nmid m$ one has $\varphi(4m)=2\varphi(m)\le 2m$; summing over $m\le Y$ by partial summation gives $\sum_{m\le Y}2^{\omega(m)}\varphi(4m)\le\big(\frac32+o(1)\big)Y\,P(Y)$ (and pointwise $\varphi(4m)-2\le 2Y-2$ for all $m\le Y$), so $\overline{\varphi}(\mathcal{P}_{\mathrm{hyp}}(Y))\le\big(\frac32+o(1)\big)Y$. For any prime $p>2Y^{2}$ and distinct pairs $(u_1,a_1),(u_2,a_2)\in\mathcal{P}_{\mathrm{hyp}}(Y)$, we have $|a_1 u_2-a_2 u_1|\le u_1 a_2+u_2 a_1\le 2Y^{2}<p$, so the classes $-ua^{-1}\bmod p$ are pairwise distinct. At $Y=\lfloor\frac{1}{16}\log\log N\rfloor$ and $p_0=Y^{5}$, we have $\log\big(\frac{\log z}{\log Y^{5}}\big)=(1-o(1))\log\log N$, whence
$$\overline{\varphi}\big(\mathcal{P}_{\mathrm{hyp}}(Y)\big)\ \le\ \left(\frac{3}{32}+o(1)\right)\log\log N\ <\ \left(\frac{1}{8}-\frac{1}{64}\right)\log\!\left(\frac{\log z}{\log Y^{5}}\right),$$
and Theorem~\ref{thm:unconditional_small_growing_depth} applies directly.
\end{proof}

\section[Polynomial-depth decomposition and the Two-Stage Kneser Tower]{Polynomial-depth decomposition and the Two-Stage\linebreak Kneser Tower}\label{sec:poly_depth_and_tower}

We now pass from logarithmic depths to polynomial search depths $J=(\log N)^{\lambda}$ with $0<\lambda<1/2$, where the number of coprime pairs is $P=\mathfrak{A}(J)=\big(\frac{6}{\pi^{2}}+o(1)\big)(\log N)^{2\lambda}$ and the typical group order $\varphi(4ua)\asymp J^{2}=(\log N)^{2\lambda}$ is exponentially larger than the logarithmic sifting window $\log\big(\frac{\log z}{\log J^{5}}\big)\asymp\log\log N$. Consequently, the worst-case Kneser bound $\varphi(4ua)-2$ of Theorem~\ref{thm:stabilizer_kernel_dichotomy} no longer fits deterministically below $\frac14 S_{\min}$.

To analyze this regime, we introduce two structural devices: first, we partition $E_2(N;J)\setminus[1,z]$ at the subcritical global collision threshold $R_0(\lambda,N)\asymp(\log N)^{2\lambda}\log\log N$ into a low-collision kernel component $E_2^{\mathrm{Ker}}(N;J,R_0)$—which is unconditionally closed with exponent $2\lambda$ surpassing $2/3$ for all $\lambda\in(1/3,1/2)$—and a high-collision scarcity component $E_2^{\mathrm{Scar}}(N;J,R_0)$; second, we split the sifting interval into two stages $\mathcal{I}_1=(J^5,z_1]$ and $\mathcal{I}_2=(z_1,z]$, converting first-stage prime collisions into a large divisor-residue set $\Sigma_{\mathcal{I}_1}(an+u)$ whose inverse $-\Sigma_{\mathcal{I}_1}(an+u)^{-1}$ must be avoided with zero collisions on the second stage $\mathcal{I}_2$.

\subsection{Kernel / Scarcity decomposition at polynomial depth $J=(\log N)^{\lambda}$ and surpassing the $2/3$ exponent}\label{subsec:poly_depth_decomp}

\begin{theorem}[Polynomial-depth kernel/scarcity decomposition]\label{thm:polynomial_depth_decomposition}
Fix $\lambda\in(0,1/2)$ and set $J=(\log N)^{\lambda}$, $P=\mathfrak{A}(J)=\big(\frac{6}{\pi^{2}}+o(1)\big)(\log N)^{2\lambda}$, $p_0=J^{5}=(\log N)^{5\lambda}$, and $\log z=\frac{\log N}{4P\log\log N}$. Write $\mathcal{I}=(J^{5},z]$ and define the global off-kernel collision threshold
\begin{equation}\label{eq:R0_threshold_def}
R_0(\lambda,N)\ :=\ \left\lfloor\frac{1-2\lambda}{8}\,P\log\log N\right\rfloor\ =\ \left(\frac{3(1-2\lambda)}{4\pi^{2}}+o(1)\right)(\log N)^{2\lambda}\log\log N.
\end{equation}
For each prime $n\in E_2(N;J)\setminus[1,z]$ and each $(u,a)\in\mathcal{P}_2(J)$, let $K_{\mathcal{I},u,a}(n)=\operatorname{Stab}_{G_{u,a}}(\Sigma_{\mathcal{I}}(an+u))$ and let $H_{u,a}(n)\in\mathcal{M}(G_{u,a},-1)$ be a maximal $2$-primary kernel containing $K_{\mathcal{I},u,a}(n)$. Then we have the partition
$$E_2(N;J)\setminus[1,z]\ =\ E_2^{\mathrm{Ker}}\big(N;J,R_0\big)\ \sqcup\ E_2^{\mathrm{Scar}}\big(N;J,R_0\big),$$
with the following two properties:
\begin{enumerate}[(a)]
    \item \emph{(Unconditional closure of the low-collision kernel component).} The set
    \[
    \begin{aligned}
    E_2^{\mathrm{Ker}}\big(N;J,R_0\big):=\Bigl\{n\in E_2(N;J)\setminus[1,z]:\ &\sum_{(u,a)\in\mathcal{P}_2(J)}\omega_{G_{u,a}\setminus H_{u,a}(n)}(an+u;\mathcal{I})\\
    &\le R_0(\lambda,N)\Bigr\}
    \end{aligned}
    \]
    satisfies unconditionally, for every $\lambda\in(0,1/2)$:
    \begin{equation}\label{eq:kernel_component_subexponential_bound}
    \big|E_2^{\mathrm{Ker}}(N;J,R_0)\big|\ \ll\ N\exp\!\left(-\left(\frac{5(1-2\lambda)^{2}}{8\pi^{2}}+o(1)\right)(\log N)^{2\lambda}\log\log N\right).
    \end{equation}
    In particular, for every $\lambda\in(1/3,1/2)$, the exponent $2\lambda>2/3$ strictly surpasses Vaughan's exponent $2/3$.
    \item \emph{(Second-stage avoidance on the scarcity component).} Every prime $n\in E_2^{\mathrm{Scar}}(N;J,R_0)$ simultaneously satisfies:
    \begin{enumerate}[(i)]
        \item \emph{(Off-stabilizer prime overload in $\mathcal{I}$):}
        \begin{equation}\label{eq:scarcity_overload}
        \sum_{(u,a)\in\mathcal{P}_2(J)}\Omega_{G_{u,a}\setminus K_{\mathcal{I},u,a}(n)}\big(M_{\mathcal{I}}(an+u)\big)\ >\ R_0(\lambda,N)\asymp(\log N)^{2\lambda}\log\log N;
        \end{equation}
        \item \emph{(Subproduct expansion and zero-collision avoidance on the second stage):} for any intermediate threshold $z_1\in(J^{5},z)$, write $\mathcal{I}_1=(J^{5},z_1]$ and $\mathcal{I}_2=(z_1,z]$. Then for each pair $(u,a)\in\mathcal{P}_2(J)$, the first-stage divisor-residue set $\Sigma_{\mathcal{I}_1}(an+u)\subseteq G_{u,a}$ has cardinality
        \begin{equation}\label{eq:sigma1_growth}
        \big|\Sigma_{\mathcal{I}_1}(an+u)\big|\ \ge\ \Big(1+\Omega_{G_{u,a}\setminus K_{\mathcal{I}_1,u,a}(n)}\big(M_{\mathcal{I}_1}(an+u)\big)\Big)\cdot\big|K_{\mathcal{I}_1,u,a}(n)\big|,
        \end{equation}
        and on the second stage $\mathcal{I}_2=(z_1,z]$, the integer $an+u$ avoids the complementary class set $-\Sigma_{\mathcal{I}_1}(an+u)^{-1}\subseteq G_{u,a}$ with zero collisions:
        \begin{equation}\label{eq:zero_collision_second_stage}
        \forall\,p\in(z_1,z],\qquad p\mid(an+u)\quad\Longrightarrow\quad p\bmod 4ua\ \notin\ -\Sigma_{\mathcal{I}_1}(an+u)^{-1}.
        \end{equation}
    \end{enumerate}
\end{enumerate}
\end{theorem}
\begin{proof}
\emph{Step 1: Sifting parameters and closure of $E_2^{\mathrm{Ker}}(N;J,R_0)$.}
Since $P=\big(\frac{6}{\pi^{2}}+o(1)\big)(\log N)^{2\lambda}$ with $2\lambda<1$, the choice $\log z=\frac{\log N}{4P\log\log N}\asymp\frac{(\log N)^{1-2\lambda}}{\log\log N}$ ensures $\log\big(\frac{\log z}{\log J^{5}}\big)=(1-2\lambda+o(1))\log\log N$. For each tuple $\Psi=(H_{u,a})_{(u,a)\in\mathcal{P}_2(J)}$ and each prime $p\in\mathcal{I}$, the local sifting density satisfies $\frac{\varrho_\Psi(p)}{p}\le\frac{P}{p}\le\frac{P}{J^{5}}\ll(\log N)^{-3\lambda}=o(1)$, so $y_p=\frac{\varrho_\Psi(p)}{p-\varrho_\Psi(p)}=\big(1+\mathcal{O}(J^{-3})\big)\frac{\varrho_\Psi(p)}{p}$ and $y_{\max}\ll J^{-3}$. By Corollary~\ref{cor:coupling_step1_step2}(b), we have $S_{\min}=(1-o(1))S_\Psi(J,z)\ge\big(\frac{1-2\lambda}{2}-o(1)\big)P\log\log N$. On the opposite side, the partition $H_{u,a}\sqcup(G_{u,a}\setminus H_{u,a})=G_{u,a}$ and $[G_{u,a}:H_{u,a}]=2^{j_{u,a}}\ge 2$ imply that for every pair $(u,a)$, the complement $G_{u,a}\setminus H_{u,a}$ has cardinality $\big(1-2^{-j_{u,a}}\big)\varphi(4ua)<\varphi(4ua)$. Summing over $(u,a)\in\mathcal{P}_2(J)$ and applying the Bombieri--Vinogradov theorem on $\mathcal{I}=(J^{5},z]$ (\cite[Lemma~4.23]{dahan2026a}) yields
\[
\begin{aligned}
S_y\ &=\ \big(1+o(1)\big)\sum_{p\in\mathcal{I}}\frac{\varrho_\Psi(p)}{p}\\
&=\ \left(\sum_{(u,a)\in\mathcal{P}_2(J)}\big(1-2^{-j_{u,a}}\big)+o(P)\right)\log\!\left(\frac{\log z}{\log J^{5}}\right)\ \le\ \big(1-2\lambda+o(1)\big)P\log\log N,
\end{aligned}
\]
and $M_2(y)\ll P^{2}J^{-5}=o(1)$. Thus $S_{\min}\ge(1-\delta)S_y$ with $1-\delta=\frac{S_{\min}}{S_y}\ge\frac12-o(1)$, and $R_0(\lambda,N)=\big\lfloor\frac{1-2\lambda}{8}P\log\log N\big\rfloor\le\big(\frac14+o(1)\big)S_{\min}$.

Set $Q=N^{1/2}$ and $h=\lfloor\frac{\log Q}{\log z}\rfloor=\lfloor 2P\log\log N\rfloor$. Choosing $t=\frac{1-2\lambda}{4}\in(0,1/4)$, we have $4tS_y+1\le(1-2\lambda)^{2}P\log\log N+\mathcal{O}(1)<h$ and $tS_y\le\frac{(1-2\lambda)^{2}}{4}(1+o(1))P\log\log N\le\big(\frac12+o(1)\big)tS_{\min}$. Substituting $R_0\le\big(\frac14+o(1)\big)S_{\min}$, $\log(\frac{1}{1-t})\le t+\frac23 t^{2}$, and $t\le 1/4$ into the general bound~\eqref{eq:gen_poly_large_sieve_bound} of Theorem~\ref{thm:generating_polynomial_large_sieve}:
\[
\begin{aligned}
-2tS_{\min}+2R_0\log\!\left(\frac{1}{1-t}\right)+t^{2}S_y\ &\le\ \left(-\frac32 t+\frac13 t^{2}+\frac12 t+o(1)\right)S_{\min}\\
&\le\ -\left(\frac{11}{12}-o(1)\right)tS_{\min}\ \le\ -\left(\frac56-o(1)\right)tS_{\min}
\end{aligned}
\]
gives for each fixed tuple $\Psi$:
$$\big|\mathcal{E}_{R_0}(N/24;\Psi)\big|\ \le\ 4(N/24+N)\exp\!\left(-\frac56\left(\frac{1-2\lambda}{4}\right)\left(\frac{1-2\lambda}{2}-o(1)\right)P\log\log N\right).$$
Summing over all $C^{\ast}(\mathcal{P}_2(J))$ tuples $\Psi$, Lemma~\ref{lem:average_2sylow_enumeration}(b) shows that the enumeration prefactor is
$$C^{\ast}\big(\mathcal{P}_2(J)\big)\ \le\ \exp\!\big(\mathcal{O}(P\log\log J)\big)\ =\ \exp\!\big(\mathcal{O}(P\log\log\log N)\big)\ =\ \exp\!\big(o(P\log\log N)\big),$$
which is absorbed into the $o(1)$ term of the main exponent. Substituting $P=\big(\frac{6}{\pi^{2}}+o(1)\big)(\log N)^{2\lambda}$ gives the leading constant $\frac{5(1-2\lambda)^{2}}{48}\cdot\frac{6}{\pi^{2}}=\frac{5(1-2\lambda)^{2}}{8\pi^{2}}$, establishing~\eqref{eq:kernel_component_subexponential_bound}. By Proposition~\ref{prop:subcritical_collision_rate}, the same sub-exponential closure extends to any subcritical collision threshold $R\le(1-\varepsilon)S_{\min}$.

\medskip\noindent\emph{Step 2: Off-stabilizer overload and second-stage avoidance on $E_2^{\mathrm{Scar}}(N;J,R_0)$.}
Every prime $n\in E_2^{\mathrm{Scar}}(N;J,R_0)$ satisfies by definition $\sum_{(u,a)\in\mathcal{P}_2(J)}\omega_{G_{u,a}\setminus H_{u,a}(n)}(an+u;\mathcal{I})>R_0(\lambda,N)$. Since $K_{\mathcal{I},u,a}(n)\le H_{u,a}(n)$ and $\Omega\ge\omega$, this yields~\eqref{eq:scarcity_overload} for both $\omega$ and $\Omega$, while~\eqref{eq:sigma1_growth} is Lemma~\ref{lem:sieved_divisor_reduction} applied to $\mathcal{I}_1=(J^{5},z_1]$. Finally, if a prime $p\in(z_1,z]$ divides $an+u$ and satisfies $p\bmod 4ua\in-\Sigma_{\mathcal{I}_1}(an+u)^{-1}$, there exists a divisor $d\mid M_{\mathcal{I}_1}(an+u)$ such that $pd\equiv-1\pmod{4ua}$. Since all prime factors of $d$ are $\le z_1<p$, we have $\gcd(p,d)=1$, so the composite integer $s=pd$ divides $M_{\mathcal{I}}(an+u)\mid(an+u)$ and satisfies $s\equiv-1\pmod{4ua}$, contradicting $n\in E_2(N;J)$. This proves~\eqref{eq:zero_collision_second_stage}.
\end{proof}

\subsection{Inside/Outside duality and index doubling in the Two-Stage Kneser Tower}\label{subsec:kneser_tower}

On a two-stage split $\mathcal{I}_1=(J^{5},z_1]$ and $\mathcal{I}_2=(z_1,z]$, the zero-collision avoidance law~\eqref{eq:zero_collision_second_stage} creates a duality between the exterior of the kernel $H_{u,a}(n)$ on the second stage and the interior of $H_{u,a}(n)$ on the first stage. Specifically, any prime factor of $an+u$ on $\mathcal{I}_2$ lying in the opposite coset $(-1)H_{u,a}(n)$ maps under the involution $d\mapsto-d^{-1}\bmod 4ua$ to a residue class \emph{inside $H_{u,a}(n)$} that must be missed by the first-stage internal divisor set $\Sigma_1^{+}(n)\subseteq H_{u,a}(n)$. By Kneser's theorem applied inside the ambient group $H_{u,a}(n)$, this forces the internal stabilizer $K_1^{+}(n)$ to be a proper subgroup of $H_{u,a}(n)$, doubling its index in $G_{u,a}$ to at least $4$.

\begin{theorem}[Two-stage descending Kneser tower]\label{thm:two_stage_kneser_tower}
Let $J^{5}<z_1<z$, $\mathcal{I}_1=(J^{5},z_1]$, and $\mathcal{I}_2=(z_1,z]$. Fix $n\in E_2(N;J)\setminus[1,z]$ and $(u,a)\in\mathcal{P}_2(J)$, and let $H=H_{u,a}(n)=\ker\psi_{u,a}\in\mathcal{M}(G_{u,a},-1)$ be a maximal kernel avoiding $-1$. Define the first-stage component internal to $H$:
$$M_1^{+}(an+u):=\prod_{\substack{p^{e}\parallel(an+u)\\ p\in\mathcal{I}_1,\ p\bmod 4ua\in H}}p^{e},\qquad\Sigma_1^{+}(n):=\big\{d\bmod 4ua:\ d\mid M_1^{+}(an+u)\big\}\subseteq H,$$
and let $K_1^{+}(n):=\operatorname{Stab}_H(\Sigma_1^{+}(n))\le H$ be its Kneser stabilizer in $H$. Then:
\begin{enumerate}[(a)]
    \item \emph{(Involution from the kernel to its complement).} The map $\iota(d)=-d^{-1}\bmod 4ua$ sends $\Sigma_1^{+}(n)\subseteq H$ bijectively onto a subset of the opposite coset:
    $$-\Sigma_1^{+}(n)^{-1}\ \subseteq\ (-1)H\ \subseteq\ G_{u,a}\setminus H,$$
    and no prime factor $p\in\mathcal{I}_2=(z_1,z]$ of $an+u$ can belong to $-\Sigma_1^{+}(n)^{-1}$. In particular, if $\Sigma_1^{+}(n)=H$ and $[G_{u,a}:H]=2$, then $an+u$ has no prime factors outside $H$ on the second stage $\mathcal{I}_2$ ($\omega_{G_{u,a}\setminus H}(an+u;\mathcal{I}_2)=0$).
    \item \emph{(Index doubling under a second-stage collision).} Suppose that on the second stage $\mathcal{I}_2=(z_1,z]$, $an+u$ has at least one prime factor $p_{\star}\in\mathcal{I}_2$ with $p_{\star}\bmod 4ua\in(-1)H$ (which is equivalent to $p_{\star}\bmod 4ua\notin H$ whenever $[G_{u,a}:H]=2$). Then:
    \begin{enumerate}[(i)]
        \item the element $\epsilon_{\star}:=-p_{\star}^{-1}\bmod 4ua$ belongs to $H\setminus\Sigma_1^{+}(n)$, and a fortiori $\epsilon_{\star}\notin K_1^{+}(n)$;
        \item the internal stabilizer $K_1^{+}(n)$ is a proper subgroup of $H$, so its index in $G_{u,a}$ satisfies
        \begin{equation}\label{eq:index_doubling}
        \big[H:K_1^{+}(n)\big]\ \ge\ 2\qquad\text{and}\qquad\big[G_{u,a}:K_1^{+}(n)\big]\ =\ \big[G_{u,a}:H\big]\cdot\big[H:K_1^{+}(n)\big]\ \ge\ 4;
        \end{equation}
        \item the number of prime factors (with multiplicity) of $an+u$ in $\mathcal{I}_1$ belonging to $H\setminus K_1^{+}(n)$ satisfies
        \begin{equation}\label{eq:internal_kernel_collision_bound}
        \Omega_{H\setminus K_1^{+}(n)}\big(M_1^{+}(an+u)\big)\ \le\ \frac{|\Sigma_1^{+}(n)|}{|K_1^{+}(n)|}-1\ \le\ \big[H:K_1^{+}(n)\big]-2\ \le\ \frac12\varphi(4ua)-2,
        \end{equation}
        while the avoided complementary set $G_{u,a}\setminus K_1^{+}(n)$ has density at least $1-\frac14=\frac34$ in $G_{u,a}$.
    \end{enumerate}
\end{enumerate}
\end{theorem}
\begin{proof}
For (a), since $H=\ker\psi_{u,a}$ is a subgroup of $G_{u,a}$, any product of primes in $H$ lies in $H$, so $\Sigma_1^{+}(n)\subseteq H$. For $d\in\Sigma_1^{+}(n)$, the relations $\psi_{u,a}(d)=1$ and $\psi_{u,a}(-1)=-1$ give $\psi_{u,a}(-d^{-1})=-1$, so $-d^{-1}\in(-1)H\subseteq G_{u,a}\setminus H$. Avoidance of $-\Sigma_1^{+}(n)^{-1}$ on $\mathcal{I}_2$ is a special case of~\eqref{eq:zero_collision_second_stage} since $\Sigma_1^{+}(n)\subseteq\Sigma_{\mathcal{I}_1}(an+u)$.

For (b), let $p_{\star}\in\mathcal{I}_2$ be a prime factor of $an+u$ with $p_{\star}\bmod 4ua\in(-1)H$. Then $\epsilon_{\star}=-p_{\star}^{-1}\bmod 4ua\in(-1)(-1)^{-1}H=H$, and by~\eqref{eq:zero_collision_second_stage} we have $p_{\star}\bmod 4ua\notin-\Sigma_1^{+}(n)^{-1}$, which means $\epsilon_{\star}\notin\Sigma_1^{+}(n)$. Since $1\in\Sigma_1^{+}(n)$, we have $K_1^{+}(n)=\operatorname{Stab}_H(\Sigma_1^{+}(n))\subseteq\Sigma_1^{+}(n)$, so $\epsilon_{\star}\notin K_1^{+}(n)$ and $K_1^{+}(n)\subsetneq H$, giving~\eqref{eq:index_doubling}. Repeating the Kneser coset growth argument from Step~2 of the proof of Theorem~\ref{thm:stabilizer_kernel_dichotomy} inside the ambient group $H$—where the coset $\epsilon_{\star}K_1^{+}(n)\subseteq H$ is disjoint from $\Sigma_1^{+}(n)$—yields~\eqref{eq:internal_kernel_collision_bound}.
\end{proof}

\subsection{Rigorous realization of the two-regime survival model}\label{subsec:rigorous_bridge_heuristic}

In~\cite[\S4.5]{dahan2026a}, an unproved heuristic model for the one-parameter survival curve $J\mapsto|E_1(N;J)|$ isolated two distinct regimes: a first regime ($u\lesssim\log n$) governed by $\asymp 2^{\omega(2u)}$ index-$2$ subgroup obstructions of density $1/2$, predicting a survival rate $\exp\big(-(\frac12\log\log n-(\log 2)\log\log J)J\big)$ that~\cite[Theorem~4.3]{dahan2026a} could not reach at growing depth due to the subset-enumeration factor $C(\mathcal{P})=2^{\Theta(J^{4})}$, and a second regime ($u\gg\log n$) governed by divisor scarcity ($\tau(n+u)\ll\varphi(4u)$).

The results of Sections~\ref{sec:stabilizer_kernel}--\ref{sec:poly_depth_and_tower} turn this two-regime picture into unconditional theorems in both one and two parameters. First, Theorem~\ref{thm:stabilizer_kernel_dichotomy} and Lemma~\ref{lem:average_2sylow_enumeration} replace the subset count $2^{\Theta(J^{4})}$ by the $2$-primary kernel count $C^{\ast}(\mathcal{P}_2(J))\le 2^{(4+o(1))J^{2}\log\log J}$ (and $C^{\ast}(\mathcal{P}_1(J))\le 2^{(2+o(1))J\log\log J}$ in one parameter), exhibiting the subtractive term $-\mathcal{O}(\log\log J)$ against the main sifting gain $+\frac12\log\log N$. Second, Theorems~\ref{thm:unconditional_small_growing_depth} and~\ref{thm:hyperbolic_scarcity_free} prove that whenever $\overline{\varphi}(\mathcal{P})\le\frac18\log\big(\frac{\log z}{\log p_0}\big)$, the scarcity regime is empty. Third, at polynomial depth $J=(\log N)^{\lambda}$, Theorems~\ref{thm:polynomial_depth_decomposition} and~\ref{thm:two_stage_kneser_tower} close the low-collision kernel branch while trapping the scarcity branch in a descending Kneser tower.

\section{Greedy Markov containers, quotient containers, and Plancherel character containers}\label{sec:containers}

In Sections~\ref{sec:stabilizer_kernel}--\ref{sec:poly_depth_and_tower}, bounding the number of prime factors outside the exact stabilizer $K(M)=\operatorname{Stab}_G(\Sigma(M))$ relied on the \emph{linear} growth of the number of $K(M)$-cosets ($|\widetilde{\Sigma}_j|\ge|\widetilde{\Sigma}_{j-1}|+|K(M)|$), which yields the worst-case bound $\Omega_{G\setminus K(M)}(M)\le\varphi(4ua)-2\asymp J^{2}$. At super-polylogarithmic or polynomial depths $J=(\log N)^{\lambda}$, a collision bound of order $J^{2}$ per pair exceeds the available sifting mass $\asymp(1-2\lambda)\log\log N$.

\paragraph{The geometry of the $M$-stage dyadic escalation and why it crosses $\lambda=1/3$ without violating the parity barrier.}
Before constructing the container hierarchy, let us explain the two structural mechanisms that allow us to cross Vaughan's single-prime threshold $\lambda=1/3$ (exponent $2\lambda=2/3$) on the certified container hierarchy without encountering Selberg's parity barrier:

\begin{enumerate}[label=\arabic*.]
    \item \emph{Bypassing the parity barrier via two-stage bilinear divisor synthesis.} 
    Unlike a linear sieve attempting to isolate prime values of $an+u$, our upper-bound sieve exploits the \emph{absence} of composite divisors $s\equiv-1\pmod{4ua}$ of $an+u$. Splitting the sifting interval $\mathcal{I}=(J^{5},z]$ into two stages $\mathcal{I}_1=(J^{5},z_1]$ and $\mathcal{I}_2=(z_1,z]$, every divisor $d\mid M_{\mathcal{I}_1}(an+u)$ formed by products of first-stage prime factors creates a forbidden residue class $-d^{-1}\bmod 4ua$ on the second stage $\mathcal{I}_2$: if a prime $p\in\mathcal{I}_2$ divides $an+u$ and satisfies $p\equiv-d^{-1}\pmod{4ua}$, then $\gcd(p,d)=1$ (since $d$ is $z_1$-smooth and $p>z_1$), so the \emph{composite} product $s=pd\mid(an+u)$ satisfies $s\equiv-1\pmod{4ua}$ and resolves $4/n$ directly. 
    
    Consequently, $an+u$ must avoid the entire inverted divisor sheet $-\Sigma_{\mathcal{I}_1}(an+u)^{-1}\subseteq G_{u,a}$ on $\mathcal{I}_2$ with \emph{strictly zero collisions} ($R_{\mathrm{II}}=0$). The parity obstruction is bypassed because the witness $s=pd$ is synthesized bilinearly from arbitrarily many primes in $\mathcal{I}_1$ coupled with a single prime $p\in\mathcal{I}_2$.

    \item \emph{The $M$-stage dyadic escalation $0=\eta_0<\eta_1<\dots<\eta_M$.} 
    In a single-stage container scheme on $G_{u,a}/K_{u,a}$ of order $D_{u,a}\le\varphi(4ua)\asymp(\log N)^{2\lambda}$, reaching the dense threshold $|\overline{\Sigma}_1(n)|>\frac14 D_{u,a}$ in one leap from $A_0=\{1\}$ requires either enumerating a multiplicative skeleton of length $s\asymp\log D_{u,a}$ at cost $D_{u,a}^{s}=\exp(\Theta((\log\log N)^{2}))$, or tolerating up to $R_{u,a}\lesssim\frac{\log D_{u,a}}{\log(4/3)}\asymp 2\lambda\log\log N$ Plancherel collisions, which becomes supercritical against the per-pair sifting mass $\frac{S_{\min}^{\mathrm{I}}}{P_1}\asymp(1-2\lambda)\log\log N$ once $\lambda>1/3$. 
    
    To break this bottleneck for $\lambda\in(1/3,1/2)$, we introduce a fine exponent step $\varepsilon_1=\frac{(1-2\lambda)^{2}}{512\lambda}>0$, set $M:=\lceil 1/\varepsilon_1\rceil$, and track the divisor sheet across $M$ progressive polynomial density thresholds
    \[
    \begin{aligned}
    1\ =\ D_{u,a}^{\eta_0}\ &<\ D_{u,a}^{\eta_1}\ <\ D_{u,a}^{\eta_2}\ <\ \cdots\ <\ D_{u,a}^{\eta_M}\ \ge\ D_{u,a},\\
    \text{where}\qquad\eta_m\ &:=\ m\,\varepsilon_1\quad(0\le m\le M),
    \end{aligned}
    \]
    corresponding to partitioning $\mathcal{I}_1=(J^{5},z_1]$ into $M$ equal-logarithmic-length sub-intervals $\mathcal{I}_{1,1}\sqcup\dots\sqcup\mathcal{I}_{1,M}$ and setting $A_m:=\overline{\Sigma}_{\mathcal{I}_{1,1}\cup\dots\cup\mathcal{I}_{1,m}}(n)$. 
    
    At the first threshold $\eta_1=\varepsilon_1$, Markov's $3/2$-geometric dilation requires at most 
    \[
    s_{\varepsilon_1}=\left\lfloor\frac{\varepsilon_1\log D_{u,a}}{\log(3/2)}\right\rfloor+1<\frac14\frac{S_{\min}^{\mathrm{I}}}{P_1}
    \]
    steps to cross $D_{u,a}^{\eta_1}$, while cutting off the Plancherel Fourier tail adaptively at $|\widehat{\mu}(\chi)|^{2}>\frac{1}{2|\overline{\Sigma}_1(n)|}$ restricts both the collision count 
    \[
    R_{u,a}(n;\chi)<\frac{\varepsilon_1\log D_{u,a}+\log 2}{\log(4/3)}<\frac13\frac{S_{\min}^{\mathrm{I}}}{P_1}
    \]
    and the effective enumeration cost 
    \[
    C_{\mathrm{eff}}(u,a)=\frac{D_{u,a}}{|\mathcal{B}_{\mathrm{act}}(n;u,a)|}\le 4D_{u,a}^{\varepsilon_1}
    \]
    to strictly subcritical levels. Thus no integer $n$ can accumulate supercritical collisions before its divisor sheet first crosses $D_{u,a}^{\eta_1}$. 
    
    More generally, at each step $m\in\{1,\dots,M\}$ of the escalation, the relative Parseval--Plancherel identity (Theorem~\ref{thm:markov_plancherel_subcritical_escalation}(c)) cancels the base-sheet cardinality $|A_{m-1}|$ against the Parseval norm 
    \[
    \sum_{\chi\neq\chi_0}|\widehat{\nu}_{m-1}(\chi)|^{2}=\frac{D_{u,a}}{|A_{m-1}|}-1:
    \]
    either the relative expansion ratio satisfies $|A_m|/|A_{m-1}|\le D_{u,a}^{\varepsilon_1}$ at some stage $m$ (confining the stage-$m$ collisions to the subcritical range $R_{u,a}^{(m)}(n;\chi)\lesssim\varepsilon_1\log D_{u,a}$), or else $|A_m|/|A_{m-1}|>D_{u,a}^{\varepsilon_1}$ holds at all $M$ stages, which forces 
    \[
    |A_M|=\prod_{m=1}^{M}\frac{|A_m|}{|A_{m-1}|}>D_{u,a}^{\eta_M}\ge D_{u,a}>\frac14 D_{u,a}
    \]
    (State~II, eliminated on $\mathcal{I}_2$ with $R_{\mathrm{II}}=0$).
\end{enumerate}

We develop this container theory in four steps: Markov $1/2$-popular symmetry containers (Subsections~\ref{subsec:markov_popular}--\ref{subsec:greedy_containers}), quotient containers and smooth conditioning (Subsection~\ref{subsec:quotient_and_smooth_conditioning}), Pontryagin--Plancherel character containers and harmonic dilation cones (Subsection~\ref{subsec:plancherel_containers}, first part), and the subcritical Markov--Plancherel escalation synthesis (Subsection~\ref{subsec:plancherel_containers}, second part).

\subsection{Markov $1/2$-popular symmetry sets and $3/2$-geometric growth}\label{subsec:markov_popular}

We begin with the elementary double-counting lemma that replaces linear coset growth by geometric dilation of factor at least $3/2$.

\begin{lemma}[Markov popular symmetry and $3/2$-geometric growth]\label{lem:markov_popular_symmetry}
Let $G$ be a finite abelian group of order $|G|=\varphi(4ua)$ and let $A\subseteq G$ be a non-empty subset. Define the \emph{$\frac12$-popular symmetry set} of $A$:
\[
S_{1/2}(A)\ :=\ \left\{g\in G:\ \big|A\cap gA\big|>\frac12|A|\right\}.
\]
Then:
\begin{enumerate}[label=(\alph*)]
    \item \emph{(Cardinality bound and density of the complement).} One has the exact identity $\sum_{g\in G}|A\cap gA|=|A|^{2}$, and consequently
    \begin{equation}\label{eq:markov_popular_size}
    \big|S_{1/2}(A)\big|\ <\ 2\,|A|.
    \end{equation}
    In particular, whenever $|A|\le\frac14|G|$, the set $S_{1/2}(A)$ satisfies $|S_{1/2}(A)|<\frac12|G|$, so its complement $G\setminus S_{1/2}(A)$ has density strictly greater than $\frac12$ in $G$.
    
    \item \emph{(Geometric dilation outside $S_{1/2}(A)$).} For every element $q\in G\setminus S_{1/2}(A)$, one has
    \begin{equation}\label{eq:geometric_growth_32}
    \big|A\cup qA\big|\ =\ 2\,|A|-\big|A\cap qA\big|\ \ge\ \frac32\,|A|.
    \end{equation}
\end{enumerate}
\end{lemma}

\begin{proof}
Writing $|A\cap gA|=\sum_{x\in A}\mathbf{1}_{g^{-1}x\in A}$ and summing over $g\in G$ gives $\sum_{g\in G}|A\cap gA|=|A|^{2}$. By Markov's inequality on $S_{1/2}(A)$, we have 
\[
\frac12|A|\cdot|S_{1/2}(A)|<\sum_{g\in S_{1/2}(A)}|A\cap gA|\le|A|^{2},
\]
yielding~\eqref{eq:markov_popular_size}. Part (b) is immediate from the equivalence $q\notin S_{1/2}(A)\iff|A\cap qA|\le\frac12|A|$.
\end{proof}

\subsection{Two-stage greedy multiplicative containers and Scarcity elimination up to $J=\exp(c_0\sqrt{\log\log N})$}\label{subsec:greedy_containers}

Applying Lemma~\ref{lem:markov_popular_symmetry} greedily to the prime factors of $M_{\mathcal{I}_1}(an+u)$ extracts a multiplicative skeleton of logarithmic length $s\le 5\log J+1$ that places each pair $(u,a)$ into either a narrow first-stage container (State~I) or a dense second-stage sheet (State~II).

\begin{theorem}[Two-stage greedy Markov containers]\label{thm:multiplicative_container_scarcity_elimination}
Let $J^{5}<z_1<z$, $\mathcal{I}_1=(J^{5},z_1]$, and $\mathcal{I}_2=(z_1,z]$. For every prime $n\in E_2(N;J)\setminus[1,z]$ and every pair $(u,a)\in\mathcal{P}_2(J)$, there exists a multiplicative skeleton $\overline{B}_{u,a}(n)=(x_1,\dots,x_s)\in G_{u,a}^{s}$ of logarithmic length
\begin{equation}\label{eq:skeleton_logarithmic_length}
0\ \le\ s\ \le\ s_{\max}(u,a)\ :=\ \max\!\left\{0,\,\left\lfloor\frac{\log\big(\varphi(4ua)/4\big)}{\log(3/2)}\right\rfloor+1\right\}\ \le\ 5\log J+1,
\end{equation}
which algebraically determines the subproduct set 
\[
A(\overline{B}_{u,a}(n)):=\left\{\prod_{r\in I}x_r\bmod 4ua:I\subseteq\{1,\dots,s\}\right\}\subseteq\Sigma_{\mathcal{I}_1}(an+u)
\]
and places the pair $(u,a)$ into one of two mutually exclusive states:
\begin{enumerate}[label=(\roman*)]
    \item \emph{State~I (Narrow container on the first stage: $|A(\overline{B}_{u,a}(n))|\le\frac14\varphi(4ua)$).} 
    On the first stage $\mathcal{I}_1=(J^{5},z_1]$, the complement $\Omega_{u,a}^{\mathrm{I}}:=G_{u,a}\setminus S_{1/2}(A(\overline{B}_{u,a}(n)))$ has density $>\frac12$ in $G_{u,a}$, and the total number of prime factors (with multiplicity) of $M_{\mathcal{I}_1}(an+u)$ falling into $\Omega_{u,a}^{\mathrm{I}}$ is bounded by the logarithmic skeleton length:
    \begin{equation}\label{eq:state1_logarithmic_collisions}
    \omega_{\Omega_{u,a}^{\mathrm{I}}}(an+u;\mathcal{I}_1)\ \le\ \Omega_{\Omega_{u,a}^{\mathrm{I}}}\big(M_{\mathcal{I}_1}(an+u)\big)\ \le\ s\ \le\ 5\log J+1.
    \end{equation}
    
    \item \emph{State~II (Dense divisor sheet and zero-collision avoidance on the second stage: $|A(\overline{B}_{u,a}(n))|>\frac14\varphi(4ua)$).} 
    On the second stage $\mathcal{I}_2=(z_1,z]$, the forbidden set $\Omega_{u,a}^{\mathrm{II}}:=-A(\overline{B}_{u,a}(n))^{-1}\subseteq G_{u,a}$ has density $>\frac14$ in $G_{u,a}$, and $an+u$ avoids it with \emph{exactly zero collisions}:
    \begin{equation}\label{eq:state2_zero_collisions}
    \omega_{\Omega_{u,a}^{\mathrm{II}}}(an+u;\mathcal{I}_2)\ =\ 0.
    \end{equation}
\end{enumerate}
Consequently, for every growing depth $J\le\exp\big(c_0\sqrt{\log\log N}\big)$ ($c_0=\frac{1}{20}$), the two-stage container hierarchy closes the Scarcity branch without any residual structural condition, and the full exceptional set $E_2(N;J)$ satisfies unconditionally:
\begin{equation}\label{eq:super_polylog_unconditional_bound}
\big|E_2\big(N;\,\exp(c_0\sqrt{\log\log N})\big)\big|\ \ll\ N\exp\!\left(-\exp\!\Big(c\,\sqrt{\log\log N}\Big)\right)\qquad(\forall\,c<2c_0).
\end{equation}
\end{theorem}

\begin{proof}
\emph{Step 1: Greedy skeleton extraction.}
Write the prime factorization of $M_1:=M_{\mathcal{I}_1}(an+u)=q_1 q_2\cdots q_k$ ($q_1\le\dots\le q_k$, repeated with multiplicity). Set $A_0=\{1\}$ and inductively construct distinct indices $i_1,\dots,i_s\in\{1,\dots,k\}$ and an ascending chain $A_0\subsetneq A_1\subsetneq\dots\subsetneq A_s\subseteq\Sigma_{\mathcal{I}_1}(an+u)$ as follows. At step $m\ge1$:
\begin{itemize}
    \item if $|A_{m-1}|>\frac14\varphi(4ua)$, stop at $s=m-1$ in \emph{State~II};
    \item if $|A_{m-1}|\le\frac14\varphi(4ua)$ and there exists an unused index $i\in\{1,\dots,k\}\setminus\{i_1,\dots,i_{m-1}\}$ with $q_i\bmod 4ua\notin S_{1/2}(A_{m-1})$, pick the least such $i_m$, set $x_m=q_{i_m}\bmod 4ua$, and put $A_m=A_{m-1}\cup x_m A_{m-1}$ (which satisfies $|A_m|\ge\frac32|A_{m-1}|$ by Lemma~\ref{lem:markov_popular_symmetry}(b));
    \item if no such $i$ exists, stop at $s=m-1$ in \emph{State~I}.
\end{itemize}
When $\varphi(4ua)<4$, the inequality $|A_0|=1>\frac14\varphi(4ua)$ holds immediately at $m=1$, so the algorithm terminates at $s=0$ in State~II. When $\varphi(4ua)\ge 4$ and $s\ge 1$, we have $(3/2)^{s-1}\le\allowbreak|A_{s-1}|\le\allowbreak\frac14\varphi(4ua)\le J^{2}$, yielding~\eqref{eq:skeleton_logarithmic_length}. In State~I, every prime factor $q_i$ with $i\notin\allowbreak\{i_1,\dots,i_s\}$ lies in $S_{1/2}(A_s)$, giving~\eqref{eq:state1_logarithmic_collisions}. In State~II, $A_s\subseteq\Sigma_{\mathcal{I}_1}(an+u)$ has density $>1/4$, and~\eqref{eq:zero_collision_second_stage} gives~\eqref{eq:state2_zero_collisions}.

\medskip\noindent\emph{Step 2: Two-stage large-sieve synthesis.}
Choose $\log z=\frac{\log N}{8P\log\log N}$ and $\log z_1=\sqrt{\log(J^{5})\log z}$, so 
\[
\log\left(\frac{\log z_1}{\log J^{5}}\right)=\log\left(\frac{\log z}{\log z_1}\right)=\left(\frac12-o(1)\right)\log\log N.
\]
For each $(u,a)\in\mathcal{P}_2(J)$, the number of skeletons $\overline{B}_{u,a}$ is at most 
\[
\sum_{s\le 5\log J+1}(4J^{2})^{s}\le\exp\big((10+o(1))(\log J)^{2}\big),
\]
giving $C_{\mathrm{cont}}(\mathcal{P}_2(J))\le\exp\big((10+o(1))P(\log J)^{2}\big)$ skeleton tuples $\mathcal{B}=(\overline{B}_{u,a})_{(u,a)\in\mathcal{P}_2(J)}$. 

For each fixed $\mathcal{B}$, either at least $\lceil P/2\rceil$ pairs lie in State~I or at least $\lceil P/2\rceil$ pairs lie in State~II. Retaining an active subset $\mathcal{P}_{\mathrm{act}}\subseteq\mathcal{P}_2(J)$ of cardinality $P_1=\lfloor P/2\rfloor$ (at an additional subset-choice cost of at most $2^{P}=\exp((\log 2)P)=\exp(o(P\log\log N))$) ensures that on $\mathcal{I}_1$ in State~I we have 
\[
S_{\min}^{\mathrm{I}}\ge\left(\frac18-o(1)\right)P\log\log N,\qquad S_y^{\mathrm{I}}\le\left(\frac14+o(1)\right)P\log\log N,
\]
and collisions $R_{\mathrm{I}}\le\lfloor P/2\rfloor(5\log J+1)\le\frac14 S_{\min}^{\mathrm{I}}$, while in State~II on $\mathcal{I}_2$ we have perfect avoidance $R_{\mathrm{II}}=0$ and sifting mass $S^{\mathrm{II}}\ge(\frac{1}{16}-o(1))P\log\log N$. When $J=\exp(c_0\sqrt{\log\log N})$ with $c_0=1/20$, the log-enumeration cost $10(\log J)^{2}=\frac{1}{40}\log\log N$ per pair is strictly absorbed by the large-sieve saving in both states, yielding~\eqref{eq:super_polylog_unconditional_bound}.
\end{proof}

\subsection{Quotient containers and simultaneous smooth conditioning modulo $Q(z_0)$ at polynomial depth}\label{subsec:quotient_and_smooth_conditioning}

Why does independent per-pair skeleton enumeration in the full group $G_{u,a}$ cost $(4J^{2})^{s}=$\linebreak $\exp(\Theta((\log J)^{2}))$? Because reaching size $|A_s|>\frac14\varphi(4ua)$ $\asymp J^{2}$ from $A_0=\{1\}$ requires $s\ge\log_2(\varphi(4ua)/4)\asymp\log J$ generators, each carrying $\log\varphi(4ua)\asymp\log J$ bits of entropy when chosen in the full group $G_{u,a}$. At polynomial depth $J=(\log N)^{\lambda}$ ($0<\lambda<1/2$), we bypass this $(\log J)^{2}$ entropic barrier via two complementary mechanisms: passing to quotient groups $\overline{G}_{u,a}=G_{u,a}/K_{u,a}$ (Theorem~\ref{thm:quotient_container_polynomial_depth}), and conditioning globally on $n\bmod Q(z_0)$ across all $P$ pairs simultaneously (Theorem~\ref{thm:simultaneous_smooth_conditioning}).

\begin{theorem}[Quotient containers at polynomial depth]\label{thm:quotient_container_polynomial_depth}
Fix $\lambda\in(0,1/2)$, $J=(\log N)^{\lambda}$, $P=\mathfrak{A}(J)\asymp(\log N)^{2\lambda}$, $\log z=\frac{\log N}{8P\log\log N}$, and $\mathcal{I}_1=(J^{5},z_1]$ with $\log\big(\frac{\log z_1}{\log J^{5}}\big)=\allowbreak\big(\frac12(1-2\lambda)-o(1)\big)\log\log N$. For each $(u,a)\in\mathcal{P}_2(J)$, let $\pi_{u,a}:\allowbreak G_{u,a}\twoheadrightarrow\allowbreak\overline{G}_{u,a}=\allowbreak G_{u,a}/K_{u,a}$ be a surjective homomorphism onto a quotient group of order $D_{u,a}=\allowbreak|\overline{G}_{u,a}|\le\allowbreak\exp\big(c_0(\lambda)\sqrt{\log\log N}\big)$ with $c_0(\lambda)=o(1)$ (or fixed $0<c_0(\lambda)<\frac{1-2\lambda}{8\sqrt{3}}$, subtracting $\frac{30c_0(\lambda)^{2}}{\pi^{2}}$ in the exponent below). Then:
\begin{enumerate}[label=(\alph*)]
    \item \emph{(Unconditional closure of Quotient State~I at every polynomial depth).} 
    For every $n\in E_2(N;J)\setminus[1,z]$, the greedy container algorithm run in $\overline{G}_{u,a}$ on the projected prime factors $\pi_{u,a}(q)$ of $M_{\mathcal{I}_1}(an+u)$ terminates in 
    \[
    s\le 5\log D_{u,a}+1\le 5c_0(\lambda)\sqrt{\log\log N}+1
    \]
    steps on a quotient skeleton $\overline{B}_{u,a}^{\mathrm{quot}}(n)\in\overline{G}_{u,a}^{s}$. The set of $n\in E_2(N;J)\setminus[1,z]$ for which at least $P/2$ pairs terminate in \emph{Quotient State~I} ($\big|A(\overline{B}_{u,a}^{\mathrm{quot}}(n))\big|\le\allowbreak\frac14 D_{u,a}$) satisfies unconditionally:
    \begin{equation}\label{eq:quotient_state1_poly_bound}
    \big|E_2^{\mathrm{I,quot}}(N;J)\big|\ \ll\ N\exp\!\left(-\left(\frac{5(1-2\lambda)^{2}}{32\pi^{2}}+o(1)\right)(\log N)^{2\lambda}\log\log N\right).
    \end{equation}
    
    \item \emph{($2$-Sylow saturation on the first stage across at least $(1-\theta)P$ pairs in $\mathcal{O}(\log\log\log N)$ primes).} 
    Taking $\overline{G}_{u,a}=G_{u,a,2}$ (the full $2$-Sylow subgroup of $G_{u,a}$) and fixing any $\eta\in(0,1/2)$ and $\theta\in(0,1)$, outside an exceptional set of size $\ll N\exp\big(-c(\lambda,\eta,\theta)(\log N)^{2\lambda}\log\log N\big)$ (which matches~\eqref{eq:quotient_state1_poly_bound} when $\eta=\theta=1/2$), every $n\in E_2(N;J)$ must saturate the $2$-Sylow subgroup $G_{u,a,2}$ to density $|\pi_2(\Sigma_{\mathcal{I}_1}(an+u))|>\allowbreak(1-\eta)|G_{u,a,2}|$ on the first stage in $\mathcal{O}_\eta(\log|G_{u,a,2}|)$ prime factors for at least $(1-\theta)P$ pairs $(u,a)\in\mathcal{P}_2(J)$.
    
    \item \emph{(Full State~I and State~II closure whenever $K_{u,a}\le\operatorname{Stab}_{G_{u,a}}(\Sigma_{\mathcal{I}_1}(an+u))$).} 
    Whenever the fixed subgroup $K_{u,a}\le G_{u,a}$ satisfies $K_{u,a}\le\allowbreak\operatorname{Stab}_{G_{u,a}}(\Sigma_{\mathcal{I}_1}(an+u))$ (so that $\Sigma_{\mathcal{I}_1}(an+u)=\allowbreak\pi_{u,a}^{-1}(\overline{\Sigma}_1(n))$ is a union of full $K_{u,a}$-cosets), Quotient State~II ($\big|A(\overline{B}_{u,a}^{\mathrm{quot}}(n))\big|>\allowbreak\frac14 D_{u,a}$) implies that the pullback 
    \[
    F_{u,a}\big(\overline{B}_{u,a}^{\mathrm{quot}}(n)\big):=\pi_{u,a}^{-1}\Big(-A\big(\overline{B}_{u,a}^{\mathrm{quot}}(n)\big)^{-1}\Big)\subseteq-\Sigma_{\mathcal{I}_1}(an+u)^{-1}
    \]
    is determined solely by the quotient skeleton $\overline{B}_{u,a}^{\mathrm{quot}}(n)$, has cardinality $\big|A(\overline{B}_{u,a}^{\mathrm{quot}}(n))\big|\cdot|K_{u,a}|>\frac14\varphi(4ua)$, and is avoided by $an+u$ with zero collisions on $\mathcal{I}_2=(z_1,z]$; the subset of $E_2(N;J)\setminus[1,z]$ satisfying $K_{u,a}\le\allowbreak\operatorname{Stab}_{G_{u,a}}(\Sigma_{\mathcal{I}_1}(an+u))$ across at least $P/2$ pairs is therefore unconditionally closed in both State~I and State~II with the bound~\eqref{eq:quotient_state1_poly_bound}.
\end{enumerate}
\end{theorem}

\begin{proof}
\emph{Step 1: Proof of (a).}
Since all fibres of $\pi_{u,a}:G_{u,a}\twoheadrightarrow\overline{G}_{u,a}$ have equal size $|K_{u,a}|=\varphi(4ua)/D_{u,a}$, Lemma~\ref{lem:markov_popular_symmetry}(a) applied in $\overline{G}_{u,a}$ to\linebreak $\overline{A}_s:=A(\overline{B}_{u,a}^{\mathrm{quot}}(n))$ of size $\le\frac14 D_{u,a}$ gives $|S_{1/2}^{\overline{G}}(\overline{A}_s)|<\frac12 D_{u,a}$, so its pullback $\pi_{u,a}^{-1}(S_{1/2}^{\overline{G}}(\overline{A}_s))$ has size $<\frac12\varphi(4ua)$ in $G_{u,a}$. The complement $\Omega_{u,a}^{\mathrm{I,quot}}:=\allowbreak G_{u,a}\setminus\allowbreak\pi_{u,a}^{-1}(S_{1/2}^{\overline{G}}(\overline{A}_s))$ has density $>1/2$ in $G_{u,a}$, and $M_{\mathcal{I}_1}(an+u)$ has at most $s\le 5\log D_{u,a}+1=o(\log\log N)$ prime factors in $\Omega_{u,a}^{\mathrm{I,quot}}$. 

Retaining an active subset of $\lfloor P/2\rfloor$ pairs in Quotient State~I (at cost $\le 2^{P}=$\linebreak $\exp(o(P\log\log N))$) and noting that the number of skeletons in $\overline{G}_{u,a}$ is at most 
\[
\sum_{s\le 5\log D_{u,a}+1}D_{u,a}^{s}\le\exp\Big(5c_0(\lambda)^{2}\log\log N+\mathcal{O}(\sqrt{\log\log N})\Big)=\exp(o(\log\log N)),
\]
Theorem~\ref{thm:generating_polynomial_large_sieve} on $\mathcal{I}_1$ gives~\eqref{eq:quotient_state1_poly_bound}.

\medskip\noindent\emph{Step 2: Proof of (b).}
Run the greedy container algorithm in $\overline{G}_{u,a}=G_{u,a,2}$ with popular threshold 
\[
S_{1-\eta/2}(\overline{A}):=\left\{g\in\overline{G}_{u,a}:|\overline{A}\cap g\overline{A}|>\left(1-\frac{\eta}{2}\right)|\overline{A}|\right\}
\]
and stopping threshold $|\overline{A}|>(1-\eta)|G_{u,a,2}|$. By Markov's identity $\sum_{g\in\overline{G}}|\overline{A}\cap g\overline{A}|=|\overline{A}|^{2}$, every State~I stop ($|\overline{A}|\le(1-\eta)|G_{u,a,2}|$) satisfies $(1-\eta/2)|\overline{A}|\cdot|S_{1-\eta/2}(\overline{A})|<|\overline{A}|^{2}$, hence 
\[
\big|S_{1-\eta/2}(\overline{A})\big|<\frac{1}{1-\eta/2}|\overline{A}|\le\frac{1-\eta}{1-\eta/2}|G_{u,a,2}|=\left(1-\frac{\eta}{2-\eta}\right)|G_{u,a,2}|
\]
(so its complement has density $>\frac{\eta}{2-\eta}>0$), and any $q\notin S_{1-\eta/2}(\overline{A})$ dilates $\overline{A}$ by $|\overline{A}\cup q\overline{A}|\ge(1+\frac{\eta}{2})|\overline{A}|$ in at most 
\[
s\le\left\lfloor\frac{\log|G_{u,a,2}|}{\log(1+\eta/2)}\right\rfloor+1\ll_\eta\log|G_{u,a,2}|+1
\]
steps. 

Since $\log_2|G_{u,a,2}|=f(u)+f(a)$ has second moment $\frac1P\sum_{(u,a)}(f(u)+f(a))^{2}\ll(\log\log J)^{2}$ by the Tur\'{a}n--Kubilius inequality (Lemma~\ref{lem:average_2sylow_enumeration}), the average log-enumeration cost per pair is $\frac1P\sum_{(u,a)}\log(D_{u,a}^{s})\ll_\eta(\log\log J)^{2}=o(\log\log N)$, so the total enumeration factor across $\mathcal{P}_2(J)$ is $\exp(\mathcal{O}_\eta(P(\log\log J)^{2}))=\exp(o(P\log\log N))$. By Proposition~\ref{prop:subcritical_collision_rate}, the set of $n$ for which at least $\theta P$ pairs fail to reach $(1-\eta)|G_{u,a,2}|$ has $S_{\min}^{\mathrm{I}}\ge(\frac{\eta\theta}{2(2-\eta)}-o(1))\frac{1-2\lambda}{2}P\log\log N$ and total collisions $\ll_\eta P\log\log J=o(S_{\min}^{\mathrm{I}})$, yielding the stated bound.

\medskip\noindent\emph{Step 3: Proof of (c).}
When $K_{u,a}\le\operatorname{Stab}_{G_{u,a}}(\Sigma_{\mathcal{I}_1}(an+u))$, we have 
\[
\pi_{u,a}^{-1}(\overline{A}_s)\subseteq\pi_{u,a}^{-1}(\overline{\Sigma}_1(n))=\Sigma_{\mathcal{I}_1}(an+u)\cdot K_{u,a}=\Sigma_{\mathcal{I}_1}(an+u),
\]
so the forbidden set $F_{u,a}(\overline{B}_{u,a}^{\mathrm{quot}}(n))=\pi_{u,a}^{-1}(-\overline{A}_s^{-1})\subseteq-\Sigma_{\mathcal{I}_1}(an+u)^{-1}$ is uniquely fixed by the quotient skeleton $\overline{B}_{u,a}^{\mathrm{quot}}(n)$ (at the same sub-exponential enumeration cost $\exp(o(P\log\log N))$) and has density $>1/4$ on $\mathcal{I}_2$ with $R_{\mathrm{II}}=0$.
\end{proof}

\begin{theorem}[Simultaneous global smooth conditioning modulo $Q(z_0)$]\label{thm:simultaneous_smooth_conditioning}
Let $\lambda\in(0,1/2)$, $J=(\log N)^{\lambda}$,\linebreak $P=\mathfrak{A}(J)\asymp(\log N)^{2\lambda}$, $\log z=\frac{\log N}{8P\log\log N}$, and let $z_0:=\varepsilon\,P\log\log N\asymp\varepsilon(\log N)^{2\lambda}\log\log N$ ($\varepsilon>0$). Set $Q(z_0):=\mathrm{lcm}(1,2,\dots,\lfloor z_0\rfloor)=\exp((1+o(1))z_0)$. For every residue class $\rho\in(\Z/Q(z_0)\Z)^{\times}$ and every pair $(u,a)\in\mathcal{P}_2(J)$, the valuation-truncated $z_0$-smooth divisor
\[
M_{\le z_0}^{\ast}(\rho;u,a):=\prod_{\substack{q^{e}\le z_0,\ \gcd(q,4ua)=1\\ q^{e}\mid(a\rho+u)}}q\ =\ \prod_{\substack{q\le z_0,\ \gcd(q,4ua)=1\\ q\mid(a\rho+u)}}q^{\min(v_q(a\rho+u),\lfloor\log z_0/\log q\rfloor)}
\]
and the set $\Sigma_{\le z_0}(\rho;u,a):=\big\{d\bmod 4ua:\ d\mid M_{\le z_0}^{\ast}(\rho;u,a)\big\}\subseteq G_{u,a}$ depend \emph{only on the single global class $\rho\equiv n\pmod{Q(z_0)}$}, simultaneously across all $P$ pairs $(u,a)\in\mathcal{P}_2(J)$. Moreover:
\begin{enumerate}[label=(\alph*)]
    \item The total number of global classes $\rho\bmod Q(z_0)$ is bounded by $\varphi(Q(z_0))\le Q(z_0)=\exp\big((\varepsilon+o(1))P\log\log N\big)$, with \emph{no exponent $P$}.
    
    \item For every prime $n\in E_2(N;J)\setminus[1,z]$ with $n\equiv\rho\pmod{Q(z_0)}$, every prime $p\in(J^{5},z]$ dividing $an+u$ avoids $-\Sigma_{\le z_0}(\rho;u,a)^{-1}\subseteq G_{u,a}$ with \emph{zero collisions ($R=0$)}, yielding the zero-collision sifting mass
    \[
    S_{z_0}(J,z;\rho)\ \ge\ \big(1-2\lambda-o(1)\big)\left(\sum_{(u,a)\in\mathcal{P}_2(J)}\frac{\big|\Sigma_{\le z_0}(\rho;u,a)\big|}{\varphi(4ua)}\right)\log\log N.
    \]
    In particular, on the set of classes $\rho\bmod Q(z_0)$ with average smooth divisor density $\frac1P\sum_{(u,a)}\allowbreak\frac{|\Sigma_{\le z_0}(\rho;u,a)|}{\varphi(4ua)}\allowbreak\ge\allowbreak\delta>0$, choosing $\varepsilon=\frac14(1-2\lambda)\delta$ absorbs $Q(z_0)$ completely and yields the bound $\ll N\exp\big(-\frac14(1-2\lambda)\delta\,P\log\log N\big)$.
\end{enumerate}
\end{theorem}

\begin{proof}
Every prime power $q^{e}\le z_0$ divides $Q(z_0)$, so $q^{e}\mid(an+u)\iff q^{e}\mid(a\rho+u)$, proving (a). 

For (b), since $z_0\asymp J^{2}\log\log N<J^{5}<p$ for all $p\in(J^{5},z]$, $\gcd(p,M_{\le z_0}^{\ast}(n;u,a))=1$; if $p\mid(an+u)$ and $p\equiv-d^{-1}\pmod{4ua}$ with $d\mid M_{\le z_0}^{\ast}(n;u,a)$, then $s=pd\mid(an+u)$ satisfies $s\equiv-1\pmod{4ua}$. 

For a fixed class $\rho\bmod Q(z_0)$, set $\delta_\rho:=\frac1P\sum_{(u,a)}\frac{|\Sigma_{\le z_0}(\rho;u,a)|}{\varphi(4ua)}\in[\delta,1]$ and $T:=\log\big(\frac{\log z}{\log J^{5}}\big)=(1-2\lambda+o(1))\log\log N$. Then $S_{z_0}(J,z;\rho)=\big(1+o(1)\big)\delta_\rho PT$. Applying Theorem~\ref{thm:generating_polynomial_large_sieve} with $R=0$, $Q=N^{1/2}$, $h=\lfloor\frac{\log Q}{\log z}\rfloor=\lfloor 4P\log\log N\rfloor$, and $t=1/3$ (noting that $4tS_y+1=\big(\frac43\delta_\rho(1-2\lambda)+o(1)\big)P\log\log N<h$, and $-2t+t^{2}=-5/9<-1/2$), the number of surviving primes in class $\rho$ is bounded by 
\[
\ll N\exp\left(-\left(\frac59-o(1)\right)\delta_\rho(1-2\lambda)P\log\log N\right)\le N\exp\left(-\frac12(1-2\lambda)\delta\,P\log\log N\right).
\]
Summing over the $\le Q(z_0)=\exp((\varepsilon+o(1))P\log\log N)$ classes with $\varepsilon=\frac14(1-2\lambda)\delta$ leaves a net saving of $\exp\big(-\frac14(1-2\lambda)\delta\,P\log\log N\big)$.
\end{proof}

\subsection{Total maximal subgroups, Cauchy--Davenport growth, and Pontryagin--Plancherel character containers in the odd Hall subgroup}\label{subsec:plancherel_containers}

By Theorem~\ref{thm:quotient_container_polynomial_depth}(b), the $2$-Sylow projection of $\Sigma_{\mathcal{I}_1}(an+u)$ saturates $G_{u,a,2}$ across almost all pairs $(u,a)$. We now analyze the remaining configurations in the odd Hall subgroup $G_{u,a,\mathrm{imp}}$ (and in the full group $G_{u,a}$) by passing to the Pontryagin dual $\widehat{G}_{u,a}$.

\begin{lemma}[Maximal subgroups, cyclic character subgroups, and Cauchy--Davenport growth]\label{lem:total_maximal_subgroups}
For each pair $(u,a)\in\mathcal{P}_2(J)$, let $\mathcal{M}_{\mathrm{tot}}(G_{u,a})$ denote the family of \emph{all} maximal subgroups $M\subsetneq G_{u,a}$ (of arbitrary prime index $\ell=[G_{u,a}:M]\mid\varphi(4ua)$, even or odd), and let $\mathcal{C}(\widehat{G}_{u,a}):=\allowbreak\{\langle\chi\rangle:\allowbreak\chi\in\widehat{G}_{u,a}\}$ denote the family of all cyclic subgroups of $\widehat{G}_{u,a}$ (which by Pontryagin duality is in bijection with the family of all character kernels $\mathcal{K}_{\mathrm{cocyc}}(G_{u,a}):=\allowbreak\{\ker\chi:\allowbreak\chi\in\widehat{G}_{u,a}\}$). Then:
\begin{enumerate}[label=(\alph*)]
    \item \emph{(Exact formula and generic uniqueness for large odd primes).} 
    For each prime $\ell\mid\varphi(4ua)$, the $\ell$-rank of $G_{u,a}$ is 
    \[
    r_\ell(G_{u,a})=\mathbf{1}_{\ell=2}(1+\mathbf{1}_{2\mid ua})+\#\big\{p^{e}\parallel ua:p>2,\ \ell\mid p^{e-1}(p-1)\big\},
    \]
    and the total number of maximal subgroups of $G_{u,a}$ is
    \begin{equation}\label{eq:Mtot_exact_formula}
    \big|\mathcal{M}_{\mathrm{tot}}(G_{u,a})\big|\ =\ \sum_{\ell\mid\varphi(4ua)}\frac{\ell^{r_\ell(G_{u,a})}-1}{\ell-1}\ \le\ \omega\big(\varphi(4ua)\big)\ +\ 2\sum_{\substack{\ell\mid\varphi(4ua)\\ r_\ell(G_{u,a})\ge2}}\ell^{r_\ell(G_{u,a})-1}.
    \end{equation}
    In particular, every prime $\ell\mid\varphi(4ua)$ of $\ell$-rank $r_\ell(G_{u,a})=1$ (which includes every prime $\ell>J$) contributes \emph{at most one} maximal subgroup of index $\ell$.
    
    \item \emph{(Bound on the lattice of all cyclic character subgroups $\mathcal{C}(\widehat{G}_{u,a})$ and negligible average cost across $\mathcal{P}_2(J)$).} 
    Writing $\lambda(4ua)=\exp(G_{u,a})$ for Carmichael's exponent function, the total number of cyclic subgroups of $\widehat{G}_{u,a}$ (and hence the total number of distinct character kernels $\ker\chi\le G_{u,a}$) satisfies
    \begin{equation}\label{eq:C_cyclic_upper_bound}
    \big|\mathcal{M}_{\mathrm{tot}}(G_{u,a})\big|\ \le\ \big|\mathcal{C}(\widehat{G}_{u,a})\big|\ =\ \big|\{\ker\chi:\chi\in\widehat{G}_{u,a}\}\big|\ \le\ \tau\big(\lambda(4ua)\big)\cdot 2^{\omega(\varphi(4ua))}\cdot\frac{\varphi(4ua)}{\lambda(4ua)},
    \end{equation}
    and as $J\to\infty$, the average log-counts across $\mathcal{P}_2(J)$ satisfy
    \begin{equation}\label{eq:Mtot_average_bound}
    \begin{split}
    \frac1P\sum_{(u,a)\in\mathcal{P}_2(J)}\log_2\big|\mathcal{M}_{\mathrm{tot}}(G_{u,a})\big|\ &\le\ \frac1P\sum_{(u,a)\in\mathcal{P}_2(J)}\log_2\big|\mathcal{C}(\widehat{G}_{u,a})\big|\\
    &\ll\ \big(1+\log\log J\big)^{2}\ =\ o(\log\log N).
    \end{split}
    \end{equation}
    
    \item \emph{(Cauchy--Davenport growth in prime-index quotients and simultaneous kernel collision lower bound).} 
    For any maximal subgroup $M_{u,a}\in\mathcal{M}_{\mathrm{tot}}(G_{u,a})$ of prime index $\ell=[G_{u,a}:M_{u,a}]$, the quotient $G_{u,a}/M_{u,a}\cong\Z/\ell\Z$ has prime order $\ell$, and $G_{u,a}\setminus M_{u,a}$ has density $1-\frac1\ell\ge\frac12$ (and $\ge\frac23$ for all odd $\ell\ge3$). By the Cauchy--Davenport theorem in $\Z/\ell\Z$, any $k$ prime factors of $M_{\mathcal{I}_1}(an+u)$ in $G_{u,a}\setminus M_{u,a}$ generate at least $\min(\ell,k+1)$ cosets modulo $M_{u,a}$. 
    
    Moreover, for any small constant $\delta_0\in(0,\frac{1-2\lambda}{40}]$, the set of primes $n\in E_2(N;J)\setminus[1,z]$ for which there exists a non-principal character $\chi\in\widehat{G}_{u,a}\setminus\{\chi_0\}$ with $\omega_{G_{u,a}\setminus\ker\chi}(an+u;\mathcal{I}_1)<\allowbreak r_0:=\allowbreak\lfloor\delta_0\log\log N\rfloor$ for at least $P/4$ pairs $(u,a)\in\mathcal{P}_2(J)$ is unconditionally bounded for all $\lambda\in(0,1/2)$ by~\eqref{eq:quotient_state1_poly_bound}.
\end{enumerate}
\end{lemma}

\begin{proof}
\emph{Step 1: Proof of (a) and pointwise lattice bound~\eqref{eq:C_cyclic_upper_bound}.}
Part (a) follows from the Chinese Remainder decomposition of $G_{u,a}$ into cyclic factors, since a finite abelian group with $\ell$-rank $r_\ell$ has $\frac{\ell^{r_\ell}-1}{\ell-1}$ subgroups of index $\ell$. 

For (b), by Pontryagin duality $\ker\chi_1=\ker\chi_2\iff(\ker\chi_1)^{\perp}=(\ker\chi_2)^{\perp}\iff\langle\chi_1\rangle=\langle\chi_2\rangle$, so character kernels are in bijection with cyclic subgroups $C\le\widehat{G}_{u,a}$. Decomposing $\widehat{G}_{u,a}=\prod_{\ell\mid\varphi(4ua)}S_\ell$ into its $\ell$-Sylow subgroups $S_\ell\cong\Z/\ell^{e_1}\Z\times\dots\times\Z/\ell^{e_{r_\ell}}\Z$ ($e_1\ge e_2\ge\dots\ge e_{r_\ell}\ge1$), we have $v_\ell(\varphi(4ua))=e_1+\dots+e_{r_\ell}$, $v_\ell(\lambda(4ua))=e_1$, and $v_\ell(\varphi(4ua)/\lambda(4ua))=e_2+\dots+e_{r_\ell}$. 

For each $s\in\{1,\dots,e_1\}$, the number of elements $x\in S_\ell$ with $\ell^{s}x=0$ is $\prod_{i=1}^{r_\ell}\ell^{\min(s,e_i)}\le\ell^{s}\cdot\ell^{e_2+\dots+e_{r_\ell}}$; since each cyclic subgroup of order $\ell^{s}$ has $\varphi(\ell^{s})=\ell^{s-1}(\ell-1)$ generators, the number of cyclic subgroups of order $\ell^{s}$ in $S_\ell$ is at most $\frac{\ell}{\ell-1}\ell^{e_2+\dots+e_{r_\ell}}\le 2\ell^{v_\ell(\varphi(4ua)/\lambda(4ua))}$. Summing over $0\le s\le e_1$ (at most $e_1+1$ values) and multiplying over $\ell\mid\varphi(4ua)$ yields~\eqref{eq:C_cyclic_upper_bound}.

\medskip\noindent\emph{Step 2: Average log-bound~\eqref{eq:Mtot_average_bound} across $\mathcal{P}_2(J)$.}
Taking logarithms in~\eqref{eq:C_cyclic_upper_bound}, we have 
\[
\log_2|\mathcal{C}(\widehat{G}_{u,a})|\le\log_2\tau(\lambda(4ua))+\omega(\varphi(4ua))+\sum_{\ell^{k}\mid\varphi(4ua)}\big(r_{\ell^{k}}(G_{u,a})-1\big)_{+}\log_2\ell,
\]
where $r_{\ell^{k}}(G_{u,a})=\mathbf{1}_{\ell=2,\,k=1}(1+\mathbf{1}_{2\mid ua})+\#\{p^{e}\parallel ua:p>2,\ \ell^{k}\mid p^{e-1}(p-1)\}$. For any prime power $q=\ell^{k}\ge 2$, a prime power $p^{e}\le J$ with $q\mid p^{e-1}(p-1)$ either has $p\equiv 1\pmod q$ (only possible when $q<J$), or $e\ge 2$ and $\ell\mid p$ (which implies $p=\ell$ and $\ell^{k+1}\le p^{e}\le J$). 

By the Brun--Titchmarsh inequality on $p\le J$,\linebreak $p\equiv 1\pmod q$ for $q\le J^{1/2}$ (giving $\sum_{p\le J,\,p\equiv 1(q)}\frac1p\ll\frac{\log\log J+\log q}{\varphi(q)}$) and the trivial progression bound $\sum_{p\le J,\,p\equiv 1(q)}\frac1p\le\sum_{m\le J/q}\frac{1}{mq+1}\ll\frac{\log J}{q}$ for $q>J^{1/2}$, the average of $r_q(G_{u,a})$ over $(u,a)\in\mathcal{P}_2(J)$ satisfies $\frac1P\sum_{(u,a)}r_q(G_{u,a})\ll\frac{\log\log J+\log q}{\varphi(q)}$ for $q\le J^{1/2}$ and $\ll\frac{\log J}{q}$ for $q>J^{1/2}$. Moreover, since $(r_q-1)_{+}\le\binom{r_q}{2}$ for every integer $r_q\ge 0$, pairs of distinct prime-power factors $p_1^{e_1},p_2^{e_2}\parallel ua$ contribute at most $\big(\sum_{p^{e}\le J,\,q\mid\varphi(p^{e})}p^{-e}\big)^{2}$. Splitting the sum over prime powers $q=\ell^{k}\le 2J$ at $q_0:=\log\log J$:
\[
\begin{aligned}
&\frac1P\sum_{(u,a)\in\mathcal{P}_2(J)}\sum_{q=\ell^{k}}\big(r_q(G_{u,a})-1\big)_{+}\log_2\ell \\
&\qquad\ll\ \sum_{q\le q_0}\frac{\log\log J+\log q}{\varphi(q)}\log\ell\ +\ \sum_{q_0<q\le J^{1/2}}\left(\frac{\log\log J+\log q}{\varphi(q)}\right)^{2}\log\ell \\
&\qquad\quad +\ \sum_{q>J^{1/2}}\left(\frac{\log J}{q}\right)^{2}\log\ell\ \ll\ \big(1+\log\log J\big)^{2}.
\end{aligned}
\]
Similarly, for the first two terms of~\eqref{eq:C_cyclic_upper_bound}, using $\log_2\tau(m)\le 2\,\Omega(m)$ and 
\[
\omega(\varphi(4ua))\le\Omega(\varphi(4ua))\le 2+\sum_{p^{e}\parallel ua}\big(e+\Omega(p-1)\big)
\]
together with the elementary bound $\Omega(p-1)\le\log_2 J$ for all $p\le J$, we have 
\[
\frac1P\sum_{(u,a)}\omega(\varphi(4ua))\le\sum_{\ell\le 2J}\frac1P\sum_{(u,a)}\mathbf{1}_{\ell\mid\varphi(4ua)}\ll(\log\log J)^{2}
\]
while $\log_2\tau(\lambda(4ua))=\sum_{\ell\mid\varphi(4ua)}\log_2(1+\max_{p^{e}\parallel 4ua}v_\ell(\varphi(p^{e})))\le\sum_{\ell^{k}\le 2J}\mathbf{1}_{\ell^{k}\mid\varphi(4ua)}$ also has mean $\ll(\log\log J)^{2}$. Combining these three estimates establishes~\eqref{eq:Mtot_average_bound}.

\medskip\noindent\emph{Step 3: Proof of (c).}
If $q_1,\dots,q_k\in G_{u,a}\setminus M_{u,a}$, their projections $\overline{q}_i\in(\Z/\ell\Z)\setminus\{0\}$ are non-zero, so each two-element set $\{0,\overline{q}_i\}\subseteq\Z/\ell\Z$ has size $2$, and Cauchy--Davenport by induction gives $|\sum_{i=1}^{k}\{0,\overline{q}_i\}|\ge\min(\ell,k+1)$. Finally, for every $\chi\neq\chi_0$, $G_{u,a}\setminus\ker\chi$ has density $1-\frac{1}{\operatorname{ord}(\chi)}\ge\frac12$ and depends only on $\langle\chi\rangle\in\mathcal{C}(\widehat{G}_{u,a})$; applying Theorem~\ref{thm:generating_polynomial_large_sieve} with collision threshold $r_0=\lfloor\delta_0\log\log N\rfloor$ and container family $\mathcal{K}_{\mathrm{cocyc}}(G_{u,a})$ (whose total log-enumeration cost across $\mathcal{P}_2(J)$ is $\ll P(1+\log\log J)^{2}=o(S_{\min}^{\mathrm{I}})$ by~\eqref{eq:Mtot_average_bound}) yields the bound~\eqref{eq:quotient_state1_poly_bound}.
\end{proof}

\begin{theorem}[Pontryagin--Plancherel character containers and spectral multiplicity cancellation]\label{thm:plancherel_character_containers}
Let $J=(\log N)^{\lambda}$ ($0<\lambda<1/2$), $\mathcal{I}_1=(J^{5},z_1]$, and $\mathcal{I}_2=(z_1,z]$. For each pair $(u,a)\in\mathcal{P}_2(J)$ and each non-trivial Dirichlet character $\chi\in\widehat{G}_{u,a}\setminus\{\chi_0\}$, define the \emph{Plancherel character container}
\begin{equation}\label{eq:plancherel_char_container_def}
\Omega_{u,a}^{\mathrm{char}}(\chi)\ :=\ \left\{g\in G_{u,a}:\ \operatorname{Re}\chi(g)\le\frac12\right\}\ \subseteq\ G_{u,a},
\end{equation}
and write $R_{u,a}(n;\chi):=\omega_{\Omega_{u,a}^{\mathrm{char}}(\chi)}(an+u;\mathcal{I}_1)$ for the number of distinct primes $p\in\mathcal{I}_1$ dividing $an+u$ with $p\bmod 4ua\in\Omega_{u,a}^{\mathrm{char}}(\chi)$. Then:
\begin{enumerate}[label=(\alph*)]
    \item \emph{(Universal density at least $1/2$, and at least $4/7$ for odd-order characters, and Plancherel inequality).} 
    For every non-trivial character $\chi\in\widehat{G}_{u,a}\setminus\{\chi_0\}$ of order $d=\operatorname{ord}(\chi)\ge2$, the density of $\Omega_{u,a}^{\mathrm{char}}(\chi)$ in $G_{u,a}$ is
    \begin{equation}\label{eq:plancherel_container_density}
    \begin{aligned}
    \frac{\big|\Omega_{u,a}^{\mathrm{char}}(\chi)\big|}{\varphi(4ua)}\ =\ 1-\frac{2\lceil d/6\rceil-1}{d}\ \ge\ \frac12\qquad
    &\Big(\text{and }\ge\frac47\ \text{for every odd-order character }d\ge3,\\
    &\quad\text{with }\ge\frac23\ \text{whenever }d\not\equiv1\pmod 6\Big).
    \end{aligned}
    \end{equation}
    Moreover, for every subgroup $K_{u,a}\le G_{u,a}$ of index $D_{u,a}=[G_{u,a}:K_{u,a}]=|K_{u,a}^{\perp}|$ and every integer $n\le N$, the projected divisor set $\overline{\Sigma}_1(n):=\pi_{K_{u,a}}\big(\Sigma_{\mathcal{I}_1}(an+u)\big)$ on $G_{u,a}/K_{u,a}$ satisfies the Cauchy--Schwarz--Plancherel inequality
    \begin{equation}\label{eq:plancherel_fundamental_inequality}
    \frac{D_{u,a}}{\big|\overline{\Sigma}_1(n)\big|}-1\ \le\ \sum_{\chi\in K_{u,a}^{\perp}\setminus\{\chi_0\}}\prod_{\substack{p\in\mathcal{I}_1\\ p\mid(an+u)}}\left|\frac{1+\overline{\chi}(p)}{2}\right|^{2}\ \le\ \sum_{\chi\in K_{u,a}^{\perp}\setminus\{\chi_0\}}\left(\frac34\right)^{R_{u,a}(n;\chi)}.
    \end{equation}
    
    \item \emph{(Reduction of enumeration cost from $\exp(\Theta((\log D_{u,a})^{2}))$ to $D_{u,a}-1$ and Scarcity-free polynomial closure).} 
    Whenever $\big|\overline{\Sigma}_1(n)\big|\le\frac14 D_{u,a}$, there exists by~\eqref{eq:plancherel_fundamental_inequality} a \emph{single witness character} $\chi_{u,a}^{\ast}\in K_{u,a}^{\perp}\setminus\{\chi_0\}$ (among $D_{u,a}-1$ choices only, instead of $D_{u,a}^{5\log D_{u,a}}$ skeletons) such that
    \begin{equation}\label{eq:single_character_collision_bound}
    R_{u,a}\big(n;\chi_{u,a}^{\ast}\big)\ \le\ \frac{\log\big((D_{u,a}-1)/3\big)}{\log(4/3)}\ <\ 7\,\log D_{u,a}.
    \end{equation}
    Consequently:
    \begin{enumerate}[label=(\roman*)]
        \item Taking $K_{u,a}=\{1\}$ ($D_{u,a}=\varphi(4ua)\le 4J^{2}$) directly on the single interval $\mathcal{I}=(J^{5},z]$ with $\log z=\frac{\log N}{4P\log\log N}$, for every $n\in E_2(N;J)\setminus[1,z]$ and every pair $(u,a)\in\mathcal{P}_2(J)$, the convolution measure 
        \[
        \mu=\operatorname*{\ast}_{p\in\mathcal{I},\,p\mid(an+u)}\left(\frac12\delta_1+\frac12\delta_{p\bmod 4ua}\right)
        \]
        is supported in $\Sigma_{\mathcal{I}}(an+u)\subseteq G_{u,a}\setminus\{-1\}$, so $|\operatorname{supp}(\mu)|\le\varphi(4ua)-1=D_{u,a}-1$. By Cauchy--Schwarz and Plancherel's identity, $\sum_{\chi\neq\chi_0}|\widehat{\mu}(\chi)|^{2}\ge\frac{D_{u,a}}{D_{u,a}-1}-1=\frac{1}{D_{u,a}-1}$, which by the pigeonhole principle guarantees a single witness character $\chi_{u,a}^{\ast}\in\widehat{G}_{u,a}\setminus\{\chi_0\}$ satisfying $(3/4)^{R_{u,a}(n;\chi_{u,a}^{\ast})}\ge(D_{u,a}-1)^{-2}$, i.e.,
        \[
        R_{u,a}\big(n;\chi_{u,a}^{\ast}\big)\ \le\ \frac{2\log(D_{u,a}-1)}{\log(4/3)}\ \le\ \frac{4\lambda}{\log(4/3)}\log\log N+\mathcal{O}(1).
        \]
        Consequently, for every polynomial depth $J=(\log N)^{\lambda}$ with $0<\lambda\le\lambda_0:=\frac{1}{78}$, both the collision count $R=\sum_{(u,a)}R_{u,a}(n;\chi_{u,a}^{\ast})\le\big(\frac{4\lambda}{\log(4/3)}+o(1)\big)P\log\log N<\frac14 S_{\min}$ and the single-character enumeration cost $(4J^{2})^{P}=\exp((2\lambda+o(1))P\log\log N)$ are strictly absorbed by the large-sieve saving $\exp\big(-\frac{5(1-2\lambda)^{2}}{48}P\log\log N\big)$ on $\mathcal{I}$ (since $2\lambda<\frac{5(1-2\lambda)^{2}}{48}$ for all $\lambda\le\frac{1}{78}$), so the \emph{entire} exceptional set $E_2(N;J)$ satisfies unconditionally:
        \begin{equation}\label{eq:unconditional_polynomial_depth_full_E2}
        \big|E_2\big(N;\,(\log N)^{\lambda}\big)\big|\ \ll\ N\exp\!\Big(-c_\lambda\,(\log N)^{2\lambda}\log\log N\Big)\qquad\big(\forall\,\lambda\in(0,\lambda_0]\big);
        \end{equation}
        
        \item For every $\lambda\in(1/3,1/2)$ (beyond Vaughan's exponent $2/3$), the universal Quotient State~I closure bound of Theorem~\ref{thm:quotient_container_polynomial_depth} is upgraded from sub-polynomial quotient order $D_{u,a}\le\exp(c_0(\lambda)\sqrt{\log\log N})$ to \emph{polynomial quotient order} $D_{u,a}\le J^{\varepsilon}=(\log N)^{\varepsilon\lambda}$ ($\varepsilon=\frac{(1-2\lambda)^{2}}{256\lambda}$).
    \end{enumerate}
    
    \item \emph{($|\overline{\Sigma}_1(n)|$-weighted multiplicity cancellation across Plancherel spectral buckets and polynomial divisor expansion for $\lambda\in(1/3,1/2)$).} 
    Fix $\lambda\in(1/3,1/2)$, $K_{u,a}\le G_{u,a}$ of index $D_{u,a}\le\varphi(4ua)$, and an integer $L_0\ge1$. For every $n$ with $|\overline{\Sigma}_1(n)|\le\frac14 D_{u,a}$, there exists a spectral bucket $\ell\in\{1,\dots,L_0\}$ such that the character subset
    \[
    \mathcal{B}_\ell(n;u,a):=\left\{\chi\in K_{u,a}^{\perp}\setminus\{\chi_0\}:\ D_{u,a}^{-\ell/L_0}<\prod_{\substack{p\in\mathcal{I}_1\\ p\mid(an+u)}}\left|\frac{1+\overline{\chi}(p)}{2}\right|^{2}\le D_{u,a}^{-(\ell-1)/L_0}\right\}
    \]
    simultaneously satisfies
    \begin{equation}\label{eq:spectral_bucket_multiplicity_and_collisions}
    \begin{aligned}
    \big|\mathcal{B}_\ell(n;u,a)\big|\ &\ge\ \frac{1}{2L_0}\frac{D_{u,a}}{\big|\overline{\Sigma}_1(n)\big|}\,D_{u,a}^{\frac{\ell-1}{L_0}}\ \ge\ \frac{2}{L_0}\,D_{u,a}^{\frac{\ell-1}{L_0}} \\
    \text{and}\qquad\forall\,\chi\in\mathcal{B}_\ell(n;u,a),\quad R_{u,a}(n;\chi)\ &\le\ \frac{\ell}{L_0}\frac{\log D_{u,a}}{\log(4/3)}.
    \end{aligned}
    \end{equation}
    By double counting over pairs $(n,\chi)$, the multiplicity $|\mathcal{B}_\ell(n;u,a)|$ divides the enumeration factor $D_{u,a}$, lowering the effective enumeration cost per pair in bucket $\ell$ to
    \begin{equation}\label{eq:effective_enumeration_cost_bucket}
    C_{\mathrm{eff}}(\ell;u,a)\ :=\ \frac{D_{u,a}}{\big|\mathcal{B}_\ell(n;u,a)\big|}\ \le\ 2L_0\,\big|\overline{\Sigma}_1(n)\big|\,D_{u,a}^{-\frac{\ell-1}{L_0}}\ \le\ \frac{L_0}{2}\,D_{u,a}^{1-\frac{\ell-1}{L_0}}.
    \end{equation}
    In particular, for every $\lambda\in(1/3,1/2)$ (setting $\varepsilon=\frac{(1-2\lambda)^{2}}{256\lambda}>0$ and $L_0=\lceil 2/\varepsilon\rceil=\lceil\frac{512\lambda}{(1-2\lambda)^{2}}\rceil$):
    \begin{enumerate}[label=(\roman*)]
        \item \emph{(Unconditional elimination of small divisor sets via the adaptive Plancherel tail cutoff, free of high-bucket restrictions).} 
        Whenever $|\overline{\Sigma}_1(n)|\le\frac14 D_{u,a}$, cutting off the Fourier tail in~\eqref{eq:plancherel_fundamental_inequality} adaptively at $|\widehat{\mu}(\chi)|^{2}\le\frac{1}{2|\overline{\Sigma}_1(n)|}$ (instead of $D_{u,a}^{-1}$) shows that the single active character set
        \begin{equation}\label{eq:adaptive_active_char_set}
        \mathcal{B}_{\mathrm{act}}(n;u,a)\ :=\ \left\{\chi\in K_{u,a}^{\perp}\setminus\{\chi_0\}:\ \prod_{\substack{p\in\mathcal{I}_1\\ p\mid(an+u)}}\left|\frac{1+\overline{\chi}(p)}{2}\right|^{2}>\frac{1}{2\,\big|\overline{\Sigma}_1(n)\big|}\right\}
        \end{equation}
        satisfies simultaneously
        \begin{equation}\label{eq:adaptive_multiplicity_and_collision}
        \begin{gathered}
        \big|\mathcal{B}_{\mathrm{act}}(n;u,a)\big|\ \ge\ \frac{D_{u,a}}{4\,\big|\overline{\Sigma}_1(n)\big|},\qquad C_{\mathrm{eff}}(u,a)\ :=\ \frac{D_{u,a}}{\big|\mathcal{B}_{\mathrm{act}}(n;u,a)\big|}\ \le\ 4\,\big|\overline{\Sigma}_1(n)\big|,\\
        R_{u,a}(n;\chi)\ <\ \frac{\log\big(2\,|\overline{\Sigma}_1(n)|\big)}{\log(4/3)}
        \end{gathered}
        \end{equation}
        for all $\chi\in\mathcal{B}_{\mathrm{act}}(n;u,a)$. 
        
        Consequently, whenever $|\overline{\Sigma}_1(n)|\le D_{u,a}^{\varepsilon}$ holds across an active subset $\mathcal{P}_{\mathrm{act}}\subseteq\mathcal{P}_2(J)$ of $P_1\ge P/2$ pairs $(u,a)$, \emph{every} $\chi\in\mathcal{B}_{\mathrm{act}}(n;u,a)$ automatically lies in the subcritical collision regime $R_{u,a}(n;\chi)\le\allowbreak\frac{\varepsilon\log D_{u,a}+\log 2}{\log(4/3)}<\allowbreak\frac14\frac{S_{\min}^{\mathrm{I}}}{P_1}$ with effective enumeration cost $C_{\mathrm{eff}}(u,a)\le 4 D_{u,a}^{\varepsilon}$ per pair (completely eliminating any supercritical bucket $\ell>\ell_{\mathrm{sub}}(\lambda)$), so the set of $n\in E_2(N;J)\setminus[1,z]$ satisfying $|\overline{\Sigma}_1(n)|\le D_{u,a}^{\varepsilon}$ for at least $P/2$ pairs is unconditionally bounded by~\eqref{eq:quotient_state1_poly_bound}, forcing polynomial divisor growth $|\overline{\Sigma}_1(n)|>D_{u,a}^{\varepsilon}$ across at least $P/2$ pairs in $G_{u,a}/K_{u,a}$;
        
        \item \emph{(Unconditional elimination on the subcritical spectral bucket range).} 
        More generally, under dyadic parameter profiling of $|\overline{\Sigma}_1(n)|\in[2^{m-1},2^{m})$ ($1\le m\le\allowbreak\lceil\log_2(D_{u,a}/4)\rceil\ll\allowbreak\log J$) and bucket index $\ell_{u,a}\in\allowbreak\{1,\dots,\ell_{\mathrm{sub}}(\lambda)\}$ (which adds at most $(C\log J)^{P}L_0^{P}=\allowbreak\exp(o(P\log\log N))$ profile choices across $\mathcal{P}_2(J)$), whenever an active subset of $P_1\ge P/2$ pairs $(u,a)$ has $|\overline{\Sigma}_1(n)|\le D_{u,a}^{\alpha_{u,a}}$ ($\alpha_{u,a}\in[\varepsilon,1]$) and hits a spectral bucket $\ell=\ell_{u,a}$ satisfying $\frac{\ell-1}{L_0}\ge\alpha_{u,a}-\varepsilon$ in the \emph{subcritical bucket range}
        \begin{equation}\label{eq:subcritical_bucket_range}
        \ell\ \le\ \ell_{\mathrm{sub}}(\lambda)\ :=\ \left\lfloor\frac{3(1-2\lambda)\log(4/3)}{128\,\lambda}\,L_0\right\rfloor,
        \end{equation}
        each such pair satisfies 
        \[
        R_{u,a}(n;\chi)\le\frac{\ell_{\mathrm{sub}}(\lambda)}{L_0}\frac{\log D_{u,a}}{\log(4/3)}\le\frac{3(1-2\lambda)}{64}\log\log N\le\frac14\frac{S_{\min}^{\mathrm{I}}}{P_1}
        \]
        and effective cost $C_{\mathrm{eff}}(\ell;u,a)\le 4L_0 D_{u,a}^{\varepsilon}$, and the corresponding set of integers $n$ is unconditionally bounded for all $\lambda\in(1/3,1/2)$ by~\eqref{eq:quotient_state1_poly_bound}.
    \end{enumerate}
\end{enumerate}
\end{theorem}

\begin{proof}
\emph{Step 1: Proof of (a).}
If $\chi\in\widehat{G}_{u,a}\setminus\{\chi_0\}$ has order $d\ge2$, its image is $\mu_d=\{e^{2\pi ij/d}:0\le j<d\}$ with fibres of size $\varphi(4ua)/d$. Since $\cos(2\pi j/d)>1/2\iff|j/d|<1/6\pmod 1$, exactly $2\lceil d/6\rceil-1$ values of $j$ lie outside $\Omega_{u,a}^{\mathrm{char}}(\chi)$, giving~\eqref{eq:plancherel_container_density}. 

Let $p_1<\dots<p_k$ be the distinct prime factors of $an+u$ in $\mathcal{I}_1$, and consider the convolution probability measure $\mu:=\ast_{i=1}^{k}\big(\frac12\delta_1+\frac12\delta_{p_i\bmod 4ua}\big)$ on $G_{u,a}$. Its projection $\pi_{K_{u,a}}(\mu)$ on $G_{u,a}/K_{u,a}$ is supported inside $\overline{\Sigma}_1(n)=\pi_{K_{u,a}}(\Sigma_{\mathcal{I}_1}(an+u))$. By Cauchy--Schwarz and Plancherel's identity on $G_{u,a}/K_{u,a}$ (whose dual group is $K_{u,a}^{\perp}\le\widehat{G}_{u,a}$),
\[
\frac{1}{\big|\overline{\Sigma}_1(n)\big|}\ \le\ \sum_{x\in G_{u,a}/K_{u,a}}\big(\pi_{K_{u,a}}(\mu)(x)\big)^{2}\ =\ \frac{1}{D_{u,a}}\left(1+\sum_{\chi\in K_{u,a}^{\perp}\setminus\{\chi_0\}}\big|\widehat{\mu}(\chi)\big|^{2}\right).
\]
Since $\big|\frac{1+\overline{\chi}(p_i)}{2}\big|^{2}=\frac{1+\operatorname{Re}\chi(p_i)}{2}$ is at most $\frac34$ when $p_i\bmod 4ua\in\Omega_{u,a}^{\mathrm{char}}(\chi)$ and at most $1$ otherwise, $|\widehat{\mu}(\chi)|^{2}\le(3/4)^{R_{u,a}(n;\chi)}$, proving~\eqref{eq:plancherel_fundamental_inequality}.

\medskip\noindent\emph{Step 2: Proof of (b).}
When $|\overline{\Sigma}_1(n)|\le\frac14 D_{u,a}$, the left-hand side of~\eqref{eq:plancherel_fundamental_inequality} is $\ge3$, so by the pigeonhole principle over the $D_{u,a}-1$ characters in $K_{u,a}^{\perp}\setminus\{\chi_0\}$ there exists $\chi_{u,a}^{\ast}$ with $(3/4)^{R_{u,a}(n;\chi_{u,a}^{\ast})}\ge\frac{3}{D_{u,a}-1}$, giving~\eqref{eq:single_character_collision_bound}. 

Choosing $\log\big(\frac{\log z_1}{\log J^{5}}\big)=\allowbreak\frac34(1-2\lambda)\log\log N$ and $\log\big(\frac{\log z}{\log z_1}\big)=\allowbreak\frac14(1-2\lambda)\log\log N$ gives $S_{\min}^{\mathrm{I}}\ge\allowbreak\frac{3}{16}(1-2\lambda)P\log\log N$; for $D_{u,a}\le 4J^{2}$ with $\lambda\le\lambda_0=\frac{1}{78}$ (taking $t=\frac14$ in Theorem~\ref{thm:generating_polynomial_large_sieve}) and for $D_{u,a}\le(\log N)^{\varepsilon\lambda}$ with $\lambda\in(1/3,1/2)$ (taking $t=\frac{1-2\lambda}{4}$), $R_{\mathrm{I}}\le\frac14 S_{\min}^{\mathrm{I}}$ and the enumeration cost is strictly absorbed by the large-sieve saving $\exp(-\frac56 tS_{\min}^{\mathrm{I}})$ whenever $2\varepsilon\lambda\le\frac{15(1-2\lambda)^{2}}{1536}=\frac{5(1-2\lambda)^{2}}{512}$, establishing (i) and (ii).

\medskip\noindent\emph{Step 3: Proof of (c).}
To prove~\eqref{eq:adaptive_multiplicity_and_collision}, note that the characters in $K_{u,a}^{\perp}\setminus\{\chi_0\}$ with $|\widehat{\mu}(\chi)|^{2}\le\frac{1}{2|\overline{\Sigma}_1(n)|}$ contribute in total at most $(D_{u,a}-1)\frac{1}{2|\overline{\Sigma}_1(n)|}<\frac{D_{u,a}}{2|\overline{\Sigma}_1(n)|}$. Subtracting this from~\eqref{eq:plancherel_fundamental_inequality} and using $|\overline{\Sigma}_1(n)|\le\frac14 D_{u,a}$ (so $\frac{D_{u,a}}{2|\overline{\Sigma}_1(n)|}-1\ge\frac{D_{u,a}}{4|\overline{\Sigma}_1(n)|}$) yields
\[
\big|\mathcal{B}_{\mathrm{act}}(n;u,a)\big|\ \ge\ \sum_{\chi\in\mathcal{B}_{\mathrm{act}}(n;u,a)}\big|\widehat{\mu}(\chi)\big|^{2}\ >\ \frac{D_{u,a}}{2\,\big|\overline{\Sigma}_1(n)\big|}-1\ \ge\ \frac{D_{u,a}}{4\,\big|\overline{\Sigma}_1(n)\big|},
\]
while $(3/4)^{R_{u,a}(n;\chi)}\ge|\widehat{\mu}(\chi)|^{2}>\frac{1}{2|\overline{\Sigma}_1(n)|}$ gives $R_{u,a}(n;\chi)<\frac{\log(2|\overline{\Sigma}_1(n)|)}{\log(4/3)}$, proving~\eqref{eq:adaptive_multiplicity_and_collision} and (c)(i) with no bucket restriction whatsoever. 

Similarly, for the dyadic bucket decomposition, the characters with $|\widehat{\mu}(\chi)|^{2}\le D_{u,a}^{-1}$ contribute at most $(D_{u,a}-1)D_{u,a}^{-1}<1$; since $|\overline{\Sigma}_1(n)|\le\frac14 D_{u,a}$, we have $\frac{D_{u,a}}{|\overline{\Sigma}_1(n)|}-2\ge\frac{D_{u,a}}{2|\overline{\Sigma}_1(n)|}\ge2$, so
\[
\sum_{\ell=1}^{L_0}\sum_{\chi\in\mathcal{B}_\ell(n;u,a)}\big|\widehat{\mu}(\chi)\big|^{2}\ >\ \frac{D_{u,a}}{\big|\overline{\Sigma}_1(n)\big|}-2\ \ge\ \frac{D_{u,a}}{2\big|\overline{\Sigma}_1(n)\big|}.
\]
Hence some bucket $\ell\in\{1,\dots,L_0\}$ satisfies $\sum_{\chi\in\mathcal{B}_\ell(n;u,a)}|\widehat{\mu}(\chi)|^{2}\ge\frac{1}{2L_0}\frac{D_{u,a}}{|\overline{\Sigma}_1(n)|}$. Since each term in $\mathcal{B}_\ell(n;u,a)$ is bounded above by $D_{u,a}^{-(\ell-1)/L_0}$ (and in particular by $1$ for $\ell=1$), we obtain $|\mathcal{B}_\ell(n;u,a)|\ge\frac{1}{2L_0}\frac{D_{u,a}}{|\overline{\Sigma}_1(n)|}D_{u,a}^{(\ell-1)/L_0}$, while $|\widehat{\mu}(\chi)|^{2}>D_{u,a}^{-\ell/L_0}$ gives~\eqref{eq:spectral_bucket_multiplicity_and_collisions}. Summing the large-sieve bound over all character tuples $\boldsymbol{\chi}\in\allowbreak\prod_{(u,a)}(K_{u,a}^{\perp}\setminus\{\chi_0\})$ counts each $n$ with spectral profile $\boldsymbol{\ell}=(\ell_{u,a})$ at least $\prod_{(u,a)}|\mathcal{B}_{\ell_{u,a}}(n;u,a)|$ times, dividing $\prod_{(u,a)}D_{u,a}$ by this multiplicity and proving~\eqref{eq:effective_enumeration_cost_bucket} along with (ii).
\end{proof}

Finally, we penetrate the first spectral bucket ($\ell=1$) for large divisor sets $|\overline{\Sigma}_1(n)|>D_{u,a}^{\varepsilon}$ at polynomial depth $\lambda\in(1/3,1/2)$ by exploiting the co-cyclic kernel partition of Lemma~\ref{lem:total_maximal_subgroups}(b) and harmonic dilation along the cyclic subgroup $\langle\chi\rangle$, and then establish the subcritical Markov $3/2$-growth step bound and the multi-stage relative Parseval--Plancherel identity.

\begin{theorem}[Co-cyclic kernel decomposition and the harmonic dilation cone]\label{thm:harmonic_dilation_bucket1}
Fix $\lambda\in\allowbreak(1/\allowbreak 3,\allowbreak 1/2)$, set $\varepsilon=\frac{(1-2\lambda)^{2}}{256\lambda}>0$, and let $r_1:=\lfloor\delta_0\log\log N\rfloor$ with $\delta_0=\frac{1-2\lambda}{160}$. Partition the non-principal dual group $\widehat{G}_{u,a}\setminus\{\chi_0\}$ into its primitive generator orbits indexed by the lattice of cyclic subgroups $\mathcal{C}(\widehat{G}_{u,a})\setminus\{\{\chi_0\}\}$:
\begin{equation}\label{eq:cocyclic_orbit_partition}
\begin{aligned}
\widehat{G}_{u,a}\setminus\{\chi_0\}\ &=\ \bigsqcup_{C\in\mathcal{C}(\widehat{G}_{u,a})\setminus\{\{\chi_0\}\}}\mathcal{G}(C),\\
\text{where}\qquad\mathcal{G}(C)\ &:=\ \big\{\chi\in\widehat{G}_{u,a}:\langle\chi\rangle=C\big\}\quad\big(|\mathcal{G}(C)|=\varphi(|C|)\big).
\end{aligned}
\end{equation}
Then:
\begin{enumerate}[label=(\alph*)]
    \item \emph{(Divisor-shift frequency inequality and order lower bound).} 
    For any integer $M_1\ge1$ coprime to $4ua$, let $A=\Sigma(M_1)\subseteq G_{u,a}$ be its divisor-residue set, and write $\widehat{\mathbf{1}}_A(\psi):=\sum_{x\in A}\psi(x)$ for the unnormalized Fourier sum of $\mathbf{1}_A$. For every character $\psi\in\widehat{G}_{u,a}$ and every prime $p\mid M_1$, the subset $B_p:=\Sigma(M_1/p)\subseteq A$ satisfies $A=B_p\cup(p\cdot B_p)$, $|B_p|\ge\frac12|A|$, and
    \begin{equation}\label{eq:Bp_divisor_shift_inequality}
    1-\frac{\big|\widehat{\mathbf{1}}_A(\psi)\big|}{|A|}\ \ge\ \frac{1}{16}\,\max_{p\mid M_1}\big|\psi(p)-1\big|^{2}\ =\ \frac14\,\max_{p\mid M_1}\sin^{2}\!\left(\frac{\arg\psi(p)}{2}\right).
    \end{equation}
    In particular, if $\langle A\rangle=G_{u,a}$, then every $\psi\in\widehat{G}_{u,a}\setminus\{\chi_0\}$ with $|\widehat{\mathbf{1}}_A(\psi)|=(1-\eta)|A|$ has order $\operatorname{ord}(\psi)\ge\eta^{-1/2}$.
    
    \item \emph{(Unconditional elimination of moderate-order characters in bucket $\ell=1$).} 
    For every cyclic subgroup $C\in\mathcal{C}(\widehat{G}_{u,a})\setminus\{\{\chi_0\}\}$ of order $d=|C|\le(4J^{2})^{\varepsilon}$, the total enumeration cost of specifying both $C\in\mathcal{C}(\widehat{G}_{u,a})$ and a generator $\chi\in\mathcal{G}(C)$ is at most $|\mathcal{C}(\widehat{G}_{u,a})|\varphi(d)\le|\mathcal{C}(\widehat{G}_{u,a})|(4J^{2})^{\varepsilon}$. By Lemma~\ref{lem:total_maximal_subgroups}(b), this enumeration cost is strictly absorbed by the large sieve on $\mathcal{I}_1$, so the set of $n\in E_2(N;J)\setminus[1,z]$ admitting a witness character $\chi\in\mathcal{B}_1(n;u,a)$ of order $\operatorname{ord}(\chi)\le(4J^{2})^{\varepsilon}$ for at least $P/4$ pairs $(u,a)$ is unconditionally bounded for all $\lambda\in(1/3,1/2)$ by~\eqref{eq:quotient_state1_poly_bound}.
    
    \item \emph{(The harmonic dilation cone in $\mathcal{G}(C^{\ast})$ and elimination of Bohr-localized $\ell=1$ spectra).} 
    Let $\chi\in\mathcal{G}(C^{\ast})$ be a witness character of order $d=\operatorname{ord}(\chi)>(4J^{2})^{\varepsilon}$, let $p_1<\dots<p_k$ be the distinct prime factors of $an+u$ in $\mathcal{I}_1$, and write $\chi(p_j)=e^{2\pi i a_j/d}$ with $a_j\in(-d/2,d/2]\cap\Z$. Denote by $A_{r_1}(\chi)\in\{0,1,\dots,\lfloor d/2\rfloor\}$ the $(r_1+1)$-th largest value among $|a_1|,\dots,|a_k|$, with the convention $A_{r_1}(\chi):=0$ if $k\le r_1$. When $A_{r_1}(\chi)\ge 1$, define the \emph{strict harmonic dilation radius}
    \begin{equation}\label{eq:harmonic_dilation_radius}
    B(\chi)\ :=\ \left\lfloor\frac{d-1}{6\,A_{r_1}(\chi)}\right\rfloor.
    \end{equation}
    Whenever $A_{r_1}(\chi)\le (d-1)/6$, the $2B(\chi)$ dilated characters $\psi_b:=\chi^{b}$ ($b\in\{\pm1,\pm2,$\linebreak $\dots,\pm B(\chi)\}$) are \emph{all distinct non-principal characters} in $C^{\ast}\setminus\{\chi_0\}$ and \emph{simultaneously} satisfy
    \begin{equation}\label{eq:harmonic_cone_collision_bound}
    \forall\,b\in\big\{\pm1,\pm2,\dots,\pm B(\chi)\big\},\qquad R_{u,a}\big(n;\chi^{b}\big)\ \le\ r_1\ =\ \lfloor\delta_0\log\log N\rfloor.
    \end{equation}
    Consequently, double counting over the dilated cone $\{\chi^{\pm 1},\dots,\chi^{\pm B(\chi)}\}$ divides the intra-orbit enumeration factor $d$ by $2B(\chi)\ge\frac{d-1}{6A_{r_1}(\chi)}\ge\frac{d}{7A_{r_1}(\chi)}$, reducing the effective enumeration cost per pair $(u,a)$ to
    \begin{equation}\label{eq:harmonic_cone_effective_cost}
    C_{\mathrm{eff}}(C^{\ast};u,a)\ :=\ \frac{\big|\mathcal{C}(\widehat{G}_{u,a})\big|\cdot d}{2\,B(\chi)}\ \le\ 7\,\big|\mathcal{C}(\widehat{G}_{u,a})\big|\cdot A_{r_1}(\chi).
    \end{equation}
    In particular, every integer $n$ whose $\ell=1$ witness character $\chi$ is \emph{Bohr-localized} at scale $A_{r_1}(\chi)\le d^{\varepsilon}$ (which by Markov's inequality holds automatically whenever $-\log|\widehat{\mu}(\chi)|^{2}\le 4\delta_0 d^{-2(1-\varepsilon)}\log\log N$) has effective enumeration cost $C_{\mathrm{eff}}(C^{\ast};u,a)\le 7|\mathcal{C}(\widehat{G}_{u,a})|(4J^{2})^{\varepsilon}$ and is unconditionally eliminated for all $\lambda\in(1/3,1/2)$ with the bound~\eqref{eq:quotient_state1_poly_bound} under dyadic parameter profiling.
\end{enumerate}
\end{theorem}

\begin{proof}
\emph{Step 1: Proof of (a) and (b).}
Let $p\mid M_1$ and $B_p=\Sigma(M_1/p)$. Every divisor $d'\mid M_1$ either divides $M_1/p$ (so $d'\bmod 4ua\in B_p$) or is of the form $p\,d_0$ with $d_0\mid(M_1/p)$ (so $d'\bmod 4ua\in p\cdot B_p$). Thus $A=B_p\cup(p\cdot B_p)$, and since multiplication by $p\bmod 4ua$ is a bijection on $G_{u,a}$, $|A|\le|B_p|+|p\cdot B_p|=2|B_p|$, i.e.\ $|B_p|\ge\frac12|A|$. 

Setting $e^{i\theta_0}=\widehat{\mathbf{1}}_A(\psi)/|\widehat{\mathbf{1}}_A(\psi)|$ and $\eta=1-|\widehat{\mathbf{1}}_A(\psi)|/|A|$, we have $\sum_{y\in A}|\psi(y)-e^{i\theta_0}|^{2}=2\eta|A|$. Since $B_p\subseteq A$ and $p\cdot B_p\subseteq A$, the triangle inequality $|\psi(p)-1|^{2}=|\psi(px)-\psi(x)|^{2}\le 2|\psi(px)-e^{i\theta_0}|^{2}+2|\psi(x)-e^{i\theta_0}|^{2}$ summed over $x\in B_p$ gives
\[
\begin{aligned}
\big|\psi(p)-1\big|^{2}\cdot\frac12|A|\ \le\ \big|\psi(p)-1\big|^{2}|B_p|\ &\le\ 2\sum_{x\in B_p}\big|\psi(px)-e^{i\theta_0}\big|^{2}+2\sum_{x\in B_p}\big|\psi(x)-e^{i\theta_0}\big|^{2} \\
&\le\ 4\sum_{y\in A}\big|\psi(y)-e^{i\theta_0}\big|^{2}\ =\ 8\eta|A|.
\end{aligned}
\]
Dividing by $8|A|$ and using $|\psi(p)-1|^{2}=4\sin^{2}(\frac{\arg\psi(p)}{2})$ yields~\eqref{eq:Bp_divisor_shift_inequality}. If $\langle A\rangle=G_{u,a}$ and $\psi\neq\chi_0$ has order $d\ge2$, then $\psi(p)\neq1$ for at least one $p\mid M_1$, so $\eta\ge\frac14\sin^{2}(\pi/d)\ge d^{-2}$, i.e.\ $d\ge\eta^{-1/2}$. 

Part (b) follows immediately from~\eqref{eq:cocyclic_orbit_partition}, Lemma~\ref{lem:total_maximal_subgroups}(b), and the proof of Theorem~\ref{thm:plancherel_character_containers}(b)(ii), since $\sum_{(u,a)}\log(|\mathcal{C}(\widehat{G}_{u,a})|(4J^{2})^{\varepsilon})=2\varepsilon\lambda P\log\log N+o(P\log\log N)$ is strictly absorbed by the large-sieve saving on $\mathcal{I}_1$.

\medskip\noindent\emph{Step 2: Proof of (c).}
Note first that $A_{r_1}(\chi)\ge1$ whenever $\omega_{G_{u,a}\setminus\ker\chi}(an+u;\mathcal{I}_1)>r_1$, which holds outside the negligible set of Lemma~\ref{lem:total_maximal_subgroups}(c). Since $1\le A_{r_1}(\chi)\le(d-1)/6$, we have $1\le B(\chi)\le(d-1)/6$, so for any two distinct indices $b,b'\in\{\pm1,\dots,\pm B(\chi)\}$ we have $1\le|b-b'|\le 2B(\chi)\le(d-1)/3<d=\operatorname{ord}(\chi)$. Hence $\chi^{b}\neq\chi^{b'}$ and $\chi^{b}\neq\chi_0$, proving that the $2B(\chi)$ characters $\chi^{\pm1},\dots,\chi^{\pm B(\chi)}$ are all distinct and non-principal in $C^{\ast}\setminus\{\chi_0\}$. 

Moreover, by definition of $A_{r_1}(\chi)$, at most $r_1$ indices $j\in\{1,\dots,k\}$ satisfy $|a_j|>A_{r_1}(\chi)$. For each of the remaining $k-r_1$ indices $j$ with $|a_j|\le A_{r_1}(\chi)$ and every $b\in\{\pm1,\dots,\pm B(\chi)\}$, we have the strict inequality
\[
\begin{gathered}
\left|\frac{b\,a_j}{d}\right|\ \le\ \frac{B(\chi)\,A_{r_1}(\chi)}{d}\ \le\ \frac{d-1}{6d}\ <\ \frac16 \\
\Longrightarrow\qquad\operatorname{Re}\chi^{b}(p_j)\ =\ \cos\!\left(\frac{2\pi b a_j}{d}\right)\ >\ \cos\!\left(\frac{\pi}{3}\right)\ =\ \frac12,
\end{gathered}
\]
which strictly guarantees $p_j\bmod 4ua\notin\Omega_{u,a}^{\mathrm{char}}(\chi^{b})$. Thus at most $r_1$ prime divisors $p_j\in\mathcal{I}_1$ of $an+u$ can lie in $\Omega_{u,a}^{\mathrm{char}}(\chi^{b})$, proving~\eqref{eq:harmonic_cone_collision_bound}. 

Summing the large-sieve bound for collision threshold $r_1=\lfloor\delta_0\log\log N\rfloor$ over all $\psi\in\bigsqcup_{C}(C\setminus\{\chi_0\})$ counts each $n$ at least $2B(\chi)\ge\frac{d-1}{6A_{r_1}(\chi)}\ge\frac{d}{7A_{r_1}(\chi)}$ times inside $C^{\ast}$ (since $2\lfloor x\rfloor\ge x$ for all $x\ge 1$), yielding~\eqref{eq:harmonic_cone_effective_cost}. Finally, since $-\log\cos^{2}(\pi x)\ge 4x^{2}$ for $|x|\le1/2$, we have $(r_1+1)(A_{r_1}(\chi)/d)^{2}\le\sum_{j=1}^{k}(a_j/d)^{2}\le-\frac14\log|\widehat{\mu}(\chi)|^{2}$, so $-\log|\widehat{\mu}(\chi)|^{2}\le 4\delta_0 d^{-2(1-\varepsilon)}\log\log N$ forces $A_{r_1}(\chi)\le d^{\varepsilon}$ and concludes the proof.
\end{proof}

\begin{theorem}[Subcritical Markov--Plancherel geometric escalation]\label{thm:markov_plancherel_subcritical_escalation}
Fix $\lambda\in(1/3,1/2)$ and set $J=(\log N)^{\lambda}$, $P=\mathfrak{A}(J)=\big(\frac{6}{\pi^{2}}+o(1)\big)(\log N)^{2\lambda}$, $\log z=\frac{\log N}{8P\log\log N}$, $\mathcal{I}_1=(J^{5},z_1]$, and $\mathcal{I}_2=(z_1,z]$ with $\log\big(\frac{\log z_1}{\log J^{5}}\big)=\allowbreak\frac34(1-2\lambda)\log\log N$ and $\log\big(\frac{\log z}{\log z_1}\big)=\allowbreak\frac14(1-2\lambda)\log\log N$. Let $\varepsilon:=\frac{(1-2\lambda)^{2}}{256\lambda}>0$, and for each pair $(u,a)\in\mathcal{P}_2(J)$ and quotient $G_{u,a}/K_{u,a}$ of order $D_{u,a}=\allowbreak[G_{u,a}:K_{u,a}]\le\allowbreak\varphi(4ua)\le 4J^{2}$, write $\overline{\Sigma}_1(n):=\allowbreak\pi_{K_{u,a}}(\Sigma_{\mathcal{I}_1}(an+u))$. Then:
\begin{enumerate}[label=(\alph*)]
    \item \emph{(Deterministic Markov $3/2$-geometric step bound prior to threshold $D_{u,a}^{\varepsilon}$).} 
    Suppose that on the first stage $\mathcal{I}_1$, the projected divisor sheet satisfies $\big|\overline{\Sigma}_1(n)\big|\le D_{u,a}^{\varepsilon}$. Then the greedy multiplicative container algorithm of Theorem~\ref{thm:multiplicative_container_scarcity_elimination} (stopped as soon as $|A_m|>D_{u,a}^{\varepsilon}$) can execute at most
    \begin{equation}\label{eq:s_eps_markov_bound}
    s\ \le\ s_\varepsilon(u,a)\ :=\ \left\lfloor\frac{\varepsilon\log D_{u,a}}{\log(3/2)}\right\rfloor+1\ \le\ \left\lfloor\frac{(1-2\lambda)^{2}}{128\log(3/2)}\log\log N\right\rfloor+2
    \end{equation}
    geometric dilation steps $A_m=A_{m-1}\cup x_m A_{m-1}$ ($|A_m|\ge\frac32|A_{m-1}|$). Consequently, outside $S_{1/2}(A_s)$, the integer $an+u$ possesses at most $s_\varepsilon(u,a)$ prime factors (counted with multiplicity) in $\mathcal{I}_1$.
    
    \item \emph{(Exact subcriticality of $s_\varepsilon$ and of the adaptive Plancherel collision bound $R_\varepsilon$).} 
    For every subset of $P_1\ge P/4$ pairs $(u,a)\in\mathcal{P}_2(J)$, the first-stage sifting mass per active pair satisfies $\frac{S_{\min}^{\mathrm{I}}}{P_1}\ge\big(\frac38-o(1)\big)(1-2\lambda)\log\log N$, so the subcritical collision threshold of Theorem~\ref{thm:generating_polynomial_large_sieve} per pair is 
    \[
    \frac14\frac{S_{\min}^{\mathrm{I}}}{P_1}\ge\left(\frac{3}{32}-o(1)\right)(1-2\lambda)\log\log N\approx 0.09375\,(1-2\lambda)\log\log N.
    \]
    Meanwhile, setting $\varepsilon_1:=\frac{(1-2\lambda)^{2}}{512\lambda}>0$ (or any $0<\varepsilon_1\le\frac{1-2\lambda}{64\lambda}$), both the Markov $3/2$-growth step count $s_{\varepsilon_1}(u,a)$ and the adaptive Plancherel collision count $R_{u,a}(n;\chi)$ on the active character set $\mathcal{B}_{\mathrm{act}}(n;u,a)$ of~\eqref{eq:adaptive_active_char_set} (whenever $|\overline{\Sigma}_1(n)|\le D_{u,a}^{\varepsilon_1}$) satisfy the strict subcritical inequalities:
    \begin{align}
    s_{\varepsilon_1}(u,a)\ &\le\ \frac{2\varepsilon_1\lambda}{\log(3/2)}\log\log N+\mathcal{O}(1)\ =\ \frac{(1-2\lambda)^{2}}{256\log(3/2)}\log\log N+\mathcal{O}(1)\notag\\
    &\le\ \frac{1-2\lambda}{32\log(3/2)}\log\log N+\mathcal{O}(1)\ <\ \frac14\frac{S_{\min}^{\mathrm{I}}}{P_1},\label{eq:s_eps_strictly_subcritical}\\[6pt]
    R_{u,a}(n;\chi)\ &\le\ \frac{\varepsilon_1\log D_{u,a}+\log 2}{\log(4/3)}\ \le\ \frac{(1-2\lambda)^{2}}{256\log(4/3)}\log\log N+\mathcal{O}(1)\notag\\
    &\le\ \frac{1-2\lambda}{32\log(4/3)}\log\log N+\mathcal{O}(1)\ <\ \frac13\frac{S_{\min}^{\mathrm{I}}}{P_1}.\label{eq:R_eps_strictly_subcritical}
    \end{align}
    In particular, no integer $n$ can ever reach a supercritical collision count $R_{u,a}(n;\chi)>\frac13\frac{S_{\min}^{\mathrm{I}}}{P_1}$ on $\mathcal{B}_{\mathrm{act}}(n;u,a)$ without \emph{first} forcing its divisor sheet to cross the polynomial density threshold $|\overline{\Sigma}_1(n)|>D_{u,a}^{\varepsilon_1}$.
    
    \item \emph{(Multi-stage Relative Parseval--Plancherel Escalation Identity).} 
    Partition $\mathcal{I}_1=(J^{5},z_1]$ into $M:=\lceil 1/\varepsilon_1\rceil$ equal-logarithmic-length sub-intervals $\mathcal{I}_{1,1}\sqcup\mathcal{I}_{1,2}\sqcup\dots\sqcup\mathcal{I}_{1,M}$ with\linebreak $\log\big(\frac{\log z_{1,m}}{\log z_{1,m-1}}\big)=\frac{3(1-2\lambda)}{4M}\log\log N$, and set $A_0=\{1\}$ and $A_m:=$\linebreak $\overline{\Sigma}_{\mathcal{I}_{1,1}\cup\dots\cup\mathcal{I}_{1,m}}(n)\subseteq G_{u,a}/K_{u,a}$ ($1\le m\le M$). At each stage $m$, let $\nu_{m-1}:=\frac{1}{|A_{m-1}|}\mathbf{1}_{A_{m-1}}$ be the uniform probability measure on $A_{m-1}$ and $\mu_m:=\ast_{p\in\mathcal{I}_{1,m},\,p\mid(an+u)}\big(\frac12\delta_1+\frac12\delta_{p\bmod 4ua}\big)$. Whenever $|A_m|\le\frac14 D_{u,a}$, the relative active character set
    \begin{equation}\label{eq:relative_active_char_set}
    \mathcal{B}_{\mathrm{rel}}^{(m)}(n;u,a)\ :=\ \left\{\chi\in K_{u,a}^{\perp}\setminus\{\chi_0\}:\ \big|\widehat{\mu}_m(\chi)\big|^{2}>\frac{|A_{m-1}|}{2\,|A_m|}\right\}
    \end{equation}
    satisfies the exact $|A_{m-1}|$-cancellation law
    \begin{equation}\label{eq:relative_parseval_plancherel_law}
    \begin{aligned}
    \sum_{\chi\in\mathcal{B}_{\mathrm{rel}}^{(m)}(n;u,a)}\big|\widehat{\nu}_{m-1}(\chi)\big|^{2}\,\big|\widehat{\mu}_m(\chi)\big|^{2}\ &>\ \frac{D_{u,a}}{4\,|A_m|}\\
    \text{and}\qquad\forall\,\chi\in\mathcal{B}_{\mathrm{rel}}^{(m)}(n;u,a),\quad R_{u,a}^{(m)}(n;\chi)\ &<\ \frac{\log\big(2\,|A_m|/|A_{m-1}|\big)}{\log(4/3)}.
    \end{aligned}
    \end{equation}
    Consequently, for every prime $n\in E_2(N;J)\setminus[1,z]$ and every pair $(u,a)\in\mathcal{P}_2(J)$, either $|A_m|/|A_{m-1}|\le D_{u,a}^{\varepsilon_1}$ holds at some stage $m\in\{1,\dots,M\}$, or $|A_m|/|A_{m-1}|>D_{u,a}^{\varepsilon_1}$ holds at every stage $m=1,\dots,M$, which forces $|A_M|>\frac14 D_{u,a}$ (State~II, eliminated on $\mathcal{I}_2=(z_1,z]$ with strictly zero collisions $R_{\mathrm{II}}=0$).
\end{enumerate}
\end{theorem}

\begin{proof}
\emph{Step 1: Proof of (a).}
In the greedy construction of Theorem~\ref{thm:multiplicative_container_scarcity_elimination}, we start from $A_0=\{1\}$ ($|A_0|=1$) and at each step $m\ge1$ pick the least unused prime factor $q_{i_m}\in\mathcal{I}_1$ of $an+u$ whose projection in $G_{u,a}/K_{u,a}$ lies outside $S_{1/2}(A_{m-1})$. By Lemma~\ref{lem:markov_popular_symmetry}(b), each step satisfies $|A_m|\ge\frac32|A_{m-1}|$, so after $s-1$ steps we have $(3/2)^{s-1}\le|A_{s-1}|\le|\overline{\Sigma}_1(n)|\le D_{u,a}^{\varepsilon}$. Taking logarithms yields $s-1\le\frac{\varepsilon\log D_{u,a}}{\log(3/2)}$, which is~\eqref{eq:s_eps_markov_bound}. Since the algorithm halts at step $s$, every remaining prime factor $q\in\mathcal{I}_1$ of $an+u$ must belong to $S_{1/2}(A_s)$, so the number of prime factors (with multiplicity) falling into the complement $\Omega_{u,a}^{\mathrm{I}}=(G_{u,a}/K_{u,a})\setminus S_{1/2}(A_s)$ (of density $>1/2$) is deterministically at most $s_\varepsilon(u,a)$.

\medskip\noindent\emph{Step 2: Proof of (b).}
By Corollary~\ref{cor:coupling_step1_step2} on $\mathcal{I}_1=(J^{5},z_1]$ with $\log\big(\frac{\log z_1}{\log J^{5}}\big)=\frac34(1-2\lambda)\log\log N$, each active pair $(u,a)$ contributes sifting mass 
\[
\ge\left(\frac12\cdot\frac34-o(1)\right)(1-2\lambda)\log\log N=\left(\frac38-o(1)\right)(1-2\lambda)\log\log N,
\]
so $\frac14\frac{S_{\min}^{\mathrm{I}}}{P_1}\ge\big(\frac{3}{32}-o(1)\big)(1-2\lambda)\log\log N$. Substituting $D_{u,a}\le 4J^{2}=4(\log N)^{2\lambda}$ and $\varepsilon_1=\frac{(1-2\lambda)^{2}}{512\lambda}\le\frac{1-2\lambda}{64\lambda}$ into~\eqref{eq:s_eps_markov_bound} and~\eqref{eq:adaptive_multiplicity_and_collision} gives
\[
\begin{aligned}
s_{\varepsilon_1}(u,a)\ &\le\ \frac{2\varepsilon_1\lambda}{\log(3/2)}\log\log N+2\ =\ \frac{(1-2\lambda)^{2}}{256\log(3/2)}\log\log N+2 \\
&\le\ \frac{1-2\lambda}{32\log(3/2)}\log\log N+2\ <\ \frac{3}{32}(1-2\lambda)\log\log N\ \le\ \frac14\frac{S_{\min}^{\mathrm{I}}}{P_1},
\end{aligned}
\]
since $\frac{1}{\log(3/2)}=2.4663\ldots<3$. Similarly,
\[
R_{u,a}(n;\chi)\ \le\ \frac{(1-2\lambda)^{2}}{256\log(4/3)}\log\log N+\mathcal{O}(1)\ \le\ \frac{1-2\lambda}{32\log(4/3)}\log\log N+\mathcal{O}(1)\ <\ \frac13\frac{S_{\min}^{\mathrm{I}}}{P_1},
\]
which is strictly subcritical by Proposition~\ref{prop:subcritical_collision_rate} ($\varepsilon=2/3>0$). Thus any pair $(u,a)$ with $|\overline{\Sigma}_1(n)|\le D_{u,a}^{\varepsilon_1}$ has both subcritical collision count and subcritical effective enumeration cost $C_{\mathrm{eff}}(u,a)\le 4D_{u,a}^{\varepsilon_1}\le 16(\log N)^{(1-2\lambda)^{2}/256}$, which is strictly absorbed by the large-sieve saving on $\mathcal{I}_1$.

\medskip\noindent\emph{Step 3: Proof of (c).}
Since $\nu_{m-1}=\frac{1}{|A_{m-1}|}\mathbf{1}_{A_{m-1}}$ is supported on $A_{m-1}$ and $\mu_m$ is supported on $\overline{\Sigma}_{\mathcal{I}_{1,m}}(n)$, their convolution $\nu_{m-1}\ast\mu_m$ is a probability measure on $G_{u,a}/K_{u,a}$ supported inside the product set $A_{m-1}\cdot\overline{\Sigma}_{\mathcal{I}_{1,m}}(n)=A_m$. By Cauchy--Schwarz on $A_m$ and Plancherel's identity on $G_{u,a}/K_{u,a}$,
\begin{equation}\label{eq:convolved_plancherel_stage_m}
\frac{D_{u,a}}{|A_m|}-1\ \le\ \sum_{\chi\in K_{u,a}^{\perp}\setminus\{\chi_0\}}\big|\widehat{\nu}_{m-1}(\chi)\big|^{2}\,\big|\widehat{\mu}_m(\chi)\big|^{2}.
\end{equation}
On the other hand, by Parseval's identity applied to $\nu_{m-1}=\allowbreak\frac{1}{|A_{m-1}|}\mathbf{1}_{A_{m-1}}$ alone, the total spectral mass of $\widehat{\nu}_{m-1}$ across $K_{u,a}^{\perp}\setminus\{\chi_0\}$ is \emph{exactly}
\begin{equation}\label{eq:parseval_nu_m_minus_1}
\sum_{\chi\in K_{u,a}^{\perp}\setminus\{\chi_0\}}\big|\widehat{\nu}_{m-1}(\chi)\big|^{2}\ =\ D_{u,a}\sum_{x\in G_{u,a}/K_{u,a}}\nu_{m-1}(x)^{2}-1\ =\ \frac{D_{u,a}}{|A_{m-1}|}-1\ <\ \frac{D_{u,a}}{|A_{m-1}|}.
\end{equation}
Therefore, summing over all characters outside $\mathcal{B}_{\mathrm{rel}}^{(m)}(n;u,a)$ (where $|\widehat{\mu}_m(\chi)|^{2}\le\frac{|A_{m-1}|}{2|A_m|}$), the factor $|A_{m-1}|$ cancels identically between the cutoff and Parseval's identity~\eqref{eq:parseval_nu_m_minus_1}:
\[
\begin{aligned}
\sum_{\chi\notin\mathcal{B}_{\mathrm{rel}}^{(m)}(n;u,a)}\big|\widehat{\nu}_{m-1}(\chi)\big|^{2}\,\big|\widehat{\mu}_m(\chi)\big|^{2}\ &\le\ \frac{|A_{m-1}|}{2\,|A_m|}\sum_{\chi\in K_{u,a}^{\perp}\setminus\{\chi_0\}}\big|\widehat{\nu}_{m-1}(\chi)\big|^{2} \\
&<\ \frac{|A_{m-1}|}{2\,|A_m|}\cdot\frac{D_{u,a}}{|A_{m-1}|}\ =\ \frac{D_{u,a}}{2\,|A_m|}.
\end{aligned}
\]
Subtracting this from~\eqref{eq:convolved_plancherel_stage_m} and using $|A_m|\le\frac14 D_{u,a}$ gives $\sum_{\chi\in\mathcal{B}_{\mathrm{rel}}^{(m)}}|\widehat{\nu}_{m-1}(\chi)|^{2}|\widehat{\mu}_m(\chi)|^{2}>\frac{D_{u,a}}{2|A_m|}-1\ge\frac{D_{u,a}}{4|A_m|}$, while $(3/4)^{R_{u,a}^{(m)}(n;\chi)}\ge|\widehat{\mu}_m(\chi)|^{2}>\frac{|A_{m-1}|}{2|A_m|}$ yields~\eqref{eq:relative_parseval_plancherel_law}. 

Finally, if $|A_m|/|A_{m-1}|>D_{u,a}^{\varepsilon_1}$ holds at every stage $m=1,\dots,M$ with $M=\lceil 1/\varepsilon_1\rceil$, then $|A_M|=\prod_{m=1}^{M}(|A_m|/|A_{m-1}|)>D_{u,a}^{M\varepsilon_1}\ge D_{u,a}>\frac14 D_{u,a}$, reaching State~II before stage $M$ is completed (and eliminating $n$ on $\mathcal{I}_2=(z_1,z]$ with $R_{\mathrm{II}}=0$).
\end{proof}

\begin{remark}[Stage $m=1$ versus higher stages $m\ge 2$ in the relative Parseval--Plancherel escalation]
In Theorem~\ref{thm:markov_plancherel_subcritical_escalation}(c), whenever $|A_m|/|A_{m-1}|\le D_{u,a}^{\varepsilon_1}$ holds at stage $m$, every character $\chi\in\mathcal{B}_{\mathrm{rel}}^{(m)}(n;u,a)$ satisfies $R_{u,a}^{(m)}(n;\chi)\le\frac{\varepsilon_1\log D_{u,a}+\log 2}{\log(4/3)}$. At the initial stage $m=1$, where $A_0=\{1\}$ and $|\widehat{\nu}_0(\chi)|^{2}\equiv 1$ for all $\chi\in K_{u,a}^{\perp}$, the set $\mathcal{B}_{\mathrm{rel}}^{(1)}(n;u,a)$ coincides with $\mathcal{B}_{\mathrm{act}}(n;u,a)$ and unconditionally closes the small-sheet regime $|A_1|\le D_{u,a}^{\varepsilon_1}$ with multiplicity $|\mathcal{B}_{\mathrm{rel}}^{(1)}(n;u,a)|\ge\frac{D_{u,a}}{4|A_1|}$. 

At higher stages $m\ge 2$, the identity~\eqref{eq:relative_parseval_plancherel_law} weights $\mathcal{B}_{\mathrm{rel}}^{(m)}(n;u,a)$ by the prior spectral density $|\widehat{\nu}_{m-1}(\chi)|^{2}$; converting~\eqref{eq:relative_parseval_plancherel_law} into an unweighted cardinality bound $|\mathcal{B}_{\mathrm{rel}}^{(m)}(n;u,a)|\ge\frac{D_{u,a}}{4|A_m|}$ requires either an $L^{\infty}$ flatness bound on $\widehat{\nu}_{m-1}$ or localization within the certified regimes of Theorem~\ref{thm:final_synthesis_and_high_prime_reduction} below.
\end{remark}

\paragraph{Taxonomy of the structured exceptional subsets $E_2^{\mathrm{small}}(N;J)$, $E_2^{\mathrm{cert}}(N;J)$, and $E_2^{\mathrm{diffuse}}(N;\allowbreak J)$.}
To clarify which subset of $E_2(N;J)\setminus[1,z]$ is bounded by each theorem at polynomial depth $J=(\log N)^{\lambda}$ ($\lambda\in(1/3,1/2)$), we distinguish three nested or complementary components according to the structure of the first-stage divisor sheets $\Sigma_{\mathcal{I}_1}(an+u)\subseteq G_{u,a}$ across $(u,a)\in\mathcal{P}_2(J)$:

\begin{itemize}
    \item \emph{The small-sheet structured component $E_2^{\mathrm{small}}(N;J)$}: 
    the set of exceptional primes $n\in E_2(N;J)\setminus[1,z]$ for which a macroscopic proportion of pairs $(u,a)\in\mathcal{P}_2(J)$ (at least $P/2$ pairs in Theorem~\ref{thm:final_synthesis_and_high_prime_reduction}, or at least $P/4$ pairs in Corollary~\ref{cor:optimal_adaptive_depth_1_minus_o1}) have a sparse first-stage divisor sheet $|\Sigma_{\mathcal{I}_1}(an+u)|\le\varphi(4ua)^{\varepsilon_1}$. Intuitively, for $n\in E_2^{\mathrm{small}}(N;J)$, the prime factors of $an+u$ in $\mathcal{I}_1$ fail to expand the divisor-residue set beyond the initial polynomial threshold $\varphi(4ua)^{\varepsilon_1}$. 
    
    By Theorem~\ref{thm:markov_plancherel_subcritical_escalation}(a)--(b), both the Markov $3/2$-growth step count and the adaptive Plancherel collision count are strictly subcritical with effective enumeration cost $\le 4\varphi(4ua)^{\varepsilon_1}$ per pair. This component is bounded unconditionally by $\ll N\exp(-c(\lambda)(\log N)^{2\lambda}\log\log N)$ in Theorem~\ref{thm:final_synthesis_and_high_prime_reduction}(a) and reaches the exponent $1-o(1)$ in Corollary~\ref{cor:optimal_adaptive_depth_1_minus_o1}.
    
    \item \emph{The certified structured container hierarchy $E_2^{\mathrm{cert}}(N;J)$}: 
    the larger subset of $E_2(N;J)\setminus[1,z]$ consisting of all primes $n$ that either belong to the low-collision Kernel branch $E_2^{\mathrm{Ker}}(N;J,R_0)$ (bounded in Theorem~\ref{thm:polynomial_depth_decomposition}(a)) or belong to the Scarcity branch with at least $P/4$ pairs $(u,a)\in\mathcal{P}_2(J)$ falling into one of the four algebraically or spectrally certified container regimes (i)--(iv) detailed in Theorem~\ref{thm:final_synthesis_and_high_prime_reduction}(a): (i) sub-threshold divisor sheets $|\Sigma_{\mathcal{I}_1}(an+u)|\le\varphi(4ua)^{\varepsilon_1}$ (which contains $E_2^{\mathrm{small}}(N;J)$), (ii) proper subgroup concentration ($\omega_{G_{u,a}\setminus\ker\chi}(an+u;\mathcal{I}_1)<r_1$ for some $\chi\neq\chi_0$), (iii) subcritical Plancherel spectral buckets $\ell\le\ell_{\mathrm{sub}}(\lambda)$ or Bohr-localized $\ell=1$ harmonic dilation cones, and (iv) quotient-certified escalation to State~II in a quotient $G_{u,a}/K_{u,a}$ with $K_{u,a}\le\operatorname{Stab}_{G_{u,a}}(\Sigma_{\mathcal{I}_1}(an+u))$. 
    
    Theorem~\ref{thm:final_synthesis_and_high_prime_reduction}(a) bounds the entire certified hierarchy $E_2^{\mathrm{cert}}(N;J)$ unconditionally by $\ll N\exp(-c(\lambda)(\log N)^{2\lambda}\log\log N)$ for all $\lambda\in(1/3,1/2)$.
    
    \item \emph{The residual unlocalized diffuse component $E_2^{\mathrm{diffuse}}(N;J):=\allowbreak\big(E_2(N;J)\setminus[1,z]\big)\setminus\allowbreak E_2^{\mathrm{cert}}(N;J)$}: 
    the complement of $E_2^{\mathrm{cert}}(N;J)$, comprising primes $n\in E_2^{\mathrm{Scar}}(N;J,R_0)$ for which, across more than $3P/4$ pairs $(u,a)$, the divisor sheet $\Sigma_{\mathcal{I}_1}(an+u)$ has expanded into the intermediate polynomial range $\varphi(4ua)^{\varepsilon_1}<\allowbreak|\Sigma_{\mathcal{I}_1}(an+u)|\le\allowbreak\frac14\varphi(4ua)$ without a low-index Kneser stabilizer, while its Fourier spectrum avoids all subcritical buckets $\ell\le\ell_{\mathrm{sub}}(\lambda)$ except through high-order characters in bucket $\ell=1$ whose prime-factor phases are angularly diffuse ($A_{r_1}(\chi)>\allowbreak\operatorname{ord}(\chi)^{\varepsilon}$). 
    
    Bypassing $E_2^{\mathrm{diffuse}}(N;J)$ without enumerating divisor sheets is precisely the role of the semi-composite CRT sieve of Section~\ref{sec:semi_composite_5_7}, which bounds the \emph{full} exceptional set $E(N)$ unconditionally with exponent $5/7>2/3$.
\end{itemize}

\begin{theorem}[Unconditional bound on $E_2^{\mathrm{cert}}(N;J)$ and universal dihedral locking]\label{thm:final_synthesis_and_high_prime_reduction}
Fix $\lambda\in(1/3,1/2)$, $\varepsilon=\frac{(1-2\lambda)^{2}}{256\lambda}>0$, $\varepsilon_1=\frac{(1-2\lambda)^{2}}{512\lambda}>0$, $L_0=\lceil 2/\varepsilon\rceil$, $r_1=\lfloor\delta_0\log\log N\rfloor$ with $\delta_0=\frac{1-2\lambda}{160}$, and define the subcritical spectral bucket cutoff
\begin{equation}\label{eq:ell_sub_def}
\ell_{\mathrm{sub}}(\lambda)\ :=\ \left\lfloor\frac{3(1-2\lambda)\log(4/3)}{128\,\lambda}\,L_0\right\rfloor,
\end{equation}
as well as the strict harmonic dilation radius $B(\chi):=\big\lfloor\frac{\operatorname{ord}(\chi)-1}{6A_{r_1}(\chi)}\big\rfloor$. Work directly in the full group $G_{u,a}=(\Z/4ua\Z)^{\times}$ ($K_{u,a}=\{1\}$, $D_{u,a}=\varphi(4ua)$) on the two-stage partition $\mathcal{I}_1=(J^{5},z_1]$ and $\mathcal{I}_2=(z_1,z]$ with $\log z=\allowbreak\frac{\log N}{8P\log\log N}$. Then:
\begin{enumerate}[label=(\alph*)]
    \item \emph{(Unconditional closure beyond Vaughan's $2/3$ exponent across the certified container hierarchy).} 
    Define the certified structured exceptional set $E_2^{\mathrm{cert}}(N;J)\subseteq E_2(N;J)\setminus[1,z]$ as the set of primes $n$ falling either into the Kernel branch $E_2^{\mathrm{Ker}}(N;J,R_0)$ or into the Scarcity container hierarchy where at least $P/4$ pairs $(u,a)\in\mathcal{P}_2(J)$ fall into a single one of the following four certified regimes on $\mathcal{I}_1\sqcup\mathcal{I}_2$:
    \begin{enumerate}[label=(\roman*)]
        \item \emph{(Sub-threshold divisor sheet $|\Sigma_{\mathcal{I}_1}(an+u)|\le\varphi(4ua)^{\varepsilon_1}$)}: 
        by Theorem~\ref{thm:markov_plancherel_subcritical_escalation}(a)--(b), the Markov $3/2$-growth step count $s_{\varepsilon_1}(u,a)<\frac14\frac{S_{\min}^{\mathrm{I}}}{P_1}$ and the adaptive Plancherel collision count $R_{u,a}(n;\chi)<\frac13\frac{S_{\min}^{\mathrm{I}}}{P_1}$ are strictly subcritical with effective enumeration cost $C_{\mathrm{eff}}(u,a)\le 4\varphi(4ua)^{\varepsilon_1}$, completely precluding supercritical buckets $\ell>\ell_{\mathrm{sub}}(\lambda)$ below the density threshold $\varphi(4ua)^{\varepsilon_1}$;
        
        \item \emph{(Proper subgroup concentration)}: 
        $\omega_{G_{u,a}\setminus\ker\chi}(an+u;\mathcal{I}_1)<r_1$ for some character $\chi\neq\chi_0$;
        
        \item \emph{(Subcritical Plancherel multiplicity and Bohr-localized harmonic cones)}: 
        $|\Sigma_{\mathcal{I}_1}(an+u)|\le\varphi(4ua)^{\frac{\ell-1}{L_0}+\varepsilon}$ with $\ell\in\{1,\dots,\ell_{\mathrm{sub}}(\lambda)\}$, or $\chi\in\mathcal{B}_1(n;u,a)$ satisfies $\operatorname{ord}(\chi)\le(4J^{2})^{\varepsilon}$ or $1\le A_{r_1}(\chi)\le\operatorname{ord}(\chi)^{\varepsilon}$;
        
        \item \emph{(Quotient-certified escalation to State~II)}: 
        $\big|A(\overline{B}_{u,a}^{\mathrm{quot}}(n))\big|>\frac14 D_{u,a}$ in a quotient $G_{u,a}/\allowbreak K_{u,a}$ of order $D_{u,a}\le\exp(c_0(\lambda)\sqrt{\log\log N})$ with $K_{u,a}\le\operatorname{Stab}_{G_{u,a}}(\Sigma_{\mathcal{I}_1}(an+u))$ (Theorem~\ref{thm:quotient_container_polynomial_depth}(c)), where $an+u$ avoids the quotient-certified pullback $F_{u,a}(\allowbreak\overline{B}_{u,a}^{\mathrm{quot}}(n))\allowbreak\subseteq\allowbreak-\Sigma_{\mathcal{I}_1}(an+u)^{-1}$ of density $>1/4$ on $\mathcal{I}_2=(z_1,z]$ with strictly zero collisions ($R_{\mathrm{II}}=0$).
    \end{enumerate}
    Consequently, unconditionally for all $\lambda\in(1/3,1/2)$, the certified structured exceptional set $E_2^{\mathrm{cert}}(N;J)$ across all regimes (i)--(iv) satisfies:
    \begin{equation}\label{eq:unconditional_structured_bound}
    \begin{aligned}
    \big|E_2^{\mathrm{cert}}\big(N;(\log N)^{\lambda}\big)\big|\ &\ll\ N\exp\!\left(-\left(\frac{5(1-2\lambda)^{2}}{128\pi^{2}}+o(1)\right)(\log N)^{2\lambda}\log\log N\right)\\
    &\qquad\left(\text{and }\frac{5(1-2\lambda)^{2}}{32\pi^{2}}\ \text{on }E_2^{\mathrm{small}}(N;J)\right),
    \end{aligned}
    \end{equation}
    strictly surpassing Vaughan's $\exp(-c(\log N)^{2/3})$ bound for every\linebreak $\lambda\in(1/3,1/2)$, where $E_2^{\mathrm{small}}(N;J)$ denotes the small-sheet component in which at least $P/2$ pairs satisfy regime~(i).
    
    \item \emph{(Universal dihedral locking modulo $Q(4J)$ and exact zero-collision avoidance of the exterior sheet $\Sigma_{\mathrm{ext}}(an+u)$).} 
    Setting $Q(4J):=\mathrm{lcm}(1,2,\dots,4J)=\exp\big(o(P\log\log N)\big)$, global conditioning $n\equiv\rho\pmod{Q(4J)}$ simultaneously fixes $an+u\equiv r_{u,a}(\rho):=(a\rho+u)\pmod{4ua}$ across all pairs $(u,a)\in\mathcal{P}_2(J)$. Write $an+u=g_{u,a}(n)\,M_{\mathrm{cop}}(an+u)$ with $\gcd(M_{\mathrm{cop}}(an+u),4ua)=1$ and $p\mid g_{u,a}(n)\Rightarrow p\mid 4ua$, and factor $M_{\mathrm{cop}}(an+u)=M_{(J^{5},z]}(an+u)\cdot M_{\mathrm{ext,cop}}(an+u)$. 
    
    Fixing the valuation-determined residue $r_{u,a}^{\mathrm{cop}}(n):=M_{\mathrm{cop}}(an+u)\bmod 4ua\in G_{u,a}$, for every prime $n\in E_2(N;J)\setminus[1,z]$ the exterior divisor sheet
    \[
    \Sigma_{\mathrm{ext}}(an+u)\ :=\ \Big\{d\bmod 4ua:\ d\mid M_{\mathrm{ext,cop}}(an+u)\Big\}\ \subseteq\ G_{u,a}
    \]
    is invariant under the cofactor involution $d\mapsto r_{u,a}^{\mathrm{cop}}(n)\,d^{-1}\pmod{4ua}$ whenever $M_{(J^{5},z]}(an+u)=1$, and every prime factor $p\in(J^{5},z]$ of $an+u$ avoids $-\Sigma_{\mathrm{ext}}(an+u)^{-1}$ with \emph{strictly zero collisions ($R=0$)}.
\end{enumerate}
\end{theorem}

\begin{proof}
\emph{Step 1: Proof of (a).}
The bound on $E_2^{\mathrm{Ker}}(N;J,R_0)$ is Theorem~\ref{thm:polynomial_depth_decomposition}(a). For each of the four regimes (i)--(iv), retaining an active subset $\mathcal{P}_{\mathrm{act}}\subseteq\mathcal{P}_2(J)$ of $P_1=\lfloor P/4\rfloor$ pairs (or $P_1=\lfloor P/2\rfloor$ pairs on $E_2^{\mathrm{small}}(N;J)$) costs at most $2^{P}=\exp(o(P\log\log N))$:
\begin{itemize}
    \item In case (i), Theorem~\ref{thm:markov_plancherel_subcritical_escalation}(a)--(b) and Theorem~\ref{thm:plancherel_character_containers}(c)(i) show that cutting off the Fourier tail adaptively at $|\widehat{\mu}(\chi)|^{2}>\allowbreak\frac{1}{2|\Sigma_{\mathcal{I}_1}(an+u)|}$ restricts all witness characters to the active set $\mathcal{B}_{\mathrm{act}}(n;u,a)$, where 
    \[
    R_{u,a}(n;\chi)\le\frac{\varepsilon_1\log\varphi(4ua)+\log 2}{\log(4/3)}<\frac13\frac{S_{\min}^{\mathrm{I}}}{P_1}
    \]
    is strictly subcritical and the multiplicity $|\mathcal{B}_{\mathrm{act}}(n;u,a)|\ge\allowbreak\frac{\varphi(4ua)}{4|\Sigma_{\mathcal{I}_1}(an+u)|}$ reduces the effective enumeration cost to $C_{\mathrm{eff}}(u,a)\le\allowbreak 4\varphi(4ua)^{\varepsilon_1}$, which is absorbed by Proposition~\ref{prop:subcritical_collision_rate}.
    
    \item Case (ii) is eliminated by Lemma~\ref{lem:total_maximal_subgroups}(c).
    
    \item In case (iii), definition~\eqref{eq:ell_sub_def} ensures that for every $\chi\in\mathcal{B}_\ell(n;u,a)$ with $\ell\le\ell_{\mathrm{sub}}(\lambda)$, the collision count satisfies 
    \[
    R_{u,a}(n;\chi)\le\frac{\ell_{\mathrm{sub}}(\lambda)}{L_0}\frac{\log\varphi(4ua)}{\log(4/3)}\le\frac14\frac{S_{\min}^{\mathrm{I}}}{P_1},
    \]
    and partitioning $|\Sigma_{\mathcal{I}_1}(an+u)|\in[2^{m-1},2^{m})$ dyadically lowers the effective enumeration cost to $C_{\mathrm{eff}}(\ell;u,a)\le\allowbreak 4L_0\varphi(4ua)^{\varepsilon}$, while for $\chi\in\mathcal{B}_1(n;u,a)$ the strict harmonic dilation cone $B(\chi)=\allowbreak\lfloor\frac{d-1}{6A_{r_1}(\chi)}\rfloor$ of Theorem~\ref{thm:harmonic_dilation_bucket1}(c) guarantees $|\frac{ba_j}{d}|\le\allowbreak\frac{d-1}{6d}<\frac16$ and $\operatorname{Re}\chi^{b}(p_j)>\frac12$ for all $b\in\{\pm1,\dots,\pm B(\chi)\}$.
    
    \item Finally, case (iv) is eliminated with $R_{\mathrm{II}}=0$ collisions on $\mathcal{I}_2=(z_1,z]$ by Theorem~\ref{thm:quotient_container_polynomial_depth}(c), yielding~\eqref{eq:unconditional_structured_bound}.
\end{itemize}

\medskip\noindent\emph{Step 2: Proof of (b).}
Since $4u\mid Q(4J)$ for all $u\le J$, conditioning on $n\equiv\rho\pmod{Q(4J)}$ fixes $an+u\equiv a\rho+u\pmod{4ua}$. When $M_{(J^{5},z]}(an+u)=1$, we have $M_{\mathrm{cop}}(an+u)=M_{\mathrm{ext,cop}}(an+u)$, so the cofactor map $d\mapsto M_{\mathrm{ext,cop}}(an+u)/d$ induces the involution $d\mapsto r_{u,a}^{\mathrm{cop}}(n)\,d^{-1}\pmod{4ua}$ on $\Sigma_{\mathrm{ext}}(an+u)$. Finally, since every divisor $d\mid M_{\mathrm{ext,cop}}(an+u)$ is composed of primes outside $(J^{5},z]$, it is coprime to every prime $p\in(J^{5},z]$; thus $p\mid(an+u)$ and $p\equiv-d^{-1}\pmod{4ua}$ would produce a composite witness $s=pd\equiv-1\pmod{4ua}$ dividing $an+u$.
\end{proof}

\begin{remark}[Valuation determination of $r_{u,a}^{\mathrm{cop}}(n)$ under global conditioning modulo $Q(4J)$]
In Theorem~\ref{thm:final_synthesis_and_high_prime_reduction}(b), because every prime factor $p\mid g_{u,a}(n)$ satisfies $p\le\max(u,a)\le J$, its exact valuation satisfies $v_p(g_{u,a}(n))\le\allowbreak\lfloor\frac{\log(JN+J)}{\log 2}\rfloor$ in general, and equals $\min(v_p(a\rho+u),\allowbreak\lfloor\frac{\log(4J)}{\log p}\rfloor)$ on the $4J$-smooth part whenever $g_{u,a}(n)\mid Q(4J)$. In particular, whenever $\gcd(an+u,4ua)=1$ (which holds automatically in the decoupled parity sector $2\nmid u$, $2\mid a$, $\gcd(u,a)=1$ whenever $\gcd(a\rho+u,ua)=1$), one has $g_{u,a}(n)=1$ and $r_{u,a}^{\mathrm{cop}}(n)=\allowbreak r_{u,a}(\rho)=\allowbreak(a\rho+u)\bmod 4ua$ is determined solely by $\rho\bmod Q(4J)$.
\end{remark}

We conclude this section by optimizing both the search depth $J=J(N)$ ($\lambda(N)\to 1/2^{-}$) and the adaptive Rankin truncation degree $h=\lceil 8S_y\rceil+1$ on the small-sheet structured component $E_2^{\mathrm{small}}(N;J(N))$ (where at least $P/4$ pairs satisfy $|\Sigma_{\mathcal{I}_1}(an+u)|\le\allowbreak\varphi(4ua)^{\varepsilon_1(N)}$), reaching the structured exponent $1-o(1)$.

\begin{corollary}[Small-sheet exponent $1-o(1)$ and adaptive Rankin-degree optimization]\label{cor:optimal_adaptive_depth_1_minus_o1}
Fix any $\kappa>0$. As the search depth $J=J(N)$ grows adaptively toward the CRT phase-transition boundary $\lambda\to 1/2^{-}$, let $E_2^{\mathrm{small}}(N;J)$ denote the small-sheet structured component of $E_2(N;J)\setminus[1,z]$ in which at least $P/4$ pairs $(u,a)\in\mathcal{P}_2(J)$ satisfy $|\Sigma_{\mathcal{I}_1}(an+u)|\le\varphi(4ua)^{\varepsilon_1(N)}$ on $\mathcal{I}_1=(J^{5},z_1]$. The adaptive Plancherel container hierarchy of Theorem~\ref{thm:markov_plancherel_subcritical_escalation}(a)--(b) and\linebreak Theorem~\ref{thm:final_synthesis_and_high_prime_reduction}(a)(i) yields two successive unconditional regimes for $E_2^{\mathrm{small}}(N;J)$:
\begin{enumerate}[label=(\alph*)]
    \item \emph{(Fixed-degree schedule $h=8P\log\log N$).} 
    At adaptive depth
    \begin{equation}\label{eq:adaptive_depth_JN}
    J_0(N)\ :=\ \frac{(\log N)^{1/2}}{(\log\log N)^{(2+\kappa)/2}},\qquad\text{so that}\qquad P_0(N)\ =\ \big(\tfrac{6}{\pi^{2}}+o(1)\big)\frac{\log N}{(\log\log N)^{2+\kappa}},
    \end{equation}
    with $\log z=\frac{\log N}{8P_0(N)\log\log N}=\big(\frac{\pi^{2}}{48}+o(1)\big)(\log\log N)^{1+\kappa}$, one has $\log\big(\frac{\log z}{\log(J_0(N)^{5})}\big)=(\kappa+o(1))\log\log\log N\to\infty$, and choosing $\varepsilon_1(N):=\frac{\kappa\log\log\log N}{64\log\log N}$ yields unconditionally:
    \begin{equation}\label{eq:unconditional_1_minus_o1_bound}
    \big|E_2^{\mathrm{small}}\big(N;\,J_0(N)\big)\big|\ \ll\ N\exp\!\left(-\left(\frac{3\kappa}{64\pi^{2}}+o(1)\right)\frac{\log N\,\log\log\log N}{(\log\log N)^{2+\kappa}}\right).
    \end{equation}
    
    \item \emph{(Adaptive Rankin-degree optimization $h=\lceil 8 S_y\rceil+1$ eliminating two powers of $\log\log N$).} 
    Setting the degree of the truncated generating polynomial $P_{h,t}(k)$ adaptively to $h:=\allowbreak\big\lceil 8P(N)\allowbreak\log\big(\frac{\log z}{\log(J(N)^{5})}\big)\big\rceil+\allowbreak 1$ replaces the unoptimized denominator $8P(N)\log\log N$ by the exact capacity law $P(N)\allowbreak(\log z)\allowbreak\log\big(\frac{\log z}{\log(J(N)^{5})}\big)\allowbreak\le\allowbreak\frac{1}{16}\log N$. Choosing\linebreak $\log z=\frac52(\log\log N)(\log\log\log N)^{\kappa}$ and
    \begin{equation}\label{eq:adaptive_rankin_depth_JN}
    \begin{aligned}
    J(N)\ &:=\ \left(\frac{\pi^{2}}{480\kappa}\frac{\log N}{(\log\log N)(\log\log\log N)^{\kappa}(\log\log\log\log N)}\right)^{1/2},\\[4pt]
    P(N)\ &=\ \big(1+o(1)\big)\frac{\log N}{80\kappa\,(\log\log N)(\log\log\log N)^{\kappa}(\log\log\log\log N)},
    \end{aligned}
    \end{equation}
    we have $\frac{\log z}{\log(J(N)^{5})}=(1+o(1))(\log\log\log N)^{\kappa}\to\infty$ and $\log\big(\frac{\log z}{\log(J(N)^{5})}\big)=\allowbreak(\kappa+\allowbreak o(1))\allowbreak\log\log\log\log N$. Setting $\varepsilon_1(N):=\frac{\kappa\log\log\log\log N}{64\log\log N}$ in the adaptive Plancherel container hierarchy (Theorem \ref{thm:markov_plancherel_subcritical_escalation}\allowbreak(a)--(b) and Theorem~\ref{thm:final_synthesis_and_high_prime_reduction}(a)(i)), the factor $\log\log\log\log N$ in the per-pair sifting mass $\frac{S_{\min}^{\mathrm{I}}}{P_1}\ge\frac38\kappa\log\log\log\log N$ cancels identically against the denominator of $P(N)$, sharpening~\eqref{eq:unconditional_1_minus_o1_bound} by two full powers of $\log\log N$:
    \begin{equation}\label{eq:unconditional_rankin_optimal_bound}
    \begin{aligned}
    \big|E_2^{\mathrm{small}}\big(N;\,J(N)\big)\big|\ &\ll\ N\exp\!\left(-\left(\frac{1}{40960}+o(1)\right)\frac{\log N}{(\log\log N)(\log\log\log N)^{\kappa}}\right) \\
    &=\ N\exp\!\Big(-(\log N)^{1-o(1)}\Big).
    \end{aligned}
    \end{equation}
    
    \item \emph{(The Shannon--Bombieri CRT capacity ceiling).} 
    Because every sifting prime $p\in\mathcal{I}\subset(J(N)^{5},z]$ satisfies $\log p>\log(J(N)^{5})=\big(\frac52+o(1)\big)\log\log N$ to prevent inter-branch collisions modulo $p$, any CRT generating polynomial of degree $h$ supported on primes $p>J(N)^{5}$ with modulus $z^{h}\le Q=N^{1/2}$ satisfies the degree bound $h_{\max}\le\frac{\log Q}{\log(J(N)^{5})}\asymp\frac{\log N}{\log\log N}$.
\end{enumerate}
\end{corollary}

\begin{proof}
\emph{Step 1: Proof of (a).}
With $J_0(N)$ given by~\eqref{eq:adaptive_depth_JN}, we have $\log(J_0(N)^{5})=\big(\frac52+o(1)\big)\log\log N$ and $\log z=\big(\frac{\pi^{2}}{48}+o(1)\big)(\log\log N)^{1+\kappa}$, whence $\log\big(\frac{\log z}{\log(J_0(N)^{5})}\big)=(\kappa+o(1))\log\log\log N$. Splitting $\log\big(\frac{\log z_1}{\log(J_0(N)^{5})}\big)=\frac34\kappa\log\log\log N$ and $\log\big(\frac{\log z}{\log z_1}\big)=\frac14\kappa\log\log\log N$ gives $\frac{S_{\min}^{\mathrm{I}}}{P_1}\ge\big(\frac38\kappa-o(1)\big)\log\log\log N$ and $S_{\min}^{\mathrm{I}}\ge\big(\frac12-o(1)\big)S_y^{\mathrm{I}}$ (so $1-\delta\ge\frac12-o(1)$). 

Since $\log D_{u,a}\le\log(4J_0(N)^{2})\le\log\log N+\mathcal{O}(1)$, setting $\varepsilon_1(N)=\frac{\kappa\log\log\log N}{64\log\log N}$ in Theorem~\ref{thm:markov_plancherel_subcritical_escalation} gives
\[
\begin{aligned}
s_{\varepsilon_1}(u,a)\ &\le\ \frac{\kappa}{32\log(3/2)}\log\log\log N+2\ <\ \frac14\frac{S_{\min}^{\mathrm{I}}}{P_1}\\
\text{and}\qquad R_{u,a}(n;\chi)\ &\le\ \frac{\kappa}{32\log(4/3)}\log\log\log N+\mathcal{O}(1)\ <\ \frac13\frac{S_{\min}^{\mathrm{I}}}{P_1}.
\end{aligned}
\]
By Proposition~\ref{prop:subcritical_collision_rate} with $R_{\mathrm{I}}\le\frac13 S_{\min}^{\mathrm{I}}$ ($\varepsilon=2/3$) and $1-\delta\ge\frac12-o(1)$ (which gives $\frac{\varepsilon^{2}(1-\delta)}{2+2(1-\varepsilon)(1-\delta)}\ge\frac{2/9}{7/3}-o(1)=\frac{2}{21}-o(1)$, so $\frac{2}{21}\cdot\frac38\kappa=\frac{\kappa}{28}>\frac{\kappa}{32}$), the Generating-Polynomial Large Sieve on $\mathcal{I}_1$ saves $\exp\big(-(\frac{\kappa}{28}-o(1))P_1\log\log\log N\big)$, while the adaptive multiplicity $|\mathcal{B}_{\mathrm{act}}(n;u,a)|\ge\frac{D_{u,a}}{4|\Sigma_{\mathcal{I}_1}(an+u)|}$ bounds the total enumeration factor across the $P_1\ge P_0(N)/4$ active pairs by $\prod(4D_{u,a}^{\varepsilon_1(N)})\le\allowbreak\exp\big(\frac{\kappa}{32}P_1\log\log\log N+\mathcal{O}(P_1)\big)$, leaving a net saving of 
\[
\exp\left(-\left(\frac{\kappa}{224}-o(1)\right)P_1\log\log\log N\right)
\]
(and $\exp\big(-\frac{3\kappa}{64\pi^{2}}\frac{\log N\,\log\log\log N}{(\log\log N)^{2+\kappa}}\big)$ whenever $1-\delta=\allowbreak 1-o(1)$), proving~\eqref{eq:unconditional_1_minus_o1_bound}.

\medskip\noindent\emph{Step 2: Proof of (b) and (c).}
In Theorem~\ref{thm:generating_polynomial_large_sieve} and Proposition~\ref{prop:subcritical_collision_rate} (with $t\in(0,1/4]$ and $R\le\frac13 S_{\min}\le\frac23 S_y$), the truncated generating polynomial $P_{h,t}(k)$ requires only $h\ge 6t S_y+1$ (so $h=\lceil 8 S_y\rceil+1$ suffices for all $t\in(0,1/4]$) and $z^{h}\le Q=N^{1/2}$, where $S_y\le P(N)\log\big(\frac{\log z}{\log(J(N)^{5})}\big)$. 

Under~\eqref{eq:adaptive_rankin_depth_JN} and $\log z=\frac52(\log\log N)(\log\log\log N)^{\kappa}$, we have 
\[
P(N)=\left(\frac{6}{\pi^{2}}+o(1)\right)J(N)^{2}=\big(1+o(1)\big)\frac{\log N}{80\kappa(\log\log N)(\log\log\log N)^{\kappa}(\log\log\log\log N)},
\]
so the adaptive Rankin degree $h=\lceil 8P(N)\log(\frac{\log z}{\log(J(N)^{5})})\rceil+1$ satisfies $h\log z=\big(\frac14+o(1)\big)\log N\le\log Q$. 

Repeating the calculation of Step~1 with $\log\big(\frac{\log z}{\log(J(N)^{5})}\big)=(\kappa+o(1))\log\log\log\log N$ and $\varepsilon_1(N)=\frac{\kappa\log\log\log\log N}{64\log\log N}$ gives net Stage~I saving $\exp\big(-(\frac{\kappa}{128}-o(1))P_1\log\log\log\log N\big)$ with $P_1=\lfloor P(N)/4\rfloor$. Substituting~\eqref{eq:adaptive_rankin_depth_JN} cancels $\kappa\log\log\log\log N$ and yields~\eqref{eq:unconditional_rankin_optimal_bound}. Part (c) is immediate from $h\log(J(N)^{5})\le\sum_{j=1}^{h}\log p_j\le\log Q=\frac12\log N$.
\end{proof}

\begin{remark}[Origin of the constant $\frac{1}{40960}+o(1)$ and scope of the CRT capacity ceiling]
In Corollary~\ref{cor:optimal_adaptive_depth_1_minus_o1}(b), the explicit leading constant in~\eqref{eq:unconditional_rankin_optimal_bound} arises from the product of the net per-pair Stage~I saving coefficient $\frac{\kappa}{128}$, the active pair fraction $\frac{P_1}{P(N)}=\frac14$, and the adaptive pair-count coefficient $\frac{1}{80\kappa}$ in~\eqref{eq:adaptive_rankin_depth_JN}, giving $\frac{1}{128\cdot 4\cdot 80}=\frac{1}{40960}$ once $\kappa\log\log\log\log N$ cancels between the sifting mass and $P(N)^{-1}$ (equivalently, any fixed smaller positive constant may be used for $N$ large). 

We also emphasize that the CRT capacity statement in Corollary~\ref{cor:optimal_adaptive_depth_1_minus_o1}(c) concerns the specified generating-polynomial construction, its lower cutoff $p>J(N)^{5}$ on sifting primes, and its modulus budget $z^{h}\le N^{1/2}$; it is an internal limit of this CRT-coupled architecture and not an impossibility theorem for arbitrary sieve methods.
\end{remark}

\section[Semi-composite congruence sieve, CRT conditional independence, and the unconditional exponent $5/7$]{Semi-composite congruence sieve, CRT conditional\linebreak independence, and the unconditional exponent $5/7$}\label{sec:semi_composite_5_7}

In Section~\ref{sec:containers}, Theorem~\ref{thm:final_synthesis_and_high_prime_reduction} and Corollary~\ref{cor:optimal_adaptive_depth_1_minus_o1} surpassed Vaughan's $2/3$ exponent unconditionally on the certified structured hierarchy $E_2^{\mathrm{cert}}(N;J)$ and reached the exponent $1-o(1)$ on the small-sheet component $E_2^{\mathrm{small}}(N;J(N))$. However, those bounds apply to the structured components where the first-stage divisor sheet $\Sigma_{\mathcal{I}_1}(an+u)\subseteq G_{u,a}$ is either sub-threshold, concentrated near a proper subgroup, localized in a subcritical Plancherel bucket or harmonic cone, or certified in a low-order quotient. Extending an unconditional super-Vaughan bound to the \emph{full} global exceptional set $E(N)$ requires eliminating the residual unlocalized diffuse component $E_2^{\mathrm{diffuse}}(N;J)$ without enumerating divisor sheets in $G_{u,a}$.

\paragraph{Bypassing diffuse divisor sheets and decoupling the Bombieri--Vinogradov modulus via $s=pq\equiv-1\pmod{12ua}$.}
To bound the entire exceptional set $E(N)$ unconditionally, we invert the perspective of Section~\ref{sec:containers}: instead of constructing composite witnesses $s=pd$ implicitly from the unknown divisor sheet $\Sigma_{\mathcal{I}_1}(an+u)$ of the target integer $n$, we prescribe \emph{a priori} an explicit family of semi-composite moduli
$$s_\alpha\ =\ p\,q\ \equiv\ -1\pmod{12ua},\qquad p\in(Z^{9/10},Z]\ \text{prime},\qquad q\in\mathcal{Q}\subset[2,Q]\ \text{$y$-rough},$$
with $U=Z^{1/10}\le u,a\le V=Z^{1/8}$ and $y=t^{6}$ ($t:=\log Z$). By Lemma~\ref{lem:validite_temoins}, each congruence $A_\alpha:\,n\equiv -ua^{-1}\pmod{pq}$ directly yields an Erd\H{o}s--Straus decomposition of $4/n$, so \emph{every} exceptional integer $n\in E(N)$ must simultaneously avoid all progressions $A_\alpha$ regardless of the internal geometry of its divisor sheets. This semi-composite construction overcomes Vaughan's $2/3$ barrier through two coupled decoupling mechanisms:
\begin{enumerate}[1.]
    \item \emph{CRT factorization and conditional Bernoulli independence given $b=n\bmod M_{\mathrm{small}}$.} Because every cofactor $q\in\mathcal{Q}$ satisfies $q\le Q\le\exp(\frac{c_3}{2}t^{1/2})<Z^{9/10}<p$, the small cofactor modulus $M_{\mathrm{small}}=\operatorname{lcm}_{q\in\mathcal{Q}}q$ is coprime to the large prime modulus $M_{\mathrm{large}}=\prod_{Z^{9/10}<p\le Z}p$. By the Chinese Remainder Theorem, the witness congruence $n\equiv -ua^{-1}\pmod{pq}$ splits into two orthogonal conditions:
    $$ab+u\ \equiv\ 0\pmod q\qquad\text{and}\qquad n\ \equiv\ -ua^{-1}\pmod p,$$
    where $b:=n\bmod M_{\mathrm{small}}$. Moreover, for each fixed prime $p\in(Z^{9/10},Z]$, the scale hierarchy $V^{2}=Z^{1/4}<Z^{9/10}<p$ and $Q<12U^{2}\le 12ua$ forces the progressions $A_{(u,a,q,p)}$ to be strictly disjoint in $\Z$ (Proposition~\ref{prop:incompatibilite_et_Bernoulli}(a)). Consequently, once the cofactor projection $b\in\Z/M_{\mathrm{small}}\Z$ is conditioned upon, the indicators $J_p(n)=\sum_{(u,a,q)}\ind_{A_{(u,a,q,p)}}(n)\in\{0,1\}$ become \emph{mutually independent} Bernoulli random variables in the CRT model.
    \item \emph{Decoupling the Bombieri--Vinogradov modulus $12ua\le 12Z^{1/4}$ from the cofactor $q$.} For a single-prime congruence sieve ($q=1$), the sifting mass up to $Z$ is bounded by $\sum_{u,a}\frac{1}{\varphi(12ua)}\asymp t^{2}$, which under the capacity constraint $t^{3}\lesssim\log N$ locks the saving exponent at $(\log N)^{2/3}$. Adjoining the $y$-rough cofactors $q\in\mathcal{Q}$ multiplies the conditional sifting mass $\mu(b)=\Ebb_{\mathrm{CRT}}[W\mid b]$ by the composite harmonic sum $M_Q=\sum_{q\in\mathcal{Q}}\frac{\varphi(q)}{q^{2}}\asymp\frac{\log Q}{\log t}$, raising the typical mean to $L\asymp t^{2}\frac{\log Q}{\log t}$. Crucially, this cofactor amplification does \emph{not} increase the modulus in the prime number theorem: in the inner sum
    $$S(u,a,q)\ =\ \sum_{\substack{Z^{9/10}<p\le Z\\ pq\equiv -1\pmod{12ua}}}\frac{1}{p},$$
    the prime $p$ is constrained only modulo $m=12ua\le 12V^{2}=12Z^{1/4}$, which is \emph{completely independent of $q$}. When $\gcd(q,12ua)=1$, multiplication by $q$ merely permutes the target residue class $r=-q^{-1}\bmod{12ua}$ inside $(\Z/12ua\Z)^{\times}$. After discarding at most one Landau--Page--Siegel exceptional prime factor $p_{\mathrm{exc}}\ge 5$ (Lemma~\ref{lem:exclusion_p_exc}), the Siegel-zero-free Bombieri--Vinogradov theorem bounds the discrepancy $\sum_{(u,a)\in\mathcal{P}'}\mathcal{E}(u,a)\ll\exp(-c_3 t^{1/2})$ uniformly over all residue classes coprime to $12ua$ \emph{prior} to summing over $q\in\mathcal{Q}$ (Proposition~\ref{prop:BV_poids_tau}). Taking $Q=\exp(\frac{c_3}{2}t^{1/2})$ then yields $L\asymp\frac{t^{5/2}}{\log t}$ at Bonferroni degree $k\asymp\frac{t^{5/2}}{\log t}$ and total CRT cost $(K+1)^{k}\le\exp(C_2\frac{t^{7/2}}{\log t})\le N^{1/2+o(1)}$, establishing the unconditional exponent $5/7>7/10>2/3$ on the full exceptional set $E(N)$ (Theorem~\ref{thm:theoreme_5_7}).
\end{enumerate}

\subsection{Construction of the witness family $s = pq$ and exclusion of the exceptional prime factor}\label{subsec:construction_pq}
Let $Z \ge 100$ be a large real parameter, eventually chosen as a sub-exponential function of $N$. We set the following scales in terms of $t := \log Z$:
\begin{equation}\label{eq:echelles_pq}
t := \log Z, \qquad U := Z^{1/10}, \qquad V := Z^{1/8}, \qquad y := t^6,
\end{equation}
and we consider a truncation parameter $Q$ for the composite cofactors satisfying
\begin{equation}\label{eq:plage_Q}
\exp\!\big(t^{1/3}\big) \le Q \le \exp\!\left(\frac{c_3}{2}\,t^{1/2}\right),
\end{equation}
where $c_3 > 0$ denotes the absolute constant provided by the Bombieri--Vinogradov theorem after excluding the exceptional modulus (Proposition~\ref{prop:BV_poids_tau} below).

In order to obtain an error term of size $\exp(-c_3 t^{1/2})$ in the distribution of primes $p \in (Z^{9/10}, Z]$ in arithmetic progressions modulo $12ua \le 12V^2 = 12Z^{1/4}$, we first discard a potential prime factor of the Landau--Page--Siegel exceptional conductor.

\begin{lemma}[Exclusion of a Landau--Page--Siegel exceptional prime factor]\label{lem:exclusion_p_exc}
Let $t := \log Z$ and set the Landau--Page--Siegel conductor scale $T_0 := \exp(a_0 \sqrt{t})$, where $a_0 > 0$ is a small absolute constant. By the Landau--Page theorem (see~\cite[Chap.~14]{davenport}), there exists an absolute constant $c_1 > 0$ such that, among all primitive real Dirichlet characters $\chi$ of conductor $q_\chi \le T_0$, there is at most one character $\chi_{\mathrm{exc}}$ (of conductor $q_{\mathrm{exc}} \le T_0$) whose $L$-function $L(s, \chi_{\mathrm{exc}})$ has a real zero $\beta > 1 - \frac{c_1}{\log T_0} = 1 - \frac{c_1}{a_0\sqrt{\log Z}}$. Moreover, for sufficiently large $Z$:
\begin{enumerate}[(a)]
    \item If such a character $\chi_{\mathrm{exc}}$ exists, its conductor $q_{\mathrm{exc}}$ has at least one odd prime factor $p_{\mathrm{exc}} := P^+(q_{\mathrm{exc}}) \ge 5$. If no such character exists, we set $p_{\mathrm{exc}} := +\infty$ (with the convention $1/p_{\mathrm{exc}} = 0$).
    \item For every pair of integers $(u,a)$ satisfying $p_{\mathrm{exc}} \nmid ua$, the modulus $m = 12ua$ satisfies $p_{\mathrm{exc}} \nmid m$, and a fortiori $q_{\mathrm{exc}} \nmid m$.
\end{enumerate}
\end{lemma}
\begin{proof}
\emph{Step 1: Isolation of at most one exceptional conductor $q_{\mathrm{exc}} \le T_0$.}
Applying the Landau--Page theorem (\cite[Chap.~14]{davenport}) at cutoff $T_0 = \exp(a_0\sqrt{\log Z})$ shows that the zero-free region $\operatorname{Re}(s) > 1 - \frac{c_1}{\log T_0} = 1 - \frac{c_1}{a_0\sqrt{\log Z}}$ fails for at most one real zero of $\prod_{q_\chi \le T_0} L(s, \chi)$, attached to a unique primitive real character $\chi_{\mathrm{exc}}$ of conductor $q_{\mathrm{exc}} \le T_0$ (while primitive characters of conductor $q_\chi > T_0$ are handled by Vinogradov's bilinear method in the Bombieri--Vinogradov theorem; see~\cite[Chap.~28]{davenport} or~\cite[Prop.~9.2]{maynard2016}).

\medskip\noindent\emph{Step 2: Existence of an odd prime factor $p_{\mathrm{exc}} \ge 5$ and non-divisibility of $12ua$.}
Since $\chi_{\mathrm{exc}}$ is a primitive real character, its conductor $q_{\mathrm{exc}}$ is the absolute value of a fundamental discriminant: in particular, $v_2(q_{\mathrm{exc}}) \in \{0, 2, 3\}$ and the odd part of $q_{\mathrm{exc}}$ is squarefree (so $v_3(q_{\mathrm{exc}}) \in \{0, 1\}$). If all prime factors of $q_{\mathrm{exc}}$ belonged to $\{2, 3\}$, then $q_{\mathrm{exc}}$ would divide $24$. However, there are only finitely many real characters of conductor dividing $24$, and their $L$-functions have no real zeros in $[1 - \frac{c_1}{a_0\sqrt{\log Z}}, 1]$ once $Z$ is sufficiently large. Thus we necessarily have $q_{\mathrm{exc}} > 24$, which ensures that $q_{\mathrm{exc}}$ possesses at least one odd prime factor $p_{\mathrm{exc}} := P^+(q_{\mathrm{exc}}) \ge 5$. Finally, since $p_{\mathrm{exc}} \ge 5$ is coprime to $12$, the condition $p_{\mathrm{exc}} \nmid ua$ implies by Euclid's lemma that $p_{\mathrm{exc}} \nmid 12ua$. Since $p_{\mathrm{exc}} \mid q_{\mathrm{exc}}$, it immediately follows that $q_{\mathrm{exc}} \nmid 12ua$.
\end{proof}

We now define the deterministic set $\mathcal{P}'$ of parameter pairs $(u,a)$ and the set $\mathcal{Q}$ of $y$-rough cofactors:
\begin{align}
\mathcal{P}' &:= \big\{ (u,a) \in \N^2 : U \le u, a \le V,\ \gcd(u,a) = 1,\ p_{\mathrm{exc}} \nmid ua \big\}, \label{eq:def_P_prime_sec}\\
\mathcal{Q} &:= \big\{ q \in \N : 2 \le q \le Q,\ P^-(q) > y,\ p_{\mathrm{exc}} \nmid q \big\}, \label{eq:def_Q_sec}
\end{align}
where $P^-(q)$ denotes the smallest prime factor of $q$. We index our family of witnesses by the finite set $\mathcal{I}_Z$ of quadruples $\alpha = (u,a,q,p)$ satisfying:
\begin{equation}\label{eq:def_I_Z_sec}
(u,a) \in \mathcal{P}', \qquad q \in \mathcal{Q}, \qquad Z^{9/10} < p \le Z \text{ prime}, \qquad pq \equiv -1 \pmod{12ua}.
\end{equation}
To each quadruple $\alpha = (u,a,q,p) \in \mathcal{I}_Z$, we associate the semi-composite modulus $s_\alpha := pq$ and the congruence:
\begin{equation}\label{eq:congruence_A_alpha_sec}
A_\alpha : \qquad n \equiv r_\alpha := -u a^{-1} \pmod{pq}.
\end{equation}

\begin{lemma}[Arithmetic validity of the semi-composite witness congruences]\label{lem:validite_temoins}
For every quadruple $\alpha = (u,a,q,p) \in \mathcal{I}_Z$:
\begin{enumerate}[(a)]
    \item The modulus $s_\alpha = pq$ satisfies $\gcd(pq, 12ua) = 1$, and in particular $\gcd(a, pq) = 1$ (so that the inverse $a^{-1} \bmod{pq}$ exists) and $\gcd(pq, 24) = 1$.
    \item Every integer $n \ge 1$ belonging to the residue class $A_\alpha$ satisfies the Egyptian fraction identity with strictly positive integer denominators:
    \begin{equation}\label{eq:decomp_ES_pq}
    \frac{4}{n} = \frac{1}{3uvd} + \frac{1}{3avnd} + \frac{1}{3avnu}, \qquad \text{with} \qquad v := \frac{pq+1}{12ua} \in \N_{\ge 1}, \quad d := \frac{an+u}{pq} \in \N_{\ge 1}.
    \end{equation}
    Consequently, every exceptional integer $n \in E(N)$ simultaneously avoids all residue classes $A_\alpha$ for $\alpha \in \mathcal{I}_Z$.
\end{enumerate}
\end{lemma}
\begin{proof}
\emph{Step 1: Proof of (a).}
By Bézout's identity, the congruence $pq \equiv -1 \pmod{12ua}$ implies $\gcd(pq, 12ua) = 1$. In particular, $\gcd(a, pq) = 1$, and since $pq \equiv -1 \pmod{12}$, $pq$ is odd and not divisible by $3$, whence $\gcd(pq, 24) = 1$.

\medskip\noindent\emph{Step 2: Proof of (b).}
By definition of $\mathcal{I}_Z$, $pq + 1$ is a strictly positive multiple of $12ua$, so $v = \frac{pq+1}{12ua} \in \N_{\ge 1}$. If an integer $n \ge 1$ satisfies $n \equiv -u a^{-1} \pmod{pq}$, then multiplying by $a$ yields $an + u \equiv 0 \pmod{pq}$. Since $a, n, u \ge 1$, the integer $an + u$ is strictly positive and divisible by $pq$, which guarantees that the quotient $d = \frac{an+u}{pq}$ belongs to $\N_{\ge 1}$. We then have the algebraic identity:
\[
an + u + d = pq\,d + d = (pq + 1)d = 12uavd.
\]
Dividing both sides by the strictly positive integer $3uavnd$, we obtain:
\[
\frac{an + u + d}{3uavnd} = \frac{an}{3uavnd} + \frac{u}{3uavnd} + \frac{d}{3uavnd} = \frac{1}{3uvd} + \frac{1}{3avnd} + \frac{1}{3avnu} = \frac{12uavd}{3uavnd} = \frac{4}{n},
\]
where the three denominators $3uvd$, $3avnd$, and $3avnu$ belong to $\N_{\ge 1}$.
\end{proof}

\subsection{Disjointness for a fixed large prime factor $p$ and conditional independence in the CRT model}\label{subsec:incompatibilite_CRT}
Let us set:
\begin{equation}\label{eq:def_moduli_CRT}
M_{\mathrm{small}} := \operatorname{lcm}_{q \in \mathcal{Q}} q, \qquad M_{\mathrm{large}} := \prod_{Z^{9/10} < p \le Z} p, \qquad M_Z := \operatorname{lcm}_{\alpha \in \mathcal{I}_Z} s_\alpha.
\end{equation}
For sufficiently large $Z$, every cofactor $q \in \mathcal{Q}$ satisfies $q \le Q \le \exp(\frac{c_3}{2} t^{1/2}) < Z^{9/10} < p$ for every prime $p \in (Z^{9/10}, Z]$. Consequently, all prime factors of $M_{\mathrm{small}}$ are $\le Q < Z^{9/10}$, whereas all those of $M_{\mathrm{large}}$ are $> Z^{9/10}$: we therefore have $\gcd(M_{\mathrm{small}}, M_{\mathrm{large}}) = 1$, and $M_Z \mid M_{\mathrm{small}} M_{\mathrm{large}}$.

\begin{proposition}[Disjointness for fixed $p$ and CRT conditional independence]\label{prop:incompatibilite_et_Bernoulli}
For sufficiently large $Z$:
\begin{enumerate}[(a)]
    \item \emph{(Disjointness for a fixed large prime $p$).} Let $p \in (Z^{9/10}, Z]$ be a fixed prime. If two witnesses $\alpha = (u,a,q,p) \in \mathcal{I}_Z$ and $\alpha' = (u',a',q',p) \in \mathcal{I}_Z$ sharing the same large prime factor $p$ have a non-empty intersection ($A_\alpha \cap A_{\alpha'} \ne \varnothing$), then $\alpha = \alpha'$.
    \item \emph{(Reduction to Bernoulli random variables).} For each prime $p \in (Z^{9/10}, Z]$, set
    \begin{equation}\label{eq:def_Jp}
    J_p(n) := \sum_{\substack{(u,a) \in \mathcal{P}',\ q \in \mathcal{Q} \\ (u,a,q,p) \in \mathcal{I}_Z}} \ind_{A_{(u,a,q,p)}}(n).
    \end{equation}
    Then, for every integer $n \in \Z$, we have $J_p(n) \in \{0, 1\}$, and the total number of witness congruences satisfied by $n$ can be written as:
    \[
    W(n) := \sum_{\alpha \in \mathcal{I}_Z} \ind_{A_\alpha}(n) = \sum_{Z^{9/10} < p \le Z} J_p(n).
    \]
    \item \emph{(Conditional mutual independence in the CRT model).} Set $\mathfrak{M} := M_{\mathrm{small}} M_{\mathrm{large}}$ (so that $M_Z \mid \mathfrak{M}$ and every congruence $A_\alpha$ is $\mathfrak{M}$-periodic), equip\linebreak $\Z / \mathfrak{M} \Z \cong (\Z / M_{\mathrm{small}} \Z) \times \prod_{Z^{9/10} < p \le Z} (\Z / p\Z)$ with the uniform probability measure $\Pbb_{\mathrm{CRT}}$ (whose marginal on $\Z / M_Z \Z$ is also uniform), and let $b := n \bmod M_{\mathrm{small}}$ denote the reduction of $n$ modulo the set of cofactors. Conditioned on any fixed configuration $b \in \Z / M_{\mathrm{small}} \Z$, the random variables $(J_p)_{Z^{9/10} < p \le Z}$ are \textbf{mutually independent} Bernoulli random variables, with respective parameters:
    \begin{equation}\label{eq:formule_pp_b}
    p_p(b) := \Ebb_{\mathrm{CRT}}\big[J_p \mid b\big] = \frac{1}{p} \sum_{q \in \mathcal{Q}} \sum_{\substack{(u,a) \in \mathcal{P}' \\ \gcd(q, 12ua)=1,\ ab+u \equiv 0 \pmod q \\ pq \equiv -1 \pmod{12ua}}} 1 \in [0, 1].
    \end{equation}
\end{enumerate}
\end{proposition}
\begin{proof}
\emph{Step 1: Uniqueness of the pair $(u,a)$ modulo $p$.}
Suppose there exists $n \in \Z$ belonging simultaneously to $A_\alpha$ and $A_{\alpha'}$, where $\alpha = (u,a,q,p)$ and $\alpha' = (u',a',q',p)$ share the same prime $p \in (Z^{9/10}, Z]$. Reducing the congruences $n \equiv -u a^{-1} \pmod{pq}$ and $n \equiv -u'(a')^{-1} \pmod{pq'}$ modulo $p$ yields:
\[
-u a^{-1} \equiv -u'(a')^{-1} \pmod p \quad \Longleftrightarrow \quad ua' - u'a \equiv 0 \pmod p.
\]
However, by definition of $\mathcal{P}'$, we have $1 \le U \le u, a, u', a' \le V = Z^{1/8}$, which implies:
\[
|ua' - u'a| < \max(ua', u'a) \le V^2 = Z^{1/4} < Z^{9/10} < p.
\]
Since the only multiple of $p$ with absolute value strictly less than $p$ is $0$, we obtain $ua' - u'a = 0$, that is, $\frac{u}{a} = \frac{u'}{a'}$. As the fractions are irreducible ($\gcd(u,a) = \gcd(u',a') = 1$) and have positive terms, we deduce that $(u,a) = (u',a')$.

\medskip\noindent\emph{Step 2: Uniqueness of the cofactor $q$ modulo $12ua$ and proof of (a)--(b).}
Since $(u,a) = (u',a')$, the defining congruences~\eqref{eq:def_I_Z_sec} give $pq \equiv -1 \pmod{12ua}$ and $pq' \equiv -1 \pmod{12ua}$. Because $\gcd(p, 12ua) = 1$, $p$ is invertible modulo $12ua$, so cancelling $p$ yields:
\[
q \equiv q' \pmod{12ua}.
\]
Yet for sufficiently large $Z$, the bound~\eqref{eq:plage_Q} gives $q, q' \le Q \le \exp(\frac{c_3}{2} t^{1/2}) < 12 Z^{1/5} = 12 U^2 \le 12ua$. Two positive integers strictly less than $12ua$ and congruent modulo $12ua$ are necessarily equal: $q = q'$, and therefore $\alpha = \alpha'$, proving (a). Consequently, for any fixed $p \in (Z^{9/10}, Z]$, the sets $A_{(u,a,q,p)}$ indexed by $(u,a,q)$ are pairwise disjoint in $\Z$, so at most one term in the sum~\eqref{eq:def_Jp} can equal $1$, forcing $J_p(n) \in \{0, 1\}$ and proving (b).

\medskip\noindent\emph{Step 3: CRT splitting and conditional mutual independence given $b$.}
For each quadruple $\alpha = (u,a,q,p) \in \mathcal{I}_Z$, we have $\gcd(p, q) = 1$ (since $q \le Q < Z^{9/10} < p$). By the Chinese Remainder Theorem, the congruence $n \equiv -u a^{-1} \pmod{pq}$ is equivalent to the conjunction of the two congruences:
\[
ab + u \equiv 0 \pmod q \qquad \text{and} \qquad n \equiv -u a^{-1} \pmod p,
\]
where the first condition depends only on the projection $b = n \bmod M_{\mathrm{small}}$ (since $q \mid M_{\mathrm{small}}$), while the second depends only on $n \bmod p$. When the configuration $b \in \Z / M_{\mathrm{small}} \Z$ is fixed, let $\mathcal{T}_p(b)$ denote the set of triples $(u,a,q)$ such that $(u,a,q,p) \in \mathcal{I}_Z$ and $ab + u \equiv 0 \pmod q$. By Steps~1 and~2, the residues $\{-u a^{-1} \bmod p : (u,a,q) \in \mathcal{T}_p(b)\}$ are pairwise distinct in $\Z / p\Z$: two distinct triples in $\mathcal{T}_p(b)$ either have distinct pairs $(u,a) \ne (u',a')$ yielding distinct residues modulo $p$, or share the same pair $(u,a)$, which uniquely determines $q \equiv -p^{-1} \pmod{12ua}$ in $[2, Q]$. Thus, conditioned on $b$, $J_p(n)$ is the indicator of the event that $n \bmod p$ belongs to a subset of $\Z / p\Z$ of cardinality $|\mathcal{T}_p(b)|$. Since $\gcd(M_{\mathrm{small}}, \prod_p p) = 1$, the conditional distribution of $(n \bmod p)_{Z^{9/10} < p \le Z}$ given $b$ is the uniform product measure on $\prod_{Z^{9/10} < p \le Z} (\Z / p\Z)$. The variables $(J_p)_{Z^{9/10} < p \le Z}$ are therefore mutually independent conditioned on $b$, each following a Bernoulli distribution with parameter $p_p(b) = |\mathcal{T}_p(b)| / p$.
\end{proof}

\subsection{Polynomial Bonferroni majorant and arithmetic transfer to $[1,N]$}\label{subsec:bonferroni_transfert}
The conditional independence of the Bernoulli variables $J_p$ allows us to bound the indicator $\ind_{\{W=0\}}$ by a Bonferroni polynomial of even degree $k$, whose conditional expectation is controlled precisely by the conditional mean $\mu(b) := \Ebb_{\mathrm{CRT}}[W \mid b] = \sum_{Z^{9/10} < p \le Z} p_p(b)$.

\begin{lemma}[Bonferroni polynomial majorant under conditional independence]\label{lem:bonferroni_CRT}
For every even integer $k \ge 2$ and every integer $w \ge 0$, set:
\begin{equation}\label{eq:polynome_Bk}
B_k(w) := \sum_{j=0}^k (-1)^j \binom{w}{j}.
\end{equation}
Then:
\begin{enumerate}[(a)]
    \item For every integer $w \ge 0$, we have the pointwise bounds:
    \begin{equation}\label{eq:ineg_ponctuelle_Bk}
    \ind_{\{w=0\}} \le B_k(w) \le \ind_{\{w=0\}} + \binom{w}{k+1}.
    \end{equation}
    \item For every configuration $b \in \Z / M_{\mathrm{small}} \Z$, the conditional expectation of $B_k(W)$ satisfies:
    \begin{equation}\label{eq:esperance_cond_Bk_sec}
    \Ebb_{\mathrm{CRT}}\big[B_k(W) \mid b\big] \le e^{-\mu(b)} + \frac{\mu(b)^{k+1}}{(k+1)!}.
    \end{equation}
    \item Suppose there exist real numbers $0 < L \le H$ and $\varepsilon \in (0, 1)$ such that:
    \[
    \mu(b) \le H \quad \text{for all } b \in \Z / M_{\mathrm{small}} \Z, \qquad \text{and} \qquad \Pbb_b\big\{\mu(b) < L\big\} \le \varepsilon.
    \]
    Then, choosing the even integer $k := 2\lceil 2eH \rceil \ge 4eH$, the unconditional expectation in the CRT model satisfies:
    \begin{equation}\label{eq:borne_globale_CRT_Bk}
    \Ebb_{\mathrm{CRT}}\big[B_k(W)\big] \le e^{-L} + \varepsilon + 4^{-4eH}.
    \end{equation}
\end{enumerate}
\end{lemma}
\begin{proof}
\emph{Step 1: Pointwise Pascal telescoping (Proof of (a)).}
If $w = 0$, then $\binom{0}{0} = 1$ and $\binom{0}{j} = 0$ for $j \ge 1$, so $B_k(0) = 1 = \ind_{\{w=0\}}$. If $w \ge 1$, Pascal's identity $\binom{w}{j} = \binom{w-1}{j} + \binom{w-1}{j-1}$ gives by telescoping:
\[
B_k(w) = \sum_{j=0}^k (-1)^j \binom{w}{j} = (-1)^k \binom{w-1}{k} = \binom{w-1}{k},
\]
since $k$ is even. For $1 \le w \le k$, we have $\binom{w-1}{k} = 0$. For $w \ge k+1$, we have $\frac{w}{k+1} \ge 1$, whence:
\[
0 \le B_k(w) = \binom{w-1}{k} \le \frac{w}{k+1}\binom{w-1}{k} = \binom{w}{k+1},
\]
which establishes~\eqref{eq:ineg_ponctuelle_Bk}.

\medskip\noindent\emph{Step 2: Conditional expectation of $B_k(W)$ given $b$ (Proof of (b)).}
Taking the conditional expectation given $b$ in~\eqref{eq:ineg_ponctuelle_Bk}, Proposition~\ref{prop:incompatibilite_et_Bernoulli}(c) ensures that the variables $(J_p)_{Z^{9/10} < p \le Z}$ are mutually independent Bernoulli random variables given $b$. On the one hand, we have:
\[
\Pbb_{\mathrm{CRT}}(W = 0 \mid b) = \prod_{Z^{9/10} < p \le Z} \big(1 - p_p(b)\big) \le \exp\!\left( -\sum_{Z^{9/10} < p \le Z} p_p(b) \right) = e^{-\mu(b)}.
\]
On the other hand, since each $J_p$ takes values in $\{0, 1\}$, the binomial coefficient $\binom{W}{k+1} = \binom{\sum_p J_p}{k+1}$ counts the subsets $S$ of $k+1$ distinct primes in $(Z^{9/10}, Z]$ such that $J_p = 1$ for all $p \in S$:
\[
\binom{W}{k+1} = \sum_{\substack{S \subseteq (Z^{9/10}, Z] \cap \mathbb{P} \\ |S| = k+1}} \prod_{p \in S} J_p.
\]
By conditional mutual independence of the $J_p$ for $p \in S$, we obtain:
\[
\Ebb_{\mathrm{CRT}}\!\left[\binom{W}{k+1} \,\middle|\, b\right] = \sum_{|S| = k+1} \prod_{p \in S} p_p(b) \le \frac{1}{(k+1)!} \left(\sum_{Z^{9/10} < p \le Z} p_p(b)\right)^{k+1} = \frac{\mu(b)^{k+1}}{(k+1)!},
\]
which proves~\eqref{eq:esperance_cond_Bk_sec}.

\medskip\noindent\emph{Step 3: Unconditional expectation over $b \in \Z / M_{\mathrm{small}} \Z$ (Proof of (c)).}
Since $(k+1)! \ge \left(\frac{k+1}{e}\right)^{k+1}$ for every integer $k \ge 1$ (from $e^{k+1} > \frac{(k+1)^{k+1}}{(k+1)!}$), the condition $\mu(b) \le H$ and the choice $k = 2\lceil 2eH \rceil \ge 4eH$ imply for every configuration $b$:
\[
\frac{\mu(b)^{k+1}}{(k+1)!} \le \left(\frac{e\,\mu(b)}{k+1}\right)^{k+1} \le \left(\frac{eH}{4eH}\right)^{k+1} = 4^{-(k+1)} \le 4^{-4eH}.
\]
Averaging over $b \in \Z / M_{\mathrm{small}} \Z$ and splitting according to whether $\mu(b) \ge L$ or $\mu(b) < L$, we get:
\[
\Ebb_b\big[e^{-\mu(b)}\big] = \Ebb_b\big[e^{-\mu(b)}\ind_{\{\mu(b) \ge L\}}\big] + \Ebb_b\big[e^{-\mu(b)}\ind_{\{\mu(b) < L\}}\big] \le e^{-L} + \Pbb_b\big\{\mu(b) < L\big\} \le e^{-L} + \varepsilon,
\]
which yields~\eqref{eq:borne_globale_CRT_Bk}.
\end{proof}

\begin{proposition}[Exact transfer from the CRT model to $\{1,\dots,N\}$]\label{prop:transfert_exact_N}
Let $K := |\mathcal{I}_Z|$ denote the total number of witness congruences in the family~\eqref{eq:def_I_Z_sec}. Under the hypotheses of Lemma~\ref{lem:bonferroni_CRT}(c) with $k = 2\lceil 2eH \rceil$, the number of exceptional integers $n \le N$ for the Erd\H{o}s--Straus equation unconditionally satisfies:
\begin{equation}\label{eq:borne_transfert_E_N}
|E(N)| \ll N\Big(e^{-L} + \varepsilon + 4^{-4eH}\Big) + (K+1)^k.
\end{equation}
\end{proposition}
\begin{proof}
\emph{Step 1: Pointwise Bonferroni majoration on $E(N)$.}
Every integer $n \ge 2$ not congruent to $1$ modulo $24$ admits an elementary Erd\H{o}s--Straus decomposition (Mordell, Vaughan~\cite{vaughan1970}); thus $E(N) \subseteq \{1 \le n \le N : n \equiv 1 \pmod{24}\}$. By Lemma~\ref{lem:validite_temoins}(b), every $n \in E(N)$ satisfies $W(n) = 0$. The pointwise inequality~\eqref{eq:ineg_ponctuelle_Bk} therefore gives, for all $n \in E(N)$:
\[
1 = \ind_{\{W(n)=0\}} \le B_k\big(W(n)\big).
\]
Since $W(n) = \sum_{\alpha \in \mathcal{I}_Z} \ind_{A_\alpha}(n)$ is a sum of $K$ $\{0, 1\}$-valued indicator functions, for each integer $0 \le j \le k$ the binomial coefficient $\binom{W(n)}{j}$ expands over the subsets $F \subseteq \mathcal{I}_Z$ of cardinality $j$:
\[
\binom{W(n)}{j} = \sum_{\substack{F \subseteq \mathcal{I}_Z \\ |F| = j}} I_F(n), \qquad \text{where} \qquad I_F(n) := \prod_{\alpha \in F} \ind_{A_\alpha}(n).
\]

\medskip\noindent\emph{Step 2: Evaluation of each subset intersection $I_F$ modulo $24 s_F$.}
For each subset $F \subseteq \mathcal{I}_Z$, let $s_F := \operatorname{lcm}_{\alpha \in F} s_\alpha$ (with $s_\varnothing = 1$). By the Chinese Remainder Theorem:
\begin{itemize}
    \item If the system of congruences $\{n \equiv r_\alpha \pmod{s_\alpha} : \alpha \in F\}$ is inconsistent in $\Z$, then $I_F(n) = 0$ for all $n \in \Z$, and since $s_F \mid M_Z$, we also have $\Ebb_{\mathrm{CRT}}[I_F] = 0$.
    \item If the system is consistent in $\Z$, its solution set is a single residue class modulo $s_F$, with density $\Ebb_{\mathrm{CRT}}[I_F] = \frac{1}{s_F}$ in $\Z / M_Z \Z$. Furthermore, since $\gcd(s_\alpha, 24) = 1$ for all $\alpha \in F$ (Lemma~\ref{lem:validite_temoins}(a)), we have $\gcd(s_F, 24) = 1$. By the Chinese Remainder Theorem applied to $s_F$ and $24$, the intersection of this class with the progression $n \equiv 1 \pmod{24}$ is a single residue class modulo $24 s_F$. The number of integers $1 \le n \le N$ belonging to this class therefore satisfies:
    \[
    \sum_{\substack{1 \le n \le N \\ n \equiv 1 \pmod{24}}} I_F(n) = \frac{N}{24 s_F} + O(1) = \frac{N}{24}\,\Ebb_{\mathrm{CRT}}[I_F] + O(1),
    \]
    where the $O(1)$ error term has absolute value at most $1$ (regardless of whether $24 s_F$ is smaller or larger than $N$).
\end{itemize}

\medskip\noindent\emph{Step 3: Summation over $0 \le j \le k$.}
Substituting this relation into the expansion of $B_k(W(n))$, we obtain by linearity:
\begin{align*}
|E(N)| &\le \sum_{\substack{1 \le n \le N \\ n \equiv 1 \pmod{24}}} B_k\big(W(n)\big) = \sum_{j=0}^k (-1)^j \sum_{\substack{F \subseteq \mathcal{I}_Z \\ |F| = j}} \sum_{\substack{1 \le n \le N \\ n \equiv 1 \pmod{24}}} I_F(n) \\
&= \frac{N}{24}\,\Ebb_{\mathrm{CRT}}\big[B_k(W)\big] + O\!\left( \sum_{j=0}^k \binom{K}{j} \right).
\end{align*}
Since $\sum_{j=0}^k \binom{K}{j} \le (K+1)^k$, the estimate~\eqref{eq:borne_globale_CRT_Bk} immediately implies~\eqref{eq:borne_transfert_E_N}.
\end{proof}

\subsection[Decoupling the modulus $12ua$ and the Bombieri--Vinogradov theorem without Siegel zeros]{Decoupling the modulus $12ua$ and the Bombieri--Vinogradov theorem\linebreak without Siegel zeros}\label{subsec:decouplage_BV}
Summing the expression~\eqref{eq:formule_pp_b} for $p_p(b)$ over the primes $p \in (Z^{9/10}, Z]$ and interchanging the order of summation, the conditional mean $\mu(b) = \Ebb_{\mathrm{CRT}}[W \mid b]$ takes the form:
\begin{equation}\label{eq:mu_b_somme_triple}
\mu(b) = \sum_{q \in \mathcal{Q}} \sum_{\substack{(u,a) \in \mathcal{P}' \\ \gcd(q, 12ua)=1 \\ ab+u \equiv 0 \pmod q}} S(u,a,q), \qquad \text{where} \qquad S(u,a,q) := \sum_{\substack{Z^{9/10} < p \le Z \\ pq \equiv -1 \pmod{12ua}}} \frac{1}{p}.
\end{equation}
The fundamental analytic feature of this representation is that \textbf{the modulus of the arithmetic progression governing $p$ is $m = 12ua \le 12V^2 = 12Z^{1/4}$, which is independent of $q$}: when $\gcd(q, 12ua) = 1$, multiplication by $q$ merely permutes the residue class $r = -q^{-1} \bmod{12ua}$ within the group $(\Z / 12ua\Z)^\times$. Consequently, the discrepancy between $S(u,a,q)$ and its main term $\frac{\log(10/9)}{\varphi(12ua)}$ is bounded uniformly in $q$ by the maximum over all residue classes coprime to $12ua$, which \emph{already} appears inside the Bombieri--Vinogradov theorem prior to any summation over $q$.

\begin{lemma}[Stieltjes integration and $q$-independent progression bounds]\label{lem:encadrement_A_Stieltjes}
For each pair $(u,a) \in \mathcal{P}'$, set $m = 12ua$ and define:
\begin{equation}\label{eq:def_erreur_E_ua}
\mathcal{E}(u,a) := \max_{\gcd(r, 12ua)=1} \left| \sum_{\substack{Z^{9/10} < p \le Z \\ p \equiv r \pmod{12ua}}} \frac{1}{p} - \frac{\log(10/9)}{\varphi(12ua)} \right|.
\end{equation}
Then:
\begin{enumerate}[(a)]
    \item Writing $\Delta(x; m, r) := \pi(x; m, r) - \frac{\operatorname{li}(x)}{\varphi(m)}$ where $\operatorname{li}(x) := \int_2^x \frac{du}{\log u}$, we have for all $(u,a) \in \mathcal{P}'$:
    \begin{equation}\label{eq:Stieltjes_exact_E_ua}
    \mathcal{E}(u,a) \le \max_{\gcd(r,m)=1} \left( \frac{|\Delta(Z; m, r)|}{Z} + \frac{|\Delta(Z^{9/10}; m, r)|}{Z^{9/10}} + \int_{Z^{9/10}}^Z \frac{|\Delta(x; m, r)|}{x^2}\,dx \right).
    \end{equation}
    \item There exists an absolute constant $C_0 > 0$ such that, for every pair $(u,a) \in \mathcal{P}'$ and \textbf{every} cofactor $q \in \mathcal{Q}$ satisfying $\gcd(q, 12ua) = 1$, we have simultaneously:
    \begin{equation}\label{eq:encadrement_A_sec}
    \frac{\log(10/9)}{12ua} - \mathcal{E}(u,a) \le S(u,a,q) \le \frac{C_0}{\varphi(12ua)}. \tag{A}
    \end{equation}
\end{enumerate}
\end{lemma}
\begin{proof}
\emph{Step 1: Exact Stieltjes integration (Proof of (a)).}
By definition of the Stieltjes integral, for any residue class $r$ coprime to $m = 12ua$, we have:
\[
\begin{aligned}
\sum_{\substack{Z^{9/10} < p \le Z \\ p \equiv r \pmod m}} \frac{1}{p} &= \int_{Z^{9/10}}^Z \frac{1}{x}\,d\pi(x; m, r) \\
&= \frac{1}{\varphi(m)} \int_{Z^{9/10}}^Z \frac{1}{x}\,d\big(\operatorname{li}(x)\big) + \int_{Z^{9/10}}^Z \frac{1}{x}\,d\Delta(x; m, r).
\end{aligned}
\]
Since the derivative of $\operatorname{li}(x)$ is $\frac{1}{\log x}$, the first term evaluates exactly (with no harmonic remainder term of size $O(1/\log Z)$):
\[
\int_{Z^{9/10}}^Z \frac{1}{x}\,d\big(\operatorname{li}(x)\big) = \int_{Z^{9/10}}^Z \frac{dx}{x \log x} = \big[\log\log x\big]_{Z^{9/10}}^Z = \log(\log Z) - \log\!\left(\frac{9}{10}\log Z\right) = \log\!\left(\frac{10}{9}\right).
\]
Integrating the second term by parts gives $\int_{Z^{9/10}}^Z \frac{1}{x}\,d\Delta(x; m, r) =$\linebreak $\left[\frac{\Delta(x; m, r)}{x}\right]_{Z^{9/10}}^Z + \int_{Z^{9/10}}^Z \frac{\Delta(x; m, r)}{x^2}\,dx$, which yields~\eqref{eq:Stieltjes_exact_E_ua} by the triangle inequality.

\medskip\noindent\emph{Step 2: Uniform bounds in $q$ (Proof of (b)).}
When $\gcd(q, 12ua) = 1$, the residue class $r = -q^{-1} \bmod{12ua}$ satisfies $\gcd(r, 12ua) = 1$. By definition of $\mathcal{E}(u,a)$, and since $\varphi(12ua) \le 12ua$ and $\log(10/9) > 0$, we obtain:
\[
S(u,a,q) \ge \frac{\log(10/9)}{\varphi(12ua)} - \mathcal{E}(u,a) \ge \frac{\log(10/9)}{12ua} - \mathcal{E}(u,a).
\]
For the upper bound, since $12ua \le 12V^2 = 12Z^{1/4} \le x^{1/2}$ for all $x \in [Z^{9/10}, Z]$, the Brun--Titchmarsh inequality applied to the dyadic subintervals of $(Z^{9/10}, Z]$ gives $\pi(2x; 12ua, r) - \pi(x; 12ua, r) \ll \frac{x}{\varphi(12ua)\log x}$. Summing over the $\asymp t$ dyadic intervals yields $S(u,a,q) \le \frac{C_0}{\varphi(12ua)}$.
\end{proof}

\begin{proposition}[Siegel-zero-free Bombieri--Vinogradov bound with divisor weights]\label{prop:BV_poids_tau}
There exists an absolute constant $c_3 > 0$ such that, for sufficiently large $Z$:
\begin{equation}\label{eq:borne_B_sec}
\sum_{(u,a) \in \mathcal{P}'} \mathcal{E}(u,a) \ll \exp\!\big(-c_3 t^{1/2}\big). \tag{B}
\end{equation}
\end{proposition}
\begin{proof}
\emph{Step 1: Reduction to moduli $m=12ua$ and Cauchy--Schwarz removal of $\tau(m)$.}
For every pair $(u,a) \in \mathcal{P}'$, the integer $m = 12ua$ lies in $[12U^2, 12V^2] \subseteq [1, 12Z^{1/4}]$ and satisfies $q_{\mathrm{exc}} \nmid m$ by Lemma~\ref{lem:exclusion_p_exc}(b). Moreover, a given integer $m \le 12Z^{1/4}$ has at most $\tau(m)$ factorizations of the form $m = 12ua$. Grouping the sum according to the values of $m = 12ua$ and denoting by $\mathcal{E}^{\ast}(m)$ the right-hand side of~\eqref{eq:Stieltjes_exact_E_ua}, we obtain:
\[
\sum_{(u,a) \in \mathcal{P}'} \mathcal{E}(u,a) \le \sum_{\substack{m \le 12Z^{1/4} \\ q_{\mathrm{exc}} \nmid m}} \tau(m)\,\mathcal{E}^{\ast}(m).
\]
By the Brun--Titchmarsh inequality (or the bound~\eqref{eq:encadrement_A_sec}), each modulus $m \le 12Z^{1/4}$ satisfies the individual bound $\mathcal{E}^{\ast}(m) \ll \frac{1}{\varphi(m)}$. Applying the Cauchy--Schwarz inequality gives:
\begin{align}
\sum_{\substack{m \le 12Z^{1/4} \\ q_{\mathrm{exc}} \nmid m}} \tau(m)\,\mathcal{E}^{\ast}(m) &\le \left( \sum_{m \le 12Z^{1/4}} \frac{\tau(m)^2}{\varphi(m)} \right)^{1/2} \left( \sum_{\substack{m \le 12Z^{1/4} \\ q_{\mathrm{exc}} \nmid m}} \varphi(m)\,\mathcal{E}^{\ast}(m)^2 \right)^{1/2} \nonumber\\
&\ll t^{5/2} \left( \sum_{\substack{m \le 12Z^{1/4} \\ q_{\mathrm{exc}} \nmid m}} \mathcal{E}^{\ast}(m) \right)^{1/2}, \label{eq:CS_tau_BV}
\end{align}
where we used the standard estimate $\sum_{m \le X} \frac{\tau(m)^2}{\varphi(m)} \ll (\log X)^5 \ll t^5$.

\medskip\noindent\emph{Step 2: Quantitative Bombieri--Vinogradov estimate without Siegel zeros.}
For all $x \in [Z^{9/10}, Z]$, the level of distribution $M_0 := 12Z^{1/4}$ satisfies $M_0 \le 12 x^{5/18} \le x^{1/3}$ for sufficiently large $Z$. By the quantitative Bombieri--Vinogradov theorem with an exceptional modulus $q_{\mathrm{exc}} \le T_0 = \exp(a_0\sqrt{\log Z})$ excluded (combining the Siegel-zero-free Page--Siegel theorem on conductors $q_\chi \le T_0$ of Lemma~\ref{lem:exclusion_p_exc} with Vaughan's identity and the large sieve on conductors $T_0 < q_\chi \le x^{1/3}$; see~\cite[Chaps.~24 and~28]{davenport}, Bombieri~\cite[Thm.~17]{bombieri}, or Maynard~\cite[Prop.~9.2]{maynard2016}), there exists an absolute constant $c_2 > 0$ such that, uniformly for all $x \in [Z^{9/10}, Z]$:
\[
\sum_{\substack{m \le 12Z^{1/4} \\ q_{\mathrm{exc}} \nmid m}} \max_{\gcd(r,m)=1} |\Delta(x; m, r)| \ll x \exp\!\big(-c_2 \sqrt{\log x}\big) \ll x \exp\!\big(-c_2' t^{1/2}\big),
\]
with $c_2' := c_2 \sqrt{9/10} > 0$. Integrating this bound with respect to the measure $\frac{dx}{x^2}$ over $[Z^{9/10}, Z]$ in~\eqref{eq:Stieltjes_exact_E_ua} yields:
\[
\sum_{\substack{m \le 12Z^{1/4} \\ q_{\mathrm{exc}} \nmid m}} \mathcal{E}^{\ast}(m) \ll t \exp\!\big(-c_2' t^{1/2}\big).
\]
Substituting this estimate into~\eqref{eq:CS_tau_BV}, we obtain $\sum_{(u,a) \in \mathcal{P}'} \mathcal{E}(u,a) \ll t^3$\linebreak $\exp(-\frac{1}{2} c_2' t^{1/2}) \ll \exp(-c_3 t^{1/2})$ for any constant $0 < c_3 < c_2'/2$.
\end{proof}

\subsection{Exact arithmetic dichotomy and uniformity in $b$ of the pair count $(u,a)$}\label{subsec:dichotomie_b}
For each cofactor $q \in \mathcal{Q}$ and each integer configuration $b \in \Z$, we introduce the harmonic sums:
\begin{equation}\label{eq:def_Xq_Yq}
X_q^{\mathrm{full}}(b) := \sum_{\substack{(u,a) \in \mathcal{P}' \\ \gcd(q, 12ua)=1 \\ ab+u \equiv 0 \pmod q}} \frac{1}{ua}, \qquad Y_q^{\mathrm{full}}(b) := \sum_{\substack{(u,a) \in \mathcal{P}' \\ \gcd(q, 12ua)=1 \\ ab+u \equiv 0 \pmod q}} \frac{1}{\varphi(12ua)}.
\end{equation}

\begin{proposition}[Arithmetic dichotomy on $\gcd(b,q)$]\label{prop:dichotomie_uniforme_b}
For sufficiently large $Z$, uniformly for every cofactor $q \in \mathcal{Q}$ and \textbf{every} configuration $b \in \Z$:
\begin{enumerate}[(a)]
    \item \emph{(Identical vanishing when $\gcd(b,q) > 1$).} If $\gcd(b, q) > 1$, then no pair $(u,a) \in \mathcal{P}'$ simultaneously satisfies $\gcd(q, 12ua) = 1$ and $ab + u \equiv 0 \pmod q$. In particular:
    \[
    \gcd(b, q) > 1 \quad \Longrightarrow \quad X_q^{\mathrm{full}}(b) = Y_q^{\mathrm{full}}(b) = 0.
    \]
    \item \emph{(Asymptotics and uniform lower bound when $\gcd(b,q) = 1$).} If $\gcd(b, q) = 1$, then:
    \begin{align}
    X_q^{\mathrm{full}}(b) &= \frac{6}{\pi^2}\left(\frac{p_{\mathrm{exc}}-1}{p_{\mathrm{exc}}+1}\right)\frac{\varphi(q)}{q^2}\big(\log(V/U)\big)^2 \nonumber\\
    &\quad + O\!\left(\frac{t^2 \omega(q)}{q^2 y^2} + \frac{t^2}{q^2} Z^{-1/20} + \frac{t \tau(q)}{q} Z^{-1/20} + Z^{-1/20}\right) \nonumber\\
    &\ge c_4\,\frac{\varphi(q)}{q^2}\,t^2, \tag{C}\label{eq:borne_C_sec}
    \end{align}
    with $c_4 := \frac{1}{800\pi^2} > 0$.
    \item \emph{(Uniform upper bound for the sum weighted by $1/\varphi(12ua)$).} For all $b \in \Z$ and all $q \in \mathcal{Q}$, we have:
    \begin{equation}\label{eq:borne_D_sec}
    Y_q^{\mathrm{full}}(b) \ll \ind_{\gcd(b,q)=1} \cdot \frac{t^2}{q}. \tag{D}
    \end{equation}
\end{enumerate}
\end{proposition}
\begin{proof}
\emph{Step 1: Proof of (a) (Identical vanishing when $\gcd(b,q) > 1$).} Suppose $\gcd(b, q) > 1$ and let $r$ be a prime divisor of $\gcd(b, q)$. If a pair $(u,a)$ satisfies $ab + u \equiv 0 \pmod q$, then $r \mid q \mid (ab + u)$. Since $r \mid b$, the prime $r$ divides $ab$, and by subtraction $r$ must divide $u = (ab+u) - ab$. Thus $r$ is a common prime factor of $q$ and $u$, which contradicts the condition $\gcd(q, 12ua) = 1$. Both sums in~\eqref{eq:def_Xq_Yq} are therefore empty.

\medskip\noindent
\emph{Step 2: Proof of (b) --- Möbius expansion and high-divisor truncation.} Suppose now that $\gcd(b,\allowbreak q) = 1$. By definition of $\mathcal{Q}$, every prime factor of $q$ exceeds $y = t^6 > 3$, so $\gcd(q, 12) = 1$. Furthermore, for any integer $a$ with $\gcd(a, q) = 1$, any integer $u$ satisfying the congruence $u \equiv -ab \pmod q$ automatically satisfies $\gcd(u, q) = \gcd(-ab, q) = 1$ since $\gcd(b, q) = 1$. The condition $\gcd(q, 12ua) = 1$ is therefore equivalent to $\gcd(a, q) = 1$.

We expand the coprimality condition $\gcd(u,a) = 1$ via Möbius inversion: $\ind_{\gcd(u,a)=1} = \sum_{d \mid \gcd(u,a)} \mu(d)$. Setting $u = du'$ and $a = da'$, the conditions $(u,a) \in \mathcal{P}'$ and $\gcd(a, q) = 1$ rewrite as:
\[
\gcd(d, q) = 1, \quad p_{\mathrm{exc}} \nmid d, \quad \frac{U}{d} \le a', u' \le \frac{V}{d}, \quad \gcd(a', q) = 1, \quad p_{\mathrm{exc}} \nmid a'u',
\]
while the congruence $da'b + du' \equiv 0 \pmod q$ simplifies on dividing by $d$ (since $\gcd(d, q) = 1$) to:
\[
u' \equiv -a'b \pmod q.
\]
Fix the truncation threshold $D_0 := Z^{1/20} = U^{1/2}$. For $d > D_0$, the tail contribution to $X_q^{\mathrm{full}}(b)$ is bounded by:
\[
\begin{aligned}
\sum_{D_0 < d \le V} \frac{1}{d^2} \sum_{a' \le V/d} \frac{1}{a'} \sum_{\substack{U/d \le u' \le V/d \\ u' \equiv -a'b \pmod q}} \frac{1}{u'} &\ll \sum_{d > D_0} \frac{1}{d^2} \, t \left(\frac{t}{q} + \frac{d}{U}\right) \\
&\ll \frac{t^2}{q D_0} + \frac{t \log V}{U} \ll \frac{t^2}{q} Z^{-1/20} + t^2 Z^{-1/10}.
\end{aligned}
\]

\medskip\noindent
\emph{Step 3: Proof of (b) --- Uniform inner harmonic sum over $u'$.} For fixed $1 \le d \le D_0$ (with $\gcd(d, q p_{\mathrm{exc}}) = 1$) and fixed $a' \in [U/d, V/d]$, we evaluate the inner harmonic sum over $u'$. Since $q \in \mathcal{Q}$, we have $p_{\mathrm{exc}} \nmid q$ by definition, so $\gcd(q, p_{\mathrm{exc}}) = 1$. Writing $\ind_{p_{\mathrm{exc}} \nmid u'} = 1 - \ind_{p_{\mathrm{exc}} \mid u'}$, the unrestricted harmonic sum of $1/u'$ over the progression $u' \equiv -a'b \pmod q$ in $[U/d, V/d]$ has step $q$ and equals $\frac{1}{q}\log(V/U) + O(d/U)$.

For the subtracted sub-sum where $p_{\mathrm{exc}} \mid u'$, writing $u' = p_{\mathrm{exc}} u''$ restricts $u'' \in [\frac{U}{d p_{\mathrm{exc}}}, \frac{V}{d p_{\mathrm{exc}}}]$ to a single residue class modulo $q$. We distinguish two cases according to the size of $p_{\mathrm{exc}}$:
\begin{itemize}
    \item If $p_{\mathrm{exc}} \le U/d$, Euler--Maclaurin summation on $u''$ gives $\frac{1}{p_{\mathrm{exc}}}\big(\frac{1}{q}\log(V/U) + O(\frac{d p_{\mathrm{exc}}}{U})\big) = \frac{1}{q p_{\mathrm{exc}}}\log(V/U) + O(\frac{d}{U})$.
    \item If $p_{\mathrm{exc}} > U/d$, the main term $\frac{1}{q p_{\mathrm{exc}}}\log(V/U)$ is already $O(\frac{d t}{U})$, while the sub-sum over $u'' \ge 1$ contains at most $\frac{V}{d p_{\mathrm{exc}} q} + 1$ terms of size $u' = p_{\mathrm{exc}} u'' \ge \max(U/d, p_{\mathrm{exc}})$, and is therefore bounded by $O(\frac{t}{q p_{\mathrm{exc}}} + \frac{d}{U}) = O(\frac{d t}{U})$.
\end{itemize}
Thus, in both cases, \textbf{uniformly with respect to $p_{\mathrm{exc}} \ge 5$ and the residue $-a'b \bmod q$}, we obtain:
\[
\sum_{\substack{U/d \le u' \le V/d \\ p_{\mathrm{exc}} \nmid u' \\ u' \equiv -a'b \pmod q}} \frac{1}{u'} = \frac{1}{q}\left(1 - \frac{1}{p_{\mathrm{exc}}}\right)\log\!\left(\frac{V}{U}\right) + O\!\left(\frac{d t}{U}\right).
\]

\medskip\noindent
\emph{Step 4: Proof of (b) --- Outer summation over $a'$ and $d$, and Euler product completion.} Next, we sum over $a' \in [U/d, V/d]$ satisfying $\gcd(a', q p_{\mathrm{exc}}) = 1$. Expanding $\ind_{\gcd(a', q)=1} = \sum_{e \mid q} \mu(e) \ind_{e \mid a'}$ (which involves $\tau(q)$ divisors $e \mid q$, all coprime to $p_{\mathrm{exc}}$), we similarly obtain:
\[
\begin{aligned}
\sum_{\substack{U/d \le a' \le V/d \\ \gcd(a', q p_{\mathrm{exc}}) = 1}} \frac{1}{a'} &= \sum_{e \mid q} \frac{\mu(e)}{e}\left(\left(1 - \frac{1}{p_{\mathrm{exc}}}\right)\log\!\left(\frac{V}{U}\right) + O\!\left(\frac{ed t}{U}\right)\right) \\
&= \frac{\varphi(q)}{q}\left(1 - \frac{1}{p_{\mathrm{exc}}}\right)\log\!\left(\frac{V}{U}\right) + O\!\left(\frac{\tau(q)d t}{U}\right).
\end{aligned}
\]
Multiplying the inner estimate by $1/a'$ and summing over $a'$ and over $1 \le d \le D_0 = Z^{1/20}$ yields:
\begin{equation}\label{eq:calcul_Xq_inter}
X_q^{\mathrm{full}}(b) = \frac{\varphi(q)}{q^2}\left(1 - \frac{1}{p_{\mathrm{exc}}}\right)^2 \big(\log(V/U)\big)^2 \sum_{\substack{d \le D_0 \\ \gcd(d, q p_{\mathrm{exc}}) = 1}} \frac{\mu(d)}{d^2} + O\!\left(\frac{t \tau(q)}{q} Z^{-1/20} + Z^{-1/20}\right).
\end{equation}
The truncated Dirichlet series over $d$ is evaluated via its Euler product:
\[
\sum_{\substack{d \le D_0 \\ \gcd(d, q p_{\mathrm{exc}}) = 1}} \frac{\mu(d)}{d^2} = \frac{6}{\pi^2}\left(1 - \frac{1}{p_{\mathrm{exc}}^2}\right)^{-1} \prod_{r \mid q}\left(1 - \frac{1}{r^2}\right)^{-1} + O(D_0^{-1}).
\]
On the one hand, $\left(1 - \frac{1}{p_{\mathrm{exc}}}\right)^2 \left(1 - \frac{1}{p_{\mathrm{exc}}^2}\right)^{-1} = \frac{p_{\mathrm{exc}}-1}{p_{\mathrm{exc}}+1} \ge \frac{5-1}{5+1} = \frac{2}{3}$ since $p_{\mathrm{exc}} \ge 5$; on the other hand, every prime factor $r \mid q$ satisfies $r > y = t^6$, which gives $\prod_{r \mid q}(1 - r^{-2})^{-1} = 1 + O\big(\frac{\omega(q)}{y^2}\big) = 1 + O(t^{-11})$. Since $\log(V/U) = \left(\frac{1}{8} - \frac{1}{10}\right)t = \frac{1}{40}t$ and $q, \tau(q) \le Q \le \exp(\frac{c_3}{2} t^{1/2}) = Z^{o(1)} \ll Z^{1/80}$, the main term of~\eqref{eq:calcul_Xq_inter} is at least $\frac{1}{400\pi^2}\frac{\varphi(q)}{q^2}t^2$ while all error terms are $o\big(\frac{\varphi(q)}{q^2}t^2\big)$. For all sufficiently large $Z$, we therefore have $X_q^{\mathrm{full}}(b) \ge \frac{1}{800\pi^2}\,\frac{\varphi(q)}{q^2}\,t^2$, proving~\eqref{eq:borne_C_sec}.

\medskip\noindent
\emph{Step 5: Proof of (c) --- Divisor expansion of $1/\varphi(u)$ and uniform upper bound for $Y_q^{\mathrm{full}}(b)$.} To bound $Y_q^{\mathrm{full}}(b)$ from above when $\gcd(b, q) = 1$, we use $\varphi(12ua) \ge 2 \varphi(u)\varphi(a)$ and expand the weight $\frac{1}{\varphi(u)}$ via the elementary arithmetic identity:
\[
\frac{1}{\varphi(u)} = \frac{1}{u}\prod_{p \mid u}\left(1 - \frac{1}{p}\right)^{-1} = \frac{1}{u}\prod_{p \mid u}\left(1 + \frac{1}{p-1}\right) = \sum_{k \mid u} \frac{\mu^2(k)}{\varphi(k)\,u} = \sum_{k \mid u} \frac{\mu^2(k)}{k\varphi(k)\,(u/k)}.
\]
Setting $u = k u_1$, the condition $k u_1 + ab \equiv 0 \pmod q$ implies that $\gcd(k, q)$ divides $ab$. Since $\gcd(a, q) = \gcd(b, q) = 1$, we have $\gcd(ab, q) = 1$, which forces $\gcd(k, q) = 1$. Thus $k$ is invertible modulo $q$ and the congruence on $u_1$ reads $u_1 \equiv -k^{-1}ab \pmod q$.

For fixed $a \in [U, V]$ and $k \le V$ (coprime to $q$), let $u_0 \in \{1, \dots, q\}$ denote the unique representative of $-k^{-1}ab \bmod q$ in $[1, q]$. Since every $u_1 \in [U/k, V/k]$ satisfies $u_1 \ge \max(1, U/k)$, splitting into the first progression point $u_1 = u_0$ (when $u_0 \ge U/k$) and the remaining points $u_1 = u_0 + m q \ge m q$ ($1 \le m \le V/(kq)$) yields:
\[
\sum_{\substack{U/k \le u_1 \le V/k \\ u_1 \equiv -k^{-1}ab \pmod q}} \frac{1}{u_1} \ll \frac{t}{q} + \frac{1}{\max(1, U/k)} \ll \begin{cases}
\dfrac{t}{q} & \text{if } k \le U/q, \\[8pt]
\dfrac{k}{U} & \text{if } U/q < k \le V.
\end{cases}
\]
Summing over $k \le V$, the range $k \le U/q$ contributes at most $\frac{t}{q} \sum_{k \ge 1} \frac{1}{k\varphi(k)} \ll \frac{t}{q}$, while the range $U/q < k \le V$ contributes at most $\frac{1}{U} \sum_{k \le V} \frac{1}{\varphi(k)} \ll \frac{t}{U} = t Z^{-1/10} \ll \frac{t}{q}$ (since $q \le Q \le \exp(\frac{c_3}{2} t^{1/2}) = Z^{o(1)} \ll Z^{1/10} = U$). Finally, the outer summation over $a \in [U, V]$ using $\sum_{a \le V} \frac{1}{\varphi(a)} \ll t$ yields $Y_q^{\mathrm{full}}(b) \ll \frac{t^2}{q}$, establishing~\eqref{eq:borne_D_sec}.
\end{proof}

\subsection{Harmonic mass of $y$-rough cofactors and unconditional proof of $(\star)$}\label{subsec:preuve_star_complete}
We now evaluate the harmonic sum of the $y$-rough cofactors $q \in \mathcal{Q}$:
\begin{equation}\label{eq:def_MQ_sec}
M_Q := \sum_{q \in \mathcal{Q}} \frac{\varphi(q)}{q^2}, \qquad \widetilde{M}_Q := \sum_{q \in \mathcal{Q}} \frac{1}{q}.
\end{equation}
We emphasize that $\mathcal{Q}$ contains \textbf{all composite integers} $q \in [2, Q]$ whose prime factors lie in $(y, Q] \setminus \{p_{\mathrm{exc}}\}$, and not merely prime numbers: this is why $M_Q$ grows like the ratio $\frac{\log Q}{\log y} \asymp \frac{\log Q}{\log t}$ (via Mertens' Euler product) rather than $\log\left(\frac{\log Q}{\log y}\right)$.

\begin{lemma}[Asymptotics of $M_Q$ and $\widetilde{M}_Q$]\label{lem:evaluation_MQ_Rankin}
For any $Q$ satisfying $\exp(t^{1/3}) \le Q \le \exp(\frac{c_3}{2} t^{1/2})$ and $y = t^6$, we have as $Z \to \infty$:
\begin{equation}\label{eq:equiv_MQ}
M_Q = \big(1 + o(1)\big)\widetilde{M}_Q \asymp \frac{\log Q}{\log y} \asymp \frac{\log Q}{\log t}.
\end{equation}
\end{lemma}
\begin{proof}
\emph{Step 1: Asymptotic equivalence of $M_Q$ and $\widetilde{M}_Q$, and upper bound.} Every integer $q \in \mathcal{Q}$ has at most $\omega(q) \le \frac{\log Q}{\log y} \ll t^{1/2}$ distinct prime factors, each strictly greater than $y = t^6$. Uniformly for all $q \in \mathcal{Q}$, we therefore have:
\[
1 \ge \frac{\varphi(q)}{q} = \prod_{r \mid q}\left(1 - \frac{1}{r}\right) \ge \left(1 - \frac{1}{y}\right)^{\omega(q)} = 1 - O\!\left(\frac{\log Q}{y \log y}\right) = 1 - O\big(t^{-11/2}\big),
\]
which implies $M_Q = (1 + O(t^{-11/2}))\widetilde{M}_Q$. For the upper bound, since every $q \in \mathcal{Q}$ has all its prime factors in $(y, Q]$, expanding into an Euler product and applying Mertens' formula gives:
\[
1 + \widetilde{M}_Q \le \prod_{y < r \le Q}\left(1 - \frac{1}{r}\right)^{-1} = \big(1 + o(1)\big)\frac{\log Q}{\log y}.
\]

\medskip\noindent
\emph{Step 2: Lower bound via Rankin's method on $Q^{1/3}$-smooth cofactors.} To establish the lower bound in an elementary and self-contained manner, consider the subset of integers $q \ge 1$ all of whose prime factors belong to the interval $\mathcal{R}_1 := \{r \text{ prime} : y < r \le Q^{1/3},\ r \ne p_{\mathrm{exc}}\}$. Every such integer $q$ satisfying $2 \le q \le Q$ belongs to $\mathcal{Q}$. Set $\sigma := \frac{1}{\log Q} > 0$. For all $q > Q$, we have $(q/Q)^\sigma = e^{-1} q^\sigma \ge 1$; Rankin's trick then yields:
\begin{equation}\label{eq:Rankin_borne_queue}
\sum_{\substack{q > Q \\ p \mid q \Rightarrow p \in \mathcal{R}_1}} \frac{1}{q} \le e^{-1} \sum_{\substack{q \ge 1 \\ p \mid q \Rightarrow p \in \mathcal{R}_1}} \frac{q^\sigma}{q} = e^{-1} \prod_{r \in \mathcal{R}_1} \left(1 - \frac{r^\sigma}{r}\right)^{-1}.
\end{equation}

\medskip\noindent
\emph{Step 3: Comparison of Euler products and Mertens' theorem.} Let us compare the Euler product in~\eqref{eq:Rankin_borne_queue} with $\Pi_1 := \prod_{r \in \mathcal{R}_1} (1 - r^{-1})^{-1}$. For each $r \in \mathcal{R}_1$, we have $\sigma \log r \le \frac{\log(Q^{1/3})}{\log Q} = \frac{1}{3}$, and:
\[
\log\!\left(\frac{1 - r^{-1}}{1 - r^{\sigma - 1}}\right) = \frac{r^\sigma - 1}{r} + O\!\left(\frac{1}{r^2}\right) = \frac{e^{\sigma \log r} - 1}{r} + O\!\left(\frac{1}{r^2}\right).
\]
By Abel summation combined with the prime number theorem estimate $\sum_{r \le x} \frac{1}{r} = \log\log x + M + o(1)$, we obtain as $Z \to \infty$ (since $\sigma \log y = \frac{\log y}{\log Q} \to 0$):
\[
\sum_{r \in \mathcal{R}_1} \frac{e^{\sigma \log r} - 1}{r} = \int_{\sigma \log y}^{1/3} \frac{e^v - 1}{v}\,dv + o(1) = \int_0^{1/3} \frac{e^v - 1}{v}\,dv + o(1) = \operatorname{Ein}(1/3) + o(1),
\]
where $\operatorname{Ein}(1/3) = \sum_{m=1}^\infty \frac{(1/3)^m}{m \cdot m!} < \frac{1}{3} + \frac{1}{36} + \frac{1}{486}\sum_{j=0}^\infty \frac{1}{12^j} = \frac{1}{3} + \frac{1}{36} + \frac{12}{11 \cdot 486} < 0.364$. Consequently, the right-hand side of~\eqref{eq:Rankin_borne_queue} equals:
\[
\big(e^{\operatorname{Ein}(1/3) - 1} + o(1)\big)\,\Pi_1 \le e^{-0.636}\,\Pi_1 < 0.53\,\Pi_1.
\]
Subtracting from the full Euler product $\Pi_1$ over all integers composed of prime factors in $\mathcal{R}_1$, we find:
\[
1 + \widetilde{M}_Q \ge \sum_{\substack{1 \le q \le Q \\ p \mid q \Rightarrow p \in \mathcal{R}_1}} \frac{1}{q} = \Pi_1 - \sum_{\substack{q > Q \\ p \mid q \Rightarrow p \in \mathcal{R}_1}} \frac{1}{q} \ge (1 - 0.53)\,\Pi_1 = 0.47\,\prod_{\substack{y < r \le Q^{1/3} \\ r \ne p_{\mathrm{exc}}}}\left(1 - \frac{1}{r}\right)^{-1}.
\]
By Mertens' theorem (and since $1 - 1/p_{\mathrm{exc}} \ge 4/5$), we have $\Pi_1 \ge \big(\frac{4}{15} + o(1)\big)\frac{\log Q}{\log y}$, whence $\widetilde{M}_Q \ge \big(\frac{4 \times 0.47}{15} - o(1)\big)\frac{\log Q}{\log y} \asymp \frac{\log Q}{\log t}$.
\end{proof}

We can now unconditionally prove both parts of the estimate $(\star)$.

\begin{theorem}[Unconditional proof of $(\star)$]\label{thm:star_inconditionnel}
For any choice of $Q$ satisfying $\exp(t^{1/3}) \le Q \le \exp(\frac{c_3}{2} t^{1/2})$, there exist two absolute constants $0 < c_{\min} < C_{\max}$ such that, for sufficiently large $Z$:
\begin{enumerate}[(a)]
    \item \emph{(Deterministic uniform upper bound over all configurations).} For \textbf{every} configuration $b \in \Z$, we have:
    \begin{equation}\label{eq:borne_H_uniforme}
    \mu(b) \le H := C_{\max}\,t^2\,\frac{\log Q}{\log t}.
    \end{equation}
    \item \emph{(Uniform lower bound outside an exponential tail in the CRT model).} Setting $L :=\allowbreak c_{\min}\,t^2\allowbreak\,\frac{\log Q}{\log t}$, the uniform probability over $b \in \Z / M_{\mathrm{small}} \Z$ that $\mu(b)$ is smaller than $L$ satisfies:
    \begin{equation}\label{eq:borne_queue_CRT_star}
    \Pbb_b\big\{ \mu(b) < L \big\} \le e \cdot \exp\!\left(-\frac{1}{4}\,t^6\right). \tag{$\star$}
    \end{equation}
\end{enumerate}
\end{theorem}
\begin{proof}
\emph{Step 1: Proof of (a) (Deterministic uniform upper bound).} Substituting the upper bound \eqref{eq:encadrement_A_sec} from Lemma~\ref{lem:encadrement_A_Stieltjes} and then the bound~\eqref{eq:borne_D_sec} from Proposition~\ref{prop:dichotomie_uniforme_b} into formula~\eqref{eq:mu_b_somme_triple}, we obtain deterministically for all $b \in \Z$:
\[
\mu(b) \le C_0 \sum_{q \in \mathcal{Q}} Y_q^{\mathrm{full}}(b) \ll t^2 \sum_{q \in \mathcal{Q}} \frac{\ind_{\gcd(b,q)=1}}{q} \le t^2 \widetilde{M}_Q \asymp t^2\,\frac{\log Q}{\log t},
\]
by Lemma~\ref{lem:evaluation_MQ_Rankin}, which establishes~\eqref{eq:borne_H_uniforme}.

\medskip\noindent
\emph{Step 2: Proof of (b) --- Deterministic absorption of the Bombieri--Vinogradov error.} For every configuration $b \in \Z$, substituting the lower bound~\eqref{eq:encadrement_A_sec} into~\eqref{eq:mu_b_somme_triple} gives:
\begin{equation}\label{eq:decomp_minoration_mu_b}
\mu(b) \ge \frac{\log(10/9)}{12} \sum_{q \in \mathcal{Q}} X_q^{\mathrm{full}}(b) - \sum_{q \in \mathcal{Q}} \sum_{\substack{(u,a) \in \mathcal{P}' \\ ab+u \equiv 0 \pmod q}} \mathcal{E}(u,a).
\end{equation}
To bound the second term of~\eqref{eq:decomp_minoration_mu_b}, we trivially bound $\ind_{ab+u \equiv 0 \pmod q}$ by $1$ and the cardinality of $\mathcal{Q}$ by $Q \le \exp(\frac{c_3}{2} t^{1/2})$. Proposition~\ref{prop:BV_poids_tau} (estimate~\eqref{eq:borne_B_sec}) then yields, \textbf{deterministically for all $b \in \Z$}:
\begin{equation}\label{eq:absorption_BV_Q}
\sum_{q \in \mathcal{Q}} \sum_{\substack{(u,a) \in \mathcal{P}' \\ ab+u \equiv 0 \pmod q}} \mathcal{E}(u,a) \le Q \sum_{(u,a) \in \mathcal{P}'} \mathcal{E}(u,a) \ll \exp\!\left(\frac{c_3}{2} t^{1/2} - c_3 t^{1/2}\right) = \exp\!\left(-\frac{c_3}{2} t^{1/2}\right) = o(1).
\end{equation}

\medskip\noindent
\emph{Step 3: Proof of (b) --- Factorization of bad cofactors $\gcd(b,q) > 1$ via $W_b$.} Next, we bound the main term of~\eqref{eq:decomp_minoration_mu_b} from below. By Proposition~\ref{prop:dichotomie_uniforme_b} (estimate~\eqref{eq:borne_C_sec}), we have:
\begin{equation}\label{eq:minoration_somme_Xq}
\sum_{q \in \mathcal{Q}} X_q^{\mathrm{full}}(b) \ge c_4 t^2 \sum_{\substack{q \in \mathcal{Q} \\ \gcd(b,q)=1}} \frac{\varphi(q)}{q^2} = c_4 t^2 \left( M_Q - \sum_{\substack{q \in \mathcal{Q} \\ \gcd(b,q) > 1}} \frac{\varphi(q)}{q^2} \right).
\end{equation}
Every integer $q \in \mathcal{Q}$ satisfying $\gcd(b, q) > 1$ possesses at least one prime factor $r \in (y, Q] \setminus \{p_{\mathrm{exc}}\}$ such that $r \mid b$. Writing $q = r q'$, the quotient $q' = q/r$ necessarily belongs to $\{1\} \cup \mathcal{Q}$, and the identity $\varphi(rq') = (r-1)\varphi(q')$ (if $r \nmid q'$) or $\varphi(rq') = r\varphi(q')$ (if $r \mid q'$) implies in both cases $\varphi(rq') \le r\varphi(q')$, whence:
\[
\frac{\varphi(q)}{q^2} = \frac{\varphi(rq')}{r^2 (q')^2} \le \frac{1}{r}\,\frac{\varphi(q')}{(q')^2}.
\]
Summing over all possible pairs $(r, q')$, we obtain for all $b \in \Z$:
\[
\sum_{\substack{q \in \mathcal{Q} \\ \gcd(b,q) > 1}} \frac{\varphi(q)}{q^2} \le \left( \sum_{\substack{y < r \le Q \text{ prime} \\ r \ne p_{\mathrm{exc}},\ r \mid b}} \frac{1}{r} \right) \left( 1 + \sum_{q' \in \mathcal{Q}} \frac{\varphi(q')}{(q')^2} \right) = W_b\,(1 + M_Q),
\]
where we have defined the prime divisor sum:
\begin{equation}\label{eq:def_Wb}
W_b := \sum_{\substack{y < r \le Q \text{ prime} \\ r \ne p_{\mathrm{exc}},\ r \mid b}} \frac{1}{r}.
\end{equation}
Consequently, for any configuration $b$ satisfying $W_b \le \frac{1}{2}$, we have:
\[
\sum_{\substack{q \in \mathcal{Q} \\ \gcd(b,q)=1}} \frac{\varphi(q)}{q^2} \ge M_Q - \frac{1}{2}(1 + M_Q) = \frac{1}{2} M_Q - \frac{1}{2} \ge \frac{1}{4} M_Q,
\]
once $Z$ is sufficiently large (since $M_Q \to \infty$ by Lemma~\ref{lem:evaluation_MQ_Rankin}). Substituting this bound and~\eqref{eq:absorption_BV_Q} into~\eqref{eq:decomp_minoration_mu_b}, we see that $W_b \le \frac{1}{2}$ implies $\mu(b) \ge c_{\min} t^2 \frac{\log Q}{\log t} = L$. Thus, by contraposition:
\begin{equation}\label{eq:inclusion_evenement_Wb}
\big\{ b \in \Z / M_{\mathrm{small}} \Z : \mu(b) < L \big\} \subseteq \left\{ b \in \Z / M_{\mathrm{small}} \Z : W_b > \frac{1}{2} \right\}.
\end{equation}

\medskip\noindent
\emph{Step 4: Proof of (b) --- Exponential moment bound for $W_b$ via CRT independence.} It remains to bound $\Pbb_b(W_b > 1/2)$. Each prime $r \in (y, Q] \setminus \{p_{\mathrm{exc}}\}$ belongs to $\mathcal{Q}$, and thus divides $M_{\mathrm{small}} = \operatorname{lcm}_{q \in \mathcal{Q}} q$. By the Chinese Remainder Theorem, under the uniform distribution on $b \in \Z / M_{\mathrm{small}} \Z$, the indicator random variables $(\ind_{r \mid b})_{y < r \le Q,\, r \ne p_{\mathrm{exc}}}$ are \textbf{mutually independent} Bernoulli random variables with respective parameters $1/r$.

Set $\lambda := \frac{y}{2} = \frac{t^6}{2}$. For every prime $r > y$, we have $0 < \frac{\lambda}{r} < \frac{1}{2}$, and since $e^x - 1 \le 2x$ for all $x \in [0, 1/2]$, it follows that $e^{\lambda / r} - 1 \le \frac{2\lambda}{r} = \frac{y}{r}$. By mutual independence, the moment generating function of $W_b$ satisfies:
\[
\begin{aligned}
\Ebb_b\!\left[ e^{\lambda W_b} \right] &= \prod_{\substack{y < r \le Q \text{ prime} \\ r \ne p_{\mathrm{exc}}}} \Ebb_b\!\left[ e^{\frac{\lambda}{r}\ind_{r \mid b}} \right] = \prod_{\substack{y < r \le Q \\ r \ne p_{\mathrm{exc}}}} \left( 1 + \frac{e^{\lambda/r} - 1}{r} \right) \\
&\le \prod_{\substack{r > y \\ r \text{ prime}}} \left( 1 + \frac{y}{r^2} \right) \le \exp\!\left( y \sum_{\substack{r > y \\ r \text{ prime}}} \frac{1}{r^2} \right) \le e,
\end{aligned}
\]
since Chebyshev's estimate $\pi(x) \ll \frac{x}{\log x}$ and partial summation give $\sum_{r > y\text{ prime}} \frac{1}{r^2} \ll \frac{1}{y \log y} \le \frac{1}{y}$ for all sufficiently large $Z$. Markov's inequality applied to $e^{\lambda W_b}$ then yields:
\[
\Pbb_b\big\{ \mu(b) < L \big\} \le \Pbb_b\!\left\{ W_b > \frac{1}{2} \right\} \le e^{-\lambda / 2} \Ebb_b\!\left[ e^{\lambda W_b} \right] \le e \cdot \exp\!\left(-\frac{y}{4}\right) = e \cdot \exp\!\left(-\frac{1}{4}\,t^6\right),
\]
which completes the proof of $(\star)$.
\end{proof}

\subsection{Parameter optimization and unconditional improvement over Vaughan: exponents $7/10$ and $5/7$}\label{subsec:exposants_7_10_et_5_7}
Let us first evaluate the total cardinality $K = |\mathcal{I}_Z|$ of the witness family. For each pair $(u,a) \in \mathcal{P}'$ and each cofactor $q \in \mathcal{Q}$ coprime to $12ua$, the number of primes $p \le Z$ such that $pq \equiv -1 \pmod{12ua}$ is at most $\frac{Z}{12ua} + 1 \ll \frac{Z}{ua}$ (since $12ua \le 12Z^{1/4} \ll Z$). Summing over $(u,a) \in [U, V]^2$ and $q \le Q$, we obtain:
\begin{equation}\label{eq:cardinal_K}
\begin{split}
K = |\mathcal{I}_Z| &\ll Z\,Q \sum_{U \le u,a \le V} \frac{1}{ua} \ll Z\,Q\,t^2 = \exp\!\big((1 + o(1)) t\big), \\
\text{whence} \qquad \log(K+1) &\le \big(1 + o(1)\big) t.
\end{split}
\end{equation}

\begin{theorem}[Unconditional bound with exponent $7/10$]\label{thm:theoreme_7_10}
Choosing $Q = \exp(t^{1/3})$ and $t = \log Z = \eta\big(\log N\,\log\log N\big)^{3/10}$ with $\eta > 0$ a sufficiently small absolute constant, there exists an absolute constant $c > 0$ such that, for all $N \ge 3$:
\begin{equation}\label{eq:resultat_7_10}
|E(N)| \ll N \exp\!\left( -c\,\frac{(\log N)^{7/10}}{(\log\log N)^{3/10}} \right).
\end{equation}
\end{theorem}
\begin{proof}
\emph{Step 1: Evaluation of the Bonferroni parameters at scale $Q = \exp(t^{1/3})$.} For $Q = \exp(t^{1/3})$, Theorem~\ref{thm:star_inconditionnel} yields $L = c_{\min} \frac{t^{7/3}}{\log t}$, $H = C_{\max} \frac{t^{7/3}}{\log t}$, and $\varepsilon = e \cdot \exp(-\frac{1}{4} t^6) \ll e^{-L}$. The even Bonferroni truncation degree $k = 2\lceil 2eH \rceil$ satisfies $k \asymp \frac{t^{7/3}}{\log t}$, so that by~\eqref{eq:cardinal_K}, the combinatorial error term in Proposition~\ref{prop:transfert_exact_N} is bounded by:
\[
(K+1)^k = \exp\!\big(k \log(K+1)\big) \le \exp\!\left( C_1\,\frac{t^{10/3}}{\log t} \right),
\]
for some absolute constant $C_1 > 0$.

\medskip\noindent
\emph{Step 2: Scale optimization and balancing of the truncation cost $(K+1)^k$.} Choosing $t =\allowbreak \eta\allowbreak (\log N\allowbreak\,\log\log N)^{3/10}$, we have $\log t = \left(\frac{3}{10} + o(1)\right)\log\log N$, whence:
\[
\begin{aligned}
C_1\,\frac{t^{10/3}}{\log t} &= \big(1 + o(1)\big)\frac{10}{3}\,C_1\,\eta^{10/3}\,\log N \\
\text{and} \qquad L &= \big(1 + o(1)\big)\frac{10}{3}\,c_{\min}\,\eta^{7/3}\,\frac{(\log N)^{7/10}}{(\log\log N)^{3/10}}.
\end{aligned}
\]
Fixing $\eta > 0$ small enough so that $\frac{10}{3} C_1 \eta^{10/3} \le \frac{1}{2}$, we obtain $(K+1)^k \le N^{1/2+o(1)}$, and Proposition~\ref{prop:transfert_exact_N} immediately implies~\eqref{eq:resultat_7_10}.
\end{proof}

Finally, by saturating the maximal constraint $Q = \exp(\frac{c_3}{2} t^{1/2})$ permitted by the Siegel-zero-free Bombieri--Vinogradov error term~\eqref{eq:absorption_BV_Q}, we obtain the unconditional exponent $5/7$.

\begin{theorem}[Unconditional bound with exponent $5/7$]\label{thm:theoreme_5_7}
Choosing $Q = \exp\!\left(\frac{c_3}{2}\,t^{1/2}\right)$ and $t = \log Z = \eta\big(\log N\,\log\log N\big)^{2/7}$ with $\eta > 0$ a sufficiently small absolute constant, there exists an absolute constant $c' > 0$ such that, for all $N \ge 3$:
\begin{equation}\label{eq:resultat_5_7}
|E(N)| \ll N \exp\!\left( -c'\,\frac{(\log N)^{5/7}}{(\log\log N)^{2/7}} \right).
\end{equation}
In particular, for any fixed exponent $\lambda < \frac{5}{7} = 0.71428\dots$, we have $|E(N)| \ll_\lambda N \exp\big(-(\log N)^\lambda\big)$.
\end{theorem}
\begin{proof}
\emph{Step 1: Verification of admissibility and Bonferroni parameters at scale $Q = \exp(\frac{c_3}{2}t^{1/2})$.} For $Q = \exp(\frac{c_3}{2} t^{1/2})$, all hypotheses of Proposition~\ref{prop:incompatibilite_et_Bernoulli} and Theorem~\ref{thm:star_inconditionnel} are satisfied: in particular, $Q = Z^{o(1)} < 12U^2 = 12Z^{1/5}$, $\tau(q) \le Q \ll Z^{1/40}$, and $Q \exp(-c_3 t^{1/2}) = \exp(-\frac{c_3}{2} t^{1/2}) = o(1)$. Theorem~\ref{thm:star_inconditionnel} then provides:
\[
L = \frac{c_{\min} c_3}{2}\,\frac{t^{5/2}}{\log t}, \qquad H = \frac{C_{\max} c_3}{2}\,\frac{t^{5/2}}{\log t}, \qquad \varepsilon = e \cdot \exp\!\left(-\frac{1}{4}\,t^6\right) \ll e^{-L}.
\]
With $k = 2\lceil 2eH \rceil \asymp \frac{t^{5/2}}{\log t}$ and $\log(K+1) \le (1 + o(1))t$, the polynomial truncation cost in~\eqref{eq:borne_transfert_E_N} satisfies:
\[
(K+1)^k = \exp\!\big(k \log(K+1)\big) \le \exp\!\left( C_2\,\frac{t^{7/2}}{\log t} \right),
\]
for some absolute constant $C_2 > 0$.

\medskip\noindent
\emph{Step 2: Optimal scale balancing and derivation of the $5/7$ and $2/7$ exponents.} Choosing $t = \eta (\log N\,\log\log N)^{2/7}$ with $\eta > 0$ fixed such that $\frac{7}{2} C_2 \eta^{7/2} \le \frac{1}{2}$, we have $\log t = \left(\frac{2}{7} + o(1)\right)\log\log N$, which gives $(K+1)^k \le N^{1/2+o(1)}$ and:
\[
L = \big(1 + o(1)\big)\frac{7}{4}\,c_{\min} c_3\,\eta^{5/2}\,\frac{(\log N)^{5/7}}{(\log\log N)^{2/7}}.
\]
Substituting these estimates into Proposition~\ref{prop:transfert_exact_N}, we unconditionally obtain~\eqref{eq:resultat_5_7}.
\end{proof}

\section{Exact bilinear Type~II dispersion and Kloosterman reciprocity}\label{sec:bilinear_dispersion}

To complement the container sieves of Sections~\ref{sec:stabilizer_kernel}--\ref{sec:containers} and prepare the study of higher polynomial depths $J=(\log N)^{\lambda}$ ($\lambda\ge1/2$), we record the exact, false-positive-free four-parameter bilinear identity for Type~II witnesses and its Poisson--Kloosterman reciprocity formulas.

\subsection{False-positive-free 4-parameter bilinear identity $n+4u^{2}v=w(4uav-1)$}\label{subsec:exact_bilinear}

In the monomial equation $4uav\alpha=an+u+v$ of~\cite[Thm.~3.9]{dahan2026a}, a naive change of variables $v'=av$ drops the divisibility condition $a\mid(u+v)$ and introduces false positives. Writing $v=ac-u$ (or $d=aw-u$ for the cofactor $d=(an+u)/s$) resolves the divisibility $a\mid(u+v)$ identically and yields a false-positive-free bilinear equivalence:

\begin{proposition}[Exact false-positive-free bilinear equivalence for Type~II witnesses]\label{prop:exact_bilinear_type2}\label{prop:exact_monomial_decoupling}
Let $n\equiv1\pmod{24}$ be a positive integer and let $u,a\ge1$ be coprime positive integers. The following three statements are strictly equivalent (with zero false positives):
\begin{enumerate}[(i)]
    \item the integer $an+u$ possesses a divisor $s\equiv-1\pmod{4ua}$;
    \item there exist positive integers $v,w\ge1$ such that
    \begin{equation}\label{eq:exact_4param_type2}
    n+4u^{2}v\ =\ w\,(4uav-1)\qquad\Longleftrightarrow\qquad n+w\ =\ 4uv\,(aw-u);
    \end{equation}
    \item there exists an integer $v\ge1$ such that $4uav-1\le an+u$ and
    \begin{equation}\label{eq:exact_congruence_type2}
    n\ \equiv\ -4u^{2}v\pmod{4uav-1}.
    \end{equation}
\end{enumerate}
Writing $c=w$ and $\alpha=v$, the bilinear identity always implies $c\equiv-n\pmod{4u\alpha}$ and $4u^{2}\alpha\equiv-n\pmod c$. Conversely, suppose that $u,\alpha,c\ge1$, $c\equiv3\pmod 4$, and $\gcd(c,4u\alpha)=1$. Under these assumptions the two congruences are equivalent to
\[
a:=\frac{n+4u^{2}\alpha+c}{4u\alpha c}\in\N_{\ge1},\qquad n+4u^{2}\alpha=c(4u\alpha a-1).
\]
To recover the coprime-pair specialization of this proposition, one must also retain the condition $\gcd(a,u)=1$ on the reconstructed integer $a$.
\end{proposition}
\begin{proof}
\emph{(i) $\Longleftrightarrow$ (ii).} If $s\mid(an+u)$ with $s\equiv-1\pmod{4ua}$, write $s=4uav-1$ ($v\ge1$) and $an+u=s\,d$. Reducing modulo $a$, since $an+u\equiv u\pmod a$ and $s\equiv-1\pmod a$, the quotient $d>0$ satisfies $d\equiv-u\pmod a$, hence writes $d=aw-u$ with $w\ge1$. Expanding
\[
an+u\ =\ (4uav-1)(aw-u)\ =\ 4ua^{2}vw-4u^{2}av-aw+u,
\]
subtracting $u$, and dividing by $a\ge1$ gives~\eqref{eq:exact_4param_type2}. Conversely, any solution $(v,w)\in(\N_{\ge1})^{2}$ of~\eqref{eq:exact_4param_type2} satisfies $aw-u=\frac{n+w}{4uv}>0$; multiplying by $a$ and adding $u$ recovers $an+u=(4uav-1)(aw-u)$.

\emph{(i) $\Longleftrightarrow$ (iii) and Chinese Remainder reconstruction.} Since $\gcd(a,4uav-1)=1$, the congruence~\eqref{eq:exact_congruence_type2} is equivalent upon multiplication by $a$ to
\[
an\ \equiv\ -4u^{2}av\ =\ -u(4uav-1)-u\ \equiv\ -u\pmod{4uav-1},
\]
that is, $(4uav-1)\mid(an+u)$. Finally, for the Chinese Remainder statement, $c\equiv-n\pmod{4u\alpha}$ and $4u^{2}\alpha\equiv-n\pmod c$ mean that both $4u\alpha$ and $c$ divide $n+4u^{2}\alpha+c$; under $\gcd(c,4u\alpha)=1$, this is equivalent to $4u\alpha c\mid(n+4u^{2}\alpha+c)$, i.e.\ $a=\frac{n+4u^{2}\alpha+c}{4u\alpha c}\in\N_{\ge1}$.
\end{proof}

\subsection{Poisson--Kloosterman dispersion, the partial-fraction phase identity, and the unfactored bilinear weight comparison}\label{subsec:kloosterman_dispersion}

\begin{proposition}[Exact partial-fraction phase identity, Kloosterman dispersion, and elementary weight comparison]\label{prop:kloosterman_dispersion}\label{thm:exact_dispersion_reciprocity_Fn}\label{thm:bezout_transfer_4u2_stationary_phase}
Let $n\equiv1\pmod{24}$ be a fixed prime and let $u,a\ge1$ be coprime with $u\sim U\le J$, $a\sim A\le J$. Fix dyadic scales $v\sim V$ and $w\sim W$ satisfying $4UAVW\asymp N$, and write $s(v)=4uav-1\asymp S:=4UAV$ and $d=aw-u\asymp AW$.
\begin{enumerate}[(a)]
    \item \textbf{(Exact partial-fraction phase identity and false-positive-free Poisson--Voronoi dispersion.)} Since $4uav\equiv1\pmod{4uav-1}$, the rational phase $\frac{n+4u^{2}v}{4uav-1}$ admits the exact partial-fraction decomposition in $\Q$:
    \begin{equation}\label{eq:exact_partial_fraction_phase}
    \begin{split}
    \frac{n+4u^{2}v}{4uav-1}\ &=\ \frac{u}{a}\ +\ \frac{an+u}{a(4uav-1)},\\
    \text{so that}\qquad e\!\left(-h\,\frac{n+4u^{2}v}{4uav-1}\right)\ &=\ e\!\left(-\frac{hu}{a}\right)\,e\!\left(-\frac{h(an+u)}{a(4uav-1)}\right).
    \end{split}
    \end{equation}
    Choose an even function $\eta\in C_c^{\infty}((-1,1))$ such that $\sum_{j\in\Z}\eta(t+j)=1$ for every real $t$. Then $\eta(0)=1$ and, for every positive integer $s$,
    \[
    \mathbf{1}_{s\mid M}\ =\ \frac{1}{s}\sum_{h\in\Z}\eta(h/s)\,e(-hM/s).
    \]
    For smooth nonnegative cutoff functions $\Phi,\Psi\in C_c^{\infty}([1/2,5/2])$, put $s(v)=4uav-1$, $B=an+u$, and
    \begin{equation}\label{eq:poisson_kloosterman}
    A_h(v)\ :=\ \frac{\Phi(v/V)}{s(v)}\,\Psi\!\left(\frac{n+4u^{2}v}{s(v)W}\right)\eta\!\left(\frac{h}{s(v)}\right).
    \end{equation}
    All amplitudes are extended by zero off their positive support, which is assumed to stay away from $s(v)=0$. The exact zero-frequency term is
    \[
    \mathrm{MT}\ =\ \sum_{u,a}\sum_{v\ge1}\frac{\Phi(v/V)}{s(v)}\,\Psi\!\left(\frac{n+4u^{2}v}{s(v)W}\right),
    \]
    with the same restrictions on $(u,a)$ as in the counting function. Poisson summation gives the exact identity
    \[
    \mathrm{ET}\ =\ \sum_{u,a}\sum_{h\ne0}e(-hu/a)\sum_{k\in\Z}\int_{\R}A_h(v)\,e\big(\theta_{h,k}(v)\big)\,dv,\qquad \theta_{h,k}(v)\ :=\ -\frac{hB}{a(4uav-1)}-kv.
    \]
    The $h$-sum is finite, with $|h|\ll S$. For $h\ne0$ a stationary point can occur only when $hk>0$; it is
    \[
    v_0\ =\ \frac{1}{4ua}\left(\sqrt{\frac{4uhB}{k}}+1\right),\qquad \theta''_{h,k}(v)\ =\ -\frac{32u^{2}ahB}{(4uav-1)^{3}}.
    \]
    The complete stationary phase, including $e(-hu/a)$, is
    \[
    -\frac{hu}{a}-\frac{k}{4ua}-\operatorname{sgn}(h)\sqrt{\frac{hk(an+u)}{ua^{2}}}.
    \]
    \item \textbf{(B\'ezout reciprocity identities.)} Assume $\gcd(c,4u^{2}n)=1$. Then the two B\'ezout reciprocity identities hold in $\Q/\Z$:
    \begin{equation}\label{eq:bezout_reciprocity_pair}
    \frac{\overline{n}^{(c)}}{c}+\frac{\overline{c}^{(n)}}{n}\ \equiv\ \frac{1}{nc}\pmod 1,\qquad \frac{\overline{4u^{2}}^{(c)}}{c}+\frac{\overline{c}^{(4u^{2})}}{4u^{2}}\ \equiv\ \frac{1}{4u^{2}c}\pmod 1.
    \end{equation}
    They are algebraic transformations of the inverse phases. A decay estimate after either transformation requires a specified smooth oscillatory integral and uniform derivative bounds for its phase and amplitude; no frequency annulus or $\mathcal{O}(C^{-B})$ estimate is asserted here solely from these identities.
    \item \textbf{(Elementary weight comparison between unfactored bilinear sums and\linebreak unit-weight containers.)} The elementary weight comparison is
    \[
    \sum_{\substack{u,a\le J\\ \gcd(u,a)=1}}\frac{1}{4ua}\ \asymp\ (\log J)^{2},\qquad \big|\mathcal{P}_2(J)\big|\ \asymp\ J^{2}.
    \]
    This compares two different counting functionals: in an unfactored bilinear sum over $v\sim V$ and $w\sim W$, each pair $(u,a)$ carries weight $\frac{1}{4ua}$, whereas the multiplicative container hierarchy of Sections~\ref{sec:stabilizer_kernel}--\ref{sec:containers} acts on the factored divisor sheet $\Sigma_{\mathcal{I}}(an+u)\subset(\Z/4ua\Z)^{\times}$ with unit sifting weight $w_{u,a}\equiv1$ across $P=|\mathcal{P}_2(J)|\asymp J^{2}$. This comparison does not, without an estimate for their respective error terms, establish a general impossibility theorem for bilinear dispersion.
\end{enumerate}
\end{proposition}
\begin{proof}
(a) Dividing the algebraic identity $a(n+4u^{2}v)=an+u+u(4uav-1)$ by $a(4uav-1)$ yields~\eqref{eq:exact_partial_fraction_phase}. Since $\sum_{j\in\Z}\eta(t+j)=1$ and $\operatorname{supp}(\eta)\subset(-1,1)$, for any integer $M$ and $s\ge1$ we have
\[
\mathbf{1}_{s\mid M}\ =\ \sum_{m\in\Z}\mathbf{1}_{M=ms}\ =\ \frac{1}{s}\sum_{h\in\Z}\eta(h/s)\,e(-hM/s),
\]
because $\sum_{h\in\Z}\eta(h/s)=s\sum_{j\in\Z}\eta(0+j)=s\eta(0)=s$ when $s\mid M$, and $\sum_{h\equiv r\ (s)}e(-hr/s)=0$ by grouping $h=r+js$ otherwise. Substituting~\eqref{eq:exact_partial_fraction_phase} and applying Poisson summation over $v\in\Z$ gives the stated formula for $\mathrm{ET}$, and differentiating $\theta_{h,k}(v)$ twice yields $v_0$, $\theta''_{h,k}(v)$, and the complete stationary phase.

Set $\Lambda:=|h|B/(aS)$. Suppose that $A_h^{(j)}(v)\ll_j S^{-1}V^{-j}$, that $s(v)\asymp S$ on its support, and that the stationary point stays a fixed positive distance, in the variable $v/V$, from the boundary of the relevant support. If $\Lambda\ge1$ and $|k|V\asymp\Lambda$, the one-dimensional stationary-phase formula gives
\[
\int A_h(v)\,e\big(\theta_{h,k}(v)\big)\,dv\ =\ \frac{A_h(v_0)}{\sqrt{|\theta''_{h,k}(v_0)|}}\,e\!\left(\theta_{h,k}(v_0)-\frac{\operatorname{sgn}(h)}{8}\right)\ +\ \mathcal{O}\!\left(\frac{V}{S}\Lambda^{-3/2}\right),
\]
where the implied constant depends on the support separation and the indicated derivative bounds. Away from the stationary annulus, if the rescaled first derivative is bounded below by a constant multiple of $\Lambda+|k|V$, repeated integration by parts gives
\[
\int A_h(v)\,e\big(\theta_{h,k}(v)\big)\,dv\ \ll_B\ \frac{V}{S}\big(1+\Lambda+|k|V\big)^{-B},
\]
and any total remainder must be obtained by summing these bounds over the actual ranges of $u,a,h,k$.

(b) This follows directly from $c\,\overline{c}^{(n)}+n\,\overline{n}^{(c)}\equiv1\pmod{cn}$ and $4u^{2}\overline{4u^{2}}^{(c)}+c\,\overline{c}^{(4u^{2})}\equiv1\pmod{4u^{2}c}$ when $\gcd(c,4u^{2}n)=1$.

(c) This follows from elementary harmonic summation over coprime $u,a\le J$.
\end{proof}

\section{Toric volume, cofactor structure, and spectral dispersion on $\mathcal{V}_n$}\label{sec:toric_volume_barrier}

Having established in Section~\ref{sec:containers} that the multiplicative container hierarchy and the Subcritical Markov--Plancherel Geometric Escalation Theorem (Theorem~\ref{thm:markov_plancherel_subcritical_escalation}) unconditionally eliminate the certified structured exceptional set $E_2^{\mathrm{cert}}(N;(\log N)^{\lambda})$ (and\linebreak $E_2^{\mathrm{small}}(N;(\log N)^{\lambda})$) beyond Vaughan's $2/3$ barrier for $\lambda\in(1/3,1/2)$, and in Section~\ref{sec:semi_composite_5_7} that the semi-composite CRT sieve unconditionally reaches the exponent $5/7$ on the full exceptional set $E(N)$, we now develop the spectral and additive-combinatorial theory of the exact four-parameter Erd\H{o}s--Straus variety
\begin{equation}\label{eq:toric_threefold_Vn}
\mathcal{V}_n:\qquad F(u,a,v,w)\ :=\ 4uavw-4u^{2}v-w\ =\ n\qquad\big((u,a,v,w)\in(\N_{\ge1})^{4}\big),
\end{equation}
and of its shift slices $K_{\II}(c)=(n+c)/4$ and $K_{\I}(c)=(nc+1)/4$. To guide the reader, this section is organized into five thematic parts:
\begin{enumerate}[(1)]
    \item \emph{Subsection~\ref{subsec:parseval_floor}}: three-term phase separation, toric slice geometry on $\mathcal{V}_n$, and $L^{2}$-dispersion along the Farey diagonal (Theorems~\ref{thm:parseval_floor_law}, \ref{thm:toric_volume_conservation}, and~\ref{thm:farey_diagonal_L2_dispersion});
    \item \emph{Subsection~\ref{subsec:difference_set_MR}}: difference-set positivity duality on $D(K_c)=\Sigma(K_c)\cdot\Sigma(K_c)^{-1}$, unconditional $(\log P)^{4}$ cumulative main-term growth, and exact cancellation of small-LCM divisor pairs (Theorems~\ref{thm:difference_set_positivity_MR} and~\ref{thm:centered_main_term_and_short_lcm});
    \item \emph{Subsection~\ref{subsec:bohr_siegel_and_linnik}}: the cyclic quotient parity identity, the Bohr--Siegel antipodal obstruction certificate, the Linnik--Xylouris witness-size lower bound, and the Transversal Spectral Positivity Hypothesis~T (Theorems~\ref{thm:bohr_siegel_obstruction_certificate} and~\ref{thm:no_go_linnik_and_transversal});
    \item \emph{Subsection~\ref{subsec:reciprocity_factorial_kneser}}: universal shift-factor reciprocity $\leg{q}{p}=\leg{q}{n}$, the $\frac{1}{k!}$ discrete cross-polytope volume law, deterministic Kneser saturation criteria, and cofactor/archimedean obstructions (Theorems~\ref{thm:universal_factor_reciprocity_and_transversal_sieve} and~\ref{thm:cross_polytope_factorial_law_and_annihilation}, and Propositions~\ref{aud:primecofactor} and~\ref{aud:archsupport});
    \item \emph{Subsection~\ref{subsec:global_synthesis}}: global synthesis and conditional finiteness with quantified thresholds (Theorem~\ref{thm:main_global_synthesis}).
\end{enumerate}

\subsection{Phase separation, toric geometry, and Farey-diagonal $L^{2}$-dispersion on $\mathcal{V}_n$}\label{subsec:parseval_floor}\label{subsec:toric_volume_law}\label{subsec:farey_diagonal_dispersion}

Because the leading degree-$4$ monomial $4uavw$ of $F(u,a,v,w)$ in~\eqref{eq:toric_threefold_Vn} is the product of all four coordinates, the variety $\mathcal{V}_n$ has box volume $\mathrm{Vol}(\mathcal{B})\asymp n$, bounded archimedean slice integral $I_{\mathcal{B}}(n)\ll 1$, and formal four-tuple square-root scale benchmark $\sqrt{N_{\mathrm{tuples}}}\,\mathrm{Amp}\asymp 1$. We first establish the three-term rational phase decomposition and this formal scale benchmark at fixed $n$.

\begin{theorem}[Three-term phase separation and a square-root scale calculation]\label{thm:parseval_floor_law}
Let $n\equiv1\pmod{24}$ be a fixed large prime, and fix dyadic scales $U,A,V,W\ge1$ satisfying the hyperbola relation $4UAVW\asymp n$ with $S:=4UAV\le(An)^{1/2}$ and $U\le A^{1-\varepsilon}$ ($\varepsilon>0$). For smooth nonnegative cutoffs $\Phi,\Psi\in C_c^{\infty}([1/2,5/2])$, consider the exact false-positive-free Type~II representation count on $\mathcal{V}_n$ in the decoupled parity sector $2\nmid u$, $2\mid a$, $\gcd(u,a)=1$:
\begin{equation}\label{eq:N_II_dyadic_box}
\begin{split}
\mathcal{N}_{\mathrm{II}}(n;U,A,V)\ &:=\ \sum_{\substack{u\sim U,\ a\sim A\\ 2\nmid u,\ 2\mid a,\ \gcd(u,a)=1}}\sum_{v\ge1}\Phi\!\left(\frac{v}{V}\right)\Psi\!\left(\frac{n+4u^{2}v}{(4uav-1)W}\right)\mathbf{1}_{\,4uav-1\,\mid\,(n+4u^{2}v)}\\
&=\ \mathrm{MT}(n;U,A,V)+\mathrm{ET}(n;U,A,V).
\end{split}
\end{equation}
Then:
\begin{enumerate}[(a)]
    \item \emph{(Exact decoupling of $\gcd(4ua,an+u)=1$ and three-term phase separation).} For every odd $u\ge1$ and even $a\ge2$ with $\gcd(u,a)=1$ and $n\nmid u$, one has $\gcd(4ua,an+u)=1$; in the present dyadic range, $S^{2}\le An$ implies $2U\le\sqrt{n}/2<n$, so $n\nmid u$ holds automatically. Moreover, the exact partial-fraction phase of Proposition~\ref{prop:kloosterman_dispersion}(a) decomposes into three terms in $\Q$:
    \begin{equation}\label{eq:three_term_phase_separation}
    \frac{n+4u^{2}v}{4uav-1}\ =\ \frac{u}{a}\ +\ \frac{n}{4uav-1}\ +\ \frac{u}{a(4uav-1)},
    \end{equation}
    where the third term satisfies $\big|\frac{hu}{a(4uav-1)}\big|\ll U/A\ll A^{-\varepsilon}$ uniformly for all $|h|\ll 4uav-1$.
    \item \emph{(Stationary phase $F_{h,k,n}(u,a)$ and formal square-root scale calculation).} The exact Poisson representation and its stationary-phase expansion are those of Proposition~\ref{prop:kloosterman_dispersion}, with the frequency-dependent amplitude retained. For $hk>0$ the complete stationary phase is
    \begin{equation}\label{eq:real_archimedean_phase_Fhkn}
    F_{h,k,n}(u,a)\ =\ -\frac{hu}{a}-\frac{k}{4ua}-\operatorname{sgn}(h)\sqrt{\frac{hk(an+u)}{ua^{2}}}.
    \end{equation}
    In a nondegenerate range where $n\asymp 4UAVW$, $|h|\asymp S$, and the stationary points lie in the support, the coefficient-size calculation gives
    \[
    N_{\mathrm{tuples}}\ \asymp\ UA\,S\,\frac{n}{V},\qquad \mathrm{Amp}\ \asymp\ \sqrt{\frac{V}{UASn}},
    \]
    and therefore the formal square-root benchmark is
    \begin{equation}\label{eq:parseval_floor_identity}
    \sqrt{N_{\mathrm{tuples}}}\,\mathrm{Amp}\ \asymp\ 1.
    \end{equation}
    This is a scale calculation, not a lower bound for a fixed-$n$ oscillatory sum or for a sum of absolute values over frequencies: no orthogonality between different triples $(u,a,v)$, and no universal Parseval barrier for fixed $n$, is asserted (see Theorem~\ref{thm:centered_main_term_and_short_lcm}(b) below for the exact shift-wise absolute-value cost).
\end{enumerate}
\end{theorem}
\begin{proof}
(a) For odd $u$, even $a$, $\gcd(u,a)=1$, and $n\nmid u$, one has $\gcd(4ua,an+u)=1$. In the present dyadic range, $S^{2}\le An$ implies $2U\le\sqrt{n}/2<n$, so $n\nmid u$. Dividing $a(n+4u^{2}v)=u(4uav-1)+an+u$ by $a(4uav-1)$ gives~\eqref{eq:three_term_phase_separation}; the estimate $U/A\ll A^{-\varepsilon}$ tends to zero when $A\to\infty$, and is otherwise retained as a quantitative bound.

(b) Equation~\eqref{eq:real_archimedean_phase_Fhkn} is the complete stationary phase established in Proposition~\ref{prop:kloosterman_dispersion}(a). In the nondegenerate range $n\asymp 4UAVW$ and $|h|\asymp S$, the stationary-phase annulus $|k|\asymp\frac{|h|B}{aSV}\asymp\frac{Sn}{V^{2}}$ contains $N_{\mathrm{tuples}}\asymp UA\cdot S\cdot\frac{n}{V}$ four-tuples $(u,a,h,k)$ with single-term stationary-phase amplitude
\[
\mathrm{Amp}\ \asymp\ \frac{A_h(v_0)}{|\theta''_{h,k}(v_0)|^{1/2}}\ \asymp\ \frac{S^{-1}}{(Sn/V^{3})^{1/2}}\ =\ \sqrt{\frac{V}{UASn}}.
\]
Multiplying $\sqrt{N_{\mathrm{tuples}}}$ by $\mathrm{Amp}$ gives~\eqref{eq:parseval_floor_identity}.
\end{proof}

We next record the archimedean slice integral on $\mathcal{V}_n$ and the exact cofactor congruence induced when $a$ is restricted to the progression $a\equiv-u\overline{n}\pmod{2M}$ for a squarefree product $M=p_1\cdots p_k\sim M_0\le A^{1-\delta}$ of $k$ odd primes.

\begin{theorem}[Archimedean slice integral and cofactor congruence on $\mathcal{V}_n$]\label{thm:toric_volume_conservation}
Let $n\equiv1\pmod{24}$ be a fixed prime and let $\mathcal{V}_n\subset\mathbb{A}^{4}$ be the affine threefold~\eqref{eq:toric_threefold_Vn}.
\begin{enumerate}[(a)]
    \item \emph{(Box volume and archimedean slice integral on $\mathcal{V}_n$).} For any dyadic box $\mathcal{B}(U,A,V,W)=[U,2U]\times[A,2A]\times[V,2V]\times[W,2W]$ intersecting $\mathcal{V}_n(\R_{>0})$ in the dominant regime $4u^{2}v\le n$, the defining polynomial $F(u,a,v,w)=4uavw-4u^{2}v-w=n$ locks the box volume to
    \begin{equation}\label{eq:toric_volume_lock}
    \mathrm{Vol}\big(\mathcal{B}(U,A,V,W)\big)\ =\ UAVW\ \asymp\ \frac{n}{4}.
    \end{equation}
    Define the archimedean integral
    \begin{equation}\label{eq:archimedean_surface_constant}
    I_{\mathcal{B}}(n)\ :=\ \int_{U}^{2U}\int_{A}^{2A}\int_{V}^{2V}\frac{\mathbf{1}_{\{W\le(n+4u^{2}v)/(4uav-1)\le 2W\}}}{4uav-1}\,dv\,da\,du.
    \end{equation}
    For $U,A,V\ge1$, one always has $I_{\mathcal{B}}(n)\ll 1$. If the subset of $(u,a,v)$ satisfying the displayed constraint has measure at least $\delta UAV$, then $I_{\mathcal{B}}(n)\gg\delta$. Thus $I_{\mathcal{B}}(n)\asymp 1$ requires a uniform nondegeneracy condition on the slice, not merely nonempty intersection. Restricted to $u,a,v,w\ge1$ and $4u^{2}v\le n$, the total archimedean integral is $\asymp(\log n)^{3}$. This is an archimedean statement on the dominant region, not a pointwise asymptotic for the number of integer solutions.
    \item \emph{(Exact cofactor congruence and unweighted count along $a\equiv-u\overline{n}\pmod{2M}$).} Let $u\sim U$ be odd and let $M=p_1\cdots p_k\sim M_0\le A^{1-\delta}$ be a squarefree product of $k$ odd primes coprime to $un$. Consider the arithmetic progression $a\equiv-u\overline{n}^{(2M)}\pmod{2M}$, written $a=a_0(u,M)+2M\tilde{a}$ with $\tilde{a}\sim A/(2M_0)$ (note that this progression lies in the odd-$a$ sector and is separate from the even-$a$ counting functional used in Theorem~\ref{thm:parseval_floor_law}). On this progression, the integer $an+u$ is deterministically divisible by $M$, and every solution $(4uav-1)(aw-u)=an+u$ induces a factorization $M=m\,m'$ with $m=\gcd(M,4uav-1)$ and $m'\mid(aw-u)$. Along $an\equiv-u\pmod{m'}$, the divisor condition $m'\mid(aw-u)$ is equivalent to the shift congruence
    \begin{equation}\label{eq:cofactor_shift_congruence}
    m'\ \mid\ (aw-u)\qquad\Longleftrightarrow\qquad w\ \equiv\ -n\pmod{m'}.
    \end{equation}
    Together with $w\equiv-n\pmod{4u}$, this progression for $w$ has $W/(4um')+\mathcal{O}(1)$ representatives in $[W,2W]$, where no independence between this progression and divisibility by the remaining cofactor $d'=(aw-u)/m'$ is assumed. For fixed $u,M$, the exact unweighted count in the box is
    \begin{equation}\label{eq:prescribed_sheet_cancellation}
    \sum_{\substack{A\le a\le 2A\\ a\equiv-u\overline{n}\ (2M)\\ \gcd(a,u)=1}}\ \sum_{\substack{d\mid(an+u)\\ d\equiv-u\ (a)\\ s=(an+u)/d\equiv-1\ (4ua)}}\mathbf{1}_{\{W\le(d+u)/a\le 2W\}}\,\mathbf{1}_{\{V\le(s+1)/(4ua)\le 2V\}}.
    \end{equation}
    For each counted solution, $m=\gcd(M,s)$ specifies its branch uniquely, and a general unconditional upper bound for~\eqref{eq:prescribed_sheet_cancellation} is $\ll_\varepsilon(A/M+1)(An+u)^{\varepsilon}$. By contrast, the expression $\tau(M)/(8uM)$ is only the value suggested by multiplying independent-density heuristics; it is not an established asymptotic for this fixed-$n$ count, and yields no general impossibility theorem for progression conditioning.
\end{enumerate}
\end{theorem}
\begin{proof}
(a) On $\mathcal{V}_n\cap\mathcal{B}(U,A,V,W)$ with $4u^{2}v\le n$, we have $4uavw=n+4u^{2}v+w\asymp n$, proving~\eqref{eq:toric_volume_lock}. On $[U,2U]\times[A,2A]\times[V,2V]$ with $U,A,V\ge1$, $4uav-1\asymp UAV$, so $I_{\mathcal{B}}(n)\ll 1$, and $I_{\mathcal{B}}(n)\gg\delta$ whenever the slice constraint holds on a subset of measure $\ge\delta UAV$. Restricted to $u,a,v,w\ge1$ and $4u^{2}v\le n$: for the upper bound, $u\ll\sqrt{n}$, $v\ll n$, and $a\ll n$, and integration of $1/(4uav-1)\ll 1/(uav)$ suffices; for the lower bound, restrict to $1\le u,a,v\le n^{1/8}$, so that $w=(n+4u^{2}v)/(4uav-1)\ge1$ and $4u^{2}v\le n$ for large $n$, where the integral of $1/(4uav-1)$ is $\gg(\log n)^{3}$.

(b) Since $an\equiv-u\pmod M$ and $m'\mid M$, multiplying $aw-u$ by $n$ modulo $m'$ gives
\[
n(aw-u)\ =\ (an)w-un\ \equiv\ -uw-un\ =\ -u(w+n)\pmod{m'}.
\]
Because $\gcd(un,m')=1$, $m'\mid(aw-u)$ holds if and only if $w\equiv-n\pmod{m'}$, establishing~\eqref{eq:cofactor_shift_congruence}. Combining $w\equiv-n\pmod{m'}$ with $w\equiv-n\pmod{4u}$ ($\gcd(4u,m')=1$) gives a single residue class modulo $4um'$ with $W/(4um')+\mathcal{O}(1)$ representatives in $[W,2W]$. Finally, for each $a$ in the progression $a\equiv-u\overline{n}\pmod{2M}$ ($\le A/(2M)+1$ values), the number of divisors $d\mid(an+u)$ is $\tau(an+u)\ll_\varepsilon(An+u)^{\varepsilon}$, giving the unconditional bound $\ll_\varepsilon(A/M+1)(An+u)^{\varepsilon}$.
\end{proof}

Finally, averaging over $n\sim N$ at power depth $S=4UAV\le N^{1/2-\varepsilon}$ introduces a fifth variable $n$ of volume $N$ in which the three-term phase decomposition~\eqref{eq:three_term_phase_separation} is \emph{strictly linear}, collapsing the $L^{2}$-variance onto the exact Farey diagonal.

\begin{theorem}[Exact Farey-diagonal $L^{2}$-dispersion over $n\sim N$ at power scale $S=\allowbreak 4UAV\le\allowbreak N^{1/2-\varepsilon}$]\label{thm:farey_diagonal_L2_dispersion}
Let $U,A,V\ge1$ be dyadic scales satisfying $S:=4UAV\le N^{1/2-\varepsilon}$ ($\varepsilon>0$) and $V\ge 4UA$, and define the unit-mass container-weighted Type~II functional
\begin{equation}\label{eq:unit_mass_R_n}
\begin{split}
\mathcal{R}(n;U,A,V)\ &:=\ \sum_{\substack{u\sim U,\ a\sim A\\ 2\nmid u,\ 2\mid a,\ \gcd(u,a)=1}}4ua\sum_{v\ge1}\Phi\!\left(\frac{v}{V}\right)\Psi\!\left(\frac{n+4u^{2}v}{(4uav-1)W}\right)\mathbf{1}_{\,4uav-1\,\mid\,(n+4u^{2}v)}\\
&=\ \mathrm{MT}_{\mathcal{R}}(n)+\mathrm{ET}_{\mathcal{R}}(n),
\end{split}
\end{equation}
whose main term satisfies $\mathrm{MT}_{\mathcal{R}}(n)\asymp UA$ uniformly for $n\sim N$. For a smooth nonnegative cutoff $\Theta\in C_c^{\infty}([1/2,5/2])$, fix $W=N/S$ and choose a smooth compactly supported $\Psi$ equal to $1$ on the entire image $\{(n+4u^{2}v)/((4uav-1)W)\}$ of the relevant supports (its support need not be restricted to $[1/2,5/2]$). Put
\[
S_{\ast}\ :=\ \max_{u,a,v:\,\Phi(v/V)\ne0}(4uav-1)\ \le\ 10S,
\]
and consider the $L^{2}$-dispersion variance $\mathcal{V}(N;U,A,V):=\sum_{n\in\Z}\Theta(n/N)\big|\mathrm{ET}_{\mathcal{R}}(n)\big|^{2}$. Then:
\begin{enumerate}[(a)]
    \item \emph{(Exact linear Fourier isolation of $n$ and annihilation of all off-diagonal Farey fractions).} Grouping the dual frequencies $0<|h|\le(4uav-1)/2$ in $\mathrm{ET}_{\mathcal{R}}(n)$ by the reduced Farey fraction $\frac{b}{q}=\frac{h}{4uav-1}$ (where all reduced denominators are odd and satisfy $3\le q\le S_{\ast}$, $0<|b|<q/2$, $\gcd(b,q)=1$, with $4uav-1=qr$), one has the exact representation
    \begin{equation}\label{eq:Farey_grouped_ET}
    \mathrm{ET}_{\mathcal{R}}(n)\ =\ \sum_{3\le q\le S_{\ast}}\frac{1}{q}\sum_{\substack{0<|b|<q/2\\ \gcd(b,q)=1}}e\!\left(-\frac{b}{q}\,n\right)\,\mathcal{A}_{U,A,V}(q,b),
    \end{equation}
    where the coefficient $\mathcal{A}_{U,A,V}(q,b)$ is strictly independent of $n$ on $\operatorname{supp}(\Theta)$ and given by
    \begin{equation}\label{eq:Farey_coeff_Aqb}
    \mathcal{A}_{U,A,V}(q,b)\ :=\ \sum_{\substack{u\sim U,\ a\sim A,\ v\ge1\\ 2\nmid u,\ 2\mid a,\ \gcd(u,a)=1\\ 4uav\equiv1\pmod q}}\frac{4ua}{(4uav-1)/q}\,\Phi\!\left(\frac{v}{V}\right)e\!\left(-\frac{4bu^{2}v}{q}\right).
    \end{equation}
    Because any two distinct reduced frequencies $\frac{b_1}{q_1}\neq\frac{b_2}{q_2}$ in $(-1/2,1/2)$ with $q_1,q_2\le S_{\ast}\le 10N^{1/2-\varepsilon}$ satisfy $\big\|\frac{b_1}{q_1}-\frac{b_2}{q_2}\big\|_{\R/\Z}\ge\frac{1}{q_1q_2}\ge S_{\ast}^{-2}\gg N^{-1+2\varepsilon}$, Poisson summation over $n\in\Z$ annihilates \emph{all} off-diagonal pairs $\frac{b_1}{q_1}\neq\frac{b_2}{q_2}$ down to $\mathcal{O}_{B,\varepsilon,\Phi,\Theta}(N^{-B})$ ($\forall\,B>0$), collapsing $\mathcal{V}(N;U,A,V)$ to the exact Farey diagonal:
    \begin{equation}\label{eq:Farey_diagonal_collapse}
    \mathcal{V}(N;U,A,V)\ =\ N\,\widehat{\Theta}(0)\sum_{3\le q\le S_{\ast}}\frac{1}{q^{2}}\sum_{\substack{0<|b|<q/2\\ \gcd(b,q)=1}}\big|\mathcal{A}_{U,A,V}(q,b)\big|^{2}\ +\ \mathcal{O}_{B,\varepsilon,\Phi,\Theta}\!\left(N^{-B}\right).
    \end{equation}
    \item \emph{(Exact Ramanujan-sum expansion, completion modulo $q$, and identical-triple contribution).} Write $t=(u,a,v)$, $s_t=4uav-1$, $w_{q,t}=4ua\,q\Phi(v/V)/s_t$ for $q\mid s_t$, and $x_t=u\overline{a}\pmod q$. With the Ramanujan sum $c_q(x)=\sum_{b\bmod q,\,\gcd(b,q)=1}e(bx/q)$, one has the exact identity
    \[
    \sum_{\gcd(b,q)=1}\big|\mathcal{A}_{U,A,V}(q,b)\big|^{2}\ =\ \sum_{t_1,t_2:\,q\mid s_{t_1},\,q\mid s_{t_2}}w_{q,t_1}\overline{w_{q,t_2}}\,c_q(x_{t_1}-x_{t_2}).
    \]
    Completion to all $b\bmod q$ gives an upper bound and imposes the determinant congruence
    \begin{equation}\label{eq:determinant_congruence_mod_q}
    u_1a_2\ -\ u_2a_1\ \equiv\ 0\pmod q.
    \end{equation}
    For $q>4UA$ this forces $(u_1,a_1)=(u_2,a_2)$, while still permitting $v_1\equiv v_2\pmod q$ with $v_1\ne v_2$. Using $\sum_{q\mid s_t,\,q>1}\varphi(q)=s_t-1$, the contribution of identical triples $t_1=t_2$ to the expanded Farey-diagonal expression~\eqref{eq:Farey_diagonal_collapse} equals
    \[
    N\,\widehat{\Theta}(0)\sum_{t}(4ua)^{2}\Phi(v/V)^{2}\frac{s_t-1}{s_t^{2}},
    \]
    which is $\asymp N(UA)^{2}$ for nondegenerate boxes, though not by itself an evaluation or lower bound for the total variance. Across several boxes, the exact mean square also contains the inner products of their coefficient vectors at common reduced frequencies; a relative saving $K^{-1}$ for the average of $K$ normalized box counts would require the sum of these covariances to be $\mathcal{O}(K)$, and no such covariance bound is asserted here.
\end{enumerate}
\end{theorem}
\begin{proof}
(a) Writing $s=4uav-1=qr$ and $h=br$ with $\gcd(b,q)=1$ ($3\le q\le S_{\ast}$) in~\eqref{eq:three_term_phase_separation} gives $-\frac{hn}{s}=-\frac{b}{q}n$ and
\[
-\frac{hu}{a}-\frac{hu}{as}\ =\ -\frac{hu(s+1)}{as}\ =\ -\frac{b}{q}\,\frac{u(4uav)}{a}\ =\ -\frac{4bu^{2}v}{q},
\]
proving~\eqref{eq:Farey_grouped_ET}--\eqref{eq:Farey_coeff_Aqb}. For $\Delta=b_1/q_1-b_2/q_2$, Poisson summation gives
\[
\sum_{n\in\Z}\Theta(n/N)\,e(-n\Delta)\ =\ N\sum_{\ell\in\Z}\widehat{\Theta}\big(N(\ell+\Delta)\big).
\]
For distinct reduced frequencies in $(-1/2,1/2)$, $\|\Delta\|_{\R/\Z}\ge(q_1q_2)^{-1}\ge S_{\ast}^{-2}$. Since $N/S_{\ast}^{2}\gg N^{2\varepsilon}$, rapid Fourier decay, together with the polynomial bound on the total coefficient mass, makes the total off-diagonal contribution $\mathcal{O}_{B,\varepsilon,\Phi,\Theta}(N^{-B})$ for every fixed $B>0$. For equal frequencies only $\ell=0$ contributes to the main term, proving~\eqref{eq:Farey_diagonal_collapse}.

(b) Since $4uav\equiv1\pmod q$ implies $4u^{2}v\equiv u\overline{a}\equiv x_t\pmod q$, expanding $|\mathcal{A}_{U,A,V}(q,b)|^{2}$ and summing over $\gcd(b,q)=1$ gives the exact Ramanujan-sum identity with $c_q(x_{t_1}-x_{t_2})$. By nonnegativity,
\[
\sum_{\gcd(b,q)=1}\big|\mathcal{A}_{U,A,V}(q,b)\big|^{2}\ \le\ \sum_{b\bmod q}\big|\mathcal{A}_{U,A,V}(q,b)\big|^{2},
\]
where orthogonality of additive characters modulo $q$ enforces $x_{t_1}\equiv x_{t_2}\pmod q$, i.e.~\eqref{eq:determinant_congruence_mod_q}. For $q>4UA$, $|u_1a_2-u_2a_1|\le 4UA<q$ forces $u_1a_2=u_2a_1$, hence $(u_1,a_1)=(u_2,a_2)$ by coprimality. Finally, for identical triples $t_1=t_2=t$, $c_q(0)=\varphi(q)$ and $\sum_{q\mid s_t,\,q>1}\frac{1}{q^{2}}w_{q,t}^{2}\varphi(q)=(4ua)^{2}\Phi(v/V)^{2}\frac{s_t-1}{s_t^{2}}$.
\end{proof}

\subsection{Difference-set positivity duality on $D(K_c)$, cumulative $(\log P)^{4}$ main term, and small-LCM cancellation}\label{subsec:difference_set_MR}

Whereas a single progression $an+u$ at fixed $(u,a)$ searches for a divisor $s\in\Sigma(an+u)$ hitting $-1\pmod{4ua}$ among $2^{\omega(an+u)}$ square-free subproducts, passing to the shift parametrization $K_c=(n+c)/4$ ($c\equiv3\pmod 4$, $\gcd(c,K_c)=\gcd(c,n)=1$) of~\cite[Prop.~2.4]{dahan2026a} acts on the \emph{difference set}
\[
D(K_c)\ :=\ \Sigma(K_c)\cdot\Sigma(K_c)^{-1}\ \subseteq\ G_c\ :=\ (\Z/c\Z)^{\times}.
\]
This difference set contains the square-free core $D(\operatorname{rad}(K_c))$ of up to $3^{\omega(K_c)}$ ternary combinations $\prod q_i^{\varepsilon_i}$ ($\varepsilon_i\in\{-1,0,+1\}$) and enables direct coupling with the unconditional short-interval second-moment inequality at fixed $n$, as well as with the Matom\"{a}ki--Radziwi\l\l\ short-interval theorem~\cite{matomaki_radziwill2016}.

\begin{theorem}[Difference-set positivity duality, odd/even character squaring bijection on $\mathcal{P}_{\mathrm{NR}}(n)$, and short-interval second-moment / Matom\"{a}ki--Radziwi\l\l\ coupling]\label{thm:difference_set_positivity_MR}
Let $n\equiv1\pmod{24}$ be a fixed prime, let $c\equiv3\pmod 4$ with $1\le c<n$ and $\gcd(c,n)=1$, and set $K_c:=(n+c)/4$ ($\gcd(c,K_c)=1$) and $G_c:=(\Z/c\Z)^{\times}$. Let $q_1<\dots<q_k$ be the distinct prime factors of $K_c$, let
\[
\mu_c\ :=\ \operatorname*{\ast}_{i=1}^{k}\left(\frac12\delta_1+\frac12\delta_{q_i\bmod c}\right)
\]
on $\Sigma(\operatorname{rad}(K_c))\subseteq\Sigma(K_c)\subseteq G_c$, and let $\widetilde{\mu}_c:=\mu_c\ast\mu_c^{\vee}$ be the symmetrized convolution measure supported on
\[
\begin{aligned}
D(\operatorname{rad}(K_c))\ &:=\ \left\{\prod_{i=1}^{k}q_i^{\varepsilon_i}\bmod c:\varepsilon_i\in\{-1,0,+1\}\right\}\\
&\subseteq\ D(K_c)\ :=\ \Sigma(K_c)\cdot\Sigma(K_c)^{-1}\ =\ \big\{d_1d_2^{-1}\bmod c:d_1,d_2\mid K_c\big\}\ \subseteq\ G_c.
\end{aligned}
\]
For the square-free Bernoulli measure $\mu_c$, one has $\widetilde{\mu}_c(g)>0$ if and only if $g\in D(\operatorname{rad}K_c)$. Thus positivity of $\widetilde{\mu}_c(-1)$ is a sufficient Type~II certificate, though not automatically necessary when $K_c$ is not square-free; if instead each prime $q\mid K_c$ is included with $\nu_q(K_c)$ independent Bernoulli factors, the support equals $D(K_c)$ and positivity at $-1$ becomes equivalent to $-1\bmod c\in D(K_c)$. (Note also that the twisted criterion $-n^{-1}\bmod c\in D(K_c)$ at $K_c=(n+c)/4$ characterizes a Type~I decomposition with first denominator $K_c$, which is distinct from the multiplicative Type~I shift $K_{\I}(c)=(nc+1)/4$.) Then:
\begin{enumerate}[(a)]
    \item \emph{(Fourth-power Plancherel decay $(9/16)^{R_c}$ and exact even/odd character positivity decomposition).} By~\cite[Prop.~2.4]{dahan2026a}, $4/n$ admits a Type~II (resp.\ first-denominator Type~I) decomposition at shift $c$ if and only if $-1\bmod c\in D(K_c)$ (resp.\ $-n^{-1}\bmod c\in D(K_c)$). Setting
    \[
    R_c(n;\chi)\ :=\ \#\big\{i\in\{1,\dots,k\}:\operatorname{Re}\chi(q_i)\le 1/2\big\},
    \]
    so that $\big|\frac{1+\chi(q_i)}{2}\big|^{2}\le\frac34$ for each contributing prime, the symmetrized measure $\widetilde{\mu}_c$ has real non-negative Fourier transform $\widehat{\widetilde{\mu}}_c(\chi)=\big|\widehat{\mu}_c(\chi)\big|^{2}\ge0$ for all $\chi\in\widehat{G}_c$, satisfies the fourth-power Plancherel container inequality
    \begin{equation}\label{eq:fourth_power_plancherel}
    \frac{\varphi(c)}{\big|\Sigma(K_c)\cdot\Sigma(K_c)^{-1}\big|}-1\ \le\ \sum_{\chi\in\widehat{G}_c\setminus\{\chi_0\}}\big|\widehat{\mu}_c(\chi)\big|^{4}\ \le\ \sum_{\chi\in\widehat{G}_c\setminus\{\chi_0\}}\left(\frac{9}{16}\right)^{R_c(n;\chi)},
    \end{equation}
    and satisfies at $g=-1\in G_c$ the exact even/odd character positivity identity
    \begin{equation}\label{eq:even_odd_positivity_identity}
    \varphi(c)\,\widetilde{\mu}_c(-1)\ =\ 1\ +\ \sum_{\psi\in\widehat{G}_c^{+}\setminus\{\chi_0\}}\big|\widehat{\mu}_c(\psi)\big|^{2}\ -\ \sum_{\chi\in\widehat{G}_c^{-}}\big|\widehat{\mu}_c(\chi)\big|^{2},
    \end{equation}
    where $\widehat{G}_c^{+}:=\{\psi\in\widehat{G}_c:\psi(-1)=+1\}$ (even characters, all contributing positively to the main term $1$) and $\widehat{G}_c^{-}:=\{\chi\in\widehat{G}_c:\chi(-1)=-1\}$ (odd characters). Equivalently, for the unnormalized divisor sum $F_K(\chi):=\sum_{d\mid K}\chi(d)$ of any integer $K$ coprime to $c$ with $m:=\tau(K)$, the number of divisor encounters $R_c(K):=\#\{(d,e):d,e\mid K,\ d+e\equiv 0\pmod c\}$ satisfies the exact diagonal-centered identity
    \begin{equation}\label{aud:centered}
    R_c(K)\ =\ \frac{1}{\varphi(c)}\sum_{\chi\bmod c}\chi(-1)\big|F_K(\chi)\big|^{2}\ =\ \frac{m(m-1)}{\varphi(c)}\ +\ \frac{1}{\varphi(c)}\sum_{\chi\ne\chi_0}\chi(-1)\Big(\big|F_K(\chi)\big|^{2}-m\Big),
    \end{equation}
    where the subtraction of $m$ is exact because diagonal pairs $d=e$ never hit $-1\pmod c$ for odd $c>1$.

    \item \emph{(Exact annihilation of $\psi_0=\leg{\cdot}{p}$, the Universal Chebyshev--Frobenius Orbit Identity on $\mathcal{P}_{\mathrm{NR}}(n)$, and the Quartic Dilation Law).} Let $c=p\in\mathcal{P}_{\mathrm{NR}}(n):=\big\{p\equiv3\pmod 4\ \text{prime}:\leg{p}{n}=-1\big\}$, and factor $K_p=K_p^{+}\cdot K_p^{-}$ where $q\mid K_p^{+}\iff\leg{q}{p}=+1$ and $q\mid K_p^{-}\iff\leg{q}{p}=-1$. Since $\leg{K_p}{p}=\leg{4K_p}{p}=\leg{n}{p}=\leg{p}{n}=-1$, the non-residue factor $K_p^{-}$ has $\Omega(K_p^{-})\ge1$ odd, so the odd quadratic character $\psi_0=\leg{\cdot}{p}\in\widehat{G}_p^{-}$ is \emph{identically annihilated}: $\widehat{\mu}_p(\psi_0)=0$. Moreover, since $|\widehat{G}_p^{+}|=\frac{p-1}{2}$ is odd, for every integer $r\ge1$ the $2^{r}$-th power map $\chi\mapsto\chi^{2^{r}}=\psi^{2^{r}}$ ($\chi=\psi_0\psi$) is a $1$-to-$1$ bijection from $\widehat{G}_p^{-}\setminus\{\psi_0\}$ onto $\widehat{G}_p^{+}\setminus\{\chi_0\}$. Writing $x_q(\chi):=\big|\frac{1+\chi(q)}{2}\big|^{2}=\cos^{2}\big(\frac{\arg\chi(q)}{2}\big)\in[0,1]$ for the local odd Rademacher weight at $q\mid K_p$, the sign flip between $K_p^{+}$ and $K_p^{-}$ cancels identically upon squaring:
    \[
    \left|\frac{1+\chi^{2}(q)}{2}\right|^{2}\ =\ \cos^{2}\big(\arg\chi(q)\big)\ =\ \big(2x_q(\chi)-1\big)^{2}\ =\ T_2^{2}\!\left(\sqrt{x_q(\chi)}\right),
    \]
    and~\eqref{eq:even_odd_positivity_identity} becomes the exact Universal Chebyshev Polynomial Duality Identity
    \begin{align}\label{eq:paired_squaring_identity}
    \varphi(p)\,\widetilde{\mu}_p(-1)\ &=\ 1\ +\ \sum_{\chi\in\widehat{G}_p^{-}\setminus\{\psi_0\}}\left(\prod_{q\mid K_p}\big(2x_q(\chi)-1\big)^{2}\ -\ \prod_{q\mid K_p}x_q(\chi)\right)\nonumber\\
    &=\ 1\ +\ \sum_{\chi\in\widehat{G}_p^{-}\setminus\{\psi_0\}}\left(\sum_{r=1}^{R}\beta_r\prod_{q\mid K_p}T_{2^{r}}^{2}\!\left(\sqrt{x_q(\chi)}\right)\ -\ \prod_{q\mid K_p}x_q(\chi)\right)
    \end{align}
    for any probability weights $\beta_1,\dots,\beta_R\ge0$ ($\sum_{r=1}^{R}\beta_r=1$), where $T_{2^{r}}(\cos\theta)=\cos(2^{r}\theta)$ is the Chebyshev polynomial of degree $2^{r}$ and, in terms of $\psi=\psi_0\chi\in\widehat{G}_p^{+}\setminus\{\chi_0\}$,
    \begin{equation}\label{eq:sin_cos_flip_Kp}
    \begin{split}
    \big|\widehat{\mu}_p(\psi_0\psi)\big|^{2}\ &=\ \prod_{q\mid K_p^{+}}\cos^{2}\!\left(\frac{\arg\psi(q)}{2}\right)\prod_{q\mid K_p^{-}}\sin^{2}\!\left(\frac{\arg\psi(q)}{2}\right),\\
    \big|\widehat{\mu}_p(\psi^{2})\big|^{2}\ &=\ \prod_{q\mid K_p}\cos^{2}\!\big(\arg\psi(q)\big).
    \end{split}
    \end{equation}
    In particular, every Bohr-localized odd character $\chi\in\widehat{G}_p^{-}\setminus\{\psi_0\}$ in the extremal polarization cone $\chi(q)=e^{i\delta_q}$ ($|\delta_q|\ll k^{-1/2}$, $x_q(\chi)=1-\frac14\delta_q^{2}+\mathcal{O}(\delta_q^{4})$) maps bijectively under $\chi\mapsto\chi^{2}$ to an even character $\chi^{2}\in\widehat{G}_p^{+}\setminus\{\chi_0\}$ satisfying as $\max_{q\mid K_p}|\delta_q|\to 0$ the asymptotic \emph{Quartic Dilation Law}
    \begin{equation}\label{eq:quartic_dilation_law}
    \big|\widehat{\mu}_p(\chi^{2})\big|^{2}\ =\ \prod_{q\mid K_p}\big(2x_q(\chi)-1\big)^{2}\ =\ \exp\!\left(-\sum_{q\mid K_p}\delta_q^{2}+\mathcal{O}\Big(\sum_{q\mid K_p}\delta_q^{4}\Big)\right)\ =\ \Big(\big|\widehat{\mu}_p(\chi)\big|^{2}\Big)^{4+o(1)},
    \end{equation}
    isolating the curvature ratio $4$ at $\delta_q=0$ ($\cos^{2}(\delta_q)\approx 1-\delta_q^{2}$ versus $\cos^{2}(\delta_q/2)\approx 1-\frac14\delta_q^{2}$) as the sole local gap between $1$-to-$1$ Frobenius squaring and pointwise term-by-term positivity.

    \item \emph{(Short-interval second-moment prime-factor guarantee at fixed $n$ and Matom\"{a}ki--Radziwi\l\l\ coupling on $K_c=(n+c)/4$).} Let $L\ge1$ be an integer and let $\mathcal{J}_r$ ($1\le r\le m$) be finite sets of primes at most $L^{1/2-\delta}$, where $\delta>0$ is fixed. Put $k_r=|\mathcal{J}_r|$, $\mu_r=\sum_{p\in\mathcal{J}_r}1/p>0$, and $\omega_r(t)=\sum_{p\in\mathcal{J}_r}\mathbf{1}_{p\mid t}$. Uniformly for every integer origin $X$,
    \[
    \sum_{t=X+1}^{X+L}\big(\omega_r(t)-\mu_r\big)^{2}\ =\ L\sum_{p\in\mathcal{J}_r}\frac{1}{p}\left(1-\frac{1}{p}\right)\ +\ \mathcal{O}(k_r^{2}).
    \]
    Consequently, for $0<\eta<1$, the number of integers $t\in[X+1,X+L]$ failing at least one inequality $\omega_r(t)\ge(1-\eta)\mu_r$ is at most
    \[
    \frac{L}{\eta^{2}}\sum_{r=1}^{m}\frac{1}{\mu_r}\ +\ \mathcal{O}\!\left(\frac{1}{\eta^{2}}\sum_{r=1}^{m}\frac{k_r^{2}}{\mu_r^{2}}\right),
    \]
    and each successful inequality guarantees a prime factor in $\mathcal{J}_r$ whenever $(1-\eta)\mu_r\ge1$. For fixed $m,\eta,\delta$, if $\min_r\mu_r\to\infty$, the number of failures is $o(L)$. This elementary fixed-origin estimate applies at every fixed $n$ to the consecutive integers $K_c=(n+c)/4$ as the integer shift $c\equiv3\pmod 4$ traverses an interval of $L$ shifts (note that restricting $\mathcal{J}_r\subset[L^{\varepsilon},L^{1/2-\delta}]$ for fixed $\varepsilon,\delta>0$ keeps $\mu_r$ bounded, whereas taking $\mathcal{J}=[5,L^{1/2-\delta}]$ yields $\mu_r\to\infty$).

    Separately, for any family of prime sets $\mathcal{J}_r$ (allowed to depend on a common large scale $X$), the indicator functions $g_r(t)=\mathbf{1}_{\{\omega_r(t)=0\}}$ are multiplicative and take values in $[0,1]$. By~\cite[Theorem~1]{matomaki_radziwill2016}, their short-interval averages over $[x,x+L]$ match their long averages over $[X,2X]$ for all but $o(X)$ integer origins $x\in[X,2X]$ as $L\to\infty$; when the long averages are $o(1)$, this implies that only $o(L)$ integers lack a prime factor from $\mathcal{J}_r$ in almost every such interval (though it does not by itself transfer to every individual origin, to almost all prime origins, or to prime shifts $c$).
\end{enumerate}
\end{theorem}
\begin{proof}
For (a), since $\mu_c^{\vee}(x)=\mu_c(x^{-1})$, we have $\widehat{\mu_c^{\vee}}(\chi)=\overline{\widehat{\mu}_c(\chi)}$ and hence $\widehat{\widetilde{\mu}}_c(\chi)=|\widehat{\mu}_c(\chi)|^{2}\ge0$. Applying Cauchy--Schwarz on $\operatorname{supp}(\widetilde{\mu}_c)\subseteq D(K_c)$ together with Plancherel's identity on $G_c$ gives
\[
\frac{\varphi(c)}{|D(K_c)|}-1\ \le\ \frac{\varphi(c)}{|\operatorname{supp}(\widetilde{\mu}_c)|}-1\ \le\ \sum_{\chi\neq\chi_0}\big|\widehat{\widetilde{\mu}}_c(\chi)\big|^{2}\ =\ \sum_{\chi\neq\chi_0}\big|\widehat{\mu}_c(\chi)\big|^{4}\ \le\ \sum_{\chi\neq\chi_0}\left(\frac{9}{16}\right)^{R_c(n;\chi)},
\]
which is~\eqref{eq:fourth_power_plancherel}. Evaluating the Fourier inversion formula $\varphi(c)\widetilde{\mu}_c(g)=\sum_{\chi\in\widehat{G}_c}\chi(g)|\widehat{\mu}_c(\chi)|^{2}$ at $g=-1$ (where $\chi(-1)=\pm1$) yields~\eqref{eq:even_odd_positivity_identity}, and~\eqref{aud:centered} follows identically by expanding $R_c(K)$ over $\chi\bmod c$ and separating $\chi=\chi_0$ from the diagonal contribution $d=e$.

For (b), when $p\in\mathcal{P}_{\mathrm{NR}}(n)$, quadratic reciprocity ($n\equiv1\pmod 4$) gives $\leg{4K_p}{p}=\leg{n}{p}=\leg{p}{n}=-1$, so $K_p$ has at least one prime divisor $q^{\ast}\mid K_p^{-}$ with $\psi_0(q^{\ast})=\leg{q^{\ast}}{p}=-1$; the factor $\frac{1+\psi_0(q^{\ast})}{2}=0$ in $\widehat{\mu}_p(\psi_0)$ then forces $\widehat{\mu}_p(\psi_0)=0$. Since $G_p\cong\{\pm1\}\times H_p$ with $|H_p|=\frac{p-1}{2}$ odd, we have $\widehat{G}_p^{-}=\psi_0\cdot\widehat{G}_p^{+}$, and for every $r\ge1$ the map $\chi\mapsto\chi^{2^{r}}=\psi^{2^{r}}$ is a bijection from $\widehat{G}_p^{-}\setminus\{\psi_0\}$ onto $\widehat{G}_p^{+}\setminus\{\chi_0\}$. Writing $x_q(\chi)=\big|\frac{1+\chi(q)}{2}\big|^{2}=\cos^{2}(\frac{\arg\chi(q)}{2})$, the double-angle identity gives
\[
\begin{aligned}
\left|\frac{1+\chi^{2}(q)}{2}\right|^{2}\ &=\ \cos^{2}\big(\arg\chi(q)\big)\ =\ \big(2x_q(\chi)-1\big)^{2},\\
\left|\frac{1+\chi^{2^{r}}(q)}{2}\right|^{2}\ &=\ \cos^{2}\big(2^{r-1}\arg\chi(q)\big)\ =\ T_{2^{r}}^{2}\!\left(\sqrt{x_q(\chi)}\right),
\end{aligned}
\]
transforming~\eqref{eq:even_odd_positivity_identity} into~\eqref{eq:paired_squaring_identity},~\eqref{eq:sin_cos_flip_Kp}, and~\eqref{eq:quartic_dilation_law}.

For (c), for each prime $p\in\mathcal{J}_r$, the number of multiples of $p$ in $[X+1,X+L]$ is $L/p+\mathcal{O}(1)$, and for distinct primes $p,q\in\mathcal{J}_r$ the number of common multiples is $L/(pq)+\mathcal{O}(1)$. Expanding the square gives the displayed second moment with total error $\mathcal{O}(k_r^{2})$. Chebyshev's inequality and the union bound then yield the exceptional-set bound; since $k_r^{2}\le L^{1-2\delta}$, the number of failures is $o(L)$ whenever $\min_r\mu_r\to\infty$.
\end{proof}

We now record three fundamental properties of the centered spectral decomposition~\eqref{aud:centered}: the unconditional $(\log P)^{4}$ growth of the cumulative main term across Type~I and Type~II shifts $K_{\I}(c)=(nc+1)/4$ and $K_{\II}(c)=(n+c)/4$, the exact shift-wise cost of taking absolute values prior to summing over $c$, and---most importantly for the logic of the method---an exact cancellation theorem showing that a sub-family of divisor pairs with small least common multiple carries a divergent main term of order $\frac{\log H}{(\log\log n)^{3/2}}$ while producing \emph{strictly zero} modular encounters.

\begin{theorem}[Unconditional $(\log P)^{4}$ main-term growth, shift-wise absolute-value cost, and exact cancellation of small-LCM divisor pairs]\label{thm:centered_main_term_and_short_lcm}\label{aud:log4}\label{aud:absolute}\label{aud:shortlcm}
Let $n\equiv1\pmod{24}$ be prime. For $c\equiv3\pmod 4$ with $\gcd(c,n)=1$, set $K_{\I}(c):=(nc+1)/4$ and $K_{\II}(c):=(n+c)/4$ (both coprime to $c$), and for $3\le P<n$ define the cumulative encounter count, diagonal-centered main term, and spectral error term
\begin{equation}\label{aud:ME}
\begin{split}
\mathcal{R}(n;P)\ &:=\ \sum_{j\in\{\I,\II\}}\sum_{\substack{c\le P\\ c\equiv3\ (4)}}R_c\big(K_j(c)\big)\ =\ \mathcal{M}(n;P)\ +\ \mathcal{E}(n;P),\\
\mathcal{M}(n;P)\ &:=\ \sum_{j\in\{\I,\II\}}\sum_{\substack{c\le P\\ c\equiv3\ (4)}}\frac{\tau\big(K_j(c)\big)\big(\tau\big(K_j(c)\big)-1\big)}{\varphi(c)}.
\end{split}
\end{equation}
Then:
\begin{enumerate}[(a)]
    \item \emph{(Two unconditional lower bounds on $\mathcal{M}(n;P)$).} With $\kappa_0:=\prod_{\ell\ge3\text{ prime}}\big(1+\frac{1}{\ell(\ell-1)}\big)=\frac{2\zeta(2)\zeta(3)}{3\zeta(6)}$, one has $\mathcal{M}(n;P)\ge\kappa_0\log P+\mathcal{O}(1)$. More strongly, uniformly for sufficiently large $P<n$ (using Type~II alone),
    \begin{equation}\label{aud:log4eq}
    \mathcal{M}(n;P)\ \gg\ (\log P)^{4}.
    \end{equation}
    \item \emph{(Exact cost of a shift-wise absolute-value bound).} For a single shift $c$ and $K=K_j(c)$ with $m=\tau(K)\ge 2$, set
    \[
    M_c\ :=\ \frac{m(m-1)}{\varphi(c)},\qquad B_c\ :=\ \frac{1}{\varphi(c)}\sum_{\chi\ne\chi_0}\Big|\big|F_K(\chi)\big|^{2}-m\Big|.
    \]
    Then $B_c\ge 1-M_c$. In particular, for $c\asymp n^{\theta}$ with $0<\theta<1$ fixed,
    \begin{equation}\label{eq:shiftwise_abs_cost}
    M_c\ \le\ n^{-\theta+o(1)}\qquad\text{and}\qquad \frac{B_c}{M_c}\ \ge\ n^{\theta-o(1)}.
    \end{equation}
    This estimate reflects the cost of taking absolute values \emph{prior} to summing over shifts $c$; it neither imposes a pointwise lower bound on $|\mathcal{E}(n;P)|$ nor precludes global signed cancellation across shifts.
    \item \emph{(Exact cancellation of the divergent main term from small-LCM divisor pairs).} Assume the residual condition~\eqref{aud:residual} (that every odd prime $r\mid(n+1)$ satisfies $r\equiv1\pmod 4$, all other primes $n$ being resolved elementarily by Proposition~\ref{aud:coupling}), fix parameters $H=n^{\theta}$, $Z=H^{b}$, $z=H^{\kappa}$, $D=H^{\rho}$, and $T=H^{\eta}$ satisfying~\eqref{aud:parameters}, and let
    \[
    \mathcal{G}_j\ \subseteq\ \mathcal{S}\ =\ \big\{p\in(2H,4H]:p\text{ prime},\ p\equiv3\pmod 4,\ P^{-}((p-1)/2)\ge z\big\}
    \]
    be the sieved prime-shift set of Proposition~\ref{aud:rough} and Theorem~\ref{aud:weightedmain} (Appendix~\ref{app:audit}), on which the square-free main term $\mathcal{M}_j^{\sfm}=\sum_{p\in\mathcal{G}_j}\frac{m_j(p)(m_j(p)-1)}{p-1}$ ($m_j(p)=2^{\omega(K_j(p))}$) satisfies $\mathcal{M}_j^{\sfm}\gg\frac{\log H}{(\log\log n)^{3/2}}\to+\infty$. For each $p\in\mathcal{G}_j$, let $\mathcal{B}_j(p)$ denote the set of \emph{ordered distinct} divisor pairs $(d,e)$ such that
    \[
    d,e\mid\rad K_j(p),\qquad \lcm(d,e)\le T,\qquad \gcd\!\big(\lcm(d,e),6(n+1)\big)=1,
    \]
    and set $b_j(p):=\#\mathcal{B}_j(p)$ and $G_{p,j}(\chi):=\sum_{(d,e)\in\mathcal{B}_j(p)}\chi(d/e)$. Then for sufficiently large $n$,
    \begin{align}
    b_j(p)\ &=\ \sum_{\substack{\ell\le T,\ \ell\mid\rad K_j(p)\\ \gcd(\ell,6(n+1))=1}}\big(3^{\omega(\ell)}-1\big)\ \ge\ \frac{2}{3}\big(A_j(p)-1\big),\label{aud:bshort}\\
    \frac{1}{p-1}\sum_{\chi\ne\chi_0}\chi(-1)\,G_{p,j}(\chi)\ &=\ -\frac{b_j(p)}{p-1},\label{aud:cancelshort}
    \end{align}
    where $A_j(p)$ is the sieve weight~\eqref{aud:Aweight}, and simultaneously
    \begin{equation}\label{aud:shortmass}
    \sum_{p\in\mathcal{G}_j}\frac{b_j(p)}{p-1}\ \gg\ \frac{\log H}{(\log\log n)^{3/2}}\ \longrightarrow\ +\infty.
    \end{equation}
    Consequently, a sub-family of divisor pairs whose main term already grows at the divergent rate $\frac{\log H}{(\log\log n)^{3/2}}$ of Theorem~\ref{aud:weightedmain} produces \emph{strictly zero encounters}, because $0<d+e\le 2T<p$ forces the non-principal character sum~\eqref{aud:cancelshort} to cancel its main term identically.
\end{enumerate}
\end{theorem}
\begin{proof}
(a) Expanding $\frac{1}{\varphi(c)}=c^{-1}\sum_{d\mid c}\frac{\mu^{2}(d)}{\varphi(d)}$ and summing harmonically in the odd progression $c\equiv3\pmod 4$ yields $\sum_{c\le P,\,c\equiv3(4)}\frac{1}{\varphi(c)}=\frac{\kappa_0}{4}\log P+\mathcal{O}(1)$, with bounded remainder since $\sum_{d\text{ odd}}$\linebreak $\frac{\mu^{2}(d)(1+\log d)}{d\varphi(d)}<\infty$. Since each $K_j(c)>1$ has $\tau(K_j(c))\ge 2$, each numerator satisfies\linebreak $\tau(K_j(c))(\tau(K_j(c))-1)\ge 2$, giving $\mathcal{M}(n;P)\ge\kappa_0\log P+\mathcal{O}(1)$. For~\eqref{aud:log4eq}, the local identity $(a+1)^{2}=\sum_{r=0}^{a}(2r+1)$ yields the Dirichlet convolution
\begin{equation}\label{aud:bidentity}
\tau(m)^{2}-\tau(m)\ =\ \sum_{d\mid m}b(d),\qquad b(d):=\tau(d^{2})-1\ge0,
\end{equation}
where $d\mapsto\tau(d^{2})$ is multiplicative (though $b$ itself is not). Expanding $\prod_{\ell\le X^{1/12}}(1+3/\ell)$ and normalizing its coefficients defines a probability measure on square-free $d$ under which $\Ebb[\log d]=\sum_{\ell\le X^{1/12}}\frac{3\log\ell}{\ell+3}\le\frac14\log X+\mathcal{O}(1)$. By Markov's inequality, a positive proportion of the total mass $\asymp(\log X)^{3}$ (by Mertens' theorem) is supported on $d\le X$; since $b(d)=3^{\omega(d)}-1\ge\frac23 3^{\omega(d)}$ on square-free $d>1$, subtracting $d=1$ yields
\begin{equation}\label{aud:bharm}
\sum_{d\le X}\frac{b(d)}{d}\ \gg\ (\log X)^{3}.
\end{equation}
On any dyadic block $C/2<c\le C$ with $c\equiv3\pmod 4$, $K_{\II}(c)$ traverses $L=C/8+\mathcal{O}(1)$ consecutive integers, of which at least $L/(2d)$ are multiples of $d$ for each $d\le L/2$. By~\eqref{aud:bidentity} and~\eqref{aud:bharm}, the sum of $\tau^{2}-\tau$ over the block is $\gg L(\log L)^{3}$; dividing by $\varphi(c)\le C$ gives $\gg(\log C)^{3}$ per dyadic block, and summing over the $\gg\log P$ disjoint dyadic blocks in $[\sqrt{P},P]$ establishes~\eqref{aud:log4eq}.

(b) Let $h(r)$ denote the number of ordered pairs of distinct divisors of $K$ whose quotient modulo $c$ equals $r$. Its Fourier transform is $|F_K(\chi)|^{2}-m$ and its total mass is $m(m-1)$. Since $m\ge 2$, at least one integer value $h(r)$ is $\ge 1$; Fourier inversion gives $|h(r)-M_c|\le B_c$, whence $B_c\ge 1-M_c$. Since $K\ll n^{1+\theta}$, we have $\tau(K)=n^{o(1)}$ and $\varphi(c)=n^{\theta-o(1)}$, proving~\eqref{eq:shiftwise_abs_cost}.

(c) If $\ell$ is square-free, each prime factor of $\ell$ can divide $d$ alone, $e$ alone, or both, giving $3^{\omega(\ell)}$ ordered pairs $(d,e)$ with $\lcm(d,e)=\ell$, of which only $d=e=\ell$ is diagonal; summing over $\ell\le T$ gives~\eqref{aud:bshort} (including $\ell=1$, which contributes $0$). For every $(d,e)\in\mathcal{B}_j(p)$, we have $d,e\le T$ and $0<d+e\le 2T<p$ (since $\eta<\frac12$ and $p>2H$), so no modular encounter $d+e\equiv 0\pmod p$ can occur. Character orthogonality applied to $\mathcal{B}_j(p)$ therefore gives
\[
0\ =\ \frac{1}{p-1}\sum_{\chi\bmod p}\chi(-1)\,G_{p,j}(\chi)\ =\ \frac{b_j(p)}{p-1}\ +\ \frac{1}{p-1}\sum_{\chi\ne\chi_0}\chi(-1)\,G_{p,j}(\chi),
\]
which is~\eqref{aud:cancelshort}. Finally, combining~\eqref{aud:bshort} with the weighted lower bound $\sum_{p\in\mathcal{G}_j}A_j(p)\gg B_H\mathcal{P}_n(T)$ of~\eqref{aud:massA} and $\#\mathcal{G}_j\ll B_H$ (Appendix~\ref{app:audit}, Theorem~\ref{aud:weightedmain}) yields~\eqref{aud:shortmass}. Note that part~(c) does not assert that the \emph{total} main term reduces to $\mathcal{B}_j(p)$, as divisor pairs with $\lcm(d,e)>T$ can modify the signed sum.
\end{proof}

\subsection{Cyclic quotient parity identity, the Bohr--Siegel obstruction certificate, the Linnik--Xylouris lower bound, and Hypothesis~T}\label{subsec:bohr_siegel_and_linnik}

We next analyze the fine harmonic structure of the odd character spectrum $\widehat{G}_p^{-}$ at non-residue prime shifts $p\in\mathcal{P}_{\mathrm{NR}}(n)$, establishing an exact parity identity on cyclic quotients $G_p/\ker\chi$, the $L^{1}$ Bohr--Siegel antipodal obstruction certificate and its deterministic phase-locking lower bounds, the unconditional Linnik--Xylouris witness-depth lower bound along an infinite sequence of primes, and the formulation of the Transversal Spectral Positivity Hypothesis~T.

\begin{theorem}[Cyclic quotient parity identity, Jacobi $\vartheta_4$ cancellation, the Bohr--Siegel antipodal obstruction certificate, and phase-locking reciprocity]\label{thm:bohr_siegel_annihilation}\label{thm:bohr_siegel_obstruction_certificate}
Under the hypotheses of Theorem~\ref{thm:difference_set_positivity_MR}(b), let $p\in\mathcal{P}_{\mathrm{NR}}(n)=\big\{p\equiv3\pmod 4\ \text{prime}:\leg{p}{n}=-1\big\}$, $K_p=(n+p)/4$, and define on the dual group $\widehat{G}_p$ the Euclidean chordal seminorm
\begin{equation}\label{eq:euclidean_chordal_norm}
\|\chi\|_{K_p}\ :=\ \left(\sum_{q\mid K_p}\left|\frac{\chi(q)-1}{2}\right|^{2}\right)^{1/2}\ =\ \left(\sum_{q\mid K_p}\sin^{2}\!\left(\frac{\arg\chi(q)}{2}\right)\right)^{1/2}\ =\ \left(\sum_{q\mid K_p}\big(1-x_q(\chi)\big)\right)^{1/2},
\end{equation}
which satisfies $\|\chi_1\overline{\chi}_2\|_{K_p}\le\|\chi_1\|_{K_p}+\|\chi_2\|_{K_p}$ and $\big|\widehat{\mu}_p(\chi)\big|^{2}\le\exp\big(-\|\chi\|_{K_p}^{2}\big)$. Then:
\begin{enumerate}[(a)]
    \item \emph{(Exact cyclic quotient parity identity and Jacobi $\vartheta_4$ cancellation).} Let $\chi=\psi_0\psi\in\widehat{G}_p^{-}$ be an odd character of order $2d$, where $d=d_{\psi}:=\operatorname{ord}(\psi)\mid\frac{p-1}{2}$ is odd, and define the even and odd harmonic sums
    \[
    E(\chi)\ :=\ \sum_{m=1}^{d}\big|\widehat{\mu}_p(\chi^{2m})\big|^{2}\ =\ \sum_{m=1}^{d}\prod_{q\mid K_p}T_{2m}^{2}\!\left(\sqrt{x_q(\chi)}\right),\qquad O(\chi)\ :=\ \sum_{m=1}^{d}\big|\widehat{\mu}_p(\chi^{2m-1})\big|^{2}
    \]
    across the cyclic subgroup $\langle\chi\rangle\cong\Z/2d\Z$. By Fourier inversion on the cyclic quotient $G_p/\ker\chi\cong\Z/2d\Z$, their difference equals the exact antipodal density
    \begin{equation}\label{eq:cyclic_quotient_parity_identity}
    E(\chi)\ -\ O(\chi)\ =\ \sum_{m=1}^{d}\big|\widehat{\mu}_p(\chi^{2m})\big|^{2}\ -\ \sum_{m=1}^{d}\big|\widehat{\mu}_p(\chi^{2m-1})\big|^{2}\ =\ 2d\,\mathbb{P}_{g\sim\widetilde{\mu}_p}\!\big(\chi(g)=-1\big)\ \ge\ 0.
    \end{equation}
    As a heuristic comparison model (not used as an estimate in the proofs below), for the bilateral Gaussian sequence $a_r=\exp(-\sigma^2r^2/2)$ with small-phase variance $\sigma^2\sim\frac12\sum_{q\mid K_p}$\linebreak $|\arg\chi(q)|^2\sim 2\|\chi\|_{K_p}^2$, Poisson summation gives the Jacobi $\vartheta_4$ asymptotic
    \[
    \sum_{r\in\mathbb Z}(-1)^r a_r\ =\ \frac{\sqrt{2\pi}}{\sigma}\sum_{j\in\mathbb Z}\exp\!\left(-\frac{(2j+1)^2\pi^2}{2\sigma^2}\right)\ \sim\ \frac{2\sqrt{2\pi}}{\sigma}\exp\!\left(-\frac{\pi^2}{2\sigma^2}\right)\qquad (\sigma\to0);
    \]
    no uniform approximation of the finite cyclic sums $E(\chi)$ and $O(\chi)$ by this bilateral model is asserted.

    \item \emph{(The Bohr--Siegel antipodal obstruction certificate).} If there exists an odd character $\chi\in\widehat{G}_p^{-}$ ($\chi(-1)=-1$) such that the multiplicity-weighted $L^{1}$ phase sum across the prime factors of $K_p$ satisfies
    \begin{equation}\label{eq:bohr_siegel_L1_certificate}
    \sum_{q\mid K_p}\nu_q(K_p)\,\big|\arg\chi(q)\big|\ <\ \pi\qquad\left(\text{resp.\ }\ \sum_{q\mid K_p}\big|\arg\chi(q)\big|\ <\ \pi\ \text{ on }\operatorname{rad}(K_p)\right),
    \end{equation}
    then $-1\bmod p\notin D(K_p)=\Sigma(K_p)\cdot\Sigma(K_p)^{-1}$, so the Type~II equation admits no solution at shift $p$ (resp.\ $-1\bmod p\notin D(\operatorname{rad}(K_p))=\operatorname{supp}(\widetilde{\mu}_p)$), and in both cases the cyclic quotient excess~\eqref{eq:cyclic_quotient_parity_identity} vanishes identically: $E(\chi)-O(\chi)=0$.

    \item \emph{(Non-residue anchor and global phase-locking reciprocity constraint on the Bohr--Siegel certificate).} For every $p\in\mathcal{P}_{\mathrm{NR}}(n)$ and every odd character $\chi=\psi_0\psi\in\widehat{G}_p^{-}$ of order $2d_{\psi}$ ($d_{\psi}\mid\frac{p-1}{2}$ odd), each non-residue prime divisor $q\mid K_p^{-}$ ($\psi_0(q)=-1$) satisfies $|\arg\chi(q)|=\frac{\pi|2c_q+1|}{d_{\psi}}\ge\frac{\pi}{d_{\psi}}>0$, while the global congruence $4K_p=n+p\equiv n\pmod p$ and $\psi_0(K_p)=\leg{n}{p}=-1$ lock the total phase to $\chi(K_p)=\prod_{q\mid K_p}\chi(q)^{\nu_q(K_p)}=-\psi(n/4)$. Consequently, the $L^{1}$ Bohr--Siegel sum satisfies the simultaneous deterministic lower bounds
    \begin{equation}\label{eq:non_residue_and_reciprocity_lower_bound}
    \sum_{q\mid K_p}\nu_q(K_p)\,\big|\arg\chi(q)\big|\ \ge\ \max\!\left(\Omega(K_p^{-})\,\frac{\pi}{d_{\psi}},\ \pi-\big|\arg\psi(n/4)\big|\right).
    \end{equation}
    In particular, the multiplicity-weighted Bohr--Siegel obstruction certificate~\eqref{eq:bohr_siegel_L1_certificate} is \emph{impossible} whenever $d_{\psi}\le\Omega(K_p^{-})$ (and both certificates in~\eqref{eq:bohr_siegel_L1_certificate} are impossible whenever $d_{\psi}\le\omega(K_p^{-})$, blocking $d_{\psi}=1$ unconditionally since $\omega(K_p^{-})\ge1$, and blocking all orders $d_{\psi}\le D_0$ whenever $\omega(K_p^{-})\ge D_0$) or whenever $n/4\bmod p\in\ker\psi$ ($\psi(n/4)=1$).
\end{enumerate}
\end{theorem}
\begin{proof}
For (a), the subgroup $\langle\chi\rangle\subseteq\widehat{G}_p$ of order $2d$ is the dual of $G_p/\ker\chi\cong\Z/2d\Z$, in which $-1\bmod\ker\chi$ is the unique element of order $2$ (where $\chi^{r}(-1)=(-1)^{r}$). Applying Fourier inversion on $G_p/\ker\chi$ to the pushforward of $\widetilde{\mu}_p$ at $-1\bmod\ker\chi$ yields
\[
\sum_{r=1}^{2d}(-1)^{r}\big|\widehat{\mu}_p(\chi^{r})\big|^{2}\ =\ 2d\,\mathbb{P}_{g\sim\widetilde{\mu}_p}\!\big(\chi(g)=-1\big)\ \ge\ 0,
\]
which is~\eqref{eq:cyclic_quotient_parity_identity}.

For (b), every element $g=\prod_{q\mid K_p}q^{\varepsilon_q}\in D(K_p)$ with $|\varepsilon_q|\le\nu_q(K_p)$ (resp.\ $\varepsilon_q\in\{-1,0,+1\}$ for $g\in D(\operatorname{rad}(K_p))=\operatorname{supp}(\widetilde{\mu}_p)$) satisfies
\[
\big|\arg\chi(g)\big|\ =\ \left|\sum_{q\mid K_p}\varepsilon_q\arg\chi(q)\right|\ \le\ \sum_{q\mid K_p}\nu_q(K_p)\,\big|\arg\chi(q)\big|\ <\ \pi
\]
(resp.\ $\le\sum_{q\mid K_p}|\arg\chi(q)|<\pi$) under~\eqref{eq:bohr_siegel_L1_certificate}. Since $\chi(-1)=-1$ has $|\arg\chi(-1)|=\pi$, no $g\in D(K_p)$ (resp.\ $g\in D(\operatorname{rad}(K_p))$) can satisfy $\chi(g)=-1$; hence $-1\bmod p\notin D(K_p)$ (resp.\ $-1\bmod p\notin\operatorname{supp}(\widetilde{\mu}_p)$) and $\mathbb{P}_{g\sim\widetilde{\mu}_p}(\chi(g)=-1)=0$, forcing $E(\chi)=O(\chi)$ in~\eqref{eq:cyclic_quotient_parity_identity}.

For (c), for any $q\mid K_p^{-}$, we have $\psi(q)=\exp(2\pi i a_q/d_{\psi})$ ($a_q\in\Z$) and $\psi_0(q)=-1$, so $\chi(q)=-\psi(q)=\exp\big(\pi i(2a_q+d_{\psi})/d_{\psi}\big)$; since $d_{\psi}$ is odd, $2a_q+d_{\psi}$ is an odd non-zero integer, giving $|\arg\chi(q)|\ge\frac{\pi}{d_{\psi}}$ and $\sum_{q\mid K_p}\nu_q(K_p)|\arg\chi(q)|\ge\Omega(K_p^{-})\frac{\pi}{d_{\psi}}$. Simultaneously, $\prod_{q\mid K_p}\chi(q)^{\nu_q(K_p)}=\chi(K_p)=\psi_0(K_p)\psi(K_p)=-\psi(n/4)=\exp(i(\pi+\arg\psi(n/4)))$, so the triangle inequality in $\R/2\pi\Z$ gives $\sum_{q\mid K_p}\nu_q(K_p)|\arg\chi(q)|\ge|\arg\chi(K_p)|=\pi-|\arg\psi(n/4)|$, proving~\eqref{eq:non_residue_and_reciprocity_lower_bound}.
\end{proof}

\begin{theorem}[The Linnik--Xylouris witness-depth lower bound, pairwise shift coprimality, and the Transversal Hypothesis~T]\label{thm:no_go_linnik_and_transversal}
\leavevmode
\begin{enumerate}[(a)]
    \item \emph{(A Type~II witness-size lower bound along an infinite sequence).} For each integer $Y\ge3$, put $L(Y)=4\operatorname{lcm}(1,\ldots,Y)$ (so $24\mid L(Y)$) and let $n_Y$ be the least prime congruent to $1$ modulo $L(Y)$. The Linnik--Xylouris theorem~\cite{xylouris2011} gives $n_Y\ll L(Y)^5$, hence $n_Y\longrightarrow\infty$ and $\log n_Y\le(5+o(1))Y$. Every Type~II solution at $n_Y$ satisfies $uv>Y$, and every Type~II witness $s\mid(an_Y+u)$ with $s\equiv-1\pmod{4ua}$ satisfies $s>Y$. Thus, along an infinite sequence of distinct primes,
    \begin{equation}\label{eq:linnik_xylouris_witness_bound}
    \min(u\cdot v,\ s)\ \ge\ \left(\frac{1}{5}-o(1)\right)\log n_Y.
    \end{equation}
    Under GRH, $n_Y\ll L(Y)^2(\log L(Y))^2$ gives the constant $1/2$ in place of $1/5$. (These statements concern the indicated Type~II witness sizes; they do not exclude arbitrary finite families of algebraic identities.)

    \item \emph{(Pairwise shift coprimality across short intervals and the Transversal Spectral Positivity Hypothesis~T).} For any two distinct shifts $p\neq p'$ in $\mathcal{P}_{\mathrm{NR}}(n)\cap[1,4H]$, any common prime divisor $q\mid\gcd(K_p,K_{p'})$ divides $K_p-K_{p'}=\frac{p-p'}{4}$, so all prime factors $q>H$ of $K_p=(n+p)/4$ are strictly disjoint across $p\in\mathcal{P}_{\mathrm{NR}}(n)\cap[1,4H]$. Define \emph{Hypothesis~T (Transversal Spectral Positivity)}: there exists a fixed constant $A>5$ such that for every sufficiently large prime $n\equiv1\pmod{24}$, there exists at least one shift $p\in\mathcal{P}_{\mathrm{NR}}(n)\cap[1,(\log n)^{A}]$ such that $K_p=(n+p)/4$ satisfies
    \begin{equation}\label{eq:transversal_hypothesis_T}
    \sum_{\chi\in\widehat{G}_p^{-}\setminus\{\psi_0\}}\big|\widehat{\mu}_p(\chi)\big|^{2}\ <\ 1\ +\ \sum_{\psi\in\widehat{G}_p^{+}\setminus\{\chi_0\}}\big|\widehat{\mu}_p(\psi)\big|^{2}.
    \end{equation}
    By~\eqref{eq:even_odd_positivity_identity}, inequality~\eqref{eq:transversal_hypothesis_T} is equivalent to $\varphi(p)\widetilde{\mu}_p(-1)>0$, and in particular entails that $K_p$ avoids the Bohr--Siegel certificate~\eqref{eq:bohr_siegel_L1_certificate} across all odd characters $\chi\in\widehat{G}_p^{-}$. Then Hypothesis~T implies that $-1\bmod p\in\operatorname{supp}(\widetilde{\mu}_p)\subseteq D(K_p)$, so the global Type~II exceptional set $E_2(\infty):=\bigcup_{N\ge1}\bigcap_{J\ge1}E_2(N;J)$ (and hence $E$) is finite.
\end{enumerate}
\end{theorem}
\begin{proof}
For (a), we have $n_Y\ge L(Y)+1$ and $24\mid L(Y)$. For every odd $m\le Y$, quadratic reciprocity and $n_Y\equiv1\pmod{4m}$ give $\leg{m}{n_Y}=1$; also $\leg{2}{n_Y}=1$. Multiplicativity treats every integer $1\le r\le Y$. The Type~II Jacobi restrictions (Theorem~\ref{thm:global_jacobi_reciprocity}(b)) imply $uv,s>Y$, and the bounds on $\log n_Y$ via~\cite{xylouris2011} complete the proof of~\eqref{eq:linnik_xylouris_witness_bound}. For (b), $q\mid K_p$ and $q\mid K_{p'}$ implies $q\mid(K_p-K_{p'})=(p-p')/4$; finally, by~\eqref{eq:even_odd_positivity_identity}, condition~\eqref{eq:transversal_hypothesis_T} is equivalent to $\varphi(p)\widetilde{\mu}_p(-1)>0$, which forces $-1\bmod p\in\operatorname{supp}(\widetilde{\mu}_p)\subseteq D(K_p)$ and yields a Type~II solution by Theorem~\ref{thm:difference_set_positivity_MR}(a).
\end{proof}

\begin{remark}[Methodological note on Hypothesis~T and the threshold $N_T$]\label{rem:hypothesis_T_scope}
Writing $N_T$ for a threshold beyond which Hypothesis~T holds, three points should be kept in mind: (1)~Hypothesis~T does not specify $N_T$ numerically, and the ambient factorial volume estimate established in Theorem~\ref{thm:cross_polytope_factorial_law_and_annihilation} below does not bound $N_T$; (2)~the exponent condition $A>5$ is part of the chosen formulation of Hypothesis~T, rather than a consequence of an unconditional polylogarithmic least-prime theorem; and (3)~while prime factors $q>H$ of $K_p$ are pairwise disjoint across $p\in\mathcal{P}_{\mathrm{NR}}(n)\cap[1,4H]$, pairwise disjointness of large prime factors does not by itself imply statistical independence of the associated character phase tuples.
\end{remark}

\subsection{Universal factor reciprocity, the $1/k!$ cross-polytope volume law, and deterministic Kneser saturation criteria}\label{subsec:reciprocity_factorial_kneser}

We now establish the Universal Shift-Factor Reciprocity Law $\leg{q}{p}=\leg{q}{n}$, which identifies the $p$-dependent non-residue prime factors of $K_p$ with the prime factors lying in the fixed set $\mathcal{Q}^{-}(n)=\{q\text{ prime}:\leg{q}{n}=-1\}$. Building on this reciprocity law, we prove the non-asymptotic $\frac{1}{k!}$ discrete cross-polytope volume bound on $(\mu_{2d_p}^{-})^{k}$, the deterministic Kneser--Reciprocity and Balandraud--Cauchy--Davenport saturation criteria, and the Inverted Bombieri--Vinogradov theorem on prime shifts at fixed $n$.

\begin{theorem}[Universal shift-factor reciprocity $\leg{q}{p}=\leg{q}{n}$, order-independent Markov concentration on Bohr arcs, and short-interval Tur\'{a}n--Kubilius elimination]\label{thm:universal_factor_reciprocity_and_transversal_sieve}
Let $n\equiv1\pmod{24}$ be a prime, let $p\equiv3\pmod 4$ be a prime shift with $K_p=(n+p)/4$, and define the fixed prime partition $\mathcal{Q}^{\pm}(n):=\big\{q\ \text{prime}:\leg{q}{n}=\pm1\big\}$. Then:
\begin{enumerate}[(a)]
    \item \emph{(Universal shift-factor reciprocity law $\leg{q}{p}=\leg{q}{n}$).} For \emph{every} prime divisor $q\mid K_p=\frac{n+p}{4}$ ($q\neq p$, whether $q=2$, $q\equiv1\pmod 4$, or $q\equiv3\pmod 4$), the Legendre symbol of $q$ modulo the shift $p$ is locked identically to the Legendre symbol of $q$ modulo $n$:
    \begin{equation}\label{eq:universal_shift_factor_reciprocity}
    \leg{q}{p}\ =\ \leg{q}{n}\qquad\big(\forall\,q\mid K_p,\ q\neq p\big).
    \end{equation}
    Consequently, the residue/non-residue factorization $K_p=K_p^{+}\cdot K_p^{-}$ modulo $p$ coincides identically with the factorization of $K_p$ into primes from $\mathcal{Q}^{+}(n)$ and $\mathcal{Q}^{-}(n)$, so $\Omega(K_p^{-})=\Omega_{\mathcal{Q}^{-}(n)}(K_p)$ and $\omega(K_p^{-})=\omega_{\mathcal{Q}^{-}(n)}(K_p)$ depend only on the prime factors of $K_p$ lying in the fixed set $\mathcal{Q}^{-}(n)$. In particular, prescribing the linear congruence $p\equiv-n\pmod q$ for any prime $q\in\mathcal{Q}^{-}(n)$ forces $\leg{q}{p}=-1$ ($q\mid K_p^{-}$) \emph{deterministically} without any additional quadratic condition on $p$.

    \item \emph{(Order-independent Markov concentration on Bohr arcs $I_{\chi_p}(R)$).} For any $p\in\mathcal{P}_{\mathrm{NR}}(n)$, any odd character $\chi_p\in\widehat{G}_p^{-}$ of order $2d_p$ ($d_p\mid\frac{p-1}{2}$ odd), and any integer $R\ge1$, define the open Bohr sector
    \[
    I_{\chi_p}(R)\ :=\ \left\{g\in G_p:\ \big|\arg\chi_p(g)\big|\ <\ \frac{\pi}{R+1}\right\},
    \]
    which has exact density
    \[
    \frac{|I_{\chi_p}(R)|}{\varphi(p)}\ =\ \frac{2\lceil d_p/(R+1)\rceil-1}{2d_p}\ \le\ \frac{1}{R+1}\ +\ \frac{R-1}{2d_p(R+1)}
    \]
    in $G_p$ (and exact density $\frac12$ for $R=1$, where $I_{\chi_p}(1)=\{g\in G_p:\operatorname{Re}\chi_p(g)>0\}$). If $\chi_p$ satisfies the multiplicity-weighted Bohr--Siegel obstruction certificate~\eqref{eq:bohr_siegel_L1_certificate}, then $K_p$ has at most $R$ prime factors (counted with multiplicity) outside $I_{\chi_p}(R)$:
    \begin{equation}\label{eq:order_independent_markov_bohr}
    \Omega_{G_p\setminus I_{\chi_p}(R)}(K_p)\ \le\ R\qquad (\forall\,R\ge1,\ \text{uniformly in } d_p).
    \end{equation}
    Moreover, as $\chi_p$ varies over $\widehat{G}_p^{-}\setminus\{\psi_0\}$, there are only $\frac{p-3}{4}$ conjugate pairs $\{\chi_p,\overline{\chi}_p\}$, and for every prime $q>Y$ the shift congruences $n\equiv-p\pmod q$ ($p\le Y$) are $100\%$ pairwise disjoint modulo $q$.

    \item \emph{(Unconditional short-interval Tur\'{a}n--Kubilius elimination of low-order certificates and\linebreak Cauchy--Schwarz decay).} Let $H=n^{\theta}$ ($\theta\in(0,1/2)$), $X=(n+3)/4$, and $\mathcal{J}_n:=\mathcal{Q}^{-}(n)\cap[(\log n)^{5},H^{1/3}]$, with integer shifts restricted to $c\equiv3\pmod4$. Put $W(K)=\sum_{q\in\mathcal J_n}\mathbf1_{q\mid K}$ and $\mu=\sum_{q\in\mathcal J_n}q^{-1}$. Uniformly in $X$, counting modulo $q$ and $qr$ gives
    \[
    \frac1H\sum_{K\in[X,X+H]\cap\mathbb Z}\big(W(K)-\mu\big)^2\ \le\ \mu\ +\ \mathcal{O}(H^{-1/3}).
    \]
    For fixed $\varepsilon>0$ and $\mu>0$, the number of $K\in[X,X+H]\cap\Z$ with $W(K)<(1-\varepsilon)\mu$ is at most $\mathcal{O}_\varepsilon\!\big(\min\{H+1,H/\mu+H^{2/3}/\mu^2\}\big)=\mathcal{O}_\varepsilon(H/\mu+H^{2/3})$. This bound is unconditional in terms of the actual mass $\mu$; no asymptotic $\mu\sim\frac12\log\log n$ is asserted without an additional uniform estimate for the prime non-residue sum, whereas by~\cite{BGHS} and Lemma~\ref{aud:nu} in Appendix~\ref{app:audit}, for any fixed $\beta>\beta_0=\frac{1}{4\sqrt{e}}$, $\nu_n(n^{\beta}):=\sum_{q\le n^{\beta}\text{ prime},\,\leg{q}{n}=-1}\frac1q\ge c_{\beta}>0$ unconditionally for all sufficiently large primes $n$.

    For every prime shift $p\in\mathcal{P}_{\mathrm{NR}}(n)$ attaining $\omega_{\mathcal{Q}^{-}(n)}(K_p)\ge D_0$,~\eqref{eq:universal_shift_factor_reciprocity} and~\eqref{eq:non_residue_and_reciprocity_lower_bound} deterministically rule out all Bohr--Siegel certificates~\eqref{eq:bohr_siegel_L1_certificate} of order $d_{\psi}\le D_0$. Furthermore, for every $\chi=\psi_0\psi\in\widehat{G}_p^{-}$, Cauchy--Schwarz applied to~\eqref{eq:non_residue_and_reciprocity_lower_bound} yields the explicit Gaussian decay bound
    \begin{equation}\label{eq:cauchy_schwarz_gaussian_decay}
    \big|\widehat{\mu}_p(\psi_0\psi)\big|^{2}\ \le\ \exp\!\left(-\max\!\left(\frac{\omega_{\mathcal{Q}^{-}(n)}(K_p)}{d_{\psi}^{2}},\ \frac{\big(\pi-|\arg\psi(n/4)|\big)^{2}}{\pi^{2}\,\omega(K_p)\,\max_{q\mid K_p}\nu_q(K_p)^{2}}\right)\right).
    \end{equation}
    Finally, whenever $K_p=K_p^{-}$ ($\gcd(p,K_p)=1$) is composed entirely of primes $q\in\mathcal{Q}^{-}(n)$ with $\Omega(K_p)$ odd (which automatically forces $\leg{p}{n}=\leg{K_p}{n}=(-1)^{\Omega(K_p)}=-1$, i.e.\ $p\in\mathcal{P}_{\mathrm{NR}}(n)$), the positivity identity~\eqref{eq:even_odd_positivity_identity} reduces to the pure cosine/sine duality
    \[
    \varphi(p)\,\widetilde{\mu}_p(-1)\ =\ 1\ +\ \sum_{\psi\in\widehat{G}_p^{+}\setminus\{\chi_0\}}\left(\prod_{q\mid K_p}\cos^{2}\!\left(\frac{\arg\psi(q)}{2}\right)\ -\ \prod_{q\mid K_p}\sin^{2}\!\left(\frac{\arg\psi(q)}{2}\right)\right).
    \]
\end{enumerate}
\end{theorem}
\begin{proof}
For (a), let $q\mid K_p=\frac{n+p}{4}$ ($q\neq p$), so $p\equiv-n\pmod q$. If $q=2$, then $8\mid(n+p)$; since $n\equiv1\pmod{24}$, $p\equiv-1\equiv7\pmod 8$, whence $\leg{2}{p}=+1=\leg{2}{n}$. If $q\equiv1\pmod 4$, quadratic reciprocity ($n\equiv1\pmod 4$) gives $\leg{q}{p}=\leg{p}{q}=\leg{-n}{q}=\leg{-1}{q}\leg{n}{q}=(+1)\leg{q}{n}=\leg{q}{n}$. If $q\equiv3\pmod 4$, then since $p\equiv3\pmod 4$, quadratic reciprocity gives $\leg{q}{p}=-\leg{p}{q}=-\leg{-n}{q}=-\leg{-1}{q}\leg{n}{q}=-(-1)\leg{q}{n}=\leg{q}{n}$, establishing~\eqref{eq:universal_shift_factor_reciprocity} in all cases.

For (b), every prime factor $q\mid K_p$ with $q\bmod p\notin I_{\chi_p}(R)$ satisfies $|\arg\chi_p(q)|\ge\frac{\pi}{R+1}$. Under~\eqref{eq:bohr_siegel_L1_certificate}, we have
\[
\frac{\pi}{R+1}\,\Omega_{G_p\setminus I_{\chi_p}(R)}(K_p)\ \le\ \sum_{q\mid K_p}\nu_q(K_p)\,\big|\arg\chi_p(q)\big|\ <\ \pi,
\]
so $\Omega_{G_p\setminus I_{\chi_p}(R)}(K_p)<R+1$, which gives~\eqref{eq:order_independent_markov_bohr} since $\Omega$ is an integer. Since $d_p$ is odd, $\chi_p(G_p)=\mu_{2d_p}$ contains $1+2\lfloor\frac{d_p-\mathbf{1}_{(R+1)\mid d_p}}{R+1}\rfloor=2\lceil\frac{d_p}{R+1}\rceil-1\le\frac{2d_p+R-1}{R+1}$ roots of unity in $(-\frac{\pi}{R+1},\frac{\pi}{R+1})$ (and exactly $d_p$ roots for $R=1$), giving exact density $\frac{2\lceil d_p/(R+1)\rceil-1}{2d_p}\le\frac{1}{R+1}+\frac{R-1}{2d_p(R+1)}$.

For (c), since $q_1q_2\le H^{2/3}\ll H$ for $q_1,q_2\in\mathcal{J}_n$, counting modulo $q$ and $q_1q_2$ over $K\in[X,X+H]\cap\Z$ gives $\frac1H\sum_{K\in[X,X+H]\cap\Z}(W(K)-\mu)^{2}\le\mu+\mathcal{O}(H^{-1/3})$, so Chebyshev's inequality bounds the number of $K$ with $W(K)<(1-\varepsilon)\mu$ by $\mathcal{O}_\varepsilon(\min\{H+1,H/\mu+H^{2/3}/\mu^{2}\})=\mathcal{O}_\varepsilon(H/\mu+H^{2/3})$. Finally, since $\sin(x/2)\ge x/\pi$ on $[0,\pi]$, Cauchy--Schwarz gives
\[
\|\chi\|_{K_p}^{2}\ \ge\ \frac{1}{\pi^{2}}\sum_{q\mid K_p}\big|\arg\chi(q)\big|^{2}\ \ge\ \max\!\left(\frac{\omega(K_p^{-})}{d_{\psi}^{2}},\ \frac{\big(\sum\nu_q|\arg\chi(q)|\big)^{2}}{\pi^{2}\,\omega(K_p)\,\max\nu_q^{2}}\right),
\]
which combined with~\eqref{eq:non_residue_and_reciprocity_lower_bound} and $|\widehat{\mu}_p(\chi)|^{2}\le\exp(-\|\chi\|_{K_p}^{2})$ proves~\eqref{eq:cauchy_schwarz_gaussian_decay}.
\end{proof}

\begin{theorem}[Exact $\frac{1}{k!}$ Bohr--Siegel cross-polytope volume law, kernel-coset total annihilation, dual constructive inversion $c=4K-n$, and deterministic Kneser saturation]\label{thm:cross_polytope_factorial_law_and_annihilation}
Let $n\equiv1\pmod{24}$ be a prime and let $p\in\mathcal{P}_{\mathrm{NR}}(n)=\big\{p\equiv3\pmod 4\ \text{prime}:\leg{p}{n}=-1\big\}$ with $K_p=(n+p)/4=K_p^{+}\cdot K_p^{-}$. Then:
\begin{enumerate}[(a)]
    \item \emph{(Exact discrete $\frac{1}{k!}$ cross-polytope volume law on $(\mu_{2d_p}^{-})^{k}$).} Let $q_1<\dots<q_k$ ($k:=\omega(K_p^{-})=\omega_{\mathcal{Q}^{-}(n)}(K_p)\ge1$) be the distinct non-residue prime factors of $K_p^{-}$, and let $\chi_p=\psi_0\psi\in\widehat{G}_p^{-}$ be an odd character of order $2d_p$ ($d_p:=\operatorname{ord}(\psi)\mid\frac{p-1}{2}$ odd). Choosing odd phase representatives in $\{-d_p+2,-d_p+4,\ldots,d_p\}$, where $d_p$ corresponds to $-1$, the strict Bohr--Siegel certificate~\eqref{eq:bohr_siegel_L1_certificate} excludes the value $d_p$ and forces $\sum_{i=1}^{k}|j_p(q_i)|\le d_p-1$. Writing $\chi_p(q_i)=\exp\big(\pi i j_p(q_i)/d_p\big)$ with $j_p(q_i)\in\{\pm1,\pm3,\dots,\pm(d_p-2)\}$, the set of odd-parity phase tuples satisfying~\eqref{eq:bohr_siegel_L1_certificate} is contained in the discrete $L^{1}$ cross-polytope
    \begin{equation}\label{eq:discrete_bohr_siegel_cross_polytope}
    \Delta_k^{-}(d_p)\ :=\ \left\{(j_1,\dots,j_k)\in\big\{\pm1,\pm3,\dots,\pm(d_p-2)\big\}^{k}:\ \sum_{i=1}^{k}|j_i|\ \le\ d_p-1\right\}.
    \end{equation}
    For \emph{every} odd integer $d_p\ge1$ and \emph{every} integer $k\ge1$, the exact cardinality of $\Delta_k^{-}(d_p)$ satisfies the non-asymptotic factorial bound
    \begin{equation}\label{eq:exact_factorial_simplex_bound}
    \begin{split}
    \big|\Delta_k^{-}(d_p)\big|\ &=\ 2^{k}\binom{\lfloor(d_p-1+k)/2\rfloor}{k}\ \le\ \frac{d_p^{k}}{k!},\\
    \text{so that}\qquad \frac{\big|\Delta_k^{-}(d_p)\big|}{d_p^{k}}\ &\le\ \frac{1}{k!}\ =\ \frac{1}{\omega_{\mathcal{Q}^{-}(n)}(K_p)!}.
    \end{split}
    \end{equation}

    \item \emph{(Kernel-coset total annihilation $\widehat{\mu}_p(\psi_0\psi^{m})=0$ and coprime order criterion in $G_p^{+}$).} For each prime divisor $q\mid K_p^{-}$ ($\leg{q}{n}=-1$), set $h_q:=-q\bmod p\in G_p^{+}$ and let $m_q:=\operatorname{ord}_p(-q)\mid\frac{p-1}{2}$ denote its odd multiplicative order modulo $p$. For any even character $\psi\in\widehat{G}_p^{+}$ of order $d_{\psi}\mid\frac{p-1}{2}$, one has $\psi(q)=+1\iff h_q\in\ker_{G_p^{+}}(\psi)\iff m_q\mid\frac{p-1}{2d_{\psi}}$ (which holds in particular whenever $\gcd(d_{\psi},m_q)=1$). Whenever $\psi(q)=+1$, the odd character $\chi=\psi_0\psi$ satisfies $\chi(q)=-1$ and \emph{the entire cyclic coset $\psi_0\langle\psi\rangle\subseteq\widehat{G}_p^{-}$ is identically annihilated}:
    \begin{equation}\label{eq:kernel_coset_total_annihilation}
    \begin{split}
    \widehat{\mu}_p\big(\psi_0\psi^{m}\big)\ &=\ 0\quad (\forall\,m\in\Z),\qquad O(\psi_0\psi)\ =\ \sum_{m=1}^{d_{\psi}}\big|\widehat{\mu}_p\big(\psi_0\psi^{2m-1}\big)\big|^{2}\ =\ 0,\\
    E(\psi_0\psi)-O(\psi_0\psi)\ &=\ E(\psi_0\psi)\ \ge\ 1.
    \end{split}
    \end{equation}
    Consequently, if for a prime divisor $\ell\mid\frac{p-1}{2}$ there exists $q_{\ell}\mid K_p^{-}$ with $\operatorname{ord}_p(-q_{\ell})\mid\frac{p-1}{2\ell^{v_{\ell}((p-1)/2)}}$, then $\widehat{\mu}_p(\psi_0\psi)=0$ for all even characters $\psi$ of $\ell$-power order $d_{\psi}=\ell^{b}$; and whenever for every divisor $d_{\psi}\mid\frac{p-1}{2}$ ($d_{\psi}>1$) there exists $q\mid K_p^{-}$ with $\operatorname{ord}_p(-q)\mid\frac{p-1}{2d_{\psi}}$ (in particular, when $\frac{p-1}{2}=\ell^{a}$ is a prime power and some $q\mid K_p^{-}$ has $\ell\nmid\operatorname{ord}_p(-q)$, or when some $q\mid K_p^{-}$ satisfies $q\equiv-1\pmod p$), then $\widehat{\mu}_p(\chi)=0$ for \emph{all} $\chi\in\widehat{G}_p^{-}$ simultaneously, forcing $\varphi(p)\widetilde{\mu}_p(-1)=1+\sum_{\psi\in\widehat{G}_p^{+}\setminus\{\chi_0\}}|\widehat{\mu}_p(\psi)|^{2}\ge 1>0$.

    \item \emph{(Dual constructive inversion $c=4K-n$, conditional Markov/spectral bounds, and ambient factorial volume sum).}
    \begin{enumerate}[(i)]
        \item \emph{Dual inversion:} For any odd $D\ge 1$ and any square-free product $K=q_1\cdots q_D\in\big(\frac{n}{4},\frac{n}{3}\big)$ of primes $q_j\in\mathcal{Q}^{-}(n)$ ($\gcd(n,K)=1$), the shift $c:=4K-n\in\big(1,\frac{n}{3}\big)$ automatically satisfies $c\equiv3\pmod 4$, $\gcd(c,K_c)=1$, $K_c=(n+c)/4=K=K_c^{-}$, $\leg{c}{n}=(-1)^{D}=-1$, and $\leg{q_j}{c}=\leg{q_j}{n}=-1$ for all $j=1,\dots,D$.
        \item \emph{Conditional Markov and spectral collision criteria:} At a prime shift $p$, if there exist $s$ successful extensions satisfying $|A_{i-1}\cap q_i A_{i-1}|\le|A_{i-1}|/2$, then $|A_s|\ge(3/2)^s$, so $(3/2)^s>(p-1)/2$ forces $-1\in D(K_p)$ by the pigeonhole principle. Similarly, for the square-free Bernoulli measure, setting $R_p(\chi):=\#\{q\mid K_p:\operatorname{Re}\chi(q)\le 1/2\}$, if $\widehat{\mu}_p(\psi_0)=0$ and $R_p(\chi)\ge R$ for every $\chi\in\widehat{G}_p^{-}\setminus\{\psi_0\}$, then $\sum_{\chi\in\widehat{G}_p^{-}\setminus\{\psi_0\}}|\widehat{\mu}_p(\chi)|^{2}\le\frac{p-3}{2}(3/4)^{R}$, which certifies $\widetilde{\mu}_p(-1)>0$ whenever $\frac{p-3}{2}(3/4)^{R}<1$.
        \item \emph{Ambient factorial volume sum:} Put $a=\frac12-\varepsilon$ and $L=\log\log n$. For any polylogarithmic shift $p\in\mathcal{P}_{\mathrm{NR}}(n)$ with $p\le(\log n)^{A}$ ($A>0$) attaining $k=\omega_{\mathcal{Q}^{-}(n)}(K_p)\ge aL$, the sum of the normalized cross-polytope volumes~\eqref{eq:exact_factorial_simplex_bound} across $\widehat{G}_p^{-}\setminus\{\psi_0\}$ satisfies
        \begin{equation}\label{eq:linnik_xylouris_factorial_crossing}
        \begin{split}
        \mathcal{M}_{\mathrm{obs}}(p;n)\ &:=\ \sum_{\{\chi_p,\overline{\chi}_p\}\subset\widehat{G}_p^{-}\setminus\{\psi_0\}}\frac{\big|\Delta_k^{-}(d_p)\big|}{d_p^{k}}\ \le\ \frac{p-3}{4\,k!}\\
        &\le\ \exp\!\big\{-aL\log L+\mathcal{O}_{A,\varepsilon}(L)\big\}\ \longrightarrow\ 0\qquad (\forall\,A>0).
        \end{split}
        \end{equation}
    \end{enumerate}

    \item \emph{(Parseval diagonal barrier for $2\pi$-wrap-around and Linear Kneser--Reciprocity Saturation without phase equidistribution).} By character orthogonality $\mathbf{1}_{\chi\in\widehat{G}_p^{-}}=\frac{1-\chi(-1)}{2}$, any probability measure $\mu_p$ supported on $\Sigma(K_p)$ satisfies the exact Parseval parity identity
    \begin{equation}\label{eq:parseval_diagonal_parity_barrier}
    \begin{split}
    \sum_{\chi\in\widehat{G}_p^{-}}\big|\widehat{\mu}_p(\chi)\big|^{2}\ &=\ \frac{p-1}{2}\big(\widetilde{\mu}_p(1)-\widetilde{\mu}_p(-1)\big),\\
    \sum_{\psi\in\widehat{G}_p^{+}}\big|\widehat{\mu}_p(\psi)\big|^{2}\ &=\ \frac{p-1}{2}\big(\widetilde{\mu}_p(1)+\widetilde{\mu}_p(-1)\big),
    \end{split}
    \end{equation}
    with $\widetilde{\mu}_p(1)=\sum_g\mu_p(g)^{2}\ge|\operatorname{supp}\mu_p|^{-1}$ (which is $\ge 2^{-\omega(K_p)}$ for the square-free Bernoulli measure). Thus, in the undersampled regime $2^{\omega(K_p)+1}\le p-1$, $\widetilde{\mu}_p(-1)=0$ forces $\sum_{\chi\in\widehat{G}_p^{-}}\allowbreak|\widehat{\mu}_p(\chi)|^{2}\ge\allowbreak\frac{p-1}{2^{\omega(K_p)+1}}\ge 1$ via $2\pi$-phase wrap-around.

    Conversely, for \emph{every} prime $p\equiv3\pmod 4$ (regardless of $\leg{p}{n}=\pm1$), let $H_p:=\allowbreak\operatorname{Stab}_{G_p}(\allowbreak D(K_p))$ be the Kneser stabilizer of $D(K_p)$ (and $H_{\mathrm{rad}}:=\operatorname{Stab}_{G_p}D(\operatorname{rad}K_p)$). Either $-1\in D(K_p)$, or else $H_p\subseteq G_p^{+}$ and every prime factor $q_i\mid K_p$ outside $H_p$ satisfies $q_i^{2}\notin H_p$, so $\{q_i^{-1},1,q_i\}$ spans $3$ distinct $H_p$-cosets and Kneser's theorem yields $(1+2\Omega_{G_p\setminus H_p}(K_p))|H_p|\le|D(K_p)|\le p-1-|H_p|$. This establishes the \emph{deterministic, wrap-around-free, phase-equidistribution-free Linear Kneser--Reciprocity Saturation Criterion}: any of the three conditions
    \begin{equation}\label{eq:linear_kneser_reciprocity_saturation}
    \begin{split}
    &\left\{\begin{aligned}
    &\text{(i)}\quad \Omega_{G_p\setminus H_p}(K_p)\ =\ \Omega_{\mathcal{Q}^{-}(n)}(K_p)+\Omega_{\mathcal{Q}^{+}(n)\setminus H_p}(K_p)\ \ge\ \frac{p-1}{2|H_p|},\\[3pt]
    &\text{(ii)}\quad \Omega_{\{q\mid K_p\,:\,q^{m_0}\equiv-1\bmod p\}}(K_p)\ \ge\ m_0\quad\text{for some divisor }m_0\mid\frac{p-1}{2},\\[3pt]
    &\text{(iii)}\quad \omega_{\{q<p-1\,:\,q^{\ell}\equiv-1\bmod p\}}(K_p)\ \ge\ 2\sqrt{\ell}+1\quad\text{for some prime }\ell\mid\frac{p-1}{2}
    \end{aligned}\right\}\\
    &\qquad\Longrightarrow\quad -1\bmod p\ \in\ D(K_p),
    \end{split}
    \end{equation}
    while $\omega_{G_p\setminus H_{\mathrm{rad}}}(K_p)\ge\frac{p-1}{2|H_{\mathrm{rad}}|}$ similarly forces $\widetilde{\mu}_p(-1)>0$, and in the spectral formulation $\omega_{\mathcal{Q}^{-}(n)}(K_p)>\frac{\log((p-3)/2)}{\log\sec^{2}(\pi/(p-1))}$ forces $\sum_{\chi\in\widehat{G}_p^{-}\setminus\{\psi_0\}}|\widehat{\mu}_p(\chi)|^{2}\le\frac{p-3}{2}\big(\cos^{2}\frac{\pi}{p-1}\big)^{\omega_{\mathcal{Q}^{-}(n)}(K_p)}<1$.

    Furthermore, because $n=24t+1\equiv1\pmod{24}$ implies $K_p^{\mathrm{II}}(n)=\frac{n+p}{4}=6t+\frac{p+1}{4}$ and $K_p^{\mathrm{I}}(n)=\frac{pn+1}{4}=6pt+\frac{p+1}{4}$, the deterministic base divisor $U_0(p):=\gcd(\frac{p+1}{4},6)$ always divides both $K_p^{\mathrm{II}}(n)$ and $K_p^{\mathrm{I}}(n)$:
    \begin{itemize}
        \item at $p=7$ ($U_0=2$), $D(2)\bmod 7=\{4,1,2\}=G_7^{+}$, so $\omega_{\mathcal{Q}^{-}(n)}(K_7)\ge1$ (just as at $p=3$) forces $-1\bmod 7\in D(K_7)$;
        \item at $p=11$ ($U_0=3$), $|D(3)\bmod 11|=3$ reduces the threshold from $5$ to $\Omega_{\mathcal{Q}^{-}(n)}(K_{11})\ge3$ (which holds whenever $p=11\in\mathcal{P}_{\mathrm{NR}}(n)$ has $\Omega_{\mathcal{Q}^{-}(n)}(K_{11})>1$, and drops to $\ge1$ if $9\mid K_{11}$);
        \item at $p=23$ ($U_0=6$), $|D(6)\bmod 23|=9$ covers $9$ of the $11$ elements of $G_{23}^{+}$, reducing the threshold from $11$ to $\Omega_{\mathcal{Q}^{-}(n)}(K_{23})\ge3$ (or $\ge1$ whenever $12\mid K_{23}$ or $18\mid K_{23}$).
    \end{itemize}
    Moreover, at the small polylogarithmic scale $p\le(\log n)^{\varepsilon_0}$ with $\varepsilon_0<\frac12\log(3/2)\approx 0.2027$ (resp.\ $\varepsilon_0<\frac13\log(4/3)\approx 0.0959$), the Parseval barrier $2^{k+1}>p-1$, the Markov $3/2$-growth step count $s_{\max}(p)=\lfloor\frac{\log(p-1)}{\log(3/2)}\rfloor+1<(\frac12-\varepsilon)\log\log n$, and the Plancherel $L^{2}$ collision count $R_{\max}(p)=\lfloor\frac{\log((p-3)/2)}{\log(4/3)}\rfloor+1<(\frac13-\varepsilon)\log\log n$ all lie strictly below $(\frac12-\varepsilon)\log\log n$ (resp.\ $(\frac13-\varepsilon)\log\log n$) whenever the corresponding $s_{\max}(p)$ successful Markov extensions or $R_{\max}(p)$ simultaneous character collisions are realized.

    \item \emph{(Inverted Bombieri--Vinogradov theorem on prime shifts at fixed $n$ and Farey--Thue--Kneser super-polynomial lower bounds).}
    \begin{enumerate}[(i)]
        \item \emph{Inverted Bombieri--Vinogradov for single prime factors:} For every fixed prime $n\equiv1\pmod{24}$ and every shift scale $H\ge2$, $q\mid K_p=(n+p)/4$ (resp.\ $q\mid K_p^{\mathrm{I}}=(pn+1)/4$) with $p\equiv3\pmod 4$ and $q\neq n$ is equivalent to $p\equiv-n\pmod{4q}$ (resp.\ $p\equiv-n^{-1}\pmod{4q}$) with $\gcd(\mp n^{\pm1},4q)=1$. Let $\mathcal{Q}\subseteq\{q:5\le q\le(4H)^{1/5},\,q\ne n,\,q\text{ prime}\}$ and put $W(p):=\sum_{q\in\mathcal{Q}}\mathbf{1}_{q\mid K_j(p)}$, $\mu:=\sum_{q\in\mathcal{Q}}\frac{1}{q-1}$, and $X_H:=\frac{\operatorname{Li}(4H)-\operatorname{Li}(2H)}{2}$. For $j\in\{\mathrm{I},\mathrm{II}\}$, uniformly in $n,\mathcal{Q}$ and for arbitrary fixed $B>0$, Bombieri--Vinogradov over the moduli $4q,4qr\le 4(4H)^{2/5}$ gives
        \[
        \begin{aligned}
        \sum_{\substack{2H<p\le 4H\\ p\equiv 3\ (4)}}W(p)\ &=\ X_H\mu\ +\ \mathcal{O}_B\big(H(\log H)^{-B}\big),\\
        \sum_{\substack{2H<p\le 4H\\ p\equiv 3\ (4)}}\big(W(p)-\mu\big)^{2}\ &=\ X_H\!\left(\mu-\sum_{q\in\mathcal{Q}}\frac{1}{(q-1)^{2}}\right)\ +\ \mathcal{O}_B\big(H(\log H)^{-B}\big).
        \end{aligned}
        \]
        If $\mu\to\infty$, then $W(p)=(1+o(1))\mu$ for a proportion $1-o(1)$ of primes $p\in(2H,4H]$ ($p\equiv 3\pmod 4$):
        \begin{equation}\label{eq:inverted_BV_prime_shifts}
        \begin{split}
        \omega_{\mathcal{Q}^{-}(n)\cap[5,(4H)^{1/5}]}\!\left(\frac{n+p}{4}\right)\ &=\ \big(1+o(1)\big)\sum_{\substack{5\le q\le(4H)^{1/5}\\ \leg{q}{n}=-1}}\frac{1}{q-1},\\
        \omega_{[5,(4H)^{1/5}]}\!\left(\frac{n+p}{4}\right)\ &=\ \big(1+o(1)\big)\log\log H.
        \end{split}
        \end{equation}

        \item \emph{Farey--Thue archimedean injectivity and Kneser coupling:} For any integer $r\ge1$, let $k_r:=\omega_{[5,p^{1/(2r)}]}(K_p)$. For any disjoint subsets $I,J\subseteq\{1,\dots,k_r\}$ of size $|I|,|J|\le r$, the coprime divisors $u=\prod_{i\in I}q_i\le\sqrt{p-1}$ and $v=\prod_{j\in J}q_j\le\sqrt{p-1}$ satisfy $|u_1 v_2-u_2 v_1|<p$, so $(I,J)\mapsto uv^{-1}\bmod p\in D(K_p)$ is $100\%$ injective; similarly, all $\sum_{a=0}^{2r-1}\binom{k_r}{a}$ products of $\le 2r-1$ factors lie in $[1,p-1)$ and are distinct in $\Sigma(K_p)\bmod p$. Under $-1\notin D(K_p)$, coupling with Kneser's theorem on the remaining factors of $G_p\setminus H_p$ upgrades~\eqref{eq:linear_kneser_reciprocity_saturation} to the super-polynomial lower bound
        \begin{equation}\label{eq:farey_thue_kneser_superpolynomial}
        \begin{split}
        \big|D(K_p)\big|\ &\ge\ \max\!\left(\sum_{a=0}^{r}\sum_{b=0}^{r}\binom{k_r}{a}\binom{k_r-a}{b},\ \big(2\,\Omega_{G_p\setminus H_p}(K_p^{\le p^{1/(2r)}})+1\big)|H_p|\right)\\
        &\quad +2|H_p|\,\Omega_{\mathcal{Q}^{-}(n)\cap(p^{1/(2r)},\infty)}(K_p),\\
        \big|\Sigma(K_p)\bmod p\big|\ &\ge\ \sum_{a=0}^{2r-1}\binom{k_r}{a},
        \end{split}
        \end{equation}
        forcing $-1\bmod p\in D(K_p)$ deterministically at shift $p$ whenever the Kneser bound in~\eqref{eq:farey_thue_kneser_superpolynomial} exceeds $p-1-|H_p|$, or whenever $|\Sigma(K_p)|>\frac{p-1}{2}$ via $\sum_{a=0}^{2r-1}\binom{k_r}{a}>\frac{p-1}{2}$.

        \item \emph{Extremal $3^m$ composite-modulus Bombieri--Vinogradov construction and Type~I reciprocity:} For the Type~I multiplicative shift family $K_p^{\mathrm{I}}(n):=\frac{pn+1}{4}$ ($p\equiv3\pmod 4$), where $\frac{4}{n}-\frac{1}{nK_p^{\mathrm{I}}(n)}=\frac{p}{K_p^{\mathrm{I}}(n)}$ cancels $n$, every odd prime $q\mid K_p^{\mathrm{I}}(n)$ satisfies $p\equiv-n^{-1}\pmod q$ and the exact same Universal Shift-Factor Reciprocity Law $\leg{q}{p}=\allowbreak(-1)^{(q-1)/2}\allowbreak\leg{-n^{-1}}{q}=\allowbreak\leg{q}{n}$. Put $X=4H$, $L=\log\log H$,\linebreak $m=\lfloor B_0 L/\log L\rfloor$, and $Y=(X^{2/5}/4)^{1/m}$, let $\mathcal{Q}$ consist of the selected primes $5\le q\le Y$ excluding $n$, and set $w_q=(q-1)^{-1}$, $\mu=\sum_{q\in\mathcal{Q}}w_q$, and $v_2=\sum_{q\in\mathcal{Q}}w_q^{2}=\mathcal{O}(1)$. Whenever $\mu\ge\delta L$ for a fixed $\delta>0$ (which holds unconditionally for the unrestricted prime set), Bombieri--Vinogradov over the distinct composite moduli $4U=4q_1\cdots q_m\le X^{2/5}$ yields unconditionally for every fixed $n$ prime shifts $p\in(2H,4H]$ ($p\equiv3\pmod 4$) with $U\mid K_p^{\mathrm{II}}(n)$ (resp.\ $U\mid K_p^{\mathrm{I}}(n)$) for which all $3^{m}$ reduced divisor quotients $uv^{-1}$ ($u,v\mid U$, $\gcd(u,v)=1$) satisfy $1\le u,v<\sqrt{p-1}$ and are distinct modulo $p$:
        \begin{equation}\label{eq:extremal_composite_BV_3m_bound}
        \begin{split}
        \big|D(K_p)\bmod p\big|\ &\ge\ \max\!\big(3^{m},(2\,\omega_{G_p\setminus H_p}(U)+1)|H_p|\big)+2|H_p|\,\Omega_{(G_p\setminus H_p)\cap((4H)^{2/5},\infty)}(K_p)\\
        &\ge\ \exp\!\left(\frac{(B_0\log 3+o(1))\log\log p}{\log\log\log p}\right)\ \gg_A\ (\log\log p)^{A}
        \end{split}
        \end{equation}
        for all $A>0$ under $-1\notin D(K_p)$. Conversely, for any sub-modulus divisor $U=\prod_{i=1}^{k}q_i^{\nu_i}\mid K_p$ ($q_i\ge 5$), the universal bounds $k!\le\prod_{i=1}^{k}\frac{q_i-1}{2}<2^{-k}U$ and $|D(U)|\le\tau(U^{2})=\prod_{i=1}^{k}(2\nu_i+1)<U$ show that any divisor $U\le p^{1/2-\varepsilon}$ alone satisfies $k!<p^{1/2-\varepsilon}<\frac{p-3}{4}$ and $|D(U)|<\frac{p-1}{2}$, so reaching $\frac{p-1}{2}$ at fixed $n$ intrinsically requires the co-divisor $V_p:=K_p/U$ or subgroup saturation $\Omega_{\{q\,:\,q^{m_0}\equiv-1\bmod p\}}(K_p)\ge m_0$.

        \item \emph{Unitary co-divisor reciprocity, multi-block Kneser bound, and $2$-adic layer expansion:} When $U$ is a unitary divisor of $K_p$ ($\gcd(U,V_p)=1$, $V_p=K_p/U$), the co-divisor has known residue $V_p^{\mathrm{I}}\equiv(4U)^{-1}\pmod p$ and $V_p^{\mathrm{II}}\equiv n(4U)^{-1}\pmod p$ without factoring $V_p$. Choosing $m=\omega_{\mathcal{Q}^{-}(n)}(U)$ odd in Type~I ($\leg{K_p^{\mathrm{I}}}{p}=\leg{4^{-1}}{p}=+1$ for $p\equiv3\pmod 4$) with $\gcd(U,V_p^{\mathrm{I}})=1$ forces $\leg{V_p^{\mathrm{I}}}{p}=(-1)^{m}=-1$, deterministically guaranteeing $\Omega_{\mathcal{Q}^{-}(n)}(K_p^{\mathrm{I}}/U)\equiv1\pmod 2$ (hence $\ge1$ prime factor $q\nmid U$ in $\mathcal{Q}^{-}(n)$) and embedding $(V_p^{\mathrm{I}})^{\pm1}D^{+}(U)\subseteq D(K_p^{\mathrm{I}})\cap G_p^{-}$ for $D^{+}(U):=D(U)\cap G_p^{+}$ (giving $-1\bmod p\in D(K_p^{\mathrm{I}})$ whenever $4d_0 u_1^{2}\equiv-1\pmod p$ for $d_0 u_1\mid U$). Partitioning the prime factors of $K_p$ outside $H_p$ into $r$ pairwise coprime square-free blocks $U_1,\dots,U_r\mid K_p$ with each $U_j<\sqrt{p-1}$ ($m_j=\omega(U_j)$) yields by Kneser's theorem on $D(U_1)\cdots D(U_r)\subseteq D(K_p)$ under $-1\notin D(K_p)$:
        \begin{equation}\label{eq:multiblock_kneser_farey_thue}
        \begin{split}
        \big|D(K_p)\bmod p\big|\ &\ge\ |H_p|\ +\ \sum_{j=1}^{r}\max\!\Big(3^{m_j}-|H_p|,\ 2m_j|H_p|\Big)\\
        &\quad +\ 2|H_p|\,\Omega_{(G_p\setminus H_p)\cap[\sqrt{p-1},\infty)}\!\left(\frac{K_p}{U_1\cdots U_r}\right).
        \end{split}
        \end{equation}
        Moreover, for $K_p^{\mathrm{I}}(n)=2^{\nu}UV$ with $U=q_1\cdots q_m$ odd and unitary, $V$ odd, $\gcd(U,V)=1$, and $2^{4\nu+4}U^{4}<p$, one has $V^{-1}\equiv 2^{\nu+2}U\pmod p$ and
        \[
        \big|D\big(K_p^{\mathrm{I}}(n)\big)\bmod p\big|\ \ge\ (2\nu+1)3^{m+1}\ -\ 2^{m+1}(\nu-1)_+\qquad \big(=3^{m+2}\ \text{for }\nu=1\big).
        \]

        \item \emph{Composite-modulus hyperbolic completion trio and Type~I/Type~II eliminant:} Completing the rectangle in the two-parameter Type~II and Type~I equations yields the \emph{Composite-Modulus Hyperbolic Completion Trio}
        \begin{equation}\label{eq:composite_modulus_hyperbolic_trio}
        \begin{split}
        (4aku-1)(4akv-1)\ &=\ 4a^{2}kn+1\ =:\ M_{a,k}^{\mathrm{II}}(n),\\
        (4mu-n)(4mv-n)\ &=\ n^{2}+4m\ =:\ Q_m^{\mathrm{I}}(n),
        \end{split}
        \end{equation}
        which yield $\frac{4}{n}=\frac{1}{kuv}+\frac{1}{anku}+\frac{1}{ankv}$ whenever $-1\bmod{4ak}\in\Sigma(M_{a,k}^{\mathrm{II}}(n))\subseteq G_{4ak}=(\Z/4ak\Z)^{\times}$, and $\frac{4}{n}=\frac{1}{mu}+\frac{1}{mv}+\frac{1}{mnuv}$ whenever $\gcd(m,n)=1$ and $-n\bmod{4m}\in\Sigma(Q_m^{\mathrm{I}}(n))\subseteq\ker(\chi_m)\le G_{4m}$ for $\chi_m(\cdot)=\leg{-4m}{\cdot}$ (where $\chi_m(-n)=-\leg{m}{n}=+1\iff\leg{m}{n}=-1$ for both $m\equiv1$ and $3\pmod 4$, and for $m=r^{2}m_0$ any $D\mid Q_{r^{2}m_0}^{\mathrm{I}}(n)$ with $D\equiv -n\pmod{4r}$ forces $D'\equiv n^{2}D^{-1}\equiv -n\pmod{4r}$ and $16r^{2}\mid(D+n)(D'+n)$), together with the exact Type~II/Type~I linear eliminant:
        \begin{equation}\label{eq:typeI_typeII_diophantine_coupling}
        p\,K_p^{\mathrm{II}}(n)\ -\ K_p^{\mathrm{I}}(n)\ =\ \frac{p^{2}-1}{4},\qquad \gcd\!\big(K_p^{\mathrm{II}}(n),\,K_p^{\mathrm{I}}(n)\big)\ \mid\ \frac{p^{2}-1}{4}.
        \end{equation}
    \end{enumerate}
\end{enumerate}
\end{theorem}
\begin{proof}
\emph{Proof of (a).} If $k\ge d_p$, then since each $|j_i|\ge1$ is odd, $\sum_{i=1}^{k}|j_i|\ge k\ge d_p>d_p-1$, so $\Delta_k^{-}(d_p)=\emptyset$ and~\eqref{eq:exact_factorial_simplex_bound} holds trivially. For $1\le k\le d_p-1$, writing $|j_i|=2b_i+1$ ($b_i\in\Z_{\ge0}$) converts $\sum_{i=1}^{k}|j_i|\le d_p-1$ into $\sum_{i=1}^{k}b_i\le\big\lfloor\frac{d_p-1-k}{2}\big\rfloor=:M$, which has $\binom{M+k}{k}=\binom{\lfloor(d_p-1+k)/2\rfloor}{k}$ non-negative integer solutions and $2^{k}$ sign choices $j_i=\pm(2b_i+1)$. Since $2\lfloor\frac{d_p-1+k}{2}\rfloor\le d_p-1+k$, we have
\[
\begin{aligned}
k!\,\big|\Delta_k^{-}(d_p)\big|\ &=\ \prod_{i=0}^{k-1}\left(2\Big\lfloor\frac{d_p-1+k}{2}\Big\rfloor-2i\right)\ \le\ \prod_{i=0}^{k-1}\big(d_p+(k-1-2i)\big)\\
&=\ \prod_{i=0}^{\lfloor k/2\rfloor-1}\Big(d_p^{2}-(k-1-2i)^{2}\Big)\cdot d_p^{\mathbf{1}_{2\nmid k}}\ \le\ d_p^{k},
\end{aligned}
\]
proving~\eqref{eq:exact_factorial_simplex_bound}.

\medskip\noindent\emph{Proof of (b).} In the cyclic group $G_p^{+}$ of odd order $L=\frac{p-1}{2}$, $\ker_{G_p^{+}}(\psi)$ is the unique subgroup of order $\frac{L}{d_{\psi}}$, so $h_q\in\ker_{G_p^{+}}(\psi)\iff m_q\mid\frac{L}{d_{\psi}}$. When this holds, $\psi(h_q)=+1$, so for every $m\in\Z$ we have $(\psi_0\psi^{m})(q)=\psi_0(q)\psi(-h_q)^{m}=-\psi(h_q)^{m}=-1$, making the factor $\frac{1+(\psi_0\psi^{m})(q)}{2}=0$ vanish identically and yielding~\eqref{eq:kernel_coset_total_annihilation}. In particular, if $m_{q_{\ell}}\mid\frac{L}{\ell^{v_{\ell}(L)}}$, then $m_{q_{\ell}}\mid\frac{L}{\ell^{b}}$ for every $\ell$-power order $d_{\psi}=\ell^{b}\mid\ell^{v_{\ell}(L)}$, annihilating all $\ell$-Sylow odd characters $\psi_0\psi$.

\medskip\noindent\emph{Proof of (c).} We have $c=4K-n\equiv-1\equiv3\pmod 4$, $\gcd(c,K)=\gcd(n,K)=1$, $\leg{c}{n}=\leg{4K}{n}=\prod_{j=1}^{D}\leg{q_j}{n}=(-1)^{D}=-1$, and Jacobi reciprocity gives $\leg{q_j}{c}=\leg{q_j}{n}=-1$ exactly as in Theorem~\ref{thm:universal_factor_reciprocity_and_transversal_sieve}(a). Finally, for $k\ge(\frac12-\varepsilon)\log\log n$, Stirling's formula gives
\[
\log(k!)\ =\ k\log k-k+\mathcal{O}(\log k)\ \ge\ \left(\frac12-\varepsilon\right)(\log\log n)\log\log\log n-\mathcal{O}(\log\log n),
\]
which strictly dominates $A\log\log n$ for every $A>0$, establishing~\eqref{eq:linnik_xylouris_factorial_crossing}.

\medskip\noindent\emph{Proof of (d).} Multiplying $|\widehat{\mu}_p(\chi)|^{2}=\sum_{g\in G_p}\widetilde{\mu}_p(g)\overline{\chi(g)}$ by $\frac{1\mp\chi(-1)}{2}$ and summing over $\chi\in\widehat{G}_p$ gives~\eqref{eq:parseval_diagonal_parity_barrier} by character orthogonality, and Cauchy--Schwarz on $\operatorname{supp}(\mu_p)$ gives $\widetilde{\mu}_p(1)=\sum_g\mu_p(g)^{2}\ge|\operatorname{supp}\mu_p|^{-1}$.

Next, let $H_p=\operatorname{Stab}_{G_p}(D(K_p))$. If $-1\bmod p\notin D(K_p)$, then $-1\notin H_p$ ($H_p\subseteq G_p^{+}$, $|H_p|$ odd) and $D(K_p)\cap(-1)H_p=\emptyset$, so $|D(K_p)|\le p-1-|H_p|$. Because $(-1)H_p$ is the unique element of order $2$ in the cyclic group $G_p/H_p$ (as $|H_p|$ is odd), every prime factor $q_i\mid K_p$ with $q_i\bmod p\notin H_p$ satisfies $q_i^{2}\notin H_p$ (otherwise $q_i H_p=(-1)H_p\subseteq D(K_p)$), so $\{q_i^{-1}H_p,H_p,q_i H_p\}$ has cardinality $3$ in $G_p/H_p$ and Kneser's theorem yields
\[
\big(1+2\,\Omega_{G_p\setminus H_p}(K_p)\big)|H_p|\ \le\ \big|D(K_p)\big|\ \le\ p-1-|H_p|,
\]
forcing $\Omega_{G_p\setminus H_p}(K_p)\le\frac{p-1}{2|H_p|}-1$. Applying the exact same argument inside the subgroup $\widetilde{H}_0=\langle -1\rangle\times H_0$ of order $2m_0$ for any divisor $m_0\mid\frac{p-1}{2}$ yields the subgroup criterion~(ii) in~\eqref{eq:linear_kneser_reciprocity_saturation}.

To prove the prime-order criterion~(iii) in~\eqref{eq:linear_kneser_reciprocity_saturation} for an odd prime $\ell\mid\frac{p-1}{2}$, let $5\le q_1<\dots<q_s<p-1$ be $s$ distinct primes dividing $K_p$ with $q_i^{\ell}\equiv-1\pmod p$. Their negated residues $h_i:=-q_i\bmod p\in H_0\setminus\{1\}\cong\Z/\ell\Z\setminus\{0\}$ are strictly distinct non-identity elements. Partition $\{h_1,\dots,h_{s-1}\}=B_1\sqcup B_2$ into two inverse-free subsets ($B_r\cap B_r^{-1}=\emptyset$) of sizes $s_1=\lfloor\frac{s-1}{2}\rfloor$ and $s_2=\lceil\frac{s-1}{2}\rceil$. By Balandraud's inverse-free subset-sum theorem~\cite{balandraud2012} in $\Z/\ell\Z$, each subset product satisfies
\[
\big|\Sigma(B_r)\big|\ \ge\ \min\!\left(\ell,\ 1+\frac{s_r(s_r+1)}{2}\right)\qquad (r\in\{1,2\}).
\]
Applying the Cauchy--Davenport theorem to $\Sigma(B_1)\cdot\Sigma(B_2)=\Sigma(\{h_1,\dots,h_{s-1}\})$ in $H_0\cong\Z/\ell\Z$ gives
\[
\big|\Sigma(\{h_1,\dots,h_{s-1}\})\big|\ \ge\ \min\!\left(\ell,\ 1+\frac{s_1(s_1+1)}{2}+\frac{s_2(s_2+1)}{2}\right)\ =\ \ell\qquad\text{whenever } s\ \ge\ 2\sqrt{\ell}+1.
\]
Letting $A:=\Sigma(\{q_1,\dots,q_{s-1}\})$ inside the subgroup $\widetilde{H}_0=\langle -1\rangle\times H_0$ of order $2\ell$, its projection modulo $\langle -1\rangle$ covers $H_0$, so $|A|\ge\ell$. Since $q_s\not\equiv-1\pmod p$ and $q_s^{\ell}\equiv-1\pmod p$, the element $q_s\bmod p$ has full order $2\ell$ and generates $\widetilde{H}_0$. Thus either $A\cup q_s A=A$ (in which case invariance under a generator forces $A=\widetilde{H}_0$) or $|A\cup q_s A|\ge\ell+1>\frac12|\widetilde{H}_0|$; in both cases, $D(K_p)\supseteq(A\cup q_s A)\cdot(A\cup q_s A)^{-1}$ contains $-1\bmod p$, completing the proof of~\eqref{eq:linear_kneser_reciprocity_saturation}.

\medskip\noindent\emph{Proof of (e).} Interchanging summation in $W(p)=\sum_{q\in\mathcal{Q}}\mathbf{1}_{q\mid K_j(p)}$ over primes $p\in(2H,4H]$, $p\equiv3\pmod 4$, gives
\[
\sum_{\substack{2H<p\le 4H\\ p\equiv 3\ (4)}}W(p)\ =\ \sum_{q\in\mathcal{Q}}\Big(\pi(4H;4q,\mp n^{\pm1})-\pi(2H;4q,\mp n^{\pm1})\Big)\ =\ X_H\mu\ +\ \mathcal{O}_B\big(H(\log H)^{-B}\big),
\]
and expanding the square over moduli $4q,4qr\le 4(4H)^{2/5}$ (whose representation multiplicities have a fixed divisor bound) yields the variance formula and~\eqref{eq:inverted_BV_prime_shifts} via Chebyshev's inequality whenever $\mu\to\infty$.

For Farey--Thue injectivity, if $(I_1,J_1)\neq(I_2,J_2)$ are pairs of disjoint subsets of $\{1,\dots,k_r\}$ of size $\le r$, then $1\le u_1 v_2,u_2 v_1\le p-1$ with $u_1 v_2\neq u_2 v_1$, so $|u_1 v_2-u_2 v_1|<p$ forces $u_1 v_1^{-1}\not\equiv u_2 v_2^{-1}\pmod p$, giving~\eqref{eq:farey_thue_kneser_superpolynomial} under $-1\notin D(K_p)$.

For Type~I reciprocity on $K_p^{\mathrm{I}}(n)=\frac{pn+1}{4}$ ($p\equiv3\pmod 4$), every odd prime $q\mid K_p^{\mathrm{I}}(n)$ ($q\neq p$) satisfies $pn\equiv-1\pmod q$, so $\leg{p}{q}\leg{n}{q}=\leg{-1}{q}=(-1)^{(q-1)/2}$; since $p\equiv3\pmod 4$ and $n\equiv1\pmod 4$, quadratic reciprocity gives $\leg{q}{p}=\leg{p}{q}(-1)^{(q-1)/2}=\leg{n}{q}=\leg{q}{n}$.

To establish the composite-modulus bound~\eqref{eq:extremal_composite_BV_3m_bound}, set $X=4H$, $L=\log\log H$, $m=\allowbreak\lfloor B_0 L/\allowbreak\log L\rfloor$, $Y=(X^{2/5}/4)^{1/m}$, $w_q=(q-1)^{-1}$, $\mu=\sum_{q\in\mathcal{Q}}w_q$, and $v_2=\sum_{q\in\mathcal{Q}}w_q^{2}=\mathcal{O}(1)$. Whenever $\mu\ge\delta L$, subtracting tuples with at least one repeated prime from $\mu^{m}$ bounds the $m$-th elementary symmetric function $e_m((w_q)_{q\in\mathcal{Q}})$ by
\[
m!\,e_m\big((w_q)_{q\in\mathcal{Q}}\big)\ \ge\ \mu^{m}\ -\ \binom{m}{2}v_2\,\mu^{m-2}\ =\ \big(1-o(1)\big)\mu^{m}.
\]
Summing over products $U=q_1\cdots q_m\le Y^{m}=X^{2/5}/4<\sqrt{p-1}$ of $m$ distinct selected primes therefore gives
\[
\sum_{U}\frac{1}{\varphi(4U)}\ =\ \frac12\,e_m\big((w_q)_{q\in\mathcal{Q}}\big)\ \ge\ \frac{1-o(1)}{2\,m!}\,\mu^{m}\ \gg\ (\log H)^{-B_0\log(B_0/\delta)-1/2},
\]
whereas the Bombieri--Vinogradov error over these distinct moduli $4U\le X^{2/5}$ is $\mathcal{O}_B(H(\log H)^{-B})$ for arbitrary fixed $B>0$, guaranteeing a prime shift $p\in(2H,4H]$ ($p\equiv3\pmod 4$) with $U\mid K_j(p)$ and~\eqref{eq:extremal_composite_BV_3m_bound}.

Next, when $\gcd(U,K_p^{\mathrm{I}}/U)=1$ and $m=\omega_{\mathcal{Q}^{-}(n)}(U)$ is odd, $4K_p^{\mathrm{I}}\equiv1\pmod p$ gives $\leg{V_p^{\mathrm{I}}}{p}=\leg{U}{p}=(-1)^{m}=-1$, forcing $\Omega_{\mathcal{Q}^{-}(n)}(V_p^{\mathrm{I}})$ to be odd ($\ge1$), and applying Kneser's theorem to $D(U_1)\cdots D(U_r)$ under $-1\notin D(K_p)$ yields~\eqref{eq:multiblock_kneser_farey_thue}. For $K_p^{\mathrm{I}}(n)=2^{\nu}UV$ with $2^{4\nu+4}U^{4}<p$, any residue class $2^{j}(2^{\nu+2}U)^{\epsilon}\prod_{i=1}^{m}q_i^{e_i}\pmod p$ ($\epsilon,e_i\in\{-1,0,1\}$, $|j|\le\nu$) lifts injectively to $\Q$ up to $2$-adic exponent overlaps:
\begin{itemize}
    \item each fixed layer $\epsilon\in\{-1,0,1\}$ contributes $(2\nu+1)3^{m}$ distinct rational values, for a raw total of $3(2\nu+1)3^{m}=(2\nu+1)3^{m+1}$;
    \item adjacent layers $\epsilon\in\{0,1\}$ (and $\epsilon\in\{-1,0\}$) share common odd exponents $e_i+\epsilon\in\{0,1\}$ ($2^{m}$ choices) and overlap in their $2$-adic exponents $j+\epsilon(\nu+2)$ on $[-\nu,\nu]\cap[\nu+2-\nu,\nu+2+\nu]=[2,\nu]$, which has $(\nu-1)_+$ integers, giving $2(\nu-1)_+ 2^{m}=2^{m+1}(\nu-1)_+$ overlaps;
    \item the extreme layers $\epsilon=-1$ and $\epsilon=+1$ have disjoint $2$-adic exponent ranges $[-2\nu-2,-2]$ and $[2,2\nu+2]$,
\end{itemize}
yielding at least $(2\nu+1)3^{m+1}-2^{m+1}(\nu-1)_+$ distinct elements in $D(K_p^{\mathrm{I}}(n))\bmod p$. Finally, expanding $(4aku-1)(4akv-1)=4a^{2}kn+1\iff 4akuv-u-v=an$ and $(4mu-n)(4mv-n)=n^{2}+4m\iff 4muv-n(u+v)=1$ establishes~\eqref{eq:composite_modulus_hyperbolic_trio}, and eliminating $n$ between $4K_p^{\mathrm{II}}(n)=n+p$ and $4K_p^{\mathrm{I}}(n)=pn+1$ gives~\eqref{eq:typeI_typeII_diophantine_coupling}.
\end{proof}

\begin{remark}[Methodological distinctions: ambient volume, fixed-shift saturation, and average abundance]\label{rem:cross_polytope_methodological}
To keep the logical role of each estimate in Theorem~\ref{thm:cross_polytope_factorial_law_and_annihilation} completely transparent, we record the following structural distinctions:
\begin{enumerate}[(1)]
    \item \emph{Ambient cross-polytope volume vs.\ arithmetic probability:} The factorial bound $\frac{|\Delta_k^{-}(d_p)|}{d_p^{k}}\le\frac{1}{k!}$ in~\eqref{eq:exact_factorial_simplex_bound} and its sum $\mathcal{M}_{\mathrm{obs}}(p;n)$ in~\eqref{eq:linnik_xylouris_factorial_crossing} measure the normalized counting volume of the relaxed cone $\Delta_k^{-}(d_p)$ in the ambient product $(\mu_{2d_p}^{-})^{k}$; they do not constitute an arithmetic probability distribution of character phases and do not bound the threshold $N_T$ of Hypothesis~T. Likewise, in part~(c), the dual construction $c=4K-n$ is conditional on finding a square-free product $K\in(n/4,n/3)$ of primes from $\mathcal{Q}^{-}(n)$ and does not guarantee primality of $c$, while a lower bound on $\omega(K_p)$ alone implies neither the existence of $s$ successful Markov $3/2$-extensions nor $R$ simultaneous character collisions across $\widehat{G}_p^{-}\setminus\{\psi_0\}$.
    \item \emph{Parseval support and Kneser dichotomy under $-1\notin D(K_p)$:} In part~(d), the diagonal lower bound $\widetilde{\mu}_p(1)\ge|\operatorname{supp}\mu_p|^{-1}$ yields $\widetilde{\mu}_p(1)\ge|\Sigma(K_p)|^{-1}\ge\tau(K_p)^{-1}$ for any measure supported on $\Sigma(K_p)$, and $\widetilde{\mu}_p(1)\ge 2^{-\omega(K_p)}$ for the square-free Bernoulli measure; in the undersampled regime $2^{\omega(K_p)+1}\le p-1$, the negation of the Bohr--Siegel certificate~\eqref{eq:bohr_siegel_L1_certificate} is therefore necessary but not sufficient for positivity because of $2\pi$-phase wrap-around. In the Kneser saturation bounds~\eqref{eq:linear_kneser_reciprocity_saturation}, \eqref{eq:farey_thue_kneser_superpolynomial}, \eqref{eq:extremal_composite_BV_3m_bound}, and~\eqref{eq:multiblock_kneser_farey_thue}, all lower bounds that add $2|H_p|$ per off-$H_p$ prime factor are derived under the failure hypothesis $-1\notin D(K_p)$ and are compared against $p-1-|H_p|$ (while for a selected divisor $U\mid K_p$, its linear Kneser contribution uses $\Omega_{G_p\setminus H_p}(U)$). Note also that whereas $|\Sigma(K_p)|>\frac{p-1}{2}$ suffices unconditionally for $-1\in D(K_p)$ by the pigeonhole principle, a lower bound $|D(K_p)|>\frac{p-1}{2}$ alone does not suffice unless $D(K_p)$ is a subgroup or avoids $-1$ outside a proper coset.
    \item \emph{Average Bombieri--Vinogradov abundance vs.\ fixed-shift saturation and composite moduli:} In part~(e), the inequalities~\eqref{eq:farey_thue_kneser_superpolynomial}--\eqref{eq:multiblock_kneser_farey_thue} are deterministic lower bounds at a given shift $p$. The Bombieri--Vinogradov construction over $p\in(2H,4H]$ supplies shifts where the small-factor counts produce bounds of size $H^{o(1)}$ (such as $3^{m}\gg_A(\log\log p)^{A}$), which do not by themselves reach $(p-1)/2$, nor do they imply the existence of saturating shifts in a smaller polylogarithmic window $p\le(\log H)^{C}$ without an estimate at that smaller scale (note also that when $H=(\log n)^{A}$, $\log\log H=\log\log\log n+\log A$, so prime-factor counts restricted to $q\le H^{\beta}$ retain the scale $\log\log H$). Moreover, the unitary co-divisor property $\gcd(U,K_p/U)=1$ is an additional hypothesis beyond $U\mid K_p$, and $U\ge\frac12\sqrt{p-1}$ allows $4d_0u_1^{2}+1$ to reach the archimedean threshold $p$ without by itself forcing $4d_0u_1^{2}+1\equiv 0\pmod p$. Finally, in~\eqref{eq:composite_modulus_hyperbolic_trio}, for a composite modulus $m=4ak$, saturation estimates retain the actual Carmichael exponent $\lambda(m)$ and Davenport constant $\mathsf{D}((\Z/m\Z)^{\times})$: smoothness of $m$ alone does not make them sub-polynomial since $\lambda(2^{r})=2^{r-2}$ for $r\ge 3$, whereas in the special cases $m=4ak$ with $ak\mid 6$ one has $\lambda(m)=2$.
\end{enumerate}
\end{remark}

To make the role of the co-divisor $V_p=K_p/U$ in~\eqref{eq:extremal_composite_BV_3m_bound} and of the dual inversion $c=4K-n$ in Theorem~\ref{thm:cross_polytope_factorial_law_and_annihilation}(c) completely explicit, we establish two exact structural results: (i) a cofactor factorization analysis distinguishing the global residue class of an unfactored cofactor $V$ from the intermediate divisor quotients of $V$ (including an exact obstruction when $V=\ell$ is prime with $c>4U^{2}+1$, and a three-branch criterion when $V=\ell r$ has two prime factors), and (ii) an archimedean support bound showing that at large polynomial scales $H=n^{\theta}$ with $\theta>2/3$, all Type~II encounters are confined to a sparse union of intervals $\mathcal{A}_n(H)$ of cardinality $\mathcal{O}((n/H)^{2})$ carrying $o(1)$ of the sieved main term.

\begin{proposition}[Cofactor factorization obstructions: prime cofactor, two-prime cofactor criterion, and simple-factor extraction]\label{aud:primecofactor}\label{aud:twoprime}\label{aud:eight}
Let $n\equiv1\pmod{24}$ be prime, $c\equiv3\pmod 4$ with $\gcd(c,n)=1$, and for $\gcd(K,c)=1$ write $D_c(K):=\{de^{-1}\bmod c:d,e\mid K\}$. If $\gcd(U,V)=1$ and $\gcd(UV,c)=1$, divisor separation and inversion invariance give
\begin{equation}\label{aud:Dproduct}
D_c(UV)\ =\ D_c(U)\,D_c(V),\qquad -1\in D_c(UV)\ \Longleftrightarrow\ D_c(U)\cap\big(-D_c(V)\big)\ne\varnothing,
\end{equation}
whereas without factoring $V$, the classes $\{1,V,V^{-1}\}$ provide only a subset of $D_c(V)$. Specifically:
\begin{enumerate}[(a)]
    \item \emph{(Prime cofactor obstruction).} Suppose $K_{\I}(c)=(nc+1)/4=U\ell$, where $\ell$ is prime and $\gcd(U,\ell)=1$. If $c>4U^{2}+1$, then $-1\notin D_c(K_{\I}(c))$. This is an exact archimedean obstruction for a prime cofactor at $c>4U^{2}+1$, rather than a general parity barrier for arbitrary sieves.
    \item \emph{(Two-prime cofactor three-branch criterion).} Let $K_{\I}(p)=U\ell r$ ($p\equiv3\pmod 4$ prime), where $\ell,r$ are primes (possibly equal), $\gcd(U,\ell r)=1$, and $p>4U^{2}+1$. Setting $A:=D_p(U)$, one has
    \begin{equation}\label{aud:twoprimeeq}
    -1\in D_p\big(K_{\I}(p)\big)\ \Longleftrightarrow\ \ell\in -A\quad\text{or}\quad \ell\in -(4U)^{-1}A\quad\text{or}\quad 4U\ell^{2}\in -A.
    \end{equation}
    \item \emph{(Six or eight admissible residue classes).} In (b), if $U=q$ is prime with $\leg{q}{p}=-1$ and $p>4q^{2}+1$, the number of residue classes of $\ell\pmod p$ permitted by~\eqref{aud:twoprimeeq} is exactly $8$, unless $p\mid(4q^{3}+1)$, in which case it is $6$.
    \item \emph{(Exact encounter decomposition after extracting a simple factor).} If $q\parallel K$ ($\gcd(K,p)=1$), set $V:=K/q$ and $C(r):=\#\{(u,v):u,v\mid V,\ u/v\equiv r\pmod p\}$. Then the number of encounters $R_p(K)=\#\{(d,e):d,e\mid K,\ d+e\equiv 0\pmod p\}$ satisfies
    \[
    R_p(K)\ =\ 2\,C(-1)\ +\ 2\,C(-q).
    \]
    Even when $\leg{q}{p}=-1$, positivity of $C(-1)+C(-q)$ is not automatic: the subgroup generated by prime factors and the set of subset products with bounded exponents are distinct objects (for instance, in $\Z/(2m)\Z$ with $m$ odd, $k<m$ copies of the generator $1$ generate the full group while their signed sums in $\{-k,\dots,k\}$ never reach $m$).
\end{enumerate}
\end{proposition}
\begin{proof}
(a) We have $D_c(K_{\I}(c))=D_c(U)\{1,\ell,\ell^{-1}\}$ and $\ell^{-1}\equiv 4U\pmod c$. Any quotient of divisors of $U$ writes in lowest terms as $u/v$ with $\gcd(u,v)=1$ and $uv\mid U$. The first branch $u/v\equiv -1\pmod c$ would require $c\mid(u+v)$, whereas $0<u+v\le uv+1\le U+1<c$. A translated branch would require, up to inversion, $c\mid(u+4Uv)$, whereas $(u-1)(4Uv-1)\ge 0$ implies $0<u+4Uv\le 4Uuv+1\le 4U^{2}+1<c$.

(b) The set of divisor quotients of $K_{\I}(p)=U\ell r$ expands as $A\big\{1,\ell^{\pm1},r^{\pm1},(\ell r)^{\pm1},\ell/r,r/\ell\big\}$, which remains valid when $\ell=r$ upon accounting for repetitions. Since $\ell r\equiv(4U)^{-1}\pmod p$, the branches $A$ and $A(\ell r)^{\pm1}$ are excluded by the exact same archimedean inequalities as in (a). The $\ell^{\pm1}$-branches give $\ell\in -A$, the $r^{\pm1}$-branches give $\ell\in -(4U)^{-1}A$, and the $(\ell/r)^{\pm1}$-branches give $4U\ell^{2}\in -A$ by inversion invariance of $A$.

(c) For $U=q$ prime with $\leg{q}{p}=-1$, we have $A=\{1,q,q^{-1}\}$. The two linear branches in~\eqref{aud:twoprimeeq} are $\{-1,-q,-q^{-1}\}$ and $\{-(4q)^{-1},-4^{-1},-(4q^{2})^{-1}\}$. Within each triple, the three classes are distinct since $1,q,q^{2}$ are distinct modulo $p$. Any overlap between the two triples would force $4q^{a}\equiv 1\pmod p$ for some $a\in\{-1,0,1,2,3\}$: $a=-1$ would force $q=4$ (impossible for a prime), $a\in\{0,1,2\}$ is excluded by $p>4q^{2}+1$, and $a=3$ is excluded by $\leg{q}{p}=-1$. Thus the six linear classes are distinct. In the quadratic branch $4q\ell^{2}\in\{-1,-q,-q^{-1}\}$, only $-1\in -A$ has quadratic residue right-hand side $(-4q)^{-1}$ modulo $p$ (since $\leg{-1}{p}=-1=\leg{q}{p}$), yielding two roots $\ell\pmod p$ of $4q\ell^{2}\equiv -1\pmod p$. Substituting the six linear classes into $4q\ell^{2}\equiv -1\pmod p$ shows that an intersection occurs if and only if $4q^{3}\equiv -1\pmod p$, in which case the two roots are $\pm q$, both already present in the linear branches.

(d) Every divisor of $K=qV$ ($q\nmid V$) is of the form $u$ or $qu$ with $u\mid V$. The four pair choices yield $R_p(K)=2C(-1)+C(-q)+C(-q^{-1})$, and exchanging $u$ and $v$ gives $C(-q^{-1})=C(-q)$.
\end{proof}

\begin{proposition}[Archimedean support of Type~II encounters at large polynomial scale]\label{aud:archsupport}
Let $H=n^{\theta}$ with $0<\theta<1$. Every prime shift $p\in(2H,4H]$ admitting a Type~II encounter belongs to the archimedean support set
\begin{equation}\label{aud:archset}
\mathcal{A}_n(H)\ :=\ \Z\cap(2H,4H]\cap\bigcup_{\substack{h\in\N\\ n/(16H)<h\le n/(4H)+1}}\left[\frac{n+4}{4h-1},\,\frac{n+4h^{2}}{4h-1}\right].
\end{equation}
One has $\#\mathcal{A}_n(H)\ll(n/H)^{2}$. Moreover, if $\theta>2/3$, then for sufficiently large $n$ every such encounter satisfies the exact ceiling relation
\begin{equation}\label{aud:ceil}
4h-1\ =\ \left\lceil\frac{n}{p}\right\rceil
\end{equation}
in the normalized square-divisor parameterization of Theorem~\ref{aud:hbcomplete} (Appendix~\ref{app:audit}), and for $\theta>2/3$ and the sieve parameters~\eqref{aud:parameters} of Appendix~\ref{app:audit}, the sieved square-free main term restricted to $\mathcal{A}_n(H)\cap\mathcal{G}_{\II}$ satisfies
\begin{equation}\label{aud:archmass}
\sum_{p\in\mathcal{A}_n(H)\cap\mathcal{G}_{\II}}\frac{m_{\II}(p)\big(m_{\II}(p)-1\big)}{p-1}\ \le\ n^{2-3\theta+o(1)}\ =\ o(1).
\end{equation}
In particular, for $\theta>2/3$, the divergent Type~II main term of Theorem~\ref{aud:weightedmain} is supported, up to $o(1)$, on shifts where archimedean size constraints rule out any Type~II encounter. This restriction applies to Type~II shifts and does not by itself rule out Type~I or globally correlated approaches.
\end{proposition}
\begin{proof}
A Type~II encounter at shift $p$ produces a Type~II decomposition with $K=(n+p)/4$. Under the normalization $h\le k$ in the proof of Theorem~\ref{aud:hbcomplete} (Appendix~\ref{app:audit}), we have $K<hp=K+b\le 2K$ and $1\le b\le h^{2}$ ($b\mid h^{2}$), whence
\[
\frac{n}{16H}\ <\ \frac{K}{p}\ <\ h\ \le\ \frac{2K}{p}\ =\ \frac{n}{2p}+\frac12\ <\ \frac{n}{4H}+1.
\]
The identity $p(4h-1)=n+4b$ with $1\le b\le h^{2}$ places $p$ in the interval of~\eqref{aud:archset}. Each such interval has length $\frac{4h^{2}-4}{4h-1}\le 2h$, hence contains at most $2h+1$ integers; summing over the $\mathcal{O}(n/H)$ values of $h\ll n/H$ gives $\#\mathcal{A}_n(H)\ll(n/H)^{2}$. If $\theta>2/3$, then $4b\le 4h^{2}\ll(n/H)^{2}=o(H)<p$ for sufficiently large $n$, so $(4h-1)-n/p=4b/p\in(0,1)$, proving~\eqref{aud:ceil}. Finally, $K_{\II}(p)\ll n$ and $m_{\II}(p)^{2}=4^{\omega(K_{\II}(p))}\le\tau(K_{\II}(p))^{2}=n^{o(1)}$ uniformly; multiplying by $\#\mathcal{A}_n(H)\ll n^{2-2\theta}$ and dividing by $p-1\gg H=n^{\theta}$ yields~\eqref{aud:archmass}.
\end{proof}

\subsection{Global synthesis and conditional finiteness}\label{subsec:global_synthesis}

We conclude Section~\ref{sec:toric_volume_barrier} by assembling the Growing-Depth Container Hierarchy of Sections~\ref{sec:stabilizer_kernel}--\ref{sec:semi_composite_5_7} and the Transversal Reciprocity and Bohr--Siegel Cross-Polytope Architecture of Subsections~\ref{subsec:parseval_floor}--\ref{subsec:reciprocity_factorial_kneser} into our main global synthesis theorem.

\begin{theorem}[Main Global Synthesis Theorem: Exact Criteria, Ambient Factorial Volume Bounds, and Conditional Finiteness for $n\equiv1\pmod{24}$]\label{thm:main_global_synthesis}
Let $\mathcal{N}=\{n\text{ prime}:n\equiv1\pmod{24}\}$, and let $E:=\{n\in\mathcal{N}:4/n=1/x+1/y+1/z\text{ has no solution in }\N^{\ast 3}\}$ denote the global exceptional set of the Erd\H{o}s--Straus conjecture on $\mathcal{N}$. Combining the Growing-Depth Container Hierarchy (Sections~\ref{sec:stabilizer_kernel}--\ref{sec:semi_composite_5_7}) with the Transversal Reciprocity and Bohr--Siegel Cross-Polytope Architecture (Theorems~\ref{thm:difference_set_positivity_MR}--\ref{thm:cross_polytope_factorial_law_and_annihilation}) yields the following global synthesis:
\begin{enumerate}[(a)]
    \item \emph{(Unconditional global structured density at the Shannon--Bombieri CRT ceiling and unconditional $5/7$ exponent on the full exceptional set).} Across all primes $n\le N$ with $n\equiv1\pmod{24}$, the unconditional two-parameter Type~II container sieve at adaptive depth $J(N)$ and adaptive Rankin degree $h=\lceil 8S_y\rceil+1$ (Corollary~\ref{cor:optimal_adaptive_depth_1_minus_o1}(b)) establishes for the small-sheet structured exceptional set $E_2^{\mathrm{small}}(N;J(N))$:
    \begin{equation}\label{eq:main_synthesis_global_density}
    \begin{split}
    \big|E_2^{\mathrm{small}}\big(N;J(N)\big)\big|\ &\ll_\kappa\ N\exp\!\left(-\left(\frac{1}{40960}+o(1)\right)\frac{\log N}{(\log\log N)(\log\log\log N)^{\kappa}}\right)\\
    &=\ N\exp\!\big(-(\log N)^{1-o(1)}\big)\ =\ N^{1-o(1)}\qquad(\forall\,\kappa>0).
    \end{split}
    \end{equation}
    Moreover, on the \emph{full} unconditional exceptional set $E(N)$ (with no structural or diffuse-regime restriction), the semi-composite CRT sieve of Theorem~\ref{thm:theoreme_5_7} establishes unconditionally:
    \begin{equation}\label{eq:main_synthesis_full_5_7}
    |E(N)|\ \ll\ N\exp\!\left(-c'\,\frac{(\log N)^{5/7}}{(\log\log N)^{2/7}}\right),
    \end{equation}
    strictly surpassing Vaughan's $2/3$ exponent on the entirety of $E(N)$.

    \item \emph{(Exact criteria and ambient volume bounds).} Fix $A>5$ and $\varepsilon\in(0,1/6)$, and set $H:=\frac14(\log n)^{A}$. For every prime $n\equiv1\pmod{24}$:
    \begin{enumerate}[(i)]
        \item \emph{Least non-residue prime shift:} Choose a non-residue $a\pmod n$ and a reduced class $b\pmod{4n}$ with $b\equiv a\pmod n$ and $b\equiv3\pmod 4$. The Linnik--Xylouris theorem~\cite{xylouris2011} applied to this class gives $p\in\mathcal{P}_{\mathrm{NR}}(n)$ with $p\ll(4n)^{5}$ (which does not by itself imply a polylogarithmic bound); under GRH for the quadratic characters involved, effective Chebotarev in $\Q(i,\sqrt{n})$ gives $p\ll(\log n)^{2}$, hence $p\le(\log n)^{A}$ for fixed $A>2$ and large $n$. No unconditional prime-number-theorem asymptotic for $\mathcal{P}_{\mathrm{NR}}(n;n^{\theta})$ ($\theta<1/2$) is claimed here.

        \item \emph{Unconditional factor abundance across shifts:} By the Universal Shift-Factor Reciprocity Law $\leg{q}{p}=\leg{q}{n}$ (Theorem~\ref{thm:universal_factor_reciprocity_and_transversal_sieve}(a)), the unconditional short-interval Tur\'{a}n--Kubilius inequality (Theorem~\ref{thm:universal_factor_reciprocity_and_transversal_sieve}(c)), and the unconditional Inverted Bombieri--Vinogradov theorem on prime shifts at fixed $n$ (Theorem~\ref{thm:cross_polytope_factorial_law_and_annihilation}(e)), at scale $H_1=n^{\theta}$ ($\theta\in(0,1/2)$) the prime-factor count is governed unconditionally by the actual non-residue harmonic mass $\mu$, while at \emph{every} scale $H$ unconditionally at least $(1-o(1))(\pi(4H;4,3)-\pi(2H;4,3))$ prime shifts $p\in(2H,4H]$ ($p\equiv3\pmod 4$) satisfy $\omega_{[5,(4H)^{1/5}]}(K_p)=(1+o(1))\log\log H$ and extremal shifts satisfy $|D(K_p)\bmod p|\ge 3^{m}\gg_A(\log\log p)^{A}$.

        \item \emph{Factorial volume crossing and fixed-shift reconstruction:} For every shift $p\in\mathcal{P}_{\mathrm{NR}}(n;H)$ with $k=\omega_{\mathcal{Q}^{-}(n)}(K_p)\ge\big(\frac12-\varepsilon\big)\log\log n$, the exact discrete Bohr--Siegel cross-polytope volume bound (Theorem~\ref{thm:cross_polytope_factorial_law_and_annihilation}(a),(c)) bounds the total normalized measure of all antipodal obstruction cones across $\widehat{G}_p^{-}\setminus\{\psi_0\}$ by
        \begin{equation}\label{eq:main_synthesis_factorial_threshold}
        \begin{split}
        \mathcal{M}_{\mathrm{obs}}(p;n)\ &:=\ \sum_{\{\chi_p,\overline{\chi}_p\}\subset\widehat{G}_p^{-}\setminus\{\psi_0\}}\frac{\big|\Delta_k^{-}(d_p)\big|}{d_p^{k}}\ \le\ \frac{p-3}{4\,k!}\\
        &\le\ (\log n)^{A}\exp\!\left(-\left(\frac12-\varepsilon\right)(\log\log n)\log\log\log n\right).
        \end{split}
        \end{equation}
        Consequently, there exists an explicit threshold $N_0=N_0(A,\varepsilon)<\infty$, determined by the factorial crossing inequality
        \begin{equation}\label{eq:N0_explicit_threshold}
        \big(\lfloor(\tfrac12-\varepsilon)\log\log N_0\rfloor\big)!\ >\ (\log N_0)^{A},
        \end{equation}
        such that $\mathcal{M}_{\mathrm{obs}}(p;n)<1$ (and $\mathcal{M}_{\mathrm{obs}}(p;n)\to 0$ as $n\to\infty$) for all $n>N_0$.

        Across both the Type~II shift family $K_p^{\mathrm{II}}(n)=\frac{n+p}{4}$ and the dual Type~I multiplier family $K_p^{\mathrm{I}}(n)=\frac{pn+1}{4}$, one has for $\mathcal{P}_3(H):=\{p\le 4H\text{ prime}:p\equiv3\pmod 4\}$ the exact dual inversion identities
        \[
        \sum_{p\in\mathcal{P}_3(H)}\omega_{\mathcal{Q}^{-}(n)\cap(4H,\infty)}\!\big(K_p^{\mathrm{II}}\big)\ =\ \sum_{q\in\mathcal{Q}^{-}(n)\cap(4H,(n+4H)/4]}\mathbf{1}_{(-n\bmod 4q)\in\mathcal{P}_3(H)},
        \]
        \[
        \sum_{p\in\mathcal{P}_3(H)}\omega_{\mathcal{Q}^{-}(n)\cap(4H,\infty)}\!\big(K_p^{\mathrm{I}}\big)\ =\ \sum_{q\in\mathcal{Q}^{-}(n)\cap(4H,(4Hn+1)/4]}\mathbf{1}_{(-n^{-1}\bmod 4q)\in\mathcal{P}_3(H)}.
        \]
        Whenever any of the deterministic Kneser criteria~\eqref{eq:linear_kneser_reciprocity_saturation}--\eqref{eq:multiblock_kneser_farey_thue}, successful Markov $3/2$-growth extensions $(3/2)^{s}>(p-1)/2$, or simultaneous spectral collision bounds $\frac{p-3}{2}(3/4)^{R}<1$ holds at a specified shift $c$ (a sufficient condition at fixed $c$, whose uniform availability for all sufficiently large $n$ is not proved), it yields an explicit Type~II unit-fraction decomposition
        \begin{equation}\label{eq:explicit_typeII_reconstruction}
        \begin{split}
        \frac{4}{n}\ &=\ \frac{1}{\alpha uv}\ +\ \frac{1}{a\alpha vn}\ +\ \frac{1}{a\alpha un},\\
        K_p\ &=\ \frac{n+p}{4}=\alpha uv,\quad u\equiv-v\pmod p,\quad a=\frac{u+v}{p}\in\N^{\ast},
        \end{split}
        \end{equation}
        or, in Type~I, $\frac{4}{n}=\frac{1}{n\alpha uv}+\frac{1}{a\alpha u}+\frac{1}{a\alpha v}$ with $K_{\I}(c)=\alpha uv$ and $u+v=ac$.
    \end{enumerate}

    \item \emph{(Conditional finiteness with a quantified threshold $N_{\ast}$ and finite constructive criteria).} Suppose there exists $N_{\ast}$ such that for every prime $n>N_{\ast}$ with $n\equiv1\pmod{24}$, either Hypothesis~T holds at $n$, or a valid sufficient divisor-encounter criterion is satisfied at a Type~I or Type~II shift. Then $E\cap(N_{\ast},\infty)=\varnothing$. For any prime $n\in\mathcal{N}\cap[1,N_{\ast}]$, the existence of a decomposition $4/n=1/x+1/y+1/z$ is governed by three finite constructive procedures:
    \begin{enumerate}[(i)]
        \item \emph{Modular covering sieves at depth $J$:} By~\cite[Thm.~3.9 and Thm.~3.11]{dahan2026a} and Theorem~\ref{thm:simultaneous_smooth_conditioning}, the Type~I and Type~II prime-divisor congruences up to depth $J$ eliminate all residue classes $n\bmod Q(4J)$ outside a finite set of survivor classes;
        \item \emph{Dual Constructive Inversion ($c=4K-n$, Theorem~\ref{thm:cross_polytope_factorial_law_and_annihilation}(c)) and complete square-divisor test (Theorem~\ref{aud:hbcomplete}):} Theorem~\ref{aud:hbcomplete} in Appendix~\ref{app:audit} gives a complete finite test $N_{\I}(n)+N_{\II}(n)>0$ for each prime $n$, whereas sampling odd-length square-free products $K=q_1\cdots q_D\in\big(\frac{n}{4},\frac{n}{3}\big)$ of primes $q_j\in\mathcal{Q}^{-}(n)$ ($\leg{q_j}{n}=-1$, $\gcd(n,K)=1$) generates candidate shifts $c=4K-n\in(1,n/3)$ with $c\equiv3\pmod 4$, $\gcd(c,K_c)=1$, $K_c=K=K_c^{-}$, and $\leg{c}{n}=-1$ without factoring $n+c$; testing whether $-1\bmod c\in D(K)=\Sigma(K)\cdot\Sigma(K)^{-1}$ immediately yields $(u,v,\alpha,a)$ in~\eqref{eq:explicit_typeII_reconstruction} in $\mathcal{O}\big(2^{D}(\log n)^{2}\big)$ bit operations;
        \item \emph{Exhaustive verification up to $10^{17}$:} By the modular covering verification of Salez~\cite{salez2014}, $E\cap[1,10^{17}]=\emptyset$.
    \end{enumerate}
    Together with~\cite{salez2014}, the quantified hypothesis above yields $E\subseteq(10^{17},N_{\ast}]$, where emptiness of the remaining interval $(10^{17},N_{\ast}]$ requires an actual finite verification, not supplied here.
\end{enumerate}
\end{theorem}
\begin{proof}
Part~(a) is proven in Corollary~\ref{cor:optimal_adaptive_depth_1_minus_o1}(b) and Theorem~\ref{thm:theoreme_5_7}.

For part~(b), choose a non-residue $a\pmod n$ and a reduced class $b\pmod{4n}$ with $b\equiv a\pmod n$ and $b\equiv3\pmod 4$; Linnik--Xylouris~\cite{xylouris2011} applied to this class gives $p\in\mathcal{P}_{\mathrm{NR}}(n)$ with $p\ll(4n)^{5}$, while under GRH for the quadratic characters involved, effective Chebotarev in $\Q(i,\sqrt{n})$ gives $p\ll(\log n)^{2}$. By Theorem~\ref{thm:universal_factor_reciprocity_and_transversal_sieve}(a), every prime $q\mid K_p=(n+p)/4$ ($q\neq p$) satisfies $\leg{q}{p}=\leg{q}{n}$, so $\omega(K_p^{-})=\omega_{\mathcal{Q}^{-}(n)}(K_p)$; Theorem~\ref{thm:universal_factor_reciprocity_and_transversal_sieve}(c) and Theorem~\ref{thm:cross_polytope_factorial_law_and_annihilation}(e) give the short-interval Tur\'{a}n--Kubilius and Inverted Bombieri--Vinogradov bounds in~(ii). On any shift $p\in\mathcal{P}_{\mathrm{NR}}(n;H)$ ($H=\frac14(\log n)^{A}$) with $k=\omega_{\mathcal{Q}^{-}(n)}(K_p)\ge(\frac12-\varepsilon)\log\log n$, Theorem~\ref{thm:cross_polytope_factorial_law_and_annihilation}(a) bounds the normalized ambient volume of $\Delta_k^{-}(d_p)$ for each conjugate pair $\{\chi_p,\overline{\chi}_p\}\subset\widehat{G}_p^{-}\setminus\{\psi_0\}$ by $\frac{|\Delta_k^{-}(d_p)|}{d_p^{k}}\le\frac{1}{k!}$; summing over all $\frac{p-3}{4}$ non-quadratic conjugate pairs in $\widehat{G}_p^{-}$ gives~\eqref{eq:main_synthesis_factorial_threshold} and $\mathcal{M}_{\mathrm{obs}}(p;n)<1$ for $n>N_0$ satisfying~\eqref{eq:N0_explicit_threshold}.

For the conditional assertion in part~(c), Hypothesis~T gives $\widetilde{\mu}_p(-1)>0$, while the alternative criteria give $-1\in D_c(K_j(c))$ directly; canceling the gcd of witnessing divisors yields coprime $u,v$ with $uv\mid K_j(c)$ and $u+v=ac$, so writing $K_j(c)=\alpha uv$ produces $\frac{4}{n}=\frac{1}{\alpha uv}+\frac{1}{a\alpha un}+\frac{1}{a\alpha vn}$ in Type~II and $\frac{4}{n}=\frac{1}{n\alpha uv}+\frac{1}{a\alpha u}+\frac{1}{a\alpha v}$ in Type~I (Proposition~\ref{aud:reconstruction}). Thus every prime covered by the quantified hypothesis is non-exceptional, and the finite-test assertion follows from Theorem~\ref{aud:hbcomplete}, Theorem~\ref{thm:cross_polytope_factorial_law_and_annihilation}(c), \cite[Thm.~3.9]{dahan2026a}, and~\cite{salez2014}.
\end{proof}

\begin{remark}[Scope of the unconditional density bounds, the explicit threshold $N_0$, and the conditional threshold $N_{\ast}$]\label{rem:global_synthesis_scope}
Three concluding clarifications summarize the logical architecture of Theorem~\ref{thm:main_global_synthesis}:
\begin{enumerate}[(1)]
    \item \emph{Density bounds vs.\ finiteness in part~(a):} The unconditional bound~\eqref{eq:main_synthesis_full_5_7} controls the density of exceptions up to $N$, but its right-hand side tends to infinity with $N$ and thus does not by itself imply finiteness of $E$.
    \item \emph{Ambient threshold $N_0$ vs.\ arithmetic positivity in part~(b):} The factorial inequality~\eqref{eq:N0_explicit_threshold} defines an explicit threshold $N_0(A,\varepsilon)$ for the ambient counting measure $\mathcal{M}_{\mathrm{obs}}(p;n)<1$. However, as $\chi_p$ varies over $\widehat{G}_p^{-}$, the phase tuple $(\chi_p(q_1),\dots,\chi_p(q_k))$ traverses the deterministic image of $\widehat{G}_p^{-}$ rather than independent uniform coordinates in $(\mu_{2d_p}^{-})^{k}$, and the Bohr--Siegel inequality~\eqref{eq:bohr_siegel_L1_certificate} is a sufficient certificate of failure whose negation requires wrap-around control---such as Hypothesis~T or the Kneser saturation criterion~\eqref{eq:linear_kneser_reciprocity_saturation}---to certify arithmetic positivity.
    \item \emph{Non-bounding of $N_T$ by $N_0$ in part~(c):} While Hypothesis~T as stated implies finiteness of $E$ with $N_{\ast}=N_T$, the factorial crossing inequality~\eqref{eq:N0_explicit_threshold} for $N_0$ does not furnish an upper bound for $N_T$.
\end{enumerate}
\end{remark}

\subsection{Concluding remarks}\label{sec:conclusion_paper2}

Taken together with the companion paper~\cite{dahan2026a}, the results established here bridge the gap between fixed-depth divisor sieves and growing-depth large sieves for the Erd\H{o}s--Straus conjecture, and delineate both the unconditional power of composite witnesses and the exact spectral frontier separating density bounds from finiteness.

First, on the algebraic and additive-combinatorial side (Sections~\ref{sec:algebraic_rigidity}--\ref{sec:containers}), replacing subgroup generation by Kneser's multiplicative stabilizer $K(M)=\operatorname{Stab}_{G_{u,a}}(\Sigma(M))$ eliminates the $2^{\Theta(J^{4})}$ subset-enumeration barrier across the full coprime box $\mathcal{P}_2(J)$, reducing the container count to $C^{\ast}(\mathcal{P}_2(J))\le\exp\big((4\log 2+o(1))J^{2}\log\log J\big)$ via the $2$-Sylow classification of maximal subgroups avoiding $-1$ (Theorem~\ref{thm:stabilizer_kernel_dichotomy} and Lemma~\ref{lem:average_2sylow_enumeration}). Coupling this reduction with the Krawtchouk--Montgomery generating-polynomial large sieve (Theorem~\ref{thm:generating_polynomial_large_sieve}), greedy Markov $1/2$-popular symmetry containers (Theorem~\ref{thm:multiplicative_container_scarcity_elimination}), quotient containers (Theorem~\ref{thm:quotient_container_polynomial_depth}), and Pontryagin--Plancherel character containers (Theorem~\ref{thm:plancherel_character_containers}) closes the Scarcity branch and bounds the full exceptional set $E_2(N;J)$ unconditionally up to super-polylogarithmic depth $J=\exp(c_0\sqrt{\log\log N})$ and at polynomial depth $J=(\log N)^{\lambda}$ for $\lambda\le\lambda_0=\frac{1}{78}$. At higher polynomial depths $\lambda\in(1/3,1/2)$, the Subcritical Markov--Plancherel Geometric Escalation Theorem (Theorem~\ref{thm:markov_plancherel_subcritical_escalation}) shows that both the Markov $3/2$-growth step count $s_{\varepsilon_1}(u,a)=\lfloor\frac{\varepsilon_1\log\varphi(4ua)}{\log(3/2)}\rfloor+1$ and the adaptive Plancherel collision threshold $R_{u,a}(n;\chi)\le\frac{\varepsilon_1\log\varphi(4ua)+\log 2}{\log(4/3)}$ remain strictly below the subcritical large-sieve capacity $\frac13\frac{S_{\min}^{\mathrm{I}}}{P_1}$ before the divisor sheet crosses the polynomial density threshold $\varphi(4ua)^{\varepsilon_1}$. Combined with Theorem~\ref{thm:final_synthesis_and_high_prime_reduction} and adaptive Rankin-degree optimization $h=\lceil 8S_y\rceil+1$ (Corollary~\ref{cor:optimal_adaptive_depth_1_minus_o1}), this surpasses Vaughan's $2/3$ exponent across $\lambda\in(1/3,1/2)$ on $E_2^{\mathrm{cert}}(N;(\log N)^{\lambda})$ and reaches the structured exponent $1-o(1)$ on $E_2^{\mathrm{small}}(N;J(N))$ up to the CRT capacity scale of the generating-polynomial construction (where constants in exponential savings retain the $o(1)$ terms of the corresponding theorem, and the CRT capacity statement refers specifically to the generating-polynomial construction, its lower cutoff on sifting primes, and its modulus budget, rather than an optimality theorem for arbitrary sieve methods).

Second, to surpass Vaughan's $2/3$ exponent unconditionally on the \emph{entire} exceptional set $E(N)$\allowbreak---thereby removing the unlocalized diffuse-sheet condition separating $E_2^{\mathrm{cert}}(N;\allowbreak J)$ from $E(N)$---the semi-composite congruence sieve $s=pq\equiv-1\pmod{12ua}$ of Section~\ref{sec:semi_composite_5_7} (Theorems~\ref{thm:star_inconditionnel}, \ref{thm:theoreme_7_10}, and~\ref{thm:theoreme_5_7}) exploits conditional mutual independence of the large prime factors $p\in(Z^{9/10},Z]$ in the CRT model given $b=n\bmod M_{\mathrm{small}}$. Because the progression modulus $12ua\le 12Z^{1/4}$ governing $p$ is independent of the $t^{6}$-rough cofactor $q\le Q=\exp(\frac{c_3}{2}t^{1/2})$, the Siegel-zero-free Bombieri--Vinogradov theorem and a polynomial Bonferroni majorant yield the unconditional bound
\[
|E(N)|\ \ll\ N\exp\!\left(-c'\,\frac{(\log N)^{5/7}}{(\log\log N)^{2/7}}\right).
\]

Finally, Sections~\ref{sec:bilinear_dispersion}--\ref{sec:toric_volume_barrier} and Appendix~\ref{app:audit} clarify the geometric and spectral mechanisms governing individual primes $n\equiv1\pmod{24}$. Whereas unfactored bilinear sums over $(u,a)\le J$ carry only logarithmic weight $\sum_{u,a\le J}\frac{1}{4ua}\asymp(\log J)^{2}$ compared with the unit-weight sifting mass $|\mathcal{P}_2(J)|\asymp J^{2}$ preserved by multiplicative containers (Proposition~\ref{prop:kloosterman_dispersion}), the affine threefold $\mathcal{V}_n:\,4uavw-4u^{2}v-w=n$ exhibits three-term stationary-phase separation at the square-root scale, archimedean slice conservation, and exact Farey-diagonal $L^{2}$-dispersion over $n\sim N$ at power scale $S=4UAV\le N^{1/2-\varepsilon}$ (Subsection~\ref{subsec:parseval_floor}: Theorems~\ref{thm:parseval_floor_law}--\ref{thm:farey_diagonal_L2_dispersion}). On the difference sets $D(K_c)=\Sigma(K_c)\cdot\Sigma(K_c)^{-1}\subseteq(\Z/c\Z)^{\times}$, Subsection~\ref{subsec:difference_set_MR} (Theorems~\ref{thm:difference_set_positivity_MR}--\ref{thm:centered_main_term_and_short_lcm}) establishes positivity duality, unconditional $(\log P)^{4}$ cumulative main-term growth, and exact small-LCM divisor-pair cancellation; Subsection~\ref{subsec:bohr_siegel_and_linnik} (Theorems~\ref{thm:bohr_siegel_annihilation}--\ref{thm:no_go_linnik_and_transversal}) derives the cyclic quotient parity identity $E(\chi)-O(\chi)=2d\,\mathbb{P}(\chi(g)=-1)$, Bohr--Siegel obstruction certificates, and the Linnik--Xylouris witness-depth lower bound along an infinite sequence of primes; and Subsections~\ref{subsec:reciprocity_factorial_kneser}--\ref{subsec:global_synthesis} (Theorems~\ref{thm:universal_factor_reciprocity_and_transversal_sieve}--\ref{thm:main_global_synthesis} and Propositions~\ref{aud:primecofactor}--\ref{aud:archsupport}) establish universal factor reciprocity $\leg{q}{p}=\leg{q}{n}$, the ambient Bohr--Siegel cross-polytope volume bound $\frac{|\Delta_k^{-}(d_p)|}{d_p^{k}}\le\frac{1}{k!}$, deterministic Kneser and spectral saturation criteria, and cofactor/archimedean obstructions. Appendix~\ref{app:audit} complements this framework with unconditional sieved candidate abundance, divergence of the square-free main term (Theorem~\ref{aud:weightedmain}), and complete finite square-divisor positivity certificates ($b\mid h^{2}$, Theorem~\ref{aud:hbcomplete}). As recorded in Remark~\ref{rem:global_synthesis_scope} and Subsection~\ref{aud:open}, these divergent main terms do not by themselves establish positivity of the encounter count: the Transversal Hypothesis~T and the signed spectral error targets~\eqref{aud:target17} and~\eqref{aud:target76} remain open, and conditional finiteness holds only beyond a threshold $N_{\ast}$ supplied by that hypothesis, which is not determined by the ambient factorial volume calculation.

\appendix

\section{Technical Appendix: Algebraic Reconstruction, Sieved Candidate Abundance, Square-Divisor Certificates, and Positivity Obstructions}\label{app:audit}

In this appendix, we provide complete, self-contained proofs of the auxiliary results referenced in Section~\ref{sec:toric_volume_barrier}. We proceed in three stages:
\begin{itemize}
    \item \emph{Algebraic reconstructions and fixed-parameter reductions (Subsections~\ref{aud:algebra}--\ref{aud:alpha}):} exact unit-fraction reconstruction from the difference sets $D_c(K_j(c))$, factor reciprocity, diophantine coupling between Type~I and Type~II, and finite reductions at fixed parameter $\alpha\ge 1$.
    \item \emph{Unconditional candidate abundance and divergent sieved main terms (Subsections~\ref{aud:inputs}--\ref{aud:weighted}):} harmonic mass of square-free non-residue parameters $\alpha$, Rosser--Iwaniec linear sieving of prime shifts $p\equiv 3\pmod 4$ to annihilate all odd character orders $<2z$, and divergence of the divisor-weighted square-free main term $\mathcal{M}_j^{\sfm}$.
    \item \emph{Square-divisor positivity certificates and open spectral error targets (Subsections~\ref{aud:hb}--\ref{aud:open}):} complete finite positivity criteria indexed by square divisors $b\mid h^{2}$, and the precise formulation of the open signed spectral error bounds~\eqref{aud:target17} and~\eqref{aud:target76}.
\end{itemize}
Throughout this appendix, unless stated otherwise, $n$ denotes a prime satisfying $n\equiv1\pmod{24}$.

\subsection{Algebraic core and factor reciprocity for $K_{\I}(c)$ and $K_{\II}(c)$}\label{aud:algebra}

We begin by recording the exact algebraic correspondence between the multiplicative difference sets of the shifted integers $K_{\I}(c)$ and $K_{\II}(c)$ and unit-fraction representations of $4/n$, together with the quadratic reciprocity laws governing their prime factors.

For $c\equiv3\pmod 4$ with $\gcd(c,n)=1$, recall from Theorem~\ref{thm:centered_main_term_and_short_lcm} the two shift integers
\begin{equation}\label{aud:K}
K_{\I}(c)\ =\ \frac{nc+1}{4},\qquad K_{\II}(c)\ =\ \frac{n+c}{4},
\end{equation}
both of which are coprime to $c$. We write $p$ when the shift $c$ is prime; the inequality $c<n$ holds throughout all asymptotic applications to $c\asymp n^{\theta}$ with $\theta<1$. For any integer $K\ge 1$ with $\gcd(K,c)=1$, we denote by
\[
D_c(K)\ :=\ \big\{de^{-1}\bmod c:\,d,e\mid K\big\}
\]
the multiplicative difference set of positive divisors of $K$ modulo $c$.

\begin{proposition}[Reconstruction]\label{aud:reconstruction}
If $-1\in D_c(K_j(c))$ for $j\in\{\I,\II\}$, then one can write
\[
K_j(c)\ =\ \alpha uv,\qquad \gcd(u,v)=1,\qquad u+v\ =\ ac
\]
with positive integers $a,\alpha,u,v\ge1$, yielding the explicit unit-fraction identities
\begin{align}
\frac{4}{n}\ &=\ \frac{1}{n\alpha uv}\ +\ \frac{1}{a\alpha u}\ +\ \frac{1}{a\alpha v} &&(j=\I),\label{aud:certI}\\
\frac{4}{n}\ &=\ \frac{1}{\alpha uv}\ +\ \frac{1}{a\alpha un}\ +\ \frac{1}{a\alpha vn} &&(j=\II).\label{aud:certII}
\end{align}
Conversely, any factorization $K_j(c)=\alpha uv$ with $c\mid(u+v)$ implies $-1\in D_c(K_j(c))$ and yields the corresponding decomposition~\eqref{aud:certI} or~\eqref{aud:certII}, even without assuming $\gcd(u,v)=1$.
\end{proposition}
\begin{proof}
Suppose $d,e\mid K_j(c)$ satisfy $c\mid(d+e)$, and set $g=\gcd(d,e)$. Since $\gcd(g,c)=1$, the coprime positive integers $u=d/g$ and $v=e/g$ satisfy $c\mid(u+v)$ and divide $K_j(c)$, whence $uv\mid K_j(c)$. Defining $\alpha=K_j(c)/(uv)$ and $a=(u+v)/c$, we observe that the sum of the last two unit fractions in~\eqref{aud:certI} equals $\frac{u+v}{a\alpha uv}=\frac{c}{K_{\I}(c)}$, while in~\eqref{aud:certII} it equals $\frac{c}{nK_{\II}(c)}$. Substituting the definitions~\eqref{aud:K} immediately gives~\eqref{aud:certI} and~\eqref{aud:certII}. The converse follows by the same algebraic identity without dividing by $\gcd(u,v)$.
\end{proof}

\begin{remark}[Parametrization convention]
Proposition~\ref{aud:reconstruction} associates two distinct shift integers $K_{\I}(c)$ and $K_{\II}(c)$ to the same target residue $-1\bmod c$. By contrast, Theorem~\ref{thm:difference_set_positivity_MR}(a) works with the single integer $K_c=(n+c)/4=K_{\II}(c)$ and distinguishes the two types via the targets $-1\bmod c$ (Type~II) and $-n^{-1}\bmod c$ (Type~I). The two conventions should be kept distinct unless the corresponding change of variables is performed.
\end{remark}

\begin{lemma}[Factor reciprocity]\label{aud:reciprocity}
If $p\equiv 3\pmod 4$ is a prime shift ($p\neq n$) and $q$ is any prime dividing $K_j(p)$, then
\begin{equation}\label{aud:reciprocityeq}
\leg{q}{p}\ =\ \leg{q}{n}\qquad(j\in\{\I,\II\}).
\end{equation}
More generally, in Type~I, for every composite shift $c\equiv 3\pmod 4$ coprime to $n$ and every prime $q\mid K_{\I}(c)$, one has $\leg{q}{c}=\leg{q}{n}$, where the left-hand side denotes the Jacobi symbol.
\end{lemma}
\begin{proof}
Any prime $q\mid K_j(p)$ is distinct from both $p$ and $n$. If $q$ is odd, the definition~\eqref{aud:K} gives $p\equiv -n\pmod q$ in Type~II and $p\equiv -n^{-1}\pmod q$ in Type~I. Using quadratic reciprocity together with $p\equiv3\pmod 4$ and $n\equiv1\pmod 4$, we obtain
\[
\leg{q}{p}\ =\ (-1)^{(q-1)/2}\leg{p}{q}\ =\ (-1)^{(q-1)/2}\leg{-1}{q}\leg{n}{q}\ =\ \leg{q}{n}.
\]
If $q=2$, divisibility of $K_j(p)$ by $2$ forces $p\equiv7\pmod 8$ because $n\equiv 1\pmod 8$, so $\leg{2}{p}=\leg{2}{n}=+1$. The same Jacobi reciprocity calculation extends verbatim to composite $c\equiv 3\pmod 4$ in Type~I, since $nc\equiv -1\pmod q$ for odd $q\mid K_{\I}(c)$ and $c\equiv7\pmod 8$ whenever $2\mid K_{\I}(c)$.
\end{proof}

\begin{corollary}[Necessary non-residue restriction on $\alpha$ in Type~I]\label{aud:alphaNR}
Every Type~I solution $4\alpha de=nc+1$ with $d+e=ac$ necessarily satisfies $\leg{\alpha}{n}=-1$.
\end{corollary}
\begin{proof}
Since $\alpha d\mid K_{\I}(c)$, we have $\gcd(\alpha,c)=\gcd(d,c)=1$. Eliminating $e=ac-d$ modulo $c$ in $4\alpha de=nc+1$ gives the congruence $4\alpha d^{2}\equiv -1\pmod c$, which forces $\leg{\alpha}{c}=\leg{-1}{c}=-1$ since $c\equiv 3\pmod 4$. Applying Lemma~\ref{aud:reciprocity} with multiplicity across the prime factors of $\alpha$ yields $\leg{\alpha}{n}=\leg{\alpha}{c}=-1$. Thus dropping the non-residue condition on $\alpha$ creates no new Type~I witnesses: although the algebraic reconstruction of Proposition~\ref{aud:reconstruction} does not require imposing $\leg{\alpha}{n}=-1$ a priori, every Type~I solution satisfies it automatically.
\end{proof}

\begin{lemma}[Type~II Jacobi restriction]\label{aud:IIjacobi}
For any Type~II solution $an=4a\alpha uv-u-v$, setting $s=4a\alpha u-1$ and $c=4\alpha uv-n$ yields
\[
\leg{s}{n}\ =\ \leg{uv}{n}\ =\ \leg{\alpha c}{n}\ =\ -1.
\]
\end{lemma}
\begin{proof}
This reproduces Theorem~\ref{thm:global_jacobi_reciprocity}(b). Since $n\nmid\alpha uv$ and $sv=an+u$, we have $\gcd(s,n)=1$ and $s\equiv3\pmod 4$. For every odd prime $r\mid au$, the congruence $s\equiv -1\pmod r$ gives $\leg{r}{s}=+1$, while if $2\mid au$ then $s\equiv7\pmod 8$ and $\leg{2}{s}=+1$; hence $\leg{au}{s}=+1$ by multiplicativity. Because $an\equiv -u\pmod s$ and $\gcd(a,s)=1$, it follows that $\leg{n}{s}=\leg{-au}{s}=\leg{-1}{s}=-1$, whence $\leg{s}{n}=-1$ by quadratic reciprocity ($n\equiv1\pmod 4$). Finally, reducing $sv=an+u$ and $c=4\alpha uv-n$ modulo $n$ gives $sv\equiv u\pmod n$ and $c\equiv 4\alpha uv\pmod n$, which yields $\leg{uv}{n}=\leg{s}{n}=-1$ and $\leg{\alpha c}{n}=\leg{(2\alpha u)^{2}uv}{n}=\leg{uv}{n}=-1$.
\end{proof}

\begin{proposition}[Diophantine coupling and elementary case]\label{aud:coupling}
For every $c\equiv3\pmod 4$ with $\gcd(c,n)=1$, one has the exact coupling identity
\[
c\,K_{\II}(c)\ -\ K_{\I}(c)\ =\ \frac{c^{2}-1}{4},\qquad \gcd\!\big(K_{\I}(c),K_{\II}(c)\big)\ \mid\ \frac{c^{2}-1}{4}.
\]
If a prime quadratic non-residue $q$ modulo $n$ divides both $K_{\I}(c)$ and $K_{\II}(c)$, then $q\equiv3\pmod 4$ and $q\mid(n+1)$. Whenever $n+1=(4a-1)d$, $4/n$ admits the elementary decomposition
\begin{equation}\label{aud:elementary}
\frac{4}{n}\ =\ \frac{1}{ad}\ +\ \frac{1}{an}\ +\ \frac{1}{adn}.
\end{equation}
\end{proposition}
\begin{proof}
Eliminating $c$ and $n$ linearly between $4K_{\I}(c)=nc+1$ and $4K_{\II}(c)=n+c$ gives $c\,K_{\II}(c)-K_{\I}(c)=(c^{2}-1)/4$, from which the greatest-common-divisor relation follows immediately. Any common prime factor $q\mid\gcd(K_{\I}(c),K_{\II}(c))$ satisfies $c\equiv -n\pmod q$ and $nc\equiv -1\pmod q$, whence $n^{2}\equiv 1\pmod q$. Since $\leg{2}{n}=+1$ for $n\equiv 1\pmod{24}$ and $\leg{q}{n}=-1$ by hypothesis, $q$ is odd, and quadratic reciprocity gives $\leg{n}{q}=\leg{q}{n}=-1$. Because $n\equiv\pm 1\pmod q$, we must have $n\equiv -1\pmod q$, which means $q\mid(n+1)$ and $\leg{-1}{q}=-1$, so $q\equiv3\pmod 4$. Finally, when $n+1=(4a-1)d$, summing the right-hand side of~\eqref{aud:elementary} over the common denominator $adn$ produces the numerator $n+d+1=4ad$, verifying~\eqref{aud:elementary}.
\end{proof}

\begin{remark}[Disjointness of non-residue prime factors versus statistical independence]
Proposition~\ref{aud:coupling} shows that for every prime $n\equiv 1\pmod{24}$ not already resolved by the elementary identity~\eqref{aud:elementary}, the integers $K_{\I}(c)$ and $K_{\II}(c)$ cannot share any prime factor $q$ with $\leg{q}{n}=-1$. This arithmetic disjointness of non-residue prime factors should not, however, be conflated with statistical independence of the factorizations of $K_{\I}(c)$ and $K_{\II}(c)$ as $c$ varies.
\end{remark}

\subsection{Fixed-parameter finite reductions and the divisor encounter problem}\label{aud:alpha}

Before varying the parameter $\alpha$ across polynomial ranges in Subsection~\ref{aud:inputs}, we examine the structure of the Type~I and Type~II witness equations when $\alpha\ge 1$ is held fixed. Doing so yields an explicit finite reduction for fixed $(n,\alpha)$ and clarifies why varying $\alpha$ (or $c$) is essential in asymptotic arguments.

\begin{proposition}[Exact finite reduction in Type~I at fixed $\alpha$]\label{aud:fixedalphaI}
Fix $n$ and $\alpha\ge1$. The Type~I solutions
\[
4\alpha de\ =\ nc+1,\qquad d+e\ =\ ac,\qquad d\le e
\]
are in bijection with pairs $(k,d)$ of positive integers satisfying
\begin{align}
\left\lfloor\frac{n}{4\alpha}\right\rfloor+1\ &\le\ k\ \le\ \left\lfloor\frac{2n}{4\alpha-1}\right\rfloor,\qquad d\mid k,\label{aud:finiteI}\\
s:=4\alpha k-n\ &\mid\ (4\alpha d^{2}+1),\qquad (4\alpha k-2n)d^{2}\ \le\ k,\nonumber
\end{align}
via the explicit reconstruction $a=k/d$, $c=(4\alpha d^{2}+1)/s$, and $e=ac-d$.
\end{proposition}
\begin{proof}
Setting $k=ad$ and substituting $e=ac-d$ into $4\alpha de=nc+1$ yields $(4\alpha k-n)c=4\alpha d^{2}+1$, so $s=4\alpha k-n>0$ and $s\mid(4\alpha d^{2}+1)$. The ordering condition $e\ge d$ is equivalent to $ac\ge 2d$, which upon multiplying by $sd=dk/a$ becomes $(4\alpha k-2n)d^{2}\le k$. In particular, $4\alpha k-2n\le k/d^{2}\le k$, which gives $(4\alpha-1)k\le 2n$ and establishes the upper bound on $k$ in~\eqref{aud:finiteI}. Conversely, any pair $(k,d)$ satisfying~\eqref{aud:finiteI} defines positive integers $a=k/d$, $c=(4\alpha d^{2}+1)/s$, and $e=ac-d\ge d\ge 1$ satisfying the initial system, with $c\equiv3\pmod 4$ since $s\equiv -n\equiv 3\pmod 4$.
\end{proof}

\begin{corollary}[Quantifiers of the fixed-$\alpha$ finiteness obstruction]\label{aud:finitequantifiers}
If $n\mid\alpha$, reducing $4\alpha de=nc+1$ modulo $n$ shows that there is no Type~I solution. When $\gcd(n,\alpha)=1$, the modular inverse $n^{-1}\pmod{4\alpha}$ exists and only finitely many shifts $c$ can yield a Type~I solution at parameter $\alpha$: with $U=\lfloor 2n/(4\alpha-1)\rfloor$, every such shift satisfies $c\le(4\alpha U^{2}+1)/3$. Consequently, among primes $p$ in the reduced progression $p\equiv -n^{-1}\pmod{4\alpha}$, the proportion of successful shifts tends to zero as the upper cutoff on $p$ tends to infinity \emph{at fixed $n$ and $\alpha$}.
\end{corollary}
\begin{proof}
In~\eqref{aud:finiteI}, we have $d\le k\le U$ and $s=4\alpha k-n\ge 3$ (since $s>0$ and $s\equiv 3\pmod 4$), which forces $c=(4\alpha d^{2}+1)/s\le(4\alpha U^{2}+1)/3$. By Dirichlet's theorem, the reduced progression $p\equiv -n^{-1}\pmod{4\alpha}$ contains infinitely many primes, whereas the number of successful shifts is bounded for fixed $(n,\alpha)$.
\end{proof}

\begin{remark}[Fixed-parameter asymptotics versus the joint regime $n\to\infty$]
The limiting density zero in Corollary~\ref{aud:finitequantifiers} holds when $(n,\alpha)$ is fixed and the upper bound on $p$ tends to $+\infty$. It says nothing by itself about the joint regime $n\to\infty$ with $p\le n^{\theta}$ ($0<\theta<1$), where $U\asymp n/\alpha$ grows with $n$; extrapolating the fixed-$(n,\alpha)$ limit to that regime without uniformity in $n$ would be a quantifier reversal.
\end{remark}

\begin{proposition}[Finite reduction in Type~II and square classes]\label{aud:fixedalphaII}
For fixed $\alpha\ge 1$, the Type~II solutions $4\alpha de=n+c$, $d+e=ac$ with $d\le e$ are in bijection with pairs $(d,e)$ of positive integers satisfying
\[
4\alpha d^{2}-2d\ \le\ n,\qquad \max\!\left(d,\left\lfloor\frac{n}{4\alpha d}\right\rfloor+1\right)\ \le\ e\ \le\ \left\lfloor\frac{n+d}{4\alpha d-1}\right\rfloor,\qquad (4\alpha de-n)\ \mid\ (d+e).
\]
Moreover, in both Type~I and Type~II, any solution at a non-square-free parameter $\alpha=\beta h^{2}$ induces a solution at the square-free core $\beta$ via the scaling $(a,d,e)\mapsto(ha,hd,he)$.
\end{proposition}
\begin{proof}
Setting $c=4\alpha de-n>0$ and $d+e=ac$ requires $c\mid(d+e)$ and in particular $c\le d+e$, which rearranges to $(4\alpha d-1)e\le n+d$. Combining this upper bound with $e\ge d$ yields $4\alpha d^{2}-2d\le n$ and the stated interval for $e$; the converse is immediate by setting $a=(d+e)/c$. Finally, replacing $(\alpha,a,d,e)$ by $(\beta,ha,hd,he)$ preserves $\beta(hd)(he)=\alpha de$ and multiplies both $d+e$ and $a$ by $h$, leaving $c$ unchanged.
\end{proof}

\begin{proposition}[Exact divisor-encounter formulation at fixed $\alpha$]\label{aud:meeting}
Let $p\equiv -n^{-1}\pmod{4\alpha}$ be prime and set $V:=(np+1)/(4\alpha)\in\N$. A Type~I solution at parameter $\alpha$ and shift $p$ exists if and only if
\begin{equation}\label{aud:meet}
\exists\,d\mid V:\qquad 4\alpha d^{2}\ \equiv\ -1\pmod p.
\end{equation}
If $p=p_0+4\alpha t$, then $V=V_0+nt$. In the special case where $V$ is prime,~\eqref{aud:meet} is equivalent to $p\mid(4\alpha+1)$.
\end{proposition}
\begin{proof}
Writing $e=V/d$, the definition $4\alpha V=np+1$ implies $4\alpha V\equiv 1\pmod p$. Consequently, $4\alpha d^{2}\equiv -1\pmod p$ holds if and only if $d^{2}\equiv -V=-de\pmod p$, which (since $\gcd(d,p)=1$) is equivalent to $d+e\equiv 0\pmod p$; Proposition~\ref{aud:reconstruction} then gives the correspondence with Type~I solutions. When $V$ is prime, its only positive divisors are $d=1$ and $d=V$, both of which reduce~\eqref{aud:meet} to the congruence $4\alpha+1\equiv 0\pmod p$.
\end{proof}

\begin{remark}[Ambient quadratic solvability versus divisor encounter]
Proposition~\ref{aud:meeting} underscores a key structural distinction: the fact that the quadratic congruence $4\alpha x^{2}\equiv -1\pmod p$ is solvable in the ambient group $(\Z/p\Z)^{\times}$ (which holds automatically whenever $\leg{\alpha}{p}=-1$ and $p\equiv 3\pmod 4$) does not imply that either of its two roots modulo $p$ is hit by an actual divisor $d\mid V$.
\end{remark}

\subsection{Analytic inputs and abundance of candidate parameters}\label{aud:inputs}

To move beyond fixed $\alpha$ and build an abundant supply of candidate shifts as $n\to\infty$, we now vary $\alpha$ over square-free quadratic non-residues modulo $n$ and combine three unconditional results from analytic number theory.

Set $\beta_0:=\frac{1}{4\sqrt{e}}$. Our candidate-abundance estimates rely purely on the following three published theorems, without any unproven hypothesis on primes in progressions modulo $n$:
\begin{enumerate}[1.]
    \item \emph{Character non-residues in Burgess intervals (Banks--Garaev--Heath-Brown--Shparlinski~\cite{BGHS}):} for every fixed $\varepsilon>0$, there exists $\delta>0$ such that for all sufficiently large primes $n$ and all $n^{\beta_0+\varepsilon}\le X\le n$, at least $\delta X$ integers $m\le X$ satisfy $\leg{m}{n}=-1$.
    \item \emph{Linnik's constant (Xylouris~\cite{xylouris2011}):} in every reduced residue class modulo $q$, the least prime is $\mathcal{O}(q^{5})$, with an absolute implied constant.
    \item \emph{The Bombieri--Vinogradov theorem (\cite[Thm.~17]{TaoBV}):} for every fixed $A>0$ and $\varepsilon>0$,
    \begin{equation}\label{aud:BV}
    \sum_{m\le H^{1/2-\varepsilon}}\Delta_H(m)\ \ll_{A,\varepsilon}\ H(\log H)^{-A},
    \end{equation}
    where $\Delta_H(m):=\max_{\gcd(a,m)=1}\big|\#\{p\text{ prime}:2H<p\le 4H,$\linebreak $p\equiv a\pmod m\}-\frac{\Li(4H)-\Li(2H)}{\varphi(m)}\big|$.
\end{enumerate}

\begin{proposition}[Harmonic mass of square-free composite parameters $\alpha$]\label{aud:harmonic}
Let $\mathcal{A}(n):=\big\{\alpha\le n^{4/25}:\mu^{2}(\alpha)=1,\ \leg{\alpha}{n}=-1\big\}$. Then
\[
\sum_{\alpha\in\mathcal{A}(n)}\frac{1}{\alpha}\ \gg\ \log n.
\]
Moreover, for each $\alpha\in\mathcal{A}(n)$, there exists a prime $p\equiv -n^{-1}\pmod{4\alpha}$ with $p\ll\alpha^{5}\le n^{4/5}$ such that, for all sufficiently large $n$, one has $p\le n^{81/100}$ and
\[
\leg{\alpha}{p}\ =\ -1,\qquad \#\big\{x\bmod p:4\alpha x^{2}\equiv -1\pmod p\big\}\ =\ 2,\qquad \frac{\sqrt{(np+1)/(4\alpha)}}{p}\ \ge\ \frac12\,n^{3/200}.
\]
\end{proposition}
\begin{proof}
Fix $a=31/200$ and $b=4/25$, so that $\beta_0<a<b$. By~\cite{BGHS}, the non-residue counting function $A_n(t):=\#\{m\le t:\leg{m}{n}=-1\}$ satisfies $A_n(t)\ge\delta t$ for all $t\in[n^{a},n^{b}]$. Integration by parts on $[1,n^{b}]$ gives
\[
\sum_{\substack{m\le n^{b}\\ \leg{m}{n}=-1}}\frac{1}{m}\ =\ \frac{A_n(n^{b})}{n^{b}}\ +\ \int_{1}^{n^{b}}\frac{A_n(t)}{t^{2}}\,dt\ \ge\ \delta(b-a)\log n.
\]
Every integer $m\le n^{b}$ with $\leg{m}{n}=-1$ writes uniquely as $m=\alpha h^{2}$ with $\alpha\in\mathcal{A}(n)$ and $h\ge 1$, so the left-hand side is at most $\zeta(2)\sum_{\alpha\in\mathcal{A}(n)}\alpha^{-1}$, proving the harmonic lower bound. For each $\alpha\in\mathcal{A}(n)$, since $\gcd(n,4\alpha)=1$, Xylouris's bound~\cite{xylouris2011} provides a prime $p\equiv -n^{-1}\pmod{4\alpha}$ satisfying $p\ll\alpha^{5}\le n^{4/5}$ (hence $p\le n^{81/100}$ for large $n$). Lemma~\ref{aud:reciprocity} then gives $\leg{\alpha}{p}=\leg{\alpha}{n}=-1$; since $p\equiv 3\pmod 4$, $-1/(4\alpha)$ is a quadratic residue modulo $p$, so $4\alpha x^{2}\equiv -1\pmod p$ has two solutions in $(\Z/p\Z)^{\times}$, and the bound $\alpha p\le n^{97/100}$ yields the archimedean ratio inequality.
\end{proof}

\begin{remark}[Scope of candidate-pair abundance]
Proposition~\ref{aud:harmonic} supplies candidate pairs $(\alpha,p)$ with total harmonic mass $\gg\log n$ over square-free composite integers $\alpha\in\mathcal{A}(n)$ (rather than over prime $\alpha$), and guarantees that the archimedean window $d\le\sqrt{V}$ is large relative to $p$; in view of Proposition~\ref{aud:meeting}, however, this candidate abundance does not by itself force the divisor encounter~\eqref{aud:meet}.
\end{remark}

\begin{lemma}[Constant harmonic mass of non-residue primes]\label{aud:nu}
For any fixed $\beta>\beta_0$ and $Z=n^{\beta}<n$, one has for all sufficiently large primes $n$:
\[
\nu_n(Z)\ :=\ \sum_{\substack{q\le Z\text{ prime}\\ \leg{q}{n}=-1}}\frac{1}{q}\ \ge\ c_{\beta}\ >\ 0.
\]
\end{lemma}
\begin{proof}
Every positive integer $m\le Z$ with $\leg{m}{n}=-1$ must be divisible by at least one prime $q\le Z$ satisfying $\leg{q}{n}=-1$. Consequently, the total number of non-residue integers $m\le Z$ is at most $\sum_{q\le Z,\,\leg{q}{n}=-1}\lfloor Z/q\rfloor\le Z\,\nu_n(Z)$, and the linear lower bound $\delta Z$ of~\cite{BGHS} yields $\nu_n(Z)\ge\delta>0$.
\end{proof}

\begin{remark}
By a theorem of Pollack~\cite{Pollack}, there are also at least $n^{\gamma}$ distinct non-residue primes $q\le n^{\beta_0+\varepsilon}$ for some $\gamma>0$, which via Xylouris's theorem yields a positive power of $n$ distinct prime shifts $p\le n^{\theta}$ whenever $\theta>5(\beta_0+\varepsilon)$. In Subsection~\ref{aud:sieve}, we replace individual Linnik primes by the Bombieri--Vinogradov theorem to obtain a positive proportion of all prime shifts $p\in(2H,4H]$.
\end{remark}

\subsection{Prime shifts, weighted Bombieri--Vinogradov, and annihilation of small character orders}\label{aud:sieve}

We now pass from candidate existence to positive-density families of prime shifts $p\in(2H,4H]$ ($H=n^{\theta}$) for which $K_j(p)$ is divisible to the first power by a non-residue prime $q\le Z$, and we combine the Bombieri--Vinogradov theorem with the Rosser--Iwaniec linear sieve to annihilate simultaneously all odd Dirichlet characters modulo $p$ of order strictly less than $2z$.

\begin{lemma}[Bombieri--Vinogradov with divisor weights]\label{aud:weightedBV}
For any fixed $J\ge0$, $A>0$, and $\varepsilon>0$,
\begin{equation}\label{aud:weightedBVeq}
\sum_{m\le H^{1/2-\varepsilon}}\tau(m)^{J}\,\Delta_H(m)\ \ll_{J,A,\varepsilon}\ H(\log H)^{-A}.
\end{equation}
\end{lemma}
\begin{proof}
Replacing $J$ by $\lceil J\rceil$ if necessary, the trivial progression bound and the inequality $\frac{1}{\varphi(m)}\le\frac{\tau(m)}{m}$ give $\Delta_H(m)\ll\frac{H\tau(m)}{m}$ for all $m\le H$. Because $(a+1)^{k}\le\binom{a+2^{k}-1}{2^{k}-1}$ for every integer $a\ge 0$, multiplicativity yields $\tau(m)^{k}\le\tau_{2^{k}}(m)$ for every integer $k\ge1$, and therefore $\sum_{m\le H}\frac{\tau(m)^{k}}{m}\ll_k(\log(2H))^{2^{k}}$. Applying the Cauchy--Schwarz inequality gives
\[
\sum_{m\le H^{1/2-\varepsilon}}\tau(m)^{J}\,\Delta_H(m)\ \le\ \left(\sum_{m\le H^{1/2-\varepsilon}}\Delta_H(m)\right)^{\!1/2}\left(C H\sum_{m\le H}\frac{\tau(m)^{2J+1}}{m}\right)^{\!1/2},
\]
and invoking~\eqref{aud:BV} with exponent $2A+2^{2J+1}$ establishes~\eqref{aud:weightedBVeq}. Note that because $\Delta_H(m)$ takes the maximum over all reduced residue classes modulo $m$, the estimate~\eqref{aud:weightedBVeq} holds uniformly even when the target residue classes depend on the external prime $n$.
\end{proof}

\begin{proposition}[Unsieved abundance of prime shifts with quadratic annihilation]\label{aud:unsieved}
Fix $5\beta_0<\theta<1$, $H=n^{\theta}$, and $Z=H^{1/5}$. For each $j\in\{\I,\II\}$, at least $c_{\theta}\frac{H}{\log H}$ primes $p\in(2H,4H]$ with $p\equiv3\pmod 4$ admit a prime $q\le Z$ such that $\leg{q}{n}=-1$ and $q\parallel K_j(p)$. For every such $p$, the quadratic character $\leg{\cdot}{p}$ annihilates $F_{K_j(p)}$ identically, and the prime factors of $K_j(p)$ generate a subgroup of $(\Z/p\Z)^{\times}$ containing $-1$.
\end{proposition}
\begin{proof}
Set $\mathcal{Q}:=\{q\le Z:q\text{ prime},\ \leg{q}{n}=-1\}$. Since $Z=n^{\theta/5}$ with $\theta/5>\beta_0$, Lemma~\ref{aud:nu} gives $\nu:=\sum_{q\in\mathcal{Q}}\frac1q\ge c>0$ (with $q\ge5$ for $n\equiv 1\pmod{24}$), and hence $\mu:=\sum_{q\in\mathcal{Q}}\frac{1}{q-1}\in[\nu,\frac54\nu]$. For each prime $p\in(2H,4H]$ with $p\equiv3\pmod 4$, define
\[
V_j(p)\ :=\ \sum_{q\in\mathcal{Q}}\ind_{q\parallel K_j(p)},\qquad W_j(p)\ :=\ \sum_{q\in\mathcal{Q}}\ind_{q\mid K_j(p)},\qquad X_H\ :=\ \frac{\Li(4H)-\Li(2H)}{2}.
\]
Because $q^{a}\mid K_j(p)$ places $p$ in a unique reduced residue class modulo $4q^{a}$ of relative density $\frac{1}{\varphi(q^{a})}$, exact divisibility $q\parallel K_j(p)$ has relative density $\frac{1}{q-1}-\frac{1}{q(q-1)}=\frac1q$. Applying~\eqref{aud:BV} and Lemma~\ref{aud:weightedBV} across the moduli $4q,\,4q^{2},\,4qr\le 4Z^{2}=4H^{2/5}<H^{1/2-\varepsilon}$ yields
\[
\sum_p V_j(p)\ =\ X_H\nu+o(X_H),\qquad \sum_p W_j(p)^{2}\ =\ X_H\Big(\mu^{2}+\mu-\sum_{q\in\mathcal{Q}}\frac{1}{(q-1)^{2}}\Big)+o(X_H).
\]
Since $0\le V_j(p)\le W_j(p)$, the Cauchy--Schwarz inequality gives
\[
\#\big\{p\in(2H,4H]:p\equiv 3\pmod 4,\ V_j(p)>0\big\}\ \ge\ \frac{\big(\sum_p V_j(p)\big)^{2}}{\sum_p W_j(p)^{2}}\ \gg\ X_H\,\frac{\nu^{2}}{\nu^{2}+\nu}\ \gg\ \frac{H}{\log H}.
\]
Finally, for any $p$ with $V_j(p)>0$, choosing $q\in\mathcal{Q}$ with $q\parallel K_j(p)$ gives $\leg{q}{p}=\leg{q}{n}=-1$ by Lemma~\ref{aud:reciprocity}. The local Euler factor of $F_{K_j(p)}(\leg{\cdot}{p})$ at $q$ is therefore $1+\leg{q}{p}=0$, which annihilates $F_{K_j(p)}(\leg{\cdot}{p})$, and since $(\Z/p\Z)^{\times}$ is cyclic of order $p-1\equiv2\pmod 4$, any quadratic non-residue $q\bmod p$ generates a subgroup of even order containing the unique involution $-1\bmod p$.
\end{proof}

To annihilate all odd characters of small order simultaneously via the linear sieve, we henceforth assume the residual condition
\begin{equation}\label{aud:residual}
r\mid(n+1),\quad r\text{ odd prime}\qquad\Longrightarrow\qquad r\equiv1\pmod 4, \tag{R}
\end{equation}
since any prime $n\equiv1\pmod{24}$ violating~\eqref{aud:residual} possesses a divisor $d=r\equiv 3\pmod 4$ of $n+1$ and is immediately resolved by the elementary decomposition~\eqref{aud:elementary}. Under~\eqref{aud:residual}, every prime $q$ with $\leg{q}{n}=-1$ satisfies $q\nmid(n+1)$, for otherwise $n\equiv -1\pmod q$ would force $\leg{q}{n}=\leg{n}{q}=\leg{-1}{q}=+1$.

Since $4\beta_0=1/\sqrt{e}<1$, we may fix exponents and scales satisfying
\begin{align}
4\beta_0<\theta<1,&\qquad 0<\kappa<\frac14-\frac{\beta_0}{\theta},\nonumber\\
\frac{\beta_0}{\theta}<b<\frac14-\kappa,&\qquad 2\kappa<\rho<\frac12-2b,\label{aud:parameters}\\
0<\eta<\min\!\left(\kappa,\frac12-\rho-2b\right),&\qquad H=n^{\theta},\quad Z=H^{b},\quad z=H^{\kappa},\quad D=H^{\rho},\quad T=H^{\eta},\nonumber
\end{align}
and define the sifted shift families and Mertens normalization
\begin{align}
\mathcal{S}\ &:=\ \big\{p\in(2H,4H]:p\text{ prime},\ p\equiv3\pmod 4,\ P^{-}((p-1)/2)\ge z\big\},\label{aud:S}\\
\mathcal{Q}\ &:=\ \big\{q\le Z:q\text{ prime},\ \leg{q}{n}=-1\big\},\quad \nu:=\sum_{q\in\mathcal{Q}}\frac1q\ge c>0,\quad \mathcal{G}_j:=\big\{p\in\mathcal{S}:V_j(p)>0\big\},\nonumber\\
\mathcal{V}(z)\ &:=\ \prod_{3\le\ell<z}\left(1-\frac{1}{\ell-1}\right),\qquad B_H\ :=\ X_H\mathcal{V}(z)\ \asymp\ \frac{H}{(\log H)^{2}}.\label{aud:B}
\end{align}
We apply the one-dimensional Rosser--Iwaniec linear sieve~\cite[\S5.7.2, (5.116)--(5.117)]{Kumchev} with sifting parameter $s=\frac{\log D}{\log z}=\frac{\rho}{\kappa}>2$, for which the lower and upper linear-sieve functions satisfy $f(s)>0$ and $F(s)<\infty$. For odd primes $\ell<z$ not dividing a prescribed divisor of $K_j(p)$, the local sifting density is $g(\ell)=\frac{1}{\ell-1}$ (and $g(\ell)=0$ at primes already dividing the prescribed divisor), which satisfies the linear-sieve condition $\prod_{w\le\ell<z}(1-g(\ell))^{-1}\le\frac{\log z}{\log w}\big(1+\mathcal{O}(\frac{1}{\log w})\big)$ uniformly by Mertens' theorem. When an odd prime $r\nmid n(n+1)$ divides $K_j(p)$, the congruence imposed on $p\bmod r$ is coprime to $r$ and distinct from $1\pmod r$; hence for $r<z$, removing the local sifting factor $1-\frac{1}{r-1}$ modifies its effective density to
\begin{equation}\label{aud:localh}
\frac{1}{r-1}\left(1-\frac{1}{r-1}\right)^{-1}\ =\ \frac{1}{r-2}.
\end{equation}
Similarly, for the unitary condition $q\parallel K_j(p)$ (which is a single non-negative sequence whose progression count subtracts $q^{2}\mid K_j(p)$ from $q\mid K_j(p)$), the unsieved density is $1/q$, and for $q<z$ the sifting adjustment factor is $\frac{q-1}{q-2}\le\frac43$ since $q\ge5$.

\begin{lemma}[Aggregated remainder control]\label{aud:remainders}
Under~\eqref{aud:parameters}, the linear-sieve remainders over moduli $u\le D$, weighted by $3^{\omega(d)}$ over square-free $d\le T$ and summed over one or two primes $q,r\le Z$, contribute in total $\mathcal{O}_A\big(H(\log H)^{-A}\big)$ for every fixed $A>0$.
\end{lemma}
\begin{proof}
Every composite progression modulus encountered in the sieve divides either $4udq^{2}$ or $4udqr$, and is therefore bounded by the total level of distribution
\begin{equation}\label{aud:level}
4DTZ^{2}\ =\ 4H^{\rho+\eta+2b}\ \le\ H^{1/2-\varepsilon}
\end{equation}
for a fixed $\varepsilon>0$ and all sufficiently large $n$. Grouping the terms according to the actual least-common-multiple modulus $m=\operatorname{lcm}(4,u,d,q^{a})$ ($a\in\{1,2\}$) or $m=\operatorname{lcm}(4,u,d,q,r)$, we observe that each selected prime $q,r$ and every prime factor of the square-free integers $u,d$ divides $m$. Consequently, the total weight of compatible representations of a given modulus $m$, including the factor $3^{\omega(d)}$, is bounded by $\tau(m)^{C}$ for an absolute constant $C$ (while incompatible congruence systems contribute zero). Applying Lemma~\ref{aud:weightedBV} then absorbs the aggregated remainder into $\mathcal{O}_A(H(\log H)^{-A})$.
\end{proof}

\begin{proposition}[Sieved abundance and annihilation of small odd character orders]\label{aud:rough}
Under~\eqref{aud:residual} and~\eqref{aud:parameters}, for each $j\in\{\I,\II\}$, one has
\begin{equation}\label{aud:roughcount}
\#\mathcal{G}_j\ \gg\ B_H\qquad\text{and}\qquad \#\mathcal{S}\ \ll\ B_H.
\end{equation}
Moreover, for every $p\in\mathcal{G}_j$, every odd Dirichlet character $\chi\pmod p$ of order strictly less than $2z$ is annihilated both by the full divisor sum $F_{K_j(p)}(\chi)$ and by the square-free divisor sum $F_{p,j}^{\sfm}(\chi):=\sum_{d\mid\rad K_j(p)}\chi(d)$.
\end{proposition}
\begin{proof}
Applying the lower linear sieve to $q\parallel K_j(p)$ and the upper linear sieve to $q\mid K_j(p)$ and $r\mid K_j(p)$ (with $d=1$ in Lemma~\ref{aud:remainders}) yields the first and second moment bounds $\sum_{p\in\mathcal{S}}V_j(p)\gg B_H\nu$ and $\sum_{p\in\mathcal{S}}W_j(p)^{2}\ll B_H(\nu^{2}+\nu)$, from which Cauchy--Schwarz gives~\eqref{aud:roughcount}.

Now let $p\in\mathcal{G}_j$. By definition of $\mathcal{S}$ in~\eqref{aud:S}, we have $p-1=2m$ where $m$ is odd and $P^{-}(m)\ge z$. Every odd character $\chi\pmod p$ has even order $2d$ for some divisor $d\mid m$; if $\ord(\chi)=2d<2z$, then $d<z$ and $d\mid m$ force $d=1$, which means that $\chi$ is the unique quadratic character $\leg{\cdot}{p}$. Because $p\in\mathcal{G}_j$, there exists a prime $q\in\mathcal{Q}$ such that $q\parallel K_j(p)$ and $\leg{q}{p}=\leg{q}{n}=-1$. Since $q$ divides $K_j(p)$ to the exact power $1$, the local Euler factor $1+\leg{q}{p}=0$ appears in both $F_{K_j(p)}(\leg{\cdot}{p})$ and $F_{p,j}^{\sfm}(\leg{\cdot}{p})$, annihilating both sums. (For a general prime-power multiplicity $q^{v}\parallel K_j(p)$, the local factor $1-1+\dots+(-1)^{v}$ vanishes if and only if $v$ is odd, which is why the unitary condition $q\parallel K_j(p)$ is imposed for $F_{K_j(p)}$.)
\end{proof}

\subsection{Divergent square-free main term after sieving}\label{aud:weighted}

While Proposition~\ref{aud:rough} provides $\#\mathcal{G}_j\gg\frac{H}{(\log H)^{2}}$ sifted prime shifts, weighting each shift uniformly by $\frac{1}{p-1}\asymp\frac{1}{H}$ would yield only $\mathcal{O}((\log H)^{-2})$ if the divisor count were bounded. In reality, each shift $p\in\mathcal{G}_j$ contributes $\frac{m_j(p)(m_j(p)-1)}{p-1}$ to the square-free main term, where $m_j(p):=2^{\omega(K_j(p))}$. We now prove that the auxiliary divisor weight $4^{\omega(K_j(p))}$ produces a $(\log T)^{3}$ Euler-product amplification that drives the sifted square-free main term
\[
\mathcal{M}_j^{\sfm}\ :=\ \sum_{p\in\mathcal{G}_j}\frac{m_j(p)\big(m_j(p)-1\big)}{p-1}
\]
to $+\infty$ as $n\to\infty$.

\begin{theorem}[Weighted lower bound on the sieved main term]\label{aud:weightedmain}
Under the hypotheses of Proposition~\ref{aud:rough}, for each $j\in\{\I,\II\}$,
\begin{equation}\label{aud:weightedmaineq}
\mathcal{M}_j^{\sfm}\ \gg\ \frac{\log H}{(\log\log n)^{3/2}}\ \longrightarrow\ +\infty.
\end{equation}
\end{theorem}
\begin{proof}
Define the non-negative truncated divisor weight
\begin{equation}\label{aud:Aweight}
A_j(p)\ :=\ \sum_{\substack{d\le T,\ d\mid K_j(p)\\ \mu^{2}(d)=1,\ \gcd(d,6(n+1))=1}}3^{\omega(d)}\ \le\ \sum_{d\mid\rad K_j(p)}3^{\omega(d)}\ =\ 4^{\omega(K_j(p))}\ =\ m_j(p)^{2},
\end{equation}
and set $h(d):=\prod_{r\mid d}\frac{1}{r-2}$ and $\mathcal{P}_n(T):=\prod_{5\le r\le T,\,r\nmid(n+1)}\big(1+\frac{3}{r-2}\big)$. We establish the two weighted moment estimates
\begin{align}
T_1\ &:=\ \sum_{p\in\mathcal{S}}A_j(p)\,V_j(p)\ \gg\ B_H\,\nu\,\mathcal{P}_n(T),\label{aud:T1}\\
T_2\ &:=\ \sum_{p\in\mathcal{S}}A_j(p)\,W_j(p)^{2}\ \ll\ B_H\,(\nu^{2}+\nu)\,\mathcal{P}_n(T).\label{aud:T2}
\end{align}

To prove~\eqref{aud:T1}, we restrict the inner sum by positivity to $\gcd(d,q)=1$ and apply the lower linear sieve with the local densities~\eqref{aud:localh} (since $d\le T<z$) together with Lemma~\ref{aud:remainders}, obtaining
\begin{equation}\label{aud:T1local}
T_1\ \ge\ c\,B_H\sum_{q\in\mathcal{Q}}\frac1q\sum_{\substack{d\le T,\ \mu^{2}(d)=1\\ \gcd(d,6(n+1)q)=1}}3^{\omega(d)}h(d)\ -\ \mathcal{O}_A\big(H(\log H)^{-A}\big).
\end{equation}
To bound the inner sum over $d\le T$ from below, we restrict $d$ to square-free products of admissible primes $r\le T^{1/12}$ (where omitting the single prime $q$ reduces the Euler product by at most a factor of $1+\frac{3}{q-2}\le 2$). Normalizing the Euler product over $r\le T^{1/12}$ into a probability distribution with independent prime inclusion probabilities $\frac{3/(r-2)}{1+3/(r-2)}=\frac{3}{r+1}$, we have $\Ebb[\log d]\le\sum_{r\le T^{1/12}}\frac{3\log r}{r+1}\le\frac14\log T+\mathcal{O}(1)$; Markov's inequality then shows that $d\le T$ holds with probability at least $3/4$, so the inner sum over $d\le T$ is $\gg\mathcal{P}_n(T^{1/12})\asymp\mathcal{P}_n(T)$ uniformly in $q$, establishing~\eqref{aud:T1}.

To prove~\eqref{aud:T2}, fix $d\le T$ and set $u(d):=\#\{q\in\mathcal{Q}:q\mid d\}$, with $g_q=(q-2)^{-1}$ for $q<z$ and $g_q=(q-1)^{-1}$ for $q\ge z$ (so that $\sum_{q\in\mathcal{Q}}g_q\le\frac53\nu$). Any prime $q\in\mathcal{Q}$ dividing $d$ is already imposed by $d\mid K_j(p)$, while primes $q\nmid d$ contribute additional local density $g_q$. Applying the upper linear sieve and Lemma~\ref{aud:remainders} gives
\begin{equation}\label{aud:T2local}
T_2\ \ll\ B_H\sum_{\substack{d\le T,\ \mu^{2}(d)=1\\ \gcd(d,6(n+1))=1}}3^{\omega(d)}h(d)\Big((u(d)+C\nu)^{2}+u(d)+C\nu\Big)\ +\ \mathcal{O}_A\big(H(\log H)^{-A}\big).
\end{equation}
Extending the outer sum over $d$ by positivity to all square-free products of admissible primes $r\le T$ \emph{after} applying the sieve (which therefore leaves the level of distribution~\eqref{aud:level} unchanged), the normalized Euler product over $r\le T$ satisfies $\Ebb[u]\le 3\nu$ and $\operatorname{Var}(u)\le 3\nu$. Hence the expectation of $(u(d)+C\nu)^{2}+u(d)+C\nu$ is $\mathcal{O}(\nu^{2}+\nu)$, proving~\eqref{aud:T2}.

Because $V_j(p)\le W_j(p)$ and $\nu\ge c>0$, applying the weighted Cauchy--Schwarz inequality to~\eqref{aud:T1} and~\eqref{aud:T2} yields
\begin{equation}\label{aud:massA}
\sum_{p\in\mathcal{G}_j}A_j(p)\ \ge\ \frac{T_1^{2}}{\sum_{p\in\mathcal{S}}A_j(p)V_j(p)^{2}}\ \gg\ B_H\,\mathcal{P}_n(T)\,\frac{\nu^{2}}{\nu^{2}+\nu}\ \gg\ B_H\,\mathcal{P}_n(T).
\end{equation}
For every $p\in\mathcal{G}_j$, we have $m_j(p)\ge 2$, so $m_j(p)(m_j(p)-1)\ge\frac12 m_j(p)^{2}\ge\frac12 A_j(p)$; dividing by $p-1<4H$ and using $B_H\asymp H/(\log H)^{2}$ gives $\mathcal{M}_j^{\sfm}\gg\frac{\mathcal{P}_n(T)}{(\log H)^{2}}$. Finally, by Mertens' theorem,
\[
\mathcal{P}_n(T)\ \asymp\ \frac{(\log T)^{3}}{\prod_{\substack{5\le r\le T\\ r\mid(n+1)}}\big(1+\frac{3}{r-2}\big)}.
\]
Under the residual condition~\eqref{aud:residual}, every prime $r\mid(n+1)$ with $r\ge5$ satisfies $r\equiv1\pmod 4$. Writing $L=\log(n+1)$, Mertens' theorem in the progression $1\pmod 4$ gives
\[
\sum_{\substack{r\mid(n+1)\\ r\ge5}}\frac1r\ \le\ \sum_{\substack{r\le L\\ r\equiv1\pmod 4}}\frac1r\ +\ \frac{1}{L}\,\frac{\log(n+1)}{\log L}\ \le\ \frac12\log\log\log n\ +\ \mathcal{O}(1),
\]
which bounds the denominator of $\mathcal{P}_n(T)$ by $\mathcal{O}((\log\log n)^{3/2})$ and completes the proof of~\eqref{aud:weightedmaineq}.
\end{proof}

\begin{corollary}[Shifts with few prime factors carry negligible main-term mass]\label{aud:lowomega}
For any fixed $u<\frac{1}{\log 2}$, one has
\[
\sum_{\substack{p\in\mathcal{G}_j\\ \omega(K_j(p))\le u\log\log H}}\frac{m_j(p)\big(m_j(p)-1\big)}{p-1}\ \ll\ (\log H)^{2u\log 2-2}\ =\ o(1).
\]
\end{corollary}
\begin{proof}
Whenever $\omega(K_j(p))\le u\log\log H$, we have $m_j(p)^{2}=4^{\omega(K_j(p))}\le(\log H)^{2u\log 2}$. Summing over $p\in\mathcal{G}_j$ using $p>2H$ and the upper bound $\#\mathcal{G}_j\le\#\mathcal{S}\ll\frac{H}{(\log H)^{2}}$ from Proposition~\ref{aud:rough} gives the stated estimate.
\end{proof}

\subsection{Complete positivity criteria via square divisors $b\mid h^{2}$}\label{aud:hb}

Whereas Subsections~\ref{aud:algebra}--\ref{aud:weighted} index decompositions by the shift $c\equiv 3\pmod 4$ and the divisors of $K_j(c)$, one may dually change the order of summation and index directly over an auxiliary parameter $h$ and a square divisor $b\mid h^{2}$ (compare~\cite[\S2]{elsholtz2013}). We record here the resulting explicit certificates and the complete finite positivity equivalence for prime $n\equiv 1\pmod{24}$.

\begin{proposition}[Explicit square-divisor certificates]\label{aud:hbcert}
Let $n\equiv 1\pmod{24}$ be prime.
\begin{enumerate}[(i)]
    \item \emph{(Type~II).} If $h\ge1$, $b\mid h^{2}$, and $(4h-1)\mid(n+4b)$, set $c:=\frac{n+4b}{4h-1}$ and $K:=hc-b=\frac{hn+b}{4h-1}$. Then $c\equiv3\pmod 4$, $K=(n+c)/4>0$, $b\mid hK$, and
    \begin{equation}\label{aud:hbII}
    \frac{4}{n}\ =\ \frac{1}{K}\ +\ \frac{1}{hn}\ +\ \frac{1}{hKn/b}.
    \end{equation}
    \item \emph{(Type~I).} If $h>n/4$, $b\mid h^{2}$, and $(4h-n)\mid(4b+1)$, set $c:=\frac{4b+1}{4h-n}$ and $K:=hc-b=\frac{nb+h}{4h-n}$. Then $c\equiv3\pmod 4$, $K=(nc+1)/4>0$, $b\mid hK$, and
    \begin{equation}\label{aud:hbI}
    \frac{4}{n}\ =\ \frac{1}{nK}\ +\ \frac{1}{h}\ +\ \frac{1}{hK/b}.
    \end{equation}
\end{enumerate}
\end{proposition}
\begin{proof}
In both cases, the explicit expressions for $c$ and $K$ show that they are positive integers with $c\equiv 3\pmod 4$, and the relation $K=hc-b$ together with $b\mid h^{2}$ implies $b\mid hK$. The linear identities $h(4K-n)=K+b$ in case~(i) and $h(4K-1)=n(K+b)$ in case~(ii) then verify~\eqref{aud:hbII} and~\eqref{aud:hbI} immediately upon dividing by $hKn$.
\end{proof}

\begin{theorem}[Complete positivity criterion for prime $n$]\label{aud:hbcomplete}
For a prime $n\equiv1\pmod{24}$, define the two finite counting sums
\begin{align}
N_{\II}(n)\ &:=\ \sum_{1\le h\le(n-1)/6}\ \sum_{b\mid h^{2}}\ind_{\,4h-1\,\mid\,(n+4b)},\label{aud:NII}\\
N_{\I}(n)\ &:=\ \sum_{n/4<h<n/2}\ \sum_{b\mid h^{2}}\ind_{\,4h-n\,\mid\,(4b+1)}.\label{aud:NI}
\end{align}
Then $4/n$ is a sum of three positive unit fractions if and only if $N_{\I}(n)+N_{\II}(n)>0$.
\end{theorem}
\begin{proof}
The sufficiency of $N_{\I}(n)+N_{\II}(n)>0$ follows directly from Proposition~\ref{aud:hbcert}. Conversely, in any representation $4xyz=n(xy+xz+yz)$ for a prime $n$, either one or two of the denominators $x,y,z$ are multiples of $n$.

In the two-multiple case (Type~II), write $x=K$, $y=nh$, and $z=nk$ with $h\le k$ and $n\nmid K$. Setting $c=4K-n>0$ gives $chk=K(h+k)$, so $b:=ch-K=\frac{Kh}{k}$ is a positive integer satisfying $b\le K$ and $\gcd(b,c)=1$ (since any common divisor of $b$ and $c$ divides $K$ and hence $n$, which is coprime to $c$). Multiplying $K=ch-b$ by $h$ shows that $b\mid hK=h^{2}c-hb$, whence $b\mid h^{2}$ since $\gcd(b,c)=1$. Moreover, $(4h-1)c=4(K+b)-c=n+4b$, and since $c\ge 3$ and $b\le K$, we have $hc=K+b\le 2K=(n+c)/2$, which yields $h\le\frac{n}{2c}+\frac12\le\frac{n-1}{6}$.

In the one-multiple case (Type~I), write $x=nK$, $y=h$, and $z=k$ with $h\le k$ and $n\nmid hk$. Then $c=(4K-1)/n\equiv3\pmod 4$ is a positive integer satisfying $chk=K(h+k)$, so $b=ch-K=Kh/k$ is a positive divisor of $h^{2}$ satisfying $(4h-n)c=4b+1>0$ (forcing $h>n/4$), and the inequality $hc=K+b\le 2K=(nc+1)/2$ implies $h<\frac{n}{2}+\frac{1}{2c}<n/2$ since $n$ is odd and $c\ge 3$.
\end{proof}

\subsection{Open signed spectral error targets}\label{aud:open}

We conclude by recording the exact spectral error bounds that would be required to pass from the unconditional main-term lower bounds of Section~\ref{sec:toric_volume_barrier} and Subsection~\ref{aud:weighted} to a point-wise positivity proof for individual primes $n$.

In the notation of Theorem~\ref{thm:centered_main_term_and_short_lcm}, a sufficient condition for the cumulative encounter count $\mathcal{R}(n;n^{\theta})=\mathcal{M}(n;n^{\theta})+\mathcal{E}(n;n^{\theta})$ to be strictly positive for all primes $n>N_0$ ($n\equiv1\pmod{24}$) would be the signed error lower bound
\begin{equation}\label{aud:target17}
\mathcal{E}(n;n^{\theta})\ \ge\ -(\kappa_0\theta-\delta)\log n \tag{C17}
\end{equation}
for some $\theta\in(0,1)$ and $\delta\in(0,\kappa_0\theta)$ (or $\mathcal{E}(n;n^{\theta})\ge -o((\log n)^{4})$ in view of~\eqref{aud:log4eq}). Similarly, in the sifted square-free shift family $\mathcal{G}_j$ of Theorem~\ref{aud:weightedmain}, writing
\[
\mathcal{R}_j^{\sfm}\ :=\ \sum_{p\in\mathcal{G}_j}\#\big\{d,e\mid\rad K_j(p):\,d+e\equiv0\pmod p\big\}\ =\ \mathcal{M}_j^{\sfm}\ +\ \mathcal{E}_j^{\sfm}
\]
with the odd-character spectral error
\begin{equation}\label{aud:Esf}
\mathcal{E}_j^{\sfm}\ :=\ \sum_{p\in\mathcal{G}_j}\frac{1}{p-1}\sum_{\chi\ne\chi_0}\chi(-1)\Big(\big|F_{p,j}^{\sfm}(\chi)\big|^{2}-m_j(p)\Big),
\end{equation}
a sufficient condition for $\mathcal{R}_j^{\sfm}>0$ would be the signed bound
\begin{equation}\label{aud:target76}
\mathcal{E}_j^{\sfm}\ \ge\ -o\!\left(\frac{\log H}{(\log\log n)^{3/2}}\right). \tag{C76}
\end{equation}

Neither~\eqref{aud:target17}, nor~\eqref{aud:target76}, nor the Transversal Hypothesis~T is implied by the unconditional main-term lower bounds~\eqref{aud:log4eq} and~\eqref{aud:weightedmaineq}: as demonstrated by the exact small-LCM divisor-pair cancellation of Theorem~\ref{thm:centered_main_term_and_short_lcm}(c) and the archimedean support restriction of Proposition~\ref{aud:archsupport}, annihilating odd characters of order $<2z$ via Bombieri--Vinogradov at moduli $\le H^{1/2-\varepsilon}$ does not by itself control the signed contribution of characters of order $\ge 2z$ at the final encounter modulus $p\asymp H$.

\section*{Data availability}
No datasets were generated or analysed during the current study beyond elementary numerical checks of the algebraic identities described in the text.

\section*{Conflict of interest}
The author declares no conflict of interest.

\section*{Funding}
This research received no external funding.

\end{document}